%% file: report.tex
\documentclass[a4paper]{amsart}

\usepackage[mathscr]{eucal}
\usepackage{amsmath,amsfonts,amssymb}
\usepackage{mathrsfs}
\usepackage{bm}
\usepackage{multirow}
\usepackage[dvips]{graphicx}
\graphicspath{{./graphs/}}
\usepackage{fancyhdr}

\newcommand{\mathsym}[1]{{}}

\theoremstyle{plain}

\newtheorem{conjecture0}{Conjecture}
\newtheorem{theorem}{Theorem}[section]
\newtheorem{proposition}[theorem]{Proposition}

\newtheorem{corollary}{Corollary}[theorem]

\theoremstyle{definition}
\newtheorem{definition}{Definition}[section]
\newtheorem{remark}{Remark}[section]

\begin{document}
\quad\vspace{-0.1in}

\title[On the zeros of certain modular functions for the normalizers of congruence subgroups]{On the zeros of certain modular functions for \\ the normalizers of congruence subgroups of low levels \\ I}
\author{Junichi Shigezumi}

\maketitle \vspace{-0.1in}
\begin{center}
Graduate School of Mathematics Kyushu University\\
Hakozaki 6-10-1 Higashi-ku, Fukuoka, 812-8581 Japan\\
{\it E-mail address} : j.shigezumi@math.kyushu-u.ac.jp \vspace{-0.05in}
\end{center} \quad

\begin{quote}
{\small\bfseries Abstract.}
We research the location of the zeros of the Eisenstein series and the modular functions from the Hecke type Faber polynomials associated with the normalizers of congruence subgroups which are of genus zero and of level at most twelve.

In Part I, we will consider the general theory of modular functions for the normalizers.\\  \vspace{-0.15in}

\vspace{-0.15in}

 \vspace{0.15in}
\end{quote}

\input{intro.tex}\\

In `Part I', we will consider the general theory of modular functions for the normalizers of the congruence subgroups $\Gamma_0(N)$ of level $N \leqslant 12$. And in `Part II', we will observe the location of the zeros of the Eisenstein series and the the modular functions from Hecke type Faber polynomials for the normalizers in Part I by numerical calculation.

\clearpage

\markboth{\sc Zeros of modular functions for the normalizers of congruence subgroups}{\sc General Theory}
\pagestyle{myheadings}
\setcounter{section}{-1}

\input{report0.tex} \clearpage

\markboth{\sc Zeros of modular functions for the normalizers of congruence subgroups}{\sc Level \thesection}

\input{report01-05.tex} \clearpage

\input{report06-08.tex} \clearpage

\input{report09-11.tex} \clearpage

\input{report12-13.tex} \quad\\

\input{bib.tex}
\markright{\sc Zeros of modular functions for the normalizers of congruence subgroups}

\end{document}

%% file: intro.tex
\section*{Introduction}

The motive of this research is to decide the location of the zeros of modular functions. The Eisenstein series and the Hecke type Faber polynomials are the most interesting and important modular forms.

F. K. C. Rankin and H. P. F. Swinnerton-Dyer considered the problem of locating the zeros of the Eisenstein series $E_k(z)$ in the standard fundamental domain $\mathbb{F}$ (See \cite{RSD}). They proved that all of the zeros of $E_k(z)$ in $\mathbb{F}$ lie on the unit circle. They also stated towards the end of their study that ``This method can equally well be applied to Eisenstein series associated with subgroups of the modular group.'' However, it seems unclear how widely this claim holds. 

Subsequently, T. Miezaki, H. Nozaki, and the present author considered the same problem for the Fricke group $\Gamma_0^{*}(p)$ (see \cite{Kr}, \cite{Q}), and proved that all of the zeros of the Eisenstein series $E_{k, p}^{*}(z)$ in a certain fundamental domain lie on a circle whose radius is equal to $1 / \sqrt{p}$,  $p = 2, 3$ (see \cite{MNS}). Furthermore, we also proved that almost all the zeros of the Eisenstein series in a certain fundamental domain lie on circles whose radius are equal to $1 / \sqrt{p}$ or $1 / (2 \sqrt{p})$,  $p = 5, 7$ (see \cite{SJ2}).

Let $\Gamma$ be a discrete subgroup of $\text{\upshape SL}_2(\mathbb{R})$, and let $h$ be the width of $\Gamma$, then we define
\begin{equation}
\mathbb{F}_{0,\Gamma} := \left\{z \in \mathbb{H} \: ; \: - h / 2 < Re(z) < h / 2 \: , \: |c z + d| > 1 \: \text{for} \: \forall \gamma = \left(\begin{smallmatrix} a & b \\ c & d \end{smallmatrix}\right) \in \Gamma \: \text{s.t.} \: c \ne 0 \right\} \label{def-fd*}.
\end{equation}
We have a fundamental domain $\mathbb{F}_{\Gamma}$ such that $\mathbb{F}_{0,\Gamma} \subset \mathbb{F}_{\Gamma} \subset \overline{\mathbb{F}_{0,\Gamma}}$. Let $\mathbb{F}_{\Gamma}$ be such a fundamental domain.

For the modular group $\text{SL}_2(\mathbb{Z})$ and the Fricke groups $\Gamma_0^{*}(p)$ ($p = 2, 3$), all the zeros of the Eisenstein series for the cusp $\infty$ lie on the arcs on the boundary of their certain fundamental domains.

H. Hahn considered that the location of the zeros of the Eisenstein series for the cusp $\infty$ for every genus zero Fucksian group $\Gamma$ of the first kind with $\infty$ as a cusp which satisfies that its hauptmodul $J_{\Gamma}$ takes real value on $\partial \mathbb{F}_{\Gamma}$, and proved that almost all the zeros of the Eisenstein series for the cusp $\infty$ for $\Gamma$ lie on $\partial \mathbb{F}_{\Gamma}$ under some more assumption (see \cite{H}).

Also, T. Asai, M. Kaneko, and H. Ninomiya considered the problem of locating the zeros of modular functions $F_m(z)$ for $\text{SL}_2(\mathbb{Z})$ which correspond to the Hecke type Faber polynomial $P_m$, that is, $F_m(z) = P_m (J(z))$ (See \cite{AKN}). They proved that all of the zeros of $F_m(z)$ in $\mathbb{F}$ lie on the unit circle for each $m \geqslant 1$. After that, E. Bannai, K. Kojima, and T. Miezaki considered the same problem for the normalizers of congruence subgroups which correspond the conjugacy classes of the Monster group (See \cite{BKM}). They observed the location of the zeros by numerical calculation, then almost all of the zeros of the modular functions from Hecke type Faber polynomial lie on the lower arcs when the group satisfy the same assumption of the theorem of H. Hahn. In particular, T. Miezaki proved that all of the zeros of the modular functions from the Hecke type Faber polynomials for the Fricke group $\Gamma_0^{*}(2)$ lie on the lower arcs of its fundamental domain in their paper.

Now, we have the following conjectures:
\begin{conjecture0}\label{conj0}
Let $\Gamma$ be a genus zero Fucksian group of the first kind with $\infty$ as a cusp. If the hauptmodul $J_{\Gamma}$ takes real value on $\partial \mathbb{F}_{\Gamma}$, all of the zeros of the Eisenstein series for the cusp $\infty$ for $\Gamma$ in $\mathbb{F}_{\Gamma}$ lie on the arcs
\begin{equation*}
\partial \mathbb{F}_{\Gamma} \setminus \{z \in \mathbb{H} \: ; \: Re(z) = \pm h / 2\}.
\end{equation*}
\end{conjecture0}

\begin{conjecture0}\label{conj01}
Let $\Gamma$ be a genus zero Fucksian group of the first kind with $\infty$ as a cusp. If the hauptmodul $J_{\Gamma}$ takes real value on $\partial \mathbb{F}_{\Gamma}$, all but at most $c_h(\Gamma)$ of the zeros of modular function from the Hecke type Faber polynomial of degree $m$ for $\Gamma$ in $\mathbb{F}_{\Gamma}$ lie on the arcs
\begin{equation*}
\partial \mathbb{F}_{\Gamma} \setminus \{z \in \mathbb{H} \: ; \: Re(z) = \pm h / 2\}
\end{equation*}
for all but finite number of $m$ and for the constant number $c_h(\Gamma)$ which does not depend on $m$.
\end{conjecture0}

In this paper, we will observe the location of the zeros of the Eisenstein series and the modular functions from Hecke type Faber polynomials for the normalizers of congruence subgroups, as a first step of a challenge for the above conjectures.\\

The normalizers of congruence subgroups of level at most $12$ which satisfies the assumption of above conjectures are 
\begin{align*}
&\text{SL}_2(\mathbb{Z}), \; \Gamma_0^{*}(2), \; \Gamma_0(2), \; \Gamma_0^{*}(3), \; \Gamma_0(3), \; \Gamma_0^{*}(4), \; \Gamma_0(4), \; \Gamma_0^{*}(5), \; \Gamma_0(6)+, \; \Gamma_0^{*}(6), \; \Gamma_0(6)+3, \; \Gamma_0(6),\\
&\Gamma_0^{*}(7), \; \Gamma_0^{*}(8), \; \Gamma_0(8), \; \Gamma_0^{*}(9), \; \Gamma_0(10)+, \; \Gamma_0^{*}(10), \; \Gamma_0(10)+5, \; \Gamma_0(12)+, \; \Gamma_0^{*}(12),\\
&\Gamma_0(12)+4, \; \text{and} \; \Gamma_0(12).
\end{align*}

For the Conjecture \ref{conj0}, $\text{SL}_2(\mathbb{Z})$, $\Gamma_0^{*}(2)$, and $\Gamma_0^{*}(3)$ verify  Conjecture \ref{conj0}. For the other cases, we can prove by numerical calculation for the Eisenstein series of weight $k \leqslant 500$.

For the Conjecture \ref{conj01}, $\text{SL}_2(\mathbb{Z})$ and $\Gamma_0^{*}(2)$ verify  Conjecture \ref{conj01} for every degree $m$, where we have $c_h(\Gamma) = 0$ for each case. Furthermore, for $\Gamma_0(2)$, $\Gamma_0^{*}(3)$, $\Gamma_0(3)$, $\Gamma_0^{*}(4)$, $\Gamma_0(4)$, $\Gamma_0(6)+$, $\Gamma_0(6)+3$, $\Gamma_0(6)$, $\Gamma_0(8)$, $\Gamma_0^{*}(9)$, $\Gamma_0(10)+$, $\Gamma_0(10)+5$, $\Gamma_0(12)+$, $\Gamma_0(12)+4$, and $\Gamma_0(12)$, we can prove all of the zeros of the modular function from the Hecke type Faber polynomial of every degee $m \leqslant 200$ in each fundamental domain lie on the lower arcs  by numerical calculation.

On the other hand, for $\Gamma_0^{*}(5)$ and $\Gamma_0^{*}(7)$, we can prove by numerical calculation for the modular function from the Hecke type Faber polynomial of every degee $m = 1$ and $3 \leqslant m \leqslant 200$, where we have $c_h(\Gamma) = 0$ for each case. When $m = 2$, there is just one zero which is on the boundary of its fundamental domain but not on the lower arcs for the each group.

For $\Gamma_0^{*}(6)$ and $\Gamma_0^{*}(8)$, we can prove by numerical calculation for the modular function from the Hecke type Faber polynomial of every degee $m \leqslant 200$ which satisfy $m \not\equiv 0 \pmod{2}$ and $m \not\equiv 2 \pmod{4}$, respectively. For the remaining degrees, there is just one zero which is on the boundary of its fundamental domain but not on the lower arcs for the each group, that is, $c_h(\Gamma) = 1$.

Finally, for $\Gamma_0^{*}(10)$ and $\Gamma_0^{*}(12)$, we have just two zeros which are not on the boundary of each fundamental domain for degrees $m =7, 9, 11$ and $m = 3, 6, 12, 13, 15$, respectively. Furthermore, there is just one zero which is on the boundary of its fundamental domain but not on the lower arcs for the case $m \equiv 0 \pmod{2}$ and $m \equiv 2, 4 \pmod{6}$, respectively. For the other cases, we can prove that all of the zeros are on the lower arcs of each fundamental domain by numerical calculation.

\begin{table}[htbp]
{\small \begin{center}
\begin{tabular}{ccc}
\hline
\hspace{-0.5in}$\Gamma$ & \hspace{-0.3in}Eisenstein series $(k \leqslant 500)$ & Hecke type Faber polynomial $(m \leqslant 200)$\\
\hline
\begin{tabular}{l}
$\text{SL}_2(\mathbb{Z})$, $\Gamma_0^{*}(2)$, $\Gamma_0(2)$, $\Gamma_0^{*}(3)$, $\Gamma_0(3)$,\\ $\Gamma_0^{*}(4)$, $\Gamma_0(4)$, $\Gamma_0(6)+$, $\Gamma_0(6)+3$,\\
$\Gamma_0(6)$, $\Gamma_0(8)$, $\Gamma_0^{*}(9)$, $\Gamma_0(10)+$, $\Gamma_0(10)+5$,\\
$\Gamma_0(12)+$, $\Gamma_0(12)+4$, $\Gamma_0(12)$.
\end{tabular} & \multirow{6}{*}{\hspace{-0.15in}$\bigcirc$} & $\bigcirc$\\
$\Gamma_0^{*}(5)$, $\Gamma_0^{*}(7)$ && $m = 2$, \quad $\langle 1 \rangle$\\
$\Gamma_0^{*}(6)$ && $m$ : even, \quad $\langle 1 \rangle$\\
$\Gamma_0^{*}(8)$ && $m \equiv 0 \pmod{4}$, \quad $\langle 1 \rangle$\\
$\Gamma_0^{*}(10)$ && $m = 7, 9, 11, \; [2], \quad m$ : even, $\langle 1 \rangle$\\
$\Gamma_0^{*}(12)$ && $m = 3, 6, 12, 13, 15, \; [2] \quad m \equiv 2, 4 \pmod{6}, \; \langle 1 \rangle$\\
\hline
\end{tabular}
\end{center}
\begin{flushleft}
\hspace{0.51in}`$\bigcirc$': all of the zeros lie on lower arcs.\\
\hspace{0.5in}$\langle \ \cdot \ \rangle$: the number of zeros which are on $\partial \mathcal{F}$ but not on lower arcs.\\
\hspace{0.51in}$[ \ \cdot \ ]$\hspace{0.02in}: the number of zeros which are not on $\partial \mathcal{F}$.
\end{flushleft}}
\caption{Result by numerical calculation}
\end{table}

\begin{figure}[hbtp]
\begin{center}
{{\small Lower arcs of $\partial \mathbb{F}_{6+2}$}\includegraphics[width=2in]{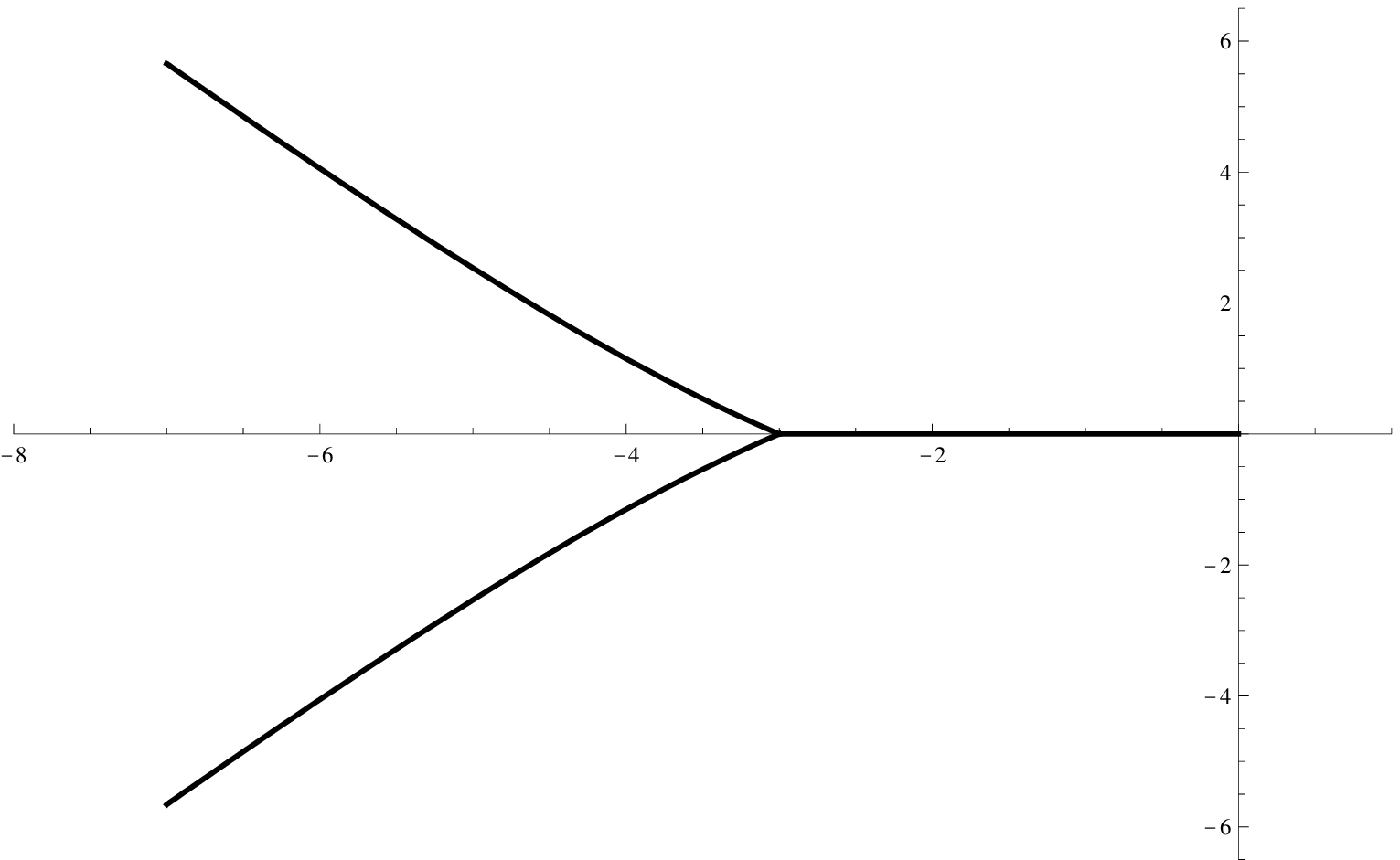}}
\end{center}
\caption{Image by $J_{6+2}$ ($\Gamma_0(6)+2$)}\label{Im-J6D-intro0}
\end{figure}

If the hauptmodul $J_{\Gamma}$ does not take real value on $\partial \mathbb{F}_{\Gamma}$ (cf. Figure \ref{Im-J6D-intro0}), it seems to be not similar. Such cases are followings;
\begin{equation*}
\Gamma_0(5), \; \Gamma_0(6)+2, \; \Gamma_0(7), \; \Gamma_0(9), \; \Gamma_0(10)+2, \; \Gamma_0(10), \; \Gamma_0^{*}(11), \; \text{and} \; \Gamma_0(12)+3.
\end{equation*}

\begin{figure}[hbtp]
\begin{center}
{{$\begin{matrix}\text{\small The zeros of $E_{k, 6+2}^{\infty}$} \\ \text{\small for $4 \leqslant k \leqslant 40$}\end{matrix}$}\includegraphics[width=4.2in]{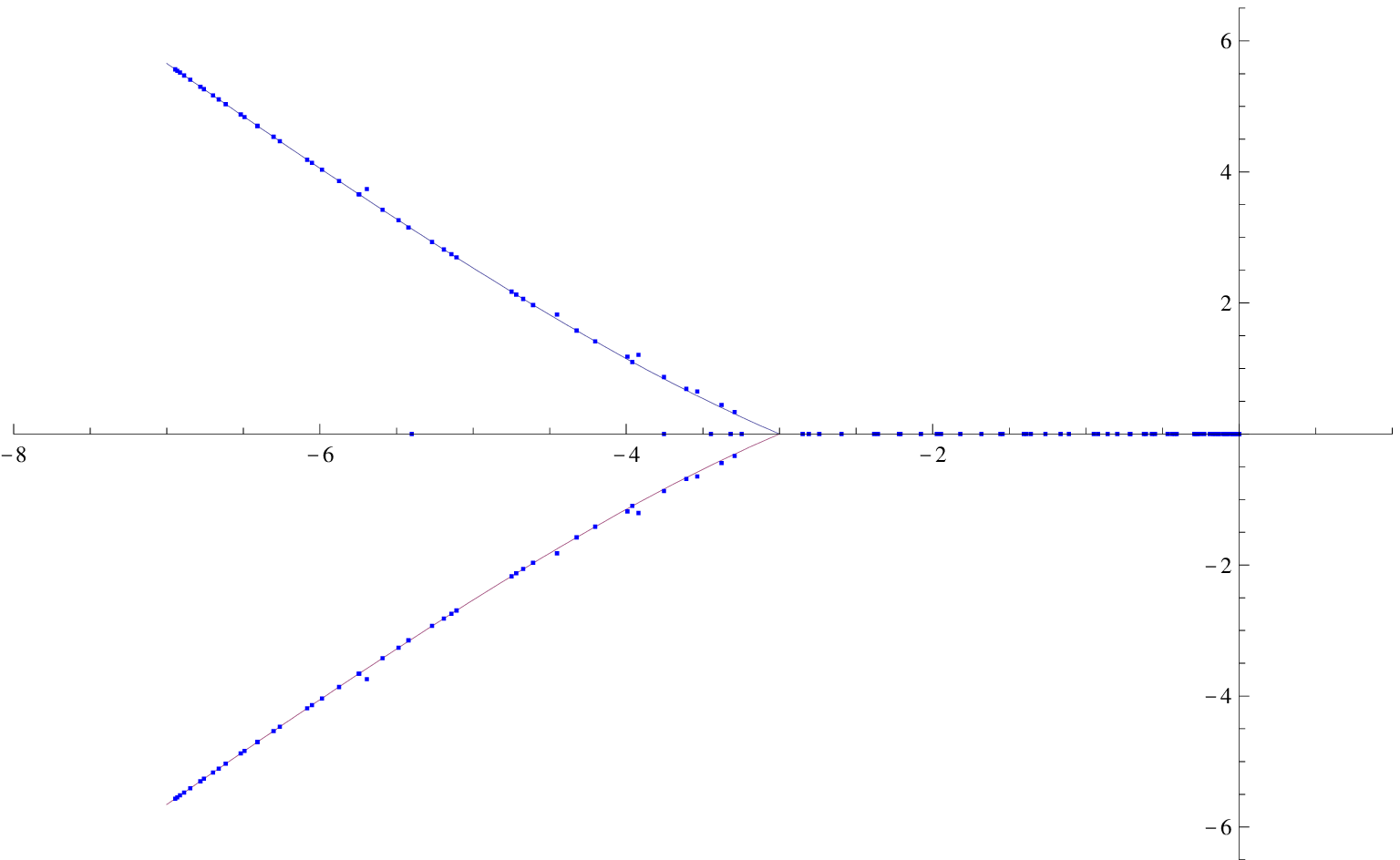}}\\
{{\small The zeros of $E_{1000, 6+2}^{\infty}$}\includegraphics[width=4.2in]{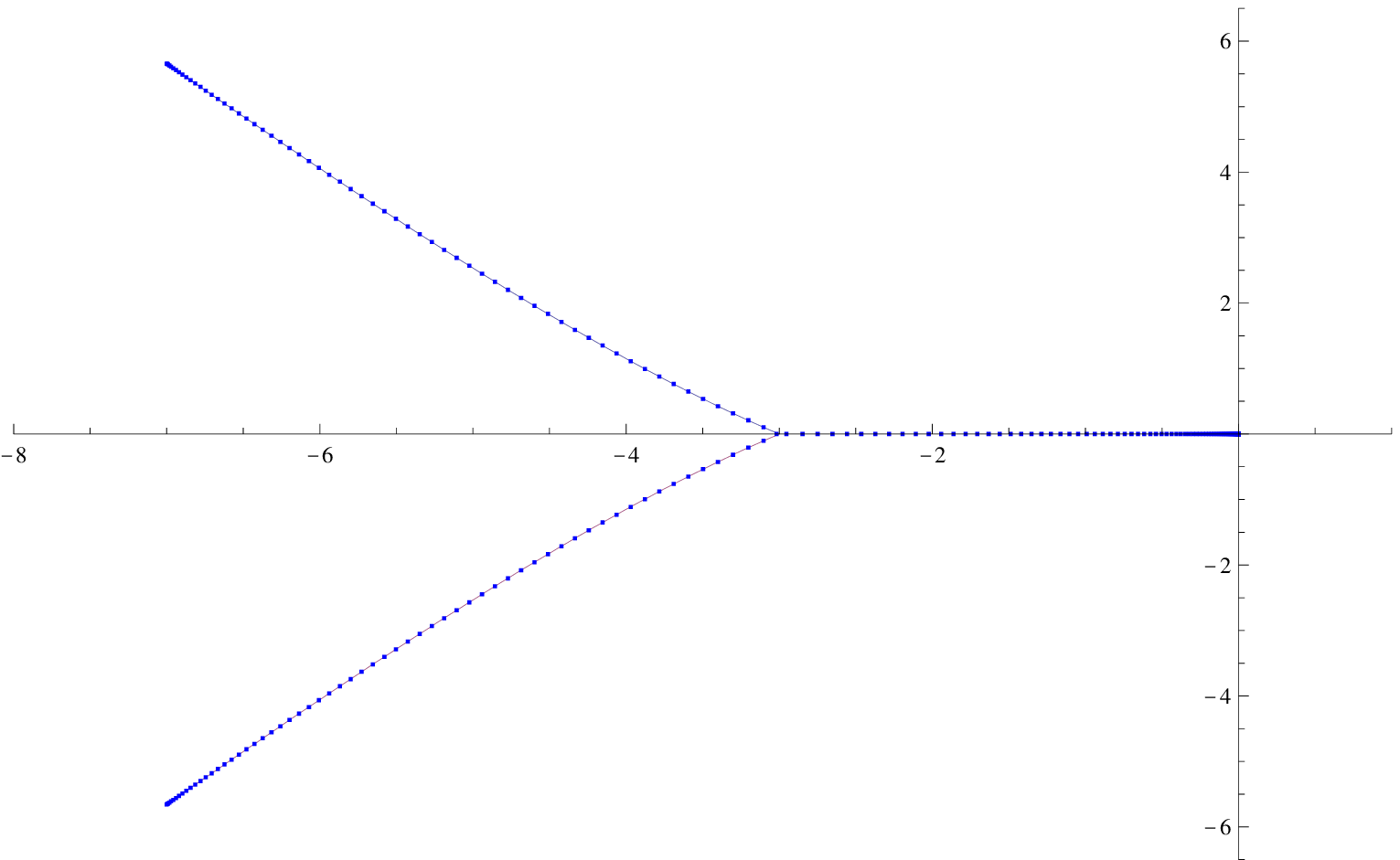}}
\end{center}
\caption{Image by $J_{6+2}$ ($\Gamma_0(6)+2$)}\label{Im-J6D-intro1}
\end{figure}

For $\Gamma_0(5)$, $\Gamma_0(6)+2$, $\Gamma_0(7)$, $\Gamma_0(10)+2$, $\Gamma_0(10)$, and $\Gamma_0^{*}(11)$, we can observe that the zeros of the Eisenstein series for cusp $\infty$ do not lie on the lower arcs of their fundamental domains by numerical calculation. However, when the weight of Eisenstein series increases, then the location of the zeros seems to approach to lower arcs. (See Figure \ref{Im-J6D-intro1})

\begin{figure}[hbtp]
\begin{center}
{{$\begin{matrix}\text{\small The zeros of $F_{m, 6+2}$} \\ \text{\small for $1 \leqslant m \leqslant 40$}\end{matrix}$}\includegraphics[width=4.2in]{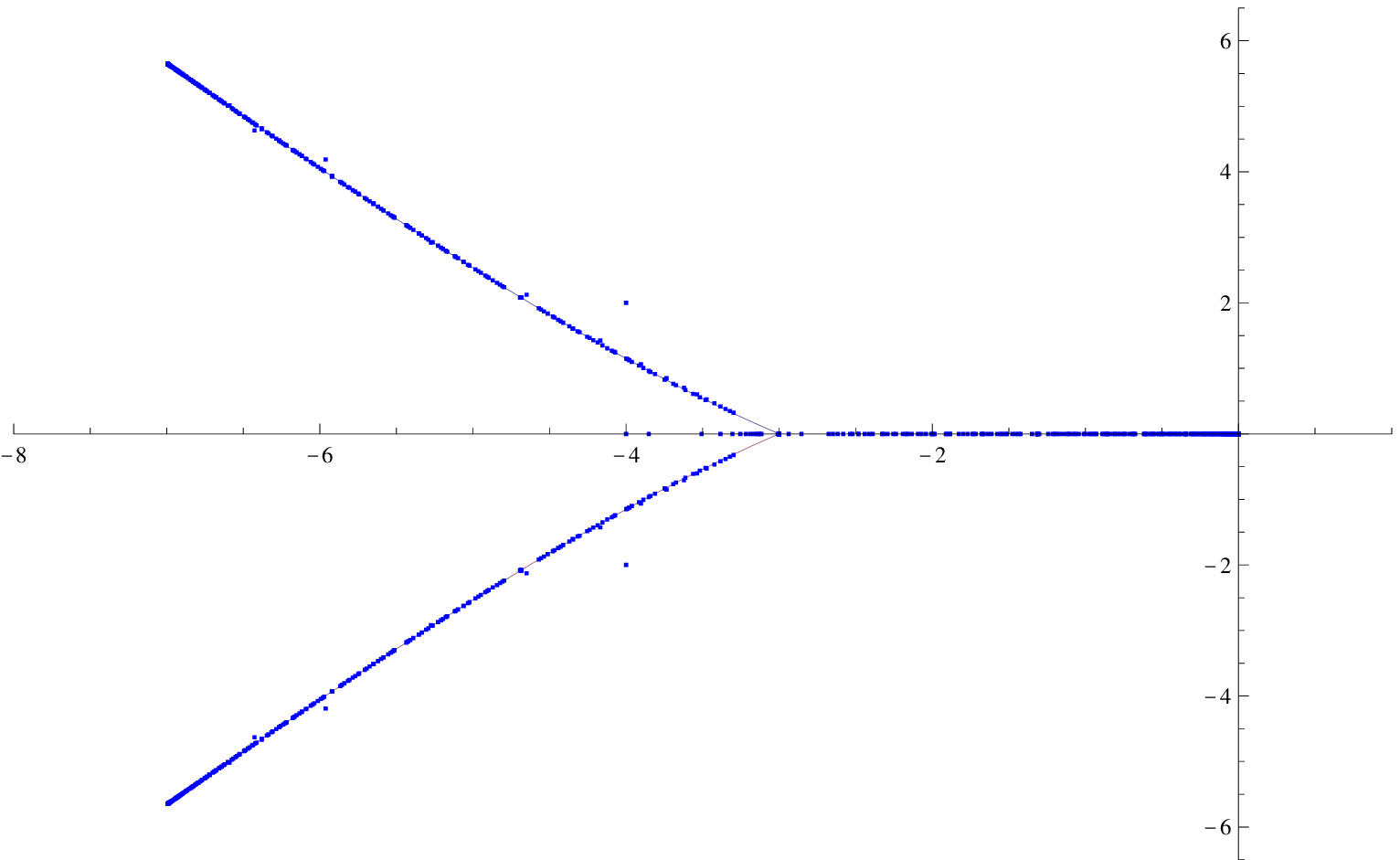}}\\
{{\small The zeros of $F_{200, 6+2}$}\includegraphics[width=4.2in]{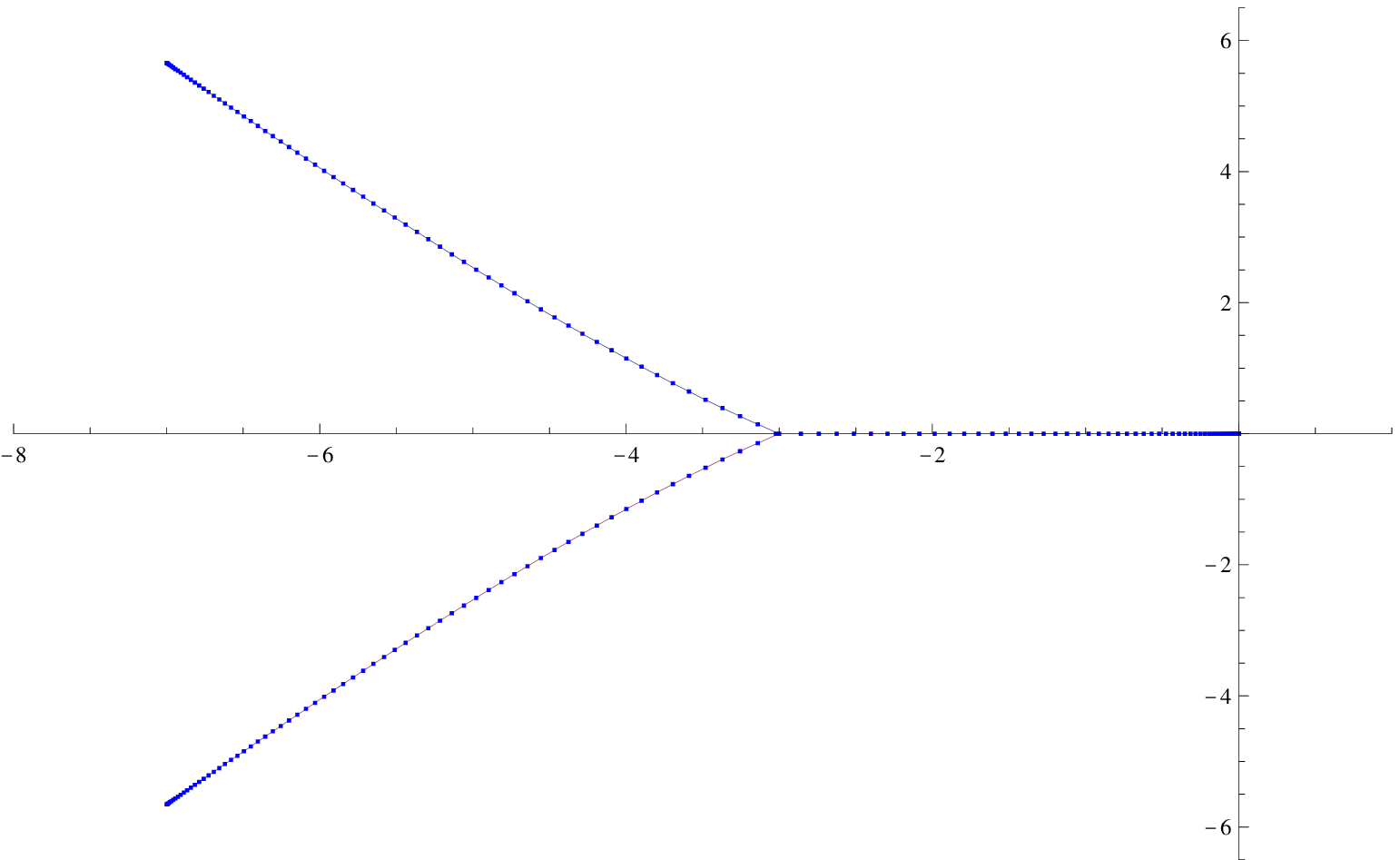}}
\end{center}
\caption{Image by $J_{6+2}$ ($\Gamma_0(6)+2$)}\label{Im-J6D-intro2}
\end{figure}

Also, for the zeros of the modular functions from the Hecke type Faber polynomials, we can observe that there are some zeros which do not lie on the lower arcs of their fundamental domains by numerical calculation. Furthermore, when the degree $m$ increases, then the location of the zeros seems to approach to lower arcs. (See Figure \ref{Im-J6D-intro2})

\begin{figure}[hbtp]
\begin{center}
{{\small $\Gamma_0(9)$}\includegraphics[width=2in]{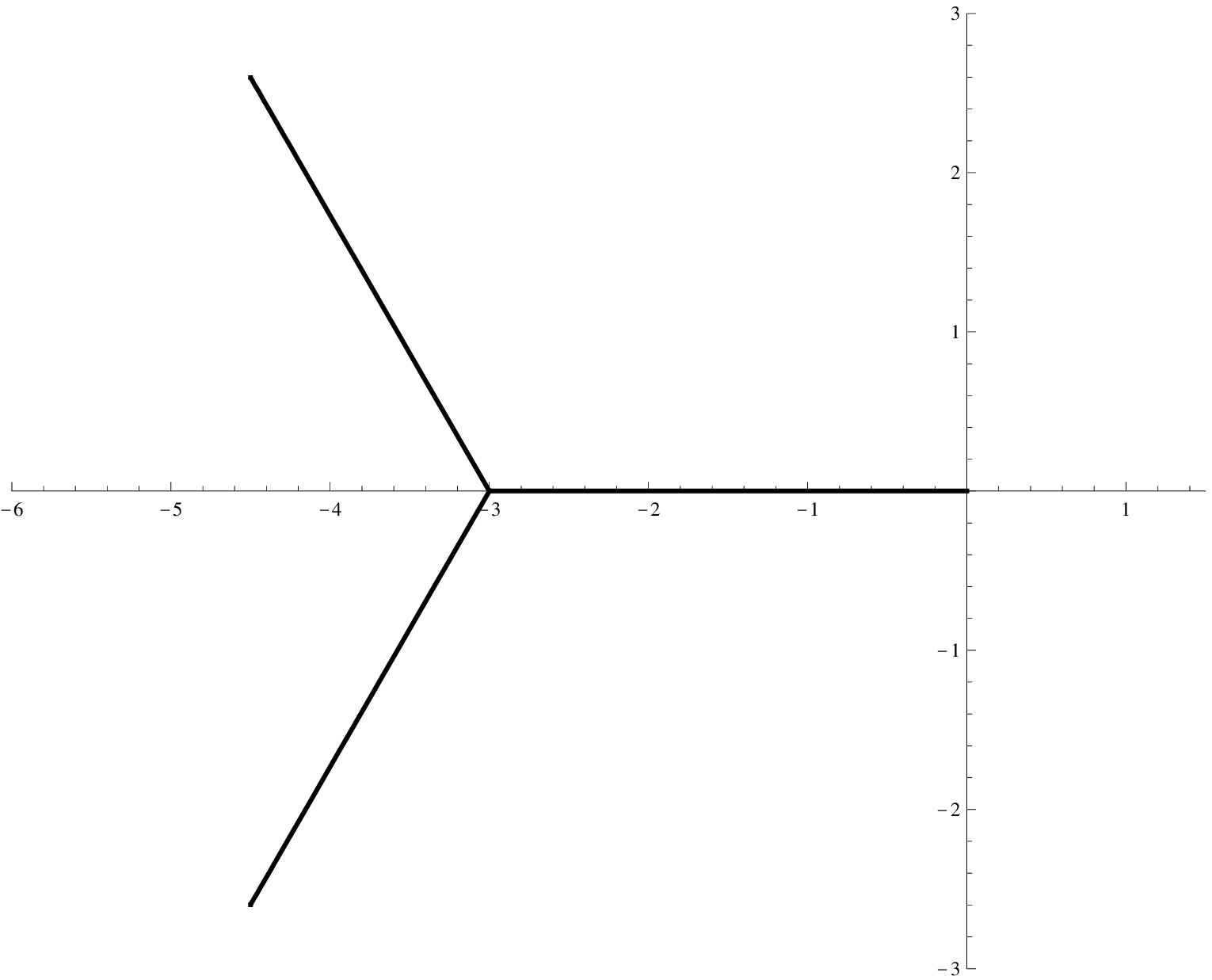}} \;
{{\small $\Gamma_0(12)+3$}\includegraphics[width=2.6in]{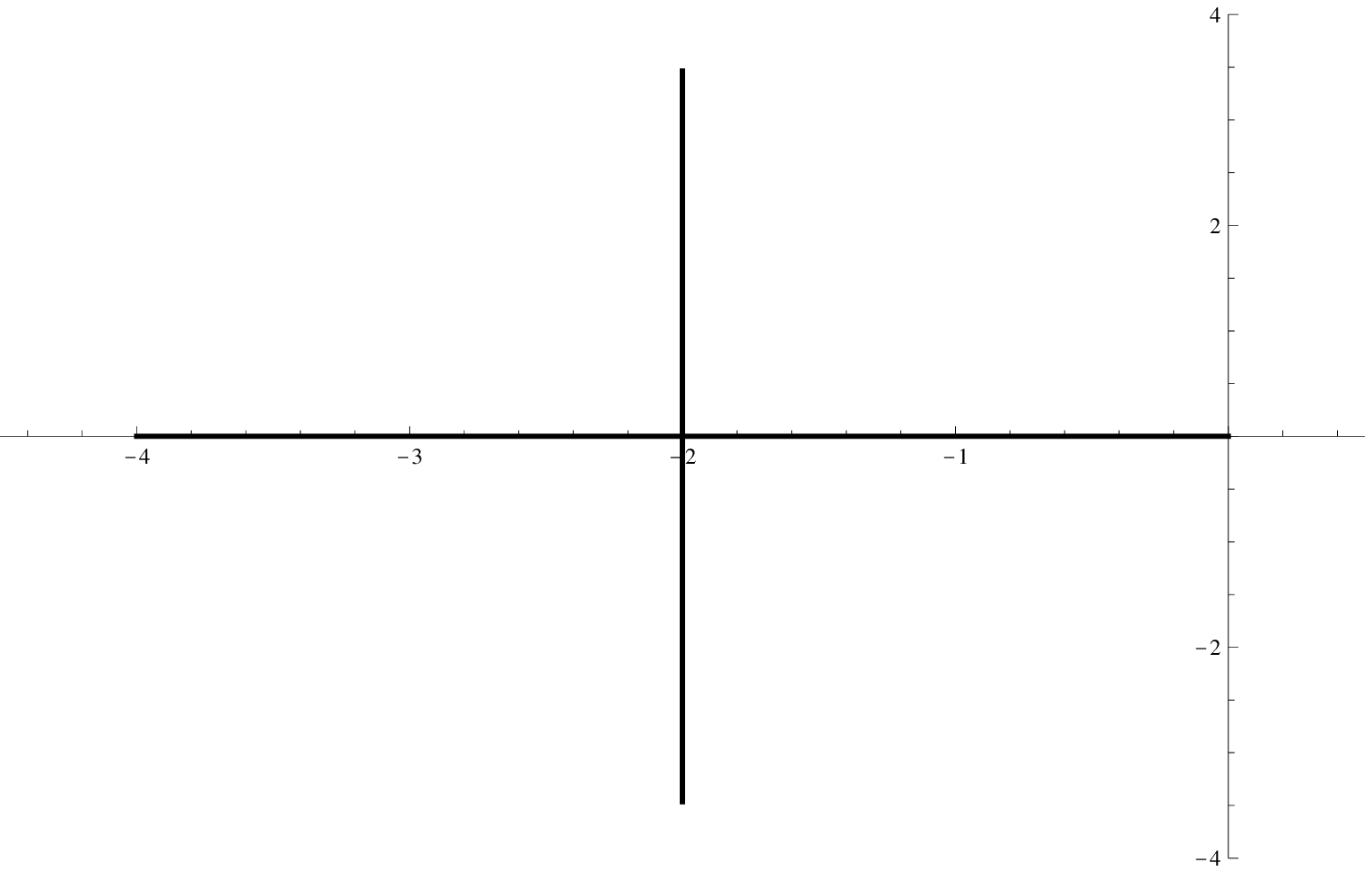}}
\end{center}
\caption{Image of the lower arcs of the fundamental domains by hauptmoduls}\label{Im-SP-intro}
\end{figure}

On the other hand, $\Gamma_0(9)$ and $\Gamma_0(12)+3$ seem to show the special cases. We can prove that all of the zeros of the Eisenstein series of weight $k \leqslant 500$ lie on the lower arcs of their fundamental domains by numerical calculation. Also, we can prove that all of the zeros of the modular function from the Hecke type Faber polynomial of degee $m \leqslant 200$ lie on the lower arcs by numerical calculation. On the other hand, they do not satisfy the assumption of Conjecture \ref{conj0} and \ref{conj01}. However, the image of lower arcs by its hauptmodul draw a interesting figure. (Figure \ref{Im-SP-intro})

\newpage
\quad

We refer to \cite{MNS}, \cite{SJ1}, and \cite{SJ2} for some groups. However, note that definitions in this paper are sometimes different from that in it.

%% file: report0.tex
\section{General theory} \label{sec-def}

Let $\Gamma$ be a Fucksian group of the first kind with $\infty$ as a cusp.

\subsection{The modular group and some groups}

\subsubsection{The modular group} (see \cite[VII.1]{Se}, \cite[I-I]{SH}, and \cite[III.1]{Ko})

We have the {\it special linear group} defined by following:
\begin{equation}
\text{SL}_2(\mathbb{R}) := \left\{ \left(\begin{smallmatrix} a & b \\ c & d \end{smallmatrix}\right) \: ; \: \forall a,b,c,d \in \mathbb{R} \: s.t. \: a d - b c = 1 \right\}.
\end{equation}

Write $\mathbb{H}:=\{z \in \mathbb{C} \: ; \: Im(z) > 0\}$, which is complex upper half-plane. We consider a action of $\text{SL}_2(\mathbb{R})$ on $\mathbb{H} \cup (\{\infty\} \cup \mathbb{R})$ in the following way:
\begin{quote}
For every $\gamma = \left(\begin{smallmatrix} a & b \\ c & d \end{smallmatrix}\right) \in \text{SL}_2(\mathbb{R})$ and every $z \in \mathbb{H}$, we put
\begin{equation}
\gamma z := \frac{a z + b}{c z + d}. \label{eq-lft}
\end{equation}
\end{quote}
Note that $- \gamma z = \gamma z$ and $Im(\gamma z) = Im(z) / |c z + d|^2$ for every $z \in \mathbb{H}$.

We also define
\begin{equation}
\text{SL}_2(\mathbb{Z}) := \left\{\begin{pmatrix} a & b \\ c & d \end{pmatrix} \in \text{SL}_2(\mathbb{R}) \: ; \: \forall a,b,c,d \in \mathbb{Z} \right\},
\end{equation}
which is called the {\it $($full$)$ modular group}.\\

\subsubsection{Congruence subgroup} (see \cite[III.1]{Ko})

For a positive integer $N$, we have
\begin{equation}
\Gamma(N) := \left\{\left(\begin{smallmatrix} a & b \\ c & d \end{smallmatrix}\right) \in \text{SL}_2(\mathbb{Z}) \: ; \: a \equiv d \equiv 1, \: b \equiv c \equiv 0 \pmod{N} \right\}.
\end{equation}
This group is a subgroup of the modular group $\text{SL}_2(\mathbb{Z})$, which is called the {\it principal congruence subgroup of level $N$}.

Also, if $\Gamma'$ is a subgroup of $\text{SL}_2(\mathbb{Z})$ such that $\Gamma' \supset \Gamma(N)$, then $\Gamma'$ is called a {\it congruence subgroup of level $N$}. Here are some examples:
\begin{gather}
\Gamma_0(N) := \left\{\left(\begin{smallmatrix} a & b \\ c & d \end{smallmatrix}\right) \in \text{SL}_2(\mathbb{Z}) \: ; \: c \equiv 0 \pmod{N} \right\},\\
\Gamma_1(N) := \left\{\left(\begin{smallmatrix} a & b \\ c & d \end{smallmatrix}\right) \in \text{SL}_2(\mathbb{Z}) \: ; \: a \equiv d \equiv 1, \: c \equiv 0 \pmod{N} \right\}.
\end{gather}
\quad

\subsubsection{Fricke group} (see \cite{Kr}, \cite{Q})

For a positive integer $N$, we consider the {\it Fricke group} $\Gamma_0^{*}(N)$. We define the following;
\begin{equation}
\Gamma_0^{*}(N) := \Gamma_0(N) \cup \Gamma_0(N) \: W_N, \quad W_N := \begin{pmatrix}0&-1 / \sqrt{N}\\ \sqrt{N}&0\end{pmatrix}.
\end{equation}
This group is a discrete subgroup of $\text{SL}_2(\mathbb{R})$, and commensurable with $\text{SL}_2(\mathbb{Z})$.\\

\subsubsection{The normalizer of $\Gamma_0(N)$} (see \cite{CN})

We introduce the normalizer of $\Gamma_0(N)$. About notations, there is something different from \cite{CN}.

We fix positive integers $n$ and $h$ such that $h | n$ and $h | 24$. Then, we define
\begin{equation}
\left\{ \begin{matrix}a & b \\ c & d\end{matrix} \right\}_e
 := \frac{1}{\sqrt{e}} \begin{pmatrix}a e & b / h \\ c n & d e\end{pmatrix}
\end{equation}
for the integers $e > 0$ and $a, b, c, d$ such that $e || (n / h)$ and $a d e - b c (n / h) / e = 1$. We denote $I_{n | h} := \{e \in \mathbb{N} \: ; \: e || (n / h) \}$, where $n || m$ means that $n | m$ and $(n, m/n) = 1$.

Then, we define
\begin{allowdisplaybreaks}
\begin{align}
\Gamma_0(n | h) &:= \left\{ \left\{ \begin{smallmatrix}a & b \\ c & d\end{smallmatrix} \right\}_ 1 \; ; \; a, b, c, d \in \mathbb{Z} \; \text{s.t.} \; a d - b c (n / h) = 1 \right\}.\\
\Gamma_0(n | h) + e_1, e_2, ... , e_m &:= \left\{ \left\{ \begin{smallmatrix}a & b \\ c & d\end{smallmatrix} \right\}_ e \; ; \;
 \begin{matrix}
  a, b, c, d \in \mathbb{Z} \; \text{s.t.} \; a d e - b c (n / h) / e = 1\\
  e \in \{ 1, e_1, ... , e_m \}
 \end{matrix} \right\},\\
 &\text{where} \; \{ 1, e_1, ... , e_m \} \subset I_{n | h}.\notag \\
\Gamma_0(n | h) + &:= \left\{ \left\{ \begin{smallmatrix}a & b \\ c & d\end{smallmatrix} \right\}_ e \; ; \;
 \begin{matrix}
  a, b, c, d \in \mathbb{Z} \; \text{s.t.} \; a d e - b c (n / h) / e = 1\\
  e \in I_{n | h} 
 \end{matrix} \right\}.
\end{align}
\end{allowdisplaybreaks}
In addition, we have $\Gamma_0(n) = \Gamma_0(n | 1)$ and denote $\Gamma_0(n | h) - := \Gamma_0(n | h)$.

For example, for a prime number $p$, we consider the case $h = 1$ and $n = p$, then we have $I_{p | 1} = \{ 1, p \}$ and
\begin{align*}
\left\{ \begin{smallmatrix}a & b \\ c & d\end{smallmatrix} \right\}_1
 &= \left( \begin{smallmatrix}a & b \\ c p & d\end{smallmatrix} \right) \in \Gamma_0(p),\\
\left\{ \begin{smallmatrix}a & b \\ c & d\end{smallmatrix} \right\}_p
 &= \left( \begin{smallmatrix}a \sqrt{p} & b / \sqrt{p} \\ c \sqrt{p} & d \sqrt{p}\end{smallmatrix} \right)
 = \left( \begin{smallmatrix}-b & a \\ -d p & c\end{smallmatrix} \right) \left( \begin{smallmatrix}0 & -1 / \sqrt{p} \\ \sqrt{p} & 0\end{smallmatrix} \right)
 \in \Gamma_0(p) W_p.
\end{align*}
Thus, we have $\Gamma_0(p)+ = \Gamma_0(p) + p = \Gamma_0^{*}(p)$ and $\Gamma_0(p)- = \Gamma_0(p)$. That is, Fricke group $\Gamma_0^{*}(p)$ is a normalizers of $\Gamma_0(p)$.\\

\subsubsection{Preliminaries}
Write $T_x := \left(\begin{smallmatrix} 1 & x \\ 0 & 1 \end{smallmatrix}\right) \: (\in \text{\upshape SL}_2(\mathbb{R}))$ and $P := \{ \pm T_x \: ; \: x \in \mathbb{R} \}$. In this paper, we may assume
\begin{equation}
\Gamma \cap P \setminus \left\{\pm I \right\} \ne \phi
\end{equation}
and $\Gamma \cap P$ is discrete. Then, we call $h := \min\{ x > 0 \: ; \: T_x \in \Gamma \}$ the {\it width of $\Gamma$}.

For a cusp $\kappa$ of $\Gamma$, we define {\it the stabilizer of the cusp $\kappa$}:
\begin{equation*}
\Gamma_{\kappa} := \{ \gamma \in \Gamma \: ; \: \gamma \kappa = \kappa \}.
\end{equation*}
In particular, we have $\Gamma_{\infty} = \Gamma \cap P = \{ T_{n h} \: ; \: n \in \mathbb{Z} \}$. Furthermore, there exist some $\gamma_{\kappa} \in \text{\upshape SL}_2(\mathbb{R})$ such that $\gamma_{\kappa} \infty = \kappa$ and
\begin{equation*}
\Gamma_{\kappa} = \gamma_{\kappa} \Gamma_{\infty} \gamma_{\kappa}^{-1}.
\end{equation*}
We call $\gamma_{\kappa}$ the {\it cusp leader} of the cusp $\kappa$. Note that, $\Gamma_{\kappa}$ and $\gamma_{\kappa}$ are depend on the group $\Gamma$.\\

\subsection{Fundamental domain}(see \cite[VII.1]{Se}, \cite[I-I.4]{SH})

In this section, we consider a fundamental domain in $\mathbb{H}$ under the action of $\Gamma$ (equation (\ref{eq-lft})).
\begin{definition}
$\mathbb{F}_\Gamma$ is a {\it fundamental domain} of $\Gamma$ if and only if it satisfies following conditions:
\begin{trivlist}
\item[(FD1)] For every $z \in \mathbb{H}$, there exists $\gamma \in \Gamma$ such that $\gamma z \in \mathbb{F}_{\Gamma}$.
\item[(FD2)] For every two distinct points $z_1$, $z_2$ $\in \mathbb{F}_{\Gamma}$, there does not exist $\gamma \in \Gamma$ such that $\gamma z_1 = z_2$.
\end{trivlist}
\label{def-fd}
\end{definition}

We define
\begin{equation*}
\mathbb{F}_{0,\Gamma} := \left\{z \in \mathbb{H} \: ; \: - h / 2 < Re(z) < h / 2 \: , \: |c z + d| > 1 \: \text{for} \: \forall \gamma = \left(\begin{smallmatrix} a & b \\ c & d \end{smallmatrix}\right) \in \Gamma \setminus P \right\}
\end{equation*}
(see \cite[I-I.4 \& 5, Theorem 1.7 and 1.15]{SH}). Now, we have the following fact:
\begin{proposition}\quad 
\def\labelenumi{(\roman{enumi})}
\begin{enumerate}
\item $\overline{\mathbb{F}_{0,\Gamma}}$ satisfies the condition {\upshape (FD1)}.
\item $\mathbb{F}_{0,\Gamma}$ satisfies the condition {\upshape (FD2)}.
\end{enumerate}
\def\labelenumi{\arabic{enumi}.}
\end{proposition}

Thus, we have a fundamental domain $\mathbb{F}_{\Gamma}$ such that $\mathbb{F}_{0,\Gamma} \subset \mathbb{F}_{\Gamma} \subset \overline{\mathbb{F}_{0,\Gamma}}$. Note that, let $z, z' \in \partial \mathbb{F}_{0, \Gamma}$. If $\gamma z = z'$ for some $\gamma = \left(\begin{smallmatrix} a & b \\ c & d \end{smallmatrix}\right) \in \Gamma \setminus P$, then $|c z + d| = 1$. Also, we have $z = (e^{i \theta} - d) / c$ and $z' = (e^{i (\pi - \theta)} + a) / c$.

We define $\Gamma^0 := \{ \gamma = \left(\begin{smallmatrix} a & b \\ c & d \end{smallmatrix}\right) \in \Gamma \: ; \: | c z + d | = 1 \: \text{for} \: \exists z \in \mathbb{F}_{0,\Gamma} \}$, then we have the following fact:
\begin{corollary}
If $\Gamma \ni -I$, then $\Gamma^0 \cup \{ T_h, \: -I \}$ generates $\Gamma$. On the other hand, if $\Gamma \not\ni -I$, $\Gamma^0 \cup \{ T_h \}$ generates $\Gamma$. \label{cor-base}
\end{corollary}\quad

\subsection{Modular forms}

\subsubsection{Preliminaries}(see \cite[III.2, 3]{Ko} and \cite[VII.1]{Se})

Let $f$ be a function on $\mathbb{H}$. For $\gamma = \left( \begin{smallmatrix} a & b \\ c & d \end{smallmatrix} \right) \in \text{SL}_2(\mathbb{R})$, we denote
\begin{equation}
f \big|_k \gamma \ (z) := (c z + d)^{-k} f(\gamma z).
\end{equation}
Then, the relation
\begin{equation}
f \big|_k \gamma \ (z) = f(z) \quad \text{for every} \: z \in \mathbb{H} \: \text{and every} \: \gamma = \left( \begin{smallmatrix} a & b \\ c & d \end{smallmatrix} \right) \in \Gamma
\end{equation}
is called the {\it transformation rule} for $\Gamma$.

Incidentally, since $\text{SL}_2(\mathbb{Z}) = \langle \left( \begin{smallmatrix} 1 & 1 \\ 0 & 1 \end{smallmatrix} \right), \: \left( \begin{smallmatrix} 0 & -1 \\ 1 & 0 \end{smallmatrix} \right) \rangle$, transformation rule for $\text{SL}_2(\mathbb{Z})$ is equivalent to the following two equations:
\begin{gather}
f(z + 1) = f(z), \label{modz+1}\\
f \left( - \frac{1}{z} \right) = z^k f(z). \label{mod-1/z}
\end{gather}

We have the Fourier expansion:
\begin{equation}
f(z) = \sum_{n \in \mathbb{Z}} a_n q_h^n, \quad \text{where} \: q_h = e^{2 \pi i z / h}.
\end{equation}
Similarly, we have the following Fourier expansion for every cusp $\kappa$ of $\Gamma$ with the cusp leader $\gamma_{\kappa}$:
\begin{equation}
f \big|_k \gamma_{\kappa} \ (z) = \sum_{n \in \mathbb{Z}} a_{\kappa, n} q_h^n, \quad \text{where} \: q_h = e^{2 \pi i z / h}. \label{eq-holo-cusp}
\end{equation}
When $h = 1$, we denote $q = q_1$. We say $f$ is {\it meromorphic at the cusp $\kappa$} if $a_{\kappa, n}$ is zero for $n$ small enough. Also, we call $f$ {\it holomorphic at the cusp $\kappa$} if $a_{\kappa, n}$ is zero for every negative integer $n$.

\begin{definition}
Let $f$ be a meromorphic function on $\mathbb{H}$. $f$ is called a {\it modular function} for $\Gamma$ if $f$ is meromorphic at every cusp and satisfies transformation rule for $\Gamma$.
\end{definition}

For a meromorphic function $f$, we assume $f(\kappa) = 0$ if and only if $a_{\kappa, n} = 0$ for every integer $n \leqslant 0$.

\begin{definition}
Let $f$ be a modular function for $\Gamma$ which is holomorphic on $\mathbb{H}$, then $f$ is called {\it weakly modular form} for $\Gamma$. In addition, if $f$ is holomorphic at evry cusp of $\Gamma$, then $f$ is called {\it modular form} for $\Gamma$. Furthermore, if $f$ is equal to $0$ at every cusp of $\Gamma$, we call $f$ {\it cusp form} for $\Gamma$.
\end{definition}

For a function $f$, let $v_p(f)$ be the order of $f$ at $p \in \mathbb{H}$. In addition, we also define the order of $f$ at a cusp $\kappa$:

\begin{equation*}
v_{\kappa}(f) := \min \{ n \in \mathbb{Z} \: ; \: a_{\kappa, n} \ne 0 \}.
\end{equation*}

Furthermore, we have following facts:
\begin{proposition}
Let $f$ be a modular form for $\Gamma$ such that every coefficient of Fourier expansion is real. Then we have $f(- \overline{z}) = \overline{f(z)}$ and $v_{\rho}(f) = v_{- \overline{\rho}}(f)$ at $\rho \in \mathbb{H}$. \label{prop-conju-z}
\end{proposition}\quad

Let $f$ be a modular function for $\Gamma$ of weight $k$. If we have $\gamma \Gamma \gamma^{-1} \subset \Gamma$ for $\gamma = \left( \begin{smallmatrix} a & b \\ c & d \end{smallmatrix} \right) \in \text{SL}_2(\mathbb{R})$, then we have
\begin{equation*}
f \big|_k \gamma \ (\gamma' \ \gamma^{-1} z) = ((c' d - c d') z + (a d' - b c'))^k f(z)
\end{equation*}
for every $\gamma' = \left( \begin{smallmatrix} a' & b' \\ c' & d' \end{smallmatrix} \right) \in \Gamma$. Then, we have
\begin{equation*}
f \big|_k \gamma \ (\gamma' z) = (c' z + d')^k (c z + d)^{-k} f(\gamma z) = (c' z + d')^k \ f \big|_k \gamma \ (z).
\end{equation*}
Thus, $f \big|_k \gamma$ is also a modular function for $\Gamma$ of weight $k$.\\

\subsubsection{Hauptmodul}\quad

Let $\Gamma$ be of genus $0$. Then, we consider the weakly modular form for $\Gamma$ which is holomorphic at every cusp but $\infty$ and has the following form of Fourier expansion:
\begin{equation}
j_{\Gamma}(z) := \frac{1}{q_h} + \sum_{n=0}^{\infty} a_n q_h^n.
\end{equation}
It is determined uniquely up to the constant term. We call $j_{\Gamma}$ {\it $($canonical$)$ hauptmodul} of $\Gamma$. Note that it is a isomoprhism from a fundamental domain $\mathcal{F}_{\Gamma}$ to the Riemann sphere $\mathbb{C} \cup \{ \infty \}$.

Similarly, we can define {\it hauptmodul for the cusp $\kappa$} $j_{\Gamma}^{\kappa}$ which has the following form of Fourier expansion for the cusp $\kappa$:
\begin{equation}
j_{\Gamma}^{\kappa} \big|_0 \gamma_{\kappa} \ (z) = \frac{1}{q_h} + \sum_{n=0}^{\infty} a_{n, \kappa} q_h^n.
\end{equation}
Note that we have $j_{\Gamma}^{\kappa} (z) = j_{\Gamma}(\gamma_{\kappa}^{-1} z)$ if $\gamma_{\kappa}^{-1} \Gamma \gamma_{\kappa} \subset \Gamma$. In this paper, we consider only hauptmodul for the cusp $\infty$.\\

\subsubsection{Eisenstein series}(see \cite{SG})

\begin{definition}
For $z \in \mathbb{H}$,
\begin{equation}
E_{k, \Gamma}^{\kappa} (z) := e \sum_{\gamma \in \Gamma_{\kappa} \setminus \Gamma} j \left( \gamma_{\kappa}^{-1} \gamma, z \right)^{-k} \quad (e \text{ : fixed number})
\end{equation}
is the {\it Eisenstein series} associated with $\Gamma$ for a cusp $\kappa$, where $j(\gamma, z) := c z + d$ for $\gamma = \left( \begin{smallmatrix} a & b \\ c & d \end{smallmatrix} \right) \in \Gamma$. $e$ is often selected so that the constant term of $E_{k, \Gamma}^{\kappa}$ is $1$.
\end{definition}

For example, let $\Gamma = \text{SL}_2(\mathbb{Z})$, then we have only $\infty$ as a cusp of $\text{SL}_2(\mathbb{Z})$. Now, for an even integer $k \geqslant 4$, we have
\begin{equation}
E_k(z) := \frac{1}{2} \sum_{(c,d)=1}(c z + d)^{-k} \label{def:e}
\end{equation}
as the Eisenstein series associated with $\text{SL}_2(\mathbb{Z})$.

For $k = 2$, we can define $E_2(z)$ as $E_k(z)$ for $k \geqslant 4$.
\begin{equation}
E_2(z) := \frac{1}{2} \sum_{(c,d)=1}(c z + d)^{-2}
 = 1 - 24 q - 72 q^2 - 96 q^3 - 168 q^4 - \cdots \label{def:e2}
\end{equation}
$E_2(z)$ is not modular form for $\text{SL}_2(\mathbb{Z})$. $E_2(z)$ is holomorphic on $\mathbb{H}$ and at $\infty$, and it satisfies transformation rule (\ref{modz+1}). On the other hand, it does not satisfy (\ref{mod-1/z}), instead, we have
\begin{equation}
E_2 \left(- \frac{1}{z}\right) = z^2 E_2(z) + \frac{12}{2 \pi i} z. \label{mod-e2}
\end{equation}\quad

If $\gamma_{\kappa}^{-1} \Gamma \gamma_{\kappa} = \Gamma$, then we have some relation between the Eisenstein series $E_{k, \Gamma}^{\infty}$ and $E_{k, \Gamma}^{\kappa}$. By the residue class $\Gamma_{\kappa} \setminus \Gamma$, we have
\begin{equation*}
\Gamma = \bigcup_{\gamma \in \Gamma_{\kappa} \setminus \Gamma} \Gamma_{\kappa} \cdot \gamma
 = \bigcup_{\gamma \in \Gamma_{\kappa} \setminus \Gamma} \left( \gamma_{\kappa} \Gamma_{\infty} \gamma_{\kappa}^{-1} \right) \cdot \gamma,\quad
\Gamma = \gamma_{\kappa}^{-1} \Gamma \gamma_{\kappa}
  = \bigcup_{\gamma \in \Gamma_{\kappa} \setminus \Gamma} \Gamma_{\infty} \cdot \left( \gamma_{\kappa}^{-1} \gamma \gamma_{\kappa}  \right).
\end{equation*}
Then, it gives us a residue class $\Gamma_{\infty} \setminus \Gamma$. Second, we have $j \left( \gamma_{\kappa}^{-1} \gamma, \gamma_{\kappa} z \right) = j \left( \gamma_{\kappa},  z \right)^{-1} j \left( \gamma_{\kappa}^{-1} \gamma \gamma_{\kappa}, z \right)$. Thus, we have
\begin{align*}
E_{k, \Gamma}^{\kappa} (\gamma_{\kappa} z)
 &= e \sum_{\gamma \in \Gamma_{\kappa} \setminus \Gamma} j \left( \gamma_{\kappa}^{-1} \gamma, \gamma_{\kappa} z \right)^{-k} \quad
 = j \left( \gamma_{\kappa},  z \right)^k \cdot e \sum_{\gamma \in \Gamma_{\kappa} \setminus \Gamma} j \left( \gamma_{\kappa}^{-1} \gamma \gamma_{\kappa}, z \right)^{-k}\\
 &= e' \; j \left( \gamma_{\kappa},  z \right)^k \cdot e'' \sum_{\gamma \in \Gamma_{\infty} \setminus \Gamma} j \left( \gamma, z \right)^{-k} \quad
  = e' \; j \left( \gamma_{\kappa},  z \right)^k E_{k, \Gamma}^{\infty} (z).
\end{align*}\quad

\subsubsection{Eta function}(see \cite[III.2]{Ko})

We put
\begin{equation}
\Delta(z) := \frac{1}{1728} \left( \left( E_4(z) \right)^3 - \left( E_6(z) \right)^2 \right).
\end{equation}
$\Delta$ is the cusp form for $\text{SL}_2(\mathbb{Z})$ of weight $12$ which satisfies $v_{\infty}(\Delta) = 1$. Now, we have

\begin{theorem}[Jacobi's product formula]
\begin{equation}
\Delta(z) = q \prod_{n = 1}^{\infty} (1 - q^n)^{24} \quad \text{where} \; q = e^{2 \pi i z}.
\end{equation}
\end{theorem}
Also, we have
\begin{equation}
\eta(z) = q^{\frac{1}{24}} \prod_{n = 1}^{\infty} (1 - q^n),
\end{equation}
which is called the {\it Dedekind $\eta$-function}. Then we have
\begin{equation}
\eta(z + 1) = e^{\frac{2 \pi i}{24}} \eta(z) \quad \text{and}  \quad \eta \left( - \frac{1}{z} \right) = \sqrt{\frac{z}{i}} \eta(z) \: \text{(see \cite{Ko}),} \label{eq-eta+1}
\end{equation}
where $\sqrt{\cdot}$ denote a square root which has nonnegative real part.
Furthermore, we have
\begin{equation}
\eta \left( \frac{a z + b}{c z + d} \right) = \epsilon \sqrt{\frac{c z + d}{i}} \eta(z) \quad \text{for} \: \left( \begin{smallmatrix} a & b \\ c & d \end{smallmatrix} \right) \in \text{SL}_2(\mathbb{Z}), \label{eq-etaabcd}
\end{equation}
where $\epsilon$ is one of the $24$th-roots of $1$ which depends on $a, b, c$, and $d$.\\

\subsubsection{Hecke type Faber polynomial}(see \cite{ACMS})

Let $\Gamma$ be a genus $0$ group. Then, we can determine the unique weakly modular form of weight $0$ for $\Gamma$ which is holomorphic at every cusp but $\infty$ and has the following form of Fourier expansion:
\begin{equation}
F_{m, \Gamma} (z) := \frac{1}{q_h^m} + \sum_{n=1}^{\infty} a_n q_h^n.
\end{equation}
Since this function is of weight $0$, we can write $F_{m, \Gamma}$ as a polynomial of the hauptmodul for $\Gamma$ as follows:
\begin{equation}
F_{m, \Gamma} (z) = P_{m, \Gamma}(j_{\Gamma} (z)),
\end{equation}
where we call $P_{m, \Gamma}(X)$ the {\it Hecke type Faber polynomial of degree $m$ for $\Gamma$}.

$F_{m, \Gamma}$ and $P_{m, \Gamma}(X)$ are also defined with the {\it twisted Hecke operator}. (see \cite{ACMS})\\

For the other cusps, we can also define the functions $F_{m, \Gamma}^{\kappa}$ and the polynomials $P_{m, \Gamma}^{\kappa}(X)$, where
\begin{align}
F_{m, \Gamma}^{\kappa} \big|_0 \gamma_{\kappa} \ (z) &:= \frac{1}{q_h^m} + \sum_{n=1}^{\infty} a_{n, \kappa} q_h^n,\\
F_{m, \Gamma}^{\kappa} (z) &= P_{m, \Gamma}^{\kappa} (j_{\Gamma}^{\kappa} (z)).
\end{align}
If $\gamma_{\kappa}^{-1} \Gamma \gamma_{\kappa} \subset \Gamma$, then we have $P_{m, \Gamma}^{\kappa} (X) = P_{m, \Gamma} (X)$ since $j_{\Gamma}^{\kappa} (z) = j_{\Gamma}(\gamma_{\kappa}^{-1} z)$.\\

\subsubsection{Semimodular form}\quad \label{subsubsec-semimod}

Let $N$ be a positive integer. If the function $f$ satisfies that $f^n$ is not a modular form (resp. function) for every integer $n < N$ and $f^N$ is a modular form (resp. function), then we would like to call $f$ the {\it $N$th semimodular form $($resp. function$)$ for $\Gamma$}.

For example, $\eta$-function is $24$th semimodular form for $\text{SL}_2(\mathbb{Z})$.

Let $f$ and $g$ be meromorphic functions on $\mathbb{H}$ and at every cusp of $\Gamma$ which satisfy transformation rule of weight $k_1$ and $k_2$ for $\Gamma$ except for $\gamma = \left( \begin{smallmatrix}a & b \\ c & d\end{smallmatrix} \right) \in \Gamma$ such that $f(\gamma z) = e^{2 \pi i / N} (c z + d)^{k_2} f(z)$ and $g(\gamma z) = e^{2 \pi i / N} (c z + d)^{k_2} g(z)$, respectively. Then, $f$ and $g$ are $N$th semimodular functions for $\Gamma$. Furthermore, $f^i g^{N - i}$ and $f / g$ are modular functions for $\Gamma$.\\

\subsection{Conjugate of the groups and the space of modular forms}\quad

We denote the space of modular forms of weight $k$ for $\Gamma$ by $M_k(\Gamma)$, and that of cusp forms by $S_k(\Gamma)$.

Some discrete subgroups of $\text{SL}_2(\mathbb{R})$ are conjugate to each other. Moreover, some spaces of modular forms for such groups are isomorphic to each other. In this section, we introduce two of such cases.

Let $\Gamma$ and $\Gamma'$ be discrete subgroups of $\text{SL}_2(\mathbb{R})$ such that they are conjugate to each other.

The first case is when $\Gamma' = V_h^{-1} \Gamma V_h$, where $V_h := \left( \begin{smallmatrix}\sqrt{h} & 0 \\ 0 & 1 / \sqrt{h}\end{smallmatrix} \right)$. Then, it is easy to show that the map
\begin{equation*}
M_k(\Gamma) \ni f(z) \mapsto f(h z) \in M_k(\Gamma')
\end{equation*}
is a isomorphism. Furthermore, we have $E_{k, \Gamma}^{\infty}(h z) = E_{k, \Gamma'}^{\infty}(z)$. For example, we have $\Gamma_0(3 | 3) = V_3^{-1} \text{SL}_2(\mathbb{Z}) V_3$.

The second case is when $\Gamma' = T_x^{-1} \Gamma T_x$. Then, we can easily show that the map
\begin{equation*}
M_k(\Gamma) \ni f(z) \mapsto f(z + x) \in M_k(\Gamma')
\end{equation*}
is a isomorphism. Furthermore, we have $E_{k, \Gamma}^{\infty}(z + x) = E_{k, \Gamma'}^{\infty}(z)$. For example, we have $\Gamma_0^{*}(4) = T_{1/2}^{-1} \Gamma_0(2) T_{1/2}$.

For both above cases, if all of the zeros of $E_{k, \Gamma}^{\infty}$ lies on the lower acrs of $\partial \mathbb{F}_{0,\Gamma}$, then all of the zeros of $E_{k, \Gamma'}^{\infty}$ lies on the lower acrs of $\partial \mathbb{F}_{0,\Gamma'}$.

\begin{figure}[hbtp]
\begin{center}
{{$\text{SL}_2(\mathbb{Z})$}\includegraphics[width=1in]{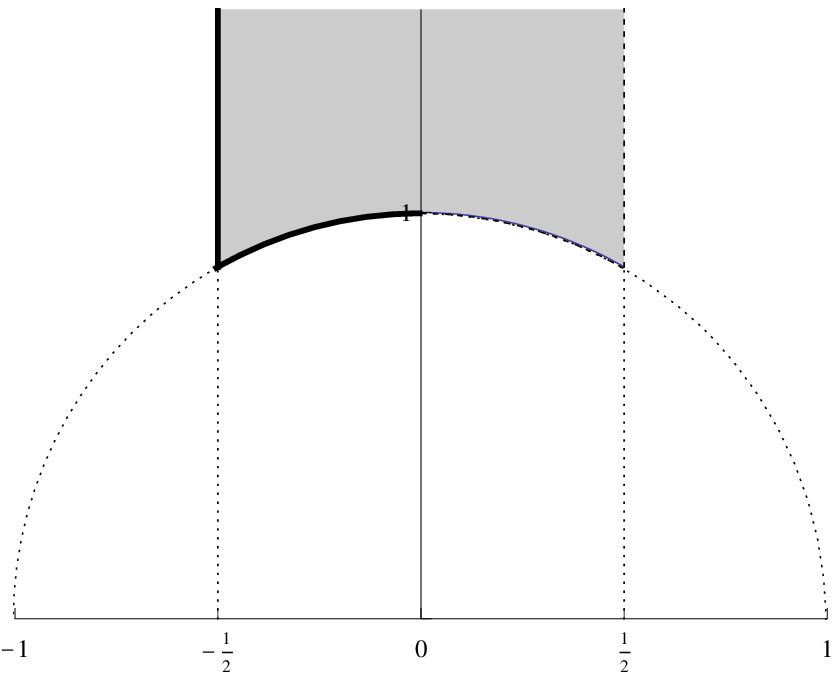}}
\qquad$\sim$\qquad
{{$\Gamma_0(3 | 3)$}\includegraphics[width=1in]{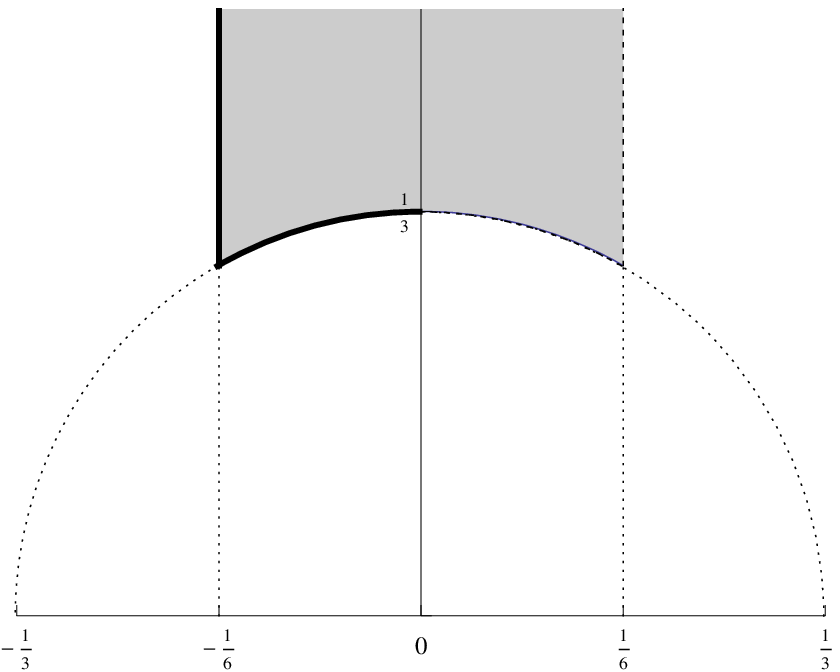}}\\
\quad\\
{{$\Gamma_0(2)$}\includegraphics[width=1in]{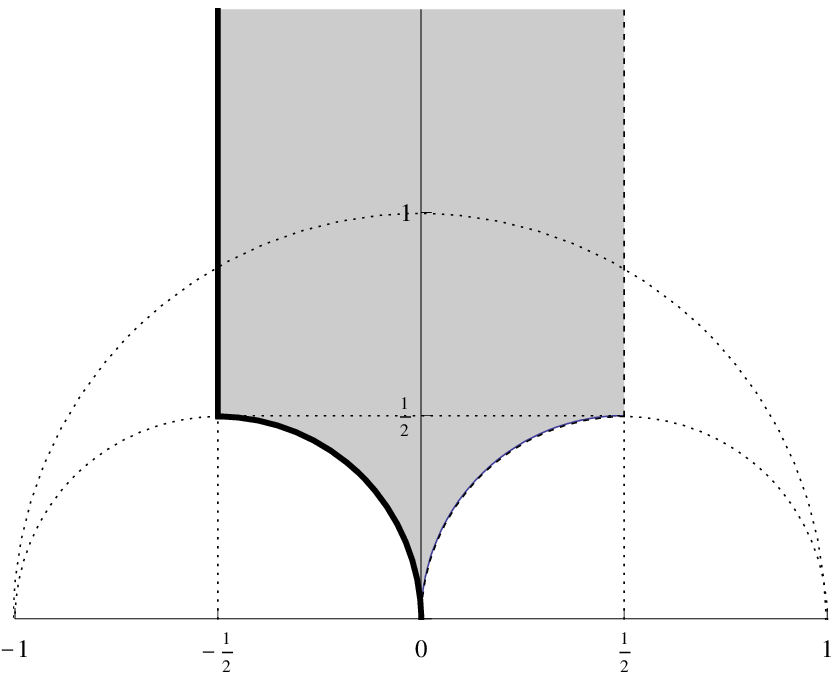}}
\qquad$\sim$\qquad
{{$\Gamma_0^{*}(4)$}\includegraphics[width=1in]{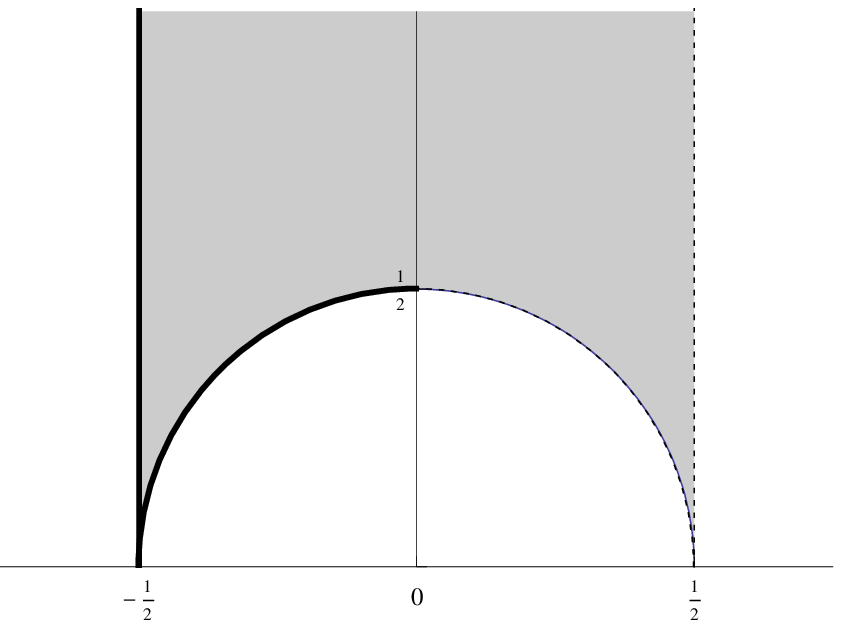}}
\end{center}
\caption{Conjugate groups}
\label{conjugate}
\end{figure}

\begin{remark}
We have
\begin{equation}
\Gamma_0(n | h) + e_1, e_2, ... , e_m = V_h^{-1} (\Gamma_0(n / h) + e_1, e_2, ... , e_m) V_h.
\end{equation}
\end{remark}

Thus, in this paper, we observe only the normalizers such that $h = 1$.

%% file: report01-05.tex
\section{Level $1$}

\subsection{$\text{SL}_2(\mathbb{Z})$}(see \cite{Se} and \cite{RSD})

\paragraph{\bf Fundamental domain}
We have a fundamental domain for $\text{\upshape SL}_2(\mathbb{Z})$ as follows:
\begin{equation}
\mathbb{F} = \left\{|z| \geqslant 1, \: - 1/2 \leqslant Re(z) \leqslant 0\right\}
 \bigcup \left\{|z| > 1, \: 0 < Re(z) < 1/2 \right\},
\end{equation}
where $\left( \begin{smallmatrix} 0 & -1 \\ 1 & 0 \end{smallmatrix} \right) : e^{i \theta} \rightarrow e^{i (\pi - \theta)}$. By corollary\ref{cor-base}, we have $\text{SL}_2(\mathbb{Z}) = \langle \left( \begin{smallmatrix} 1 & 1 \\ 0 & 1 \end{smallmatrix} \right), \: \left( \begin{smallmatrix} 0 & -1 \\ 1 & 0 \end{smallmatrix} \right) \rangle$.
\begin{figure}[hbtp]
\begin{center}
{{$\mathbb{F}$}\includegraphics[width=1.5in]{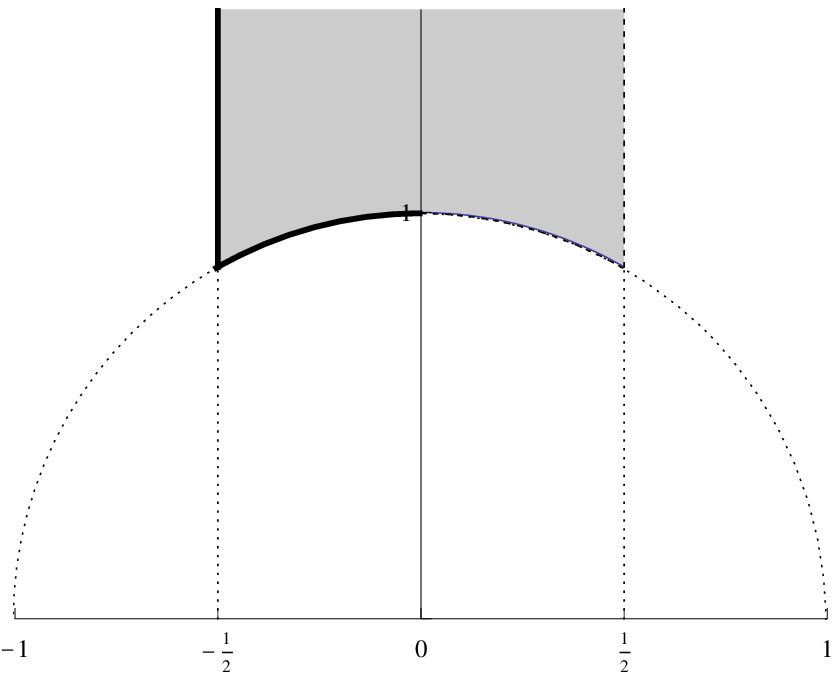}}
\end{center}
\caption{$\text{SL}_2(\mathbb{Z})$}
\end{figure}

\paragraph{\bf Valence formula}
The cusp of $\text{SL}_2(\mathbb{Z})$ is only $\infty$, and the elliptic points are $i$ and $\rho := e^{2 \pi / 3}$. Let $f$ be a modular function of weight $k$ for $\text{SL}_2(\mathbb{Z})$, which is not identically zero. We have
\begin{equation}
v_{\infty}(f) + \frac{1}{2} v_{i}(f) + \frac{1}{3} v_{\rho} (f) + \sum_{\begin{subarray}{c} p \in \text{SL}_2(\mathbb{Z}) \setminus \mathbb{H} \\ p \ne i, \; \rho\end{subarray}} v_p(f) = \frac{k}{12}.
\end{equation}

Furthermore, the stabilizer of the elliptic point $i$ is $\left\{ \pm I, \pm \left( \begin{smallmatrix}0 & -1 \\ 1 & 0\end{smallmatrix} \right) \right\}$, and that of $\rho$ is\\
$\left\{ \pm I, \pm \left( \begin{smallmatrix}-1 & -1 \\ 1 & 0\end{smallmatrix} \right), \pm \left( \begin{smallmatrix}0 & -1 \\ 1 & 1\end{smallmatrix} \right) \right\}$.\\

\paragraph{\bf For the cusp $\infty$}
We have $\Gamma_{\infty} = \left\{ \pm \left( \begin{smallmatrix}1 & n \\ 0 & 1\end{smallmatrix} \right) \: ; \: n \in \mathbb{Z} \right\}$, and we have the Eisenstein series associated with $\text{SL}_2(\mathbb{Z})$:
\begin{equation*}
E_k(z) = \frac{1}{2} \sum_{(c,d)=1}(c z + d)^{-k} \quad \text{for} \; k \geqslant 4.
\end{equation*}
We also have its Fourier Expansion:
\begin{equation}
E_k(z) = 1 - \frac{2 k}{B_k} \sum_{n = 1}^{\infty} \sigma_{k-1}(n) q^n,
\end{equation}
where $q := e^{2 \pi i z}$, $\sigma_k(n) := \sum_{d \mid n} d^k$ which is called {\it divisor function}, and $B_k$ are {\it Bernoulli number}.\\

\paragraph{\bf The space of modular forms}
We have $M_k(\text{SL}_2(\mathbb{Z})) = \mathbb{C} E_k \oplus S_k(\text{SL}_2(\mathbb{Z}))$ and $S_k(\text{SL}_2(\mathbb{Z})) = \Delta M_{k - 12}(\text{SL}_2(\mathbb{Z}))$ for every even integer $k \geqslant 4$. Then, we have
\begin{equation*}
M_k(\text{SL}_2(\mathbb{Z})) = E_{k - 12 n} (\mathbb{C} ((E_4)^3)^n \oplus \mathbb{C} ((E_4)^3)^{n-1} \Delta \oplus \cdots \oplus \mathbb{C} (\Delta)^n),
\end{equation*}
where $n = \dim(M_k(\text{SL}_2(\mathbb{Z}))) - 1 = \lfloor k/12 - (k/4 - \lfloor k/4 \rfloor)\rfloor$. Furethermore, on the zeros of $E_{k - 12 n}$, we have $v_{\rho}(E_4) = 1$, $v_{i}(E_6) = 1$, $E_8 = (E_4)^2$, $E_{10} = E_4 E_6$, and $E_{14} = (E_4)^2 E_6$.\\

\paragraph{\bf Hauptmodul}
We define the {\it hauptmodul} of $\text{SL}_2(\mathbb{Z})$:
\begin{equation}
J := (E_4)^3 / \Delta = \frac{1}{q} + 744 + 196884 q + 21493760 q^2 + 864299970 q^3 + \cdots,
\end{equation}
where $v_{\infty}(J) = -1$ and $v_{\rho}(J) = 3$. Then, we have
\begin{equation}
J : \partial \mathbb{F} \setminus \{z \in \mathbb{H} \: ; \: Re(z) = \pm 1/2\} \to [0, 1728] \subset \mathbb{R}.
\end{equation}

\clearpage

\section{Level $2$}

We have $\Gamma_0(2)+=\Gamma_0^{*}(2)$ and $\Gamma_0(2)-=\Gamma_0(2)$.

We have $W_2 = \left(\begin{smallmatrix}0&-1 / \sqrt{2}\\ \sqrt{2}&0\end{smallmatrix}\right)$, and denote $\rho_2 := - 1/2 + i / 2$. We define
\begin{equation}
\begin{split}
&\Delta_2^{\infty}(z) := \eta^{16}(2 z) / \eta^8(z), \quad \Delta_2^0(z) := \eta^{16}(z) / \eta^8(2 z),\\
&\Delta_2(z) := \Delta_2^{\infty}(z) \Delta_2^0(z) = \eta^8(z) \eta^8(2 z),
\end{split}
\end{equation}
where $\Delta_2^{\infty}$ and $\Delta_2^0$ are modular forms for $\Gamma_0(2)$ of weight $4$ such that $v_{\infty}(\Delta_2^{\infty}) = v_0(\Delta_2^0) = 1$, and $\Delta_2$ is a cusp form for $\Gamma_0(2)$ and $\Gamma_0^{*}(2)$ of weight $8$. Furthermore, we define
\begin{equation}
{E_{2,2}}'(z) := 2 E_2(2 z) - E_2(z),
\end{equation}
where $E_2$ is the Eisenstein series for $\text{SL}_2(\mathbb{Z})$, and ${E_{2,2}}'$ is not a Eisenstein series but a modular form for $\Gamma_0(2)$ of weight $2$ such that $v_{\rho_2}({E_{2,2}}') = 1$.\\

\subsection{$\Gamma_0^{*}(2)$} (see \cite{MNS})

\paragraph{\bf Fundamental domain}
We have a fundamental domain for $\Gamma_0^{*}(2)$ as follows:
\begin{equation}
\mathbb{F}_{2+} = \left\{|z| \geqslant 1 / \sqrt{2}, \: - 1/2 \leqslant Re(z) \leqslant 0\right\}
 \bigcup \left\{|z| > 1 / \sqrt{2}, \: 0 < Re(z) < 1/2 \right\},
\end{equation}
where $W_2 : e^{i \theta} / \sqrt{2} \rightarrow e^{i (\pi - \theta)} / \sqrt{2}$. Then, we have
\begin{equation}
\Gamma_0^{*}(2) = \langle \left( \begin{smallmatrix} 1 & 1 \\ 0 & 1 \end{smallmatrix} \right), \: W_2 \rangle.
\end{equation}
\begin{figure}[hbtp]
\begin{center}
{{$\mathbb{F}_{2+}$}\includegraphics[width=1.5in]{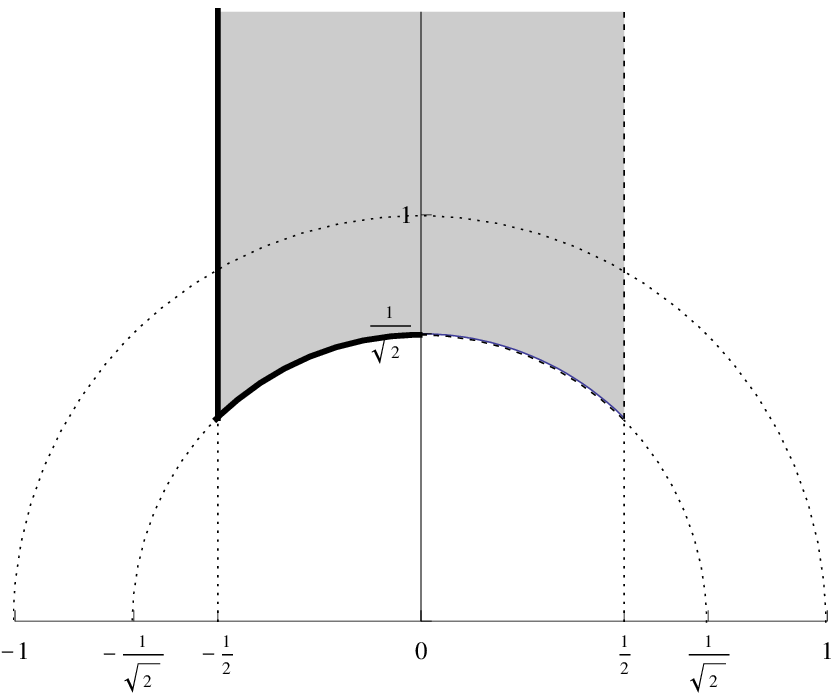}}
\end{center}
\caption{$\Gamma_0^{*}(2)$}
\end{figure}

\paragraph{\bf Valence formula}
The cusp of $\Gamma_0^{*}(2)$ is $\infty$, and the elliptic points are $i / \sqrt{2}$ and $\rho_2 = - 1/2 + i / 2$. Let $f$ be a modular function of weight $k$ for $\Gamma_0^{*}(2)$, which is not identically zero. We have
\begin{equation}
v_{\infty}(f) + \frac{1}{2} v_{i / \sqrt{2}}(f) + \frac{1}{4} v_{\rho_2} (f) + \sum_{\begin{subarray}{c} p \in \Gamma_0^{*}(2) \setminus \mathbb{H} \\ p \ne i / \sqrt{2}, \; \rho_2\end{subarray}} v_p(f) = \frac{k}{8}.
\end{equation}

Furthermore, the stabilizer of the elliptic point $i / \sqrt{2}$ is $\left\{ \pm I, \pm W_2 \right\}$, and that of $\rho_2$ is\\
$\left\{ \pm I, \pm \left( \begin{smallmatrix}-1 & -1 \\ 2 & 1\end{smallmatrix} \right), \pm \left( \begin{smallmatrix}1 & -1 \\ 0 & 1\end{smallmatrix} \right) W_2, \pm \left( \begin{smallmatrix}1 & 0 \\ -2 & 1\end{smallmatrix} \right) W_2 \right\}$.\\

\paragraph{\bf For the cusp $\infty$}
We have $\Gamma_{\infty} = \left\{ \pm \left( \begin{smallmatrix}1 & n \\ 0 & 1\end{smallmatrix} \right) \: ; \: n \in \mathbb{Z} \right\}$, and we have the Eisenstein series associated with $\Gamma_0^{*}(2)$:
\begin{equation}
E_{k, 2+}(z) := \frac{2^{k/2} E_k(2 z) + E_k(z)}{2^{k/2} + 1} \quad \text{for} \; k \geqslant 4.
\end{equation}\quad

\paragraph{\bf The space of modular forms}
We have $M_k(\Gamma_0^{*}(2)) = \mathbb{C} E_{k, 2+} \oplus S_k(\Gamma_0^{*}(2))$ and $S_k(\Gamma_0^{*}(2)) = \Delta_2 M_{k - 8}(\Gamma_0^{*}(2))$ for every even integer $k \geqslant 4$. Then, we have
\begin{equation*}
M_k(\Gamma_0^{*}(2)) = E_{k - 8 n, 2+} (\mathbb{C} ((E_{4, 2+})^2)^n \oplus \mathbb{C} ((E_{4, 2+})^2)^{n-1} \Delta_2 \oplus \cdots \oplus \mathbb{C} (\Delta_2)^n),
\end{equation*}
where $n = \dim(M_k(\Gamma_0^{*}(2))) - 1 = \lfloor k/8 - (k/4 - \lfloor k/4 \rfloor)\rfloor$. Furethermore, on the zeros of $E_{k - 8 n, 2+}$, we have $v_{\rho_2}(E_{4, 2+}) = 2$, $v_{i / \sqrt{2}}(E_{6, 2+}) = v_{\rho_2}(E_{6, 2+}) = 1$, and $E_{10, 2+}= E_{4, 2+} E_{6, 2+}$.\\

\paragraph{\bf Hauptmodul}
We define the {\it hauptmodul} of $\Gamma_0^{*}(2)$:
\begin{equation}
J_{2+} := (E_{4, 2+})^2 / \Delta_2 = \frac{1}{q} + 104 + 4372 q + 96256 q^2 + 1240002 q^3 + \cdots,
\end{equation}
where $v_{\infty}(J_{2+}) = -1$ and $v_{\rho_2}(J_{2+}) = 4$. Then, we have
\begin{equation}
J_{2+} : \partial \mathbb{F}_{2+} \setminus \{z \in \mathbb{H} \: ; \: Re(z) = \pm 1/2\} \to [0, 256] \subset \mathbb{R}.
\end{equation}\quad

\subsection{$\Gamma_0(2)$} (see \cite{SJ1})

\paragraph{\bf Fundamental domain}
We have a fundamental domain for $\Gamma_0(2)$ as follows:
\begin{equation}
\mathbb{F}_2 = \left\{|z + 1/2| \geqslant 1/2, \: - 1/2 \leqslant Re(z) \leqslant 0\right\}
 \bigcup \left\{|z - 1/2| > 1/2, \: 0 < Re(z) < 1/2 \right\},
\end{equation}
where $\left( \begin{smallmatrix} -1 & 0 \\ 2 & -1 \end{smallmatrix} \right) : (e^{i \theta} + 1) / 2 \rightarrow (e^{i (\pi - \theta)} - 1) / 2$. Then we have
\begin{equation}
\Gamma_0(2) = \langle \left( \begin{smallmatrix} 1 & 1 \\ 0 & 1 \end{smallmatrix} \right), \: \left( \begin{smallmatrix} 1 & 0 \\ 2 & 1 \end{smallmatrix} \right) \rangle.
\end{equation}
\begin{figure}[hbtp]
\begin{center}
\includegraphics[width=1.5in]{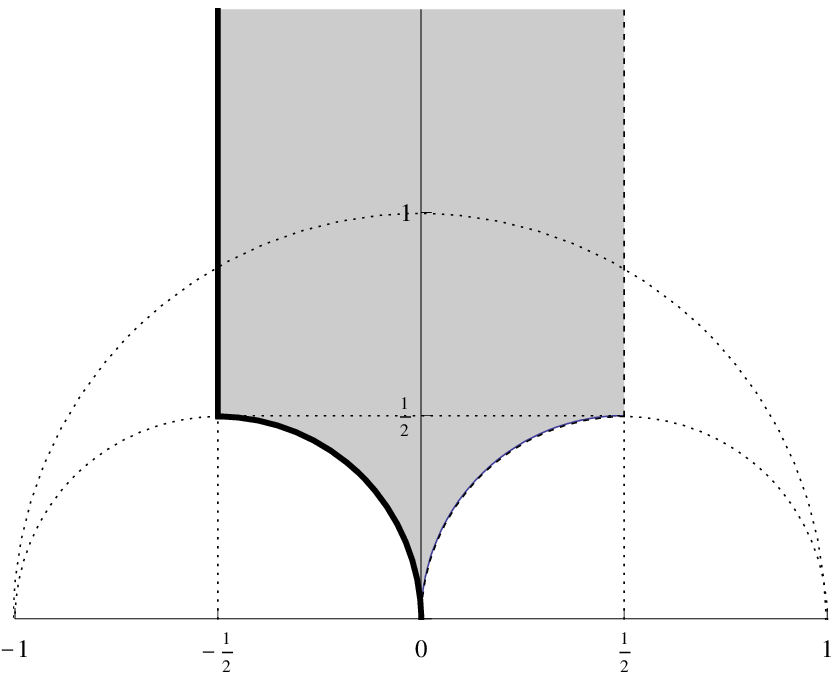}
\end{center}
\caption{$\Gamma_0(2)$}
\end{figure}

\paragraph{\bf Valence formula}
The cusps of $\Gamma_0(2)$ are $\infty$ and $0$, and the elliptic point is $\rho_2$. Let $f$ be a modular function of weight $k$ for $\Gamma_0(2)$, which is not identically zero. We have
\begin{equation}
v_{\infty}(f) + v_0(f) + \frac{1}{2} v_{\rho_2} (f) + \sum_{\begin{subarray}{c} p \in \Gamma_0(2) \setminus \mathbb{H} \\ p \ne \rho_2\end{subarray}} v_p(f) = \frac{k}{4}.
\end{equation}

Furthermore, the stabilizer of the elliptic point $\rho_2$ is $\left\{ \pm I, \pm \left( \begin{smallmatrix}-1 & -1 \\ 2 & 1\end{smallmatrix} \right) \right\}$.\\

\paragraph{\bf For the cusp $\infty$}
We have $\Gamma_{\infty} = \left\{ \pm \left( \begin{smallmatrix}1 & n \\ 0 & 1\end{smallmatrix} \right) \: ; \: n \in \mathbb{Z} \right\}$, and we have the Eisenstein series for the cusp $\infty$ associated with $\Gamma_0(2)$:
\begin{equation}
E_{k, 2}^{\infty}(z) := \frac{2^k E_k(2 z) - E_k(z)}{2^k - 1} \quad \text{for} \; k \geqslant 4.
\end{equation}\quad

\paragraph{\bf For the cusp $0$}
We have $\Gamma_0 = \left\{ \pm \left( \begin{smallmatrix}1 & 0 \\ 2 n & 1\end{smallmatrix} \right) \: ; \: n \in \mathbb{Z} \right\}$ and $\gamma_0 = W_2$, and we have the Eisenstein series for the cusp $0$ associated with $\Gamma_0(2)$:
\begin{equation}
E_{k, 2}^0(z) := \frac{- 2^{k/2} (E_k(2 z) - E_k(z))}{2^k - 1} \quad \text{for} \; k \geqslant 4.
\end{equation}
We also have $\gamma_0^{-1} \ \Gamma_0(2) \ \gamma_0 = \Gamma_0(2)$.\\

\paragraph{\bf The space of modular forms}
We have $M_k(\Gamma_0(2)) = \mathbb{C} E_{k, 2}^{\infty} \oplus \mathbb{C} E_{k, 2}^0 \oplus S_k(\Gamma_0(2))$ and $S_k(\Gamma_0(2)) = \Delta_2 M_{k - 8}(\Gamma_0(2))$ for every even integer $k \geqslant 4$. Then, we have $M_{4 n + 2}(\Gamma_0(2)) = {E_{2,2}}' M_{4 n}(\Gamma_0(2))$ and
\begin{allowdisplaybreaks}
\begin{align*}
M_{8 n}(\Gamma_0(2)) &= \mathbb{C} ((E_{4, 2}^{\infty})^2)^n \oplus \mathbb{C} ((E_{4, 2}^{\infty})^2)^{n-1} \Delta_2 \oplus \cdots \oplus \mathbb{C} (E_{4, 2}^{\infty})^2 (\Delta_2)^{n-1}\\
 &\oplus \mathbb{C} ((E_{4, 2}^0)^2)^n \oplus \mathbb{C} ((E_{4, 2}^0)^2)^{n-1} \Delta_2 \oplus \cdots \oplus \mathbb{C} (E_{4, 2}^0)^2 (\Delta_2)^{n-1} \oplus \mathbb{C} (\Delta_2)^n,\\
M_{8 n + 4}(\Gamma_0(2)) &= E_{4, 2}^{\infty} (\mathbb{C} ((E_{4, 2}^{\infty})^2)^n \oplus \mathbb{C} ((E_{4, 2}^{\infty})^2)^{n-1} \Delta_2 \oplus \cdots \oplus \mathbb{C} (\Delta_2)^n)\\
 &\oplus E_{4, 2}^0 (\mathbb{C} ((E_{4, 2}^0)^2)^n \oplus \mathbb{C} ((E_{4, 2}^0)^2)^{n-1} \Delta_2 \oplus \cdots \oplus \mathbb{C} (\Delta_2)^n).
\end{align*}\end{allowdisplaybreaks}

Here, since we have $E_{4, 2}^{\infty} = \Delta_2^0$ and $E_{4, 2}^0 = \Delta_2^{\infty}$, we can write
\begin{equation*}
M_{4 n}(\Gamma_0(2)) = \mathbb{C} (\Delta_2^{\infty})^n \oplus \mathbb{C} (\Delta_2^{\infty})^{n-1} \Delta_2^0 \oplus \cdots \oplus \mathbb{C} (\Delta_2^0)^n.
\end{equation*}\quad

\paragraph{\bf Hauptmodul}
We define the {\it hauptmodul} of $\Gamma_0(2)$:
\begin{equation}
J_2 := \Delta_2^0 / \Delta_2^{\infty} \: (= \eta^{24}(z) / \eta^{24}(2 z)) = \frac{1}{q} - 24 + 276 q - 2048 q^2 + 11202 q^3 + \cdots,
\end{equation}
where $v_{\infty}(J_2) = -1$ and $v_0(J_2) = 1$. Then, we have
\begin{equation}
J_2 : \partial \mathbb{F}_2 \setminus \{z \in \mathbb{H} \: ; \: Re(z) = \pm 1/2\} \to [-64, 0] \subset \mathbb{R}.
\end{equation}

\clearpage

\section{Level $3$}

We have $\Gamma_0(3)+=\Gamma_0^{*}(3)$ and $\Gamma_0(3)-=\Gamma_0(3)$.

We have $W_3 = \left(\begin{smallmatrix}0&-1 / \sqrt{3}\\ \sqrt{3}&0\end{smallmatrix}\right)$, and denote $\rho_3 := - 1/2 + i / (2 \sqrt{3})$. We define
\begin{equation}
\begin{split}
&\Delta_3^{\infty}(z) := \eta^9(3 z) / \eta^3(z), \quad \Delta_3^0(z) := \eta^9(z) / \eta^3(3 z),\\
&\Delta_3(z) := \Delta_3^{\infty}(z) \Delta_3^0(z) = \eta^6(z) \eta^6(3 z),
\end{split}
\end{equation}
where $\Delta_3^{\infty}$ and $\Delta_3^0$ are $2$nd semimodular forms (cf. Section \ref{subsubsec-semimod}) for $\Gamma_0(3)$ of weight $3$ such that $v_{\infty}(\Delta_3^{\infty}) = v_0(\Delta_3^0) = 1$, and $\Delta_3$ is a cusp form for $\Gamma_0(3)$ and a $2$nd semimodular form for $\Gamma_0^{*}(3)$ of weight $6$. Furthermore, we define
\begin{equation}
{E_{2, 3}}'(z) := (3 E_2(3 z) - E_2(z)) / 2,
\end{equation}
which is a modular form for $\Gamma_0(3)$ of weight $2$ such that $v_{\rho_3}({E_{2, 3}}') = 2$.\\

\subsection{$\Gamma_0^{*}(3)$} (see \cite{MNS})

\paragraph{\bf Fundamental domain}
We have a fundamental domain for $\Gamma_0^{*}(3)$ as follows:
\begin{equation}
\mathbb{F}_{3+} = \left\{|z| \geqslant 1 / \sqrt{3}, \: - 1/2 \leqslant Re(z) \leqslant 0\right\}
 \bigcup \left\{|z| > 1 / \sqrt{3}, \: 0 < Re(z) < 1/2 \right\},
\end{equation}
where $W_3 : e^{i \theta} / \sqrt{3} \rightarrow e^{i (\pi - \theta)} / \sqrt{3}$. Then, we have
\begin{equation}
\Gamma_0^{*}(3) = \langle \left( \begin{smallmatrix} 1 & 1 \\ 0 & 1 \end{smallmatrix} \right), \: W_3 \rangle.
\end{equation}
\begin{figure}[hbtp]
\begin{center}
{{$\mathbb{F}_{3+}$}\includegraphics[width=1.5in]{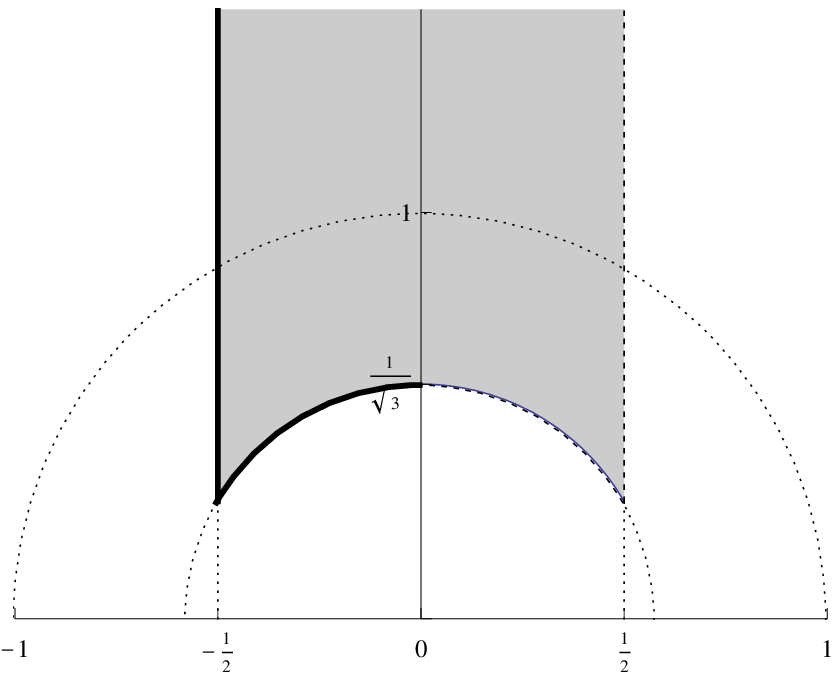}}
\end{center}
\caption{$\Gamma_0^{*}(3)$}
\end{figure}

\paragraph{\bf Valence formula}
The cusp of $\Gamma_0^{*}(3)$ is $\infty$, and the elliptic points are $i / \sqrt{3}$ and $\rho_3 = - 1/2 + i / (2 \sqrt{3})$. Let $f$ be a modular function of weight $k$ for $\Gamma_0^{*}(3)$, which is not identically zero. We have
\begin{equation}
v_{\infty}(f) + \frac{1}{2} v_{i / \sqrt{3}}(f) + \frac{1}{6} v_{\rho_3} (f) + \sum_{\begin{subarray}{c} p \in \Gamma_0^{*}(3) \setminus \mathbb{H} \\ p \ne i / \sqrt{3}, \; \rho_3\end{subarray}} v_p(f) = \frac{k}{6}.
\end{equation}

Furthermore, the stabilizer of the elliptic point $i / \sqrt{3}$ is $\left\{ \pm I, \pm W_3 \right\}$, and that of $\rho_3$ is\\
$\left\{ \pm I, \pm \left( \begin{smallmatrix}-2 & -1 \\ 3 & 1\end{smallmatrix} \right), \pm \left( \begin{smallmatrix}-1 & -1 \\ 3 & 2\end{smallmatrix} \right), \pm \left( \begin{smallmatrix}1 & -1 \\ 0 & 1\end{smallmatrix} \right) W_3, \pm \left( \begin{smallmatrix}1 & 0 \\ -3 & 1\end{smallmatrix} \right) W_3, \pm \left( \begin{smallmatrix}-2 & 1 \\ 3 & -2\end{smallmatrix} \right) W_3 \right\}$.\\

\paragraph{\bf For the cusp $\infty$}
We have $\Gamma_{\infty} = \left\{ \pm \left( \begin{smallmatrix}1 & n \\ 0 & 1\end{smallmatrix} \right) \: ; \: n \in \mathbb{Z} \right\}$, and we have the Eisenstein series associated with $\Gamma_0^{*}(3)$:
\begin{equation}
E_{k, 3+}(z) := \frac{3^{k/2} E_k(3 z) + E_k(z)}{3^{k/2} + 1} \quad \text{for} \; k \geqslant 4.
\end{equation}\quad

\paragraph{\bf The space of modular forms}
We define the following functions:
\begin{align*}
{\Delta_{8, 3}}' &:= (41/1728) ((E_{4, 3+})^2 - E_{8, 3+}), & {\Delta_{10, 3}}' &:= (61/432) (E_{4, 3+} E_{6, 3+} - E_{10, 3+}),\\
{\Delta_{12, 3}}' &:= E_{4, 3+} {\Delta_{8, 3}}', & {\Delta_{14, 3}}' &:= E_{4, 3+} {\Delta_{10, 3}}'.
\end{align*}

Now, we have $M_k(\Gamma_0^{*}(3)) = \mathbb{C} E_{k, 3+} \oplus S_k(\Gamma_0^{*}(3))$ and $S_k(\Gamma_0^{*}(3)) = (\mathbb{C} {\Delta_{12, 3}}' \oplus \mathbb{C} (\Delta_3)^2) M_{k - 12}(\Gamma_0^{*}(3))$ for every even integer $k \geqslant 4$. Then, we have $M_{12 n + l} = E_{l, 3+} M_{12 n}$ for $l = 4, 6$, $M_{12 n + l} = E_{l, 3+} M_{12 n} \oplus \mathbb{C} {\Delta_{l, 3}}' (\Delta_3)^{2 n}$ for $l = 8, 10, 14$, and
\begin{equation*}
M_{12 n}(\Gamma_0^{*}(3)) = \mathbb{C} (E_{12, 3+})^n \oplus \mathbb{C} (E_{12, 3+})^{n-1} {\Delta_{12, 3}}' \oplus \mathbb{C} (E_{12, 3+})^{n-1} (\Delta_3)^2 \oplus \cdots \oplus \mathbb{C} (\Delta_3)^{2 n}.
\end{equation*}

Here, we define ${E_{4, 3}}' := E_{6, 3+} / {E_{2, 3}}'$, which is a $2$nd semimodular form such that $v_{i / \sqrt{3}}({E_{4, 3}}') = v_{\rho_3}({E_{4, 3}}') = 1$. Then, we can write $E_{4, 3+} = ({E_{2, 3}}')^2$, $E_{6, 3+} = {E_{2, 3}}' {E_{4, 3}}'$, ${\Delta_{8, 3}}' = {E_{2, 3}}' \Delta_3$, ${\Delta_{10, 3}}' = {E_{4, 3}}' \Delta_3$, ${\Delta_{12, 3}}' = ({E_{2, 3}}')^3 \Delta_3$, and ${\Delta_{14, 3}}' = ({E_{2, 3}}')^2 {E_{4, 3}}' \Delta_3$.

Now, we have
\begin{equation*}
M_k(\Gamma_0^{*}(3)) = {E_{\overline{k}, 3}}' (\mathbb{C} (({E_{2, 3}}')^3)^n \oplus \mathbb{C} (({E_{2, 3}}')^3)^{n-1} \Delta_3 \oplus \cdots \oplus \mathbb{C} (\Delta_3)^n),
\end{equation*}
where $n = \dim(M_k(\Gamma_0^{*}(3))) - 1 = \lfloor k/6 - (k/4 - \lfloor k/4 \rfloor)\rfloor$, and where ${E_{\overline{k}, 3}}' := 1$, $({E_{2, 3}}')^2 {E_{4, 3}}'$, $({E_{2, 3}}')^2$, ${E_{2, 3}}' {E_{4, 3}}'$, ${E_{2, 3}}'$, and ${E_{4, 3}}'$, when $k \equiv 0$, $2$, $4$, $6$, $8$, and $10 \pmod{12}$, respectively.\\

\paragraph{\bf Hauptmodul}
We define the {\it hauptmodul} of $\Gamma_0^{*}(3)$:
\begin{equation}
J_{3+} := ({E_{2, 3}}')^3 / \Delta_3 = \frac{1}{q} + 42 + 783 q + 8672 q^2 + 65367 q^3 + \cdots,
\end{equation}
where $v_{\infty}(J_{3+}) = -1$ and $v_{\rho_3}(J_{3+}) = 6$. Then, we have
\begin{equation}
J_{3+} : \partial \mathbb{F}_{3+} \setminus \{z \in \mathbb{H} \: ; \: Re(z) = \pm 1/2\} \to [0, 108] \subset \mathbb{R}.
\end{equation}\quad

\subsection{$\Gamma_0(3)$} (see \cite{SJ1})

\paragraph{\bf Fundamental domain}
We have a fundamental domain for $\Gamma_0(3)$ as follows:
\begin{equation}
\mathbb{F}_3 = \left\{|z + 1/3| \geqslant 1/3, \: - 1/2 \leqslant Re(z) \leqslant 0\right\}
 \bigcup \left\{|z - 1/3| > 1/3, \: 0 < Re(z) < 1/2 \right\},
\end{equation}
where $\left( \begin{smallmatrix} -1 & 0 \\ 3 & -1 \end{smallmatrix} \right) : (e^{i \theta} + 1) / 3 \rightarrow (e^{i (\pi - \theta)} - 1) / 3$. Then, we have
\begin{equation}
\Gamma_0(3) = \langle \left( \begin{smallmatrix} 1 & 1 \\ 0 & 1 \end{smallmatrix} \right), \: - \left( \begin{smallmatrix} 1 & 0 \\ 3 & 1 \end{smallmatrix} \right) \rangle.
\end{equation}
\begin{figure}[hbtp]
\begin{center}
\includegraphics[width=1.5in]{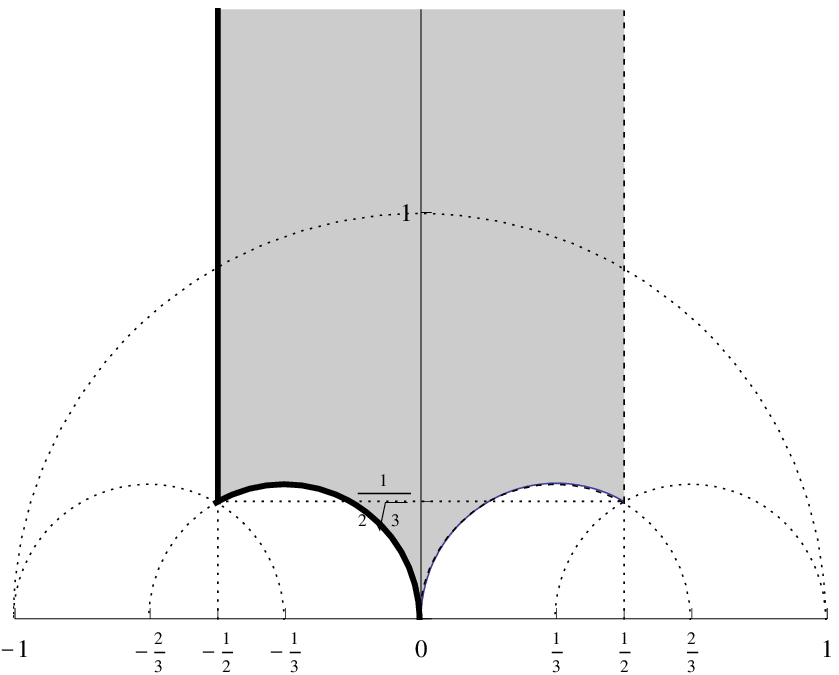}
\end{center}
\caption{$\Gamma_0(3)$}
\end{figure}

\paragraph{\bf Valence formula}
The cusps of $\Gamma_0(3)$ are $\infty$ and $0$, and the elliptic point is $\rho_3$. Let $f$ be a modular function of weight $k$ for $\Gamma_0(3)$, which is not identically zero. We have
\begin{equation}
v_{\infty}(f) + v_0(f) + \frac{1}{3} v_{\rho_3} (f) + \sum_{\begin{subarray}{c} p \in \Gamma_0(3) \setminus \mathbb{H} \\ p \ne \rho_3\end{subarray}} v_p(f) = \frac{k}{3}.
\end{equation}

Furthermore, the stabilizer of the elliptic point $\rho_3$ is $\left\{ \pm I, \pm \left( \begin{smallmatrix}-2 & -1 \\ 3 & 1\end{smallmatrix} \right), \pm \left( \begin{smallmatrix}-1 & -1 \\ 3 & 2\end{smallmatrix} \right) \right\}$.\\

\paragraph{\bf For the cusp $\infty$}
We have $\Gamma_{\infty} = \left\{ \pm \left( \begin{smallmatrix}1 & n \\ 0 & 1\end{smallmatrix} \right) \: ; \: n \in \mathbb{Z} \right\}$, and we have the Eisenstein series for the cusp $\infty$ associated with $\Gamma_0(3)$:
\begin{equation}
E_{k, 3}^{\infty}(z) := \frac{3^k E_k(3 z) - E_k(z)}{3^k - 1} \quad \text{for} \; k \geqslant 4.
\end{equation}\quad

\paragraph{\bf For the cusp $0$}
We have $\Gamma_0 = \left\{ \pm \left( \begin{smallmatrix}1 & 0 \\ 3 n & 1\end{smallmatrix} \right) \: ; \: n \in \mathbb{Z} \right\}$ and $\gamma_0 = W_3$, and we have the Eisenstein series for the cusp $0$ associated with $\Gamma_0(3)$:
\begin{equation}
E_{k, 3}^0(z) := \frac{- 3^{k/2} (E_k(3 z) - E_k(z))}{3^k - 1} \quad \text{for} \; k \geqslant 4.
\end{equation}
We also have $\gamma_0^{-1} \ \Gamma_0(3) \ \gamma_0 = \Gamma_0(3)$.\\

\paragraph{\bf The space of modular forms}
We have $M_k(\Gamma_0(3)) = \mathbb{C} E_{k, 3}^{\infty} \oplus \mathbb{C} E_{k, 3}^0 \oplus S_k(\Gamma_0(3))$ and $S_k(\Gamma_0(3)) = \Delta_3 M_{k - 6}(\Gamma_0(3))$ for every even integer $k \geqslant 4$. Then, we have $M_{6 n + 2}(\Gamma_0(3)) = {E_{2, 3}}' M_{6 n}(\Gamma_0(3))$ and
\begin{align*}
M_{6 n}(\Gamma_0(3)) &= \mathbb{C} (E_{6, 3}^{\infty})^n \oplus \mathbb{C} (E_{6, 3}^{\infty})^{n-1} \Delta_3 \oplus \cdots \oplus \mathbb{C} E_{6, 3}^{\infty} (\Delta_3)^{n-1}\\
 &\oplus \mathbb{C} (E_{6, 3}^0)^n \oplus \mathbb{C} (E_{6, 3}^0)^{n-1} \Delta_3 \oplus \cdots \oplus \mathbb{C} E_{6, 3}^0 (\Delta_3)^{n-1} \oplus \mathbb{C} (\Delta_3)^n,\\
M_{6 n + 4}(\Gamma_0(3)) &= E_{4, 3}^{\infty} (\mathbb{C} (E_{6, 3}^{\infty})^n \oplus \mathbb{C} (E_{6, 3}^{\infty})^{n-1} \Delta_3 \oplus \cdots \oplus \mathbb{C} (\Delta_3)^n)\\
 &\oplus E_{4, 3}^0 (\mathbb{C} (E_{6, 3}^0)^n \oplus \mathbb{C} (E_{6, 3}^0)^{n-1} \Delta_3 \oplus \cdots \oplus \mathbb{C} (\Delta_3)^n).
\end{align*}

Here, we define ${E_{1, 3}}' := \sqrt{{E_{2, 3}}'}$, which satisfies $v_{\rho_3}({E_{1, 3}}') = 1$. Then, we can write $E_{4, 3}^{\infty} = {E_{1, 3}}' \Delta_3^0$ and $(1/27) E_{4, 3}^0 = {E_{1, 3}}' \Delta_3^{\infty}$. Furthermore, we have $E_{6, 3}^{\infty} = (\Delta_3^{\infty})^2 + (243/13) \Delta_3^{\infty} \Delta_3^0$ and $(-13/243) E_{6, 3}^0 = \Delta_3^{\infty} \Delta_3^0 + 39 (\Delta_3^0)^2$. Now, we have
\begin{equation*}
M_k(\Gamma_0(3)) = {E_{k - 3n, 3}}' (\mathbb{C} (\Delta_3^{\infty})^n \oplus \mathbb{C} (\Delta_3^{\infty})^{n-1} \Delta_3^0 \oplus \cdots \oplus \mathbb{C} (\Delta_3^0)^n),
\end{equation*}
where $n = \dim(M_k(\Gamma_0(3))) - 1 = \lfloor k/3 \rfloor$ and where ${E_{0, 3}}' := 1$.\\

\paragraph{\bf Hauptmodul}
We define the {\it hauptmodul} of $\Gamma_0(3)$:
\begin{equation}
J_3 := \Delta_3^0 / \Delta_3^{\infty} \: (= \eta^{12}(z) / \eta^{12}(3 z))) = \frac{1}{q} -12 + 54 q - 76 q^2 - 243 q^3 + \cdots,
\end{equation}
where $v_{\infty}(J_3) = -1$ and $v_0(J_3) = 1$. Then, we have
\begin{equation}
J_3 : \partial \mathbb{F}_3 \setminus \{z \in \mathbb{H} \: ; \: Re(z) = \pm 1/2\} \to [-27, 0] \subset \mathbb{R}.
\end{equation}

\clearpage

\section{Level $4$}

We have $\Gamma_0(4)+=\Gamma_0(4)+4=\Gamma_0^{*}(4)$ and $\Gamma_0(4)-=\Gamma_0(4)$.

We have $W_4 = \left(\begin{smallmatrix} 0 & -1/2 \\ 2 & 0 \end{smallmatrix}\right)$ and define $W_{4-, 2} = \left(\begin{smallmatrix} -1 & -1 \\ 2 & 1 \end{smallmatrix}\right)$ and $W_{4+, 2} = \left(\begin{smallmatrix} -1/\sqrt{2} & -3/(2 \sqrt{2}) \\ \sqrt{2} & 1/\sqrt{2} \end{smallmatrix}\right)$. We define
\begin{equation}
\begin{split}
&\Delta_4^{\infty}(z) := \eta^8(4 z) / \eta^4(2 z),
 \quad \Delta_4^0(z) := \eta^8(z) / \eta^4(2 z),
 \quad \Delta_4^{-1/2}(z) := \eta^{20}(2 z) / (\eta^8(z) \eta^8(4 z)),\\
&\Delta_4(z) := \Delta_4^{\infty}(z) \Delta_4^0(z) \Delta_4^{-1/2}(z) = \eta^{12}(2 z),\\
&\Delta_{4+4}(z) := \Delta_4^{\infty}(z) \Delta_4^0(z) (\Delta_4^{-1/2}(z))^2 = \eta^{32}(2 z) / (\eta^8(z) \eta^8(4 z)),
\end{split}
\end{equation}
where $\Delta_4^{\infty}$, $\Delta_4^0$, and $\Delta_4^{-1/2}$ are modular forms for $\Gamma_0(4)$ of weight $2$ such that $v_{\infty}(\Delta_4^{\infty}) = v_0(\Delta_4^0) = v_{-1/2}(\Delta_4^{-1/2}) = 1$. Then, $\Delta_4$ is a cusp form for $\Gamma_0(4)$ of weight $6$, $\Delta_{4+4}$ is a cusp form for $\Gamma_0^{*}(4)$ of weight $8$. Furthermore, we define
\begin{equation}
\begin{split}
{E_{2, 4}}'(z) &:= {E_{2, 2}}'(2 z) = 2 E_2(4 z) - E_2(2 z),\\
{E_{2, 4+4}}'(z) &:= 2 {E_{2, 4}}'(z) - {E_{2, 2}}'(z) = 4 E_2(4 z) - 4 E_2(2 z) + E_2(z),
\end{split}
\end{equation}
which are modular forms for $\Gamma_0(4)$, $\Gamma_0^{*}(4)$ of weight $2$, respectively. Then, we have $v_{-1/4 + i/4}({E_{2, 4}}') =1$ and $v_{i/2}({E_{2, 4+4}}') =1$.\\

\subsection{$\Gamma_0^{*}(4)$} (see \cite{SJ1})

We have $\Gamma_0^{*}(4) = T_{1/2}^{-1} \: \Gamma_0(2) \: T_{1/2}$.\\

\paragraph{\bf Fundamental domain}
We have a fundamental domain for $\Gamma_0^{*}(4)$ as follows:
\begin{equation}
\mathbb{F}_{4+4} = \left\{|z| \geqslant 1/2, \: - 1/2 \leqslant Re(z) \leqslant 0\right\}
 \bigcup \left\{|z| > 1/2, \: 0 < Re(z) < 1/2 \right\},
\end{equation}
where $W_4 : e^{i \theta} / 2 \rightarrow e^{i (\pi - \theta)} / 2$. Then, we have
\begin{equation}
\Gamma_0^{*}(4) = \langle \left( \begin{smallmatrix} 1 & 1 \\ 0 & 1 \end{smallmatrix} \right), \: W_4 \rangle.
\end{equation}
\begin{figure}[hbtp]
\begin{center}
\includegraphics[width=1.5in]{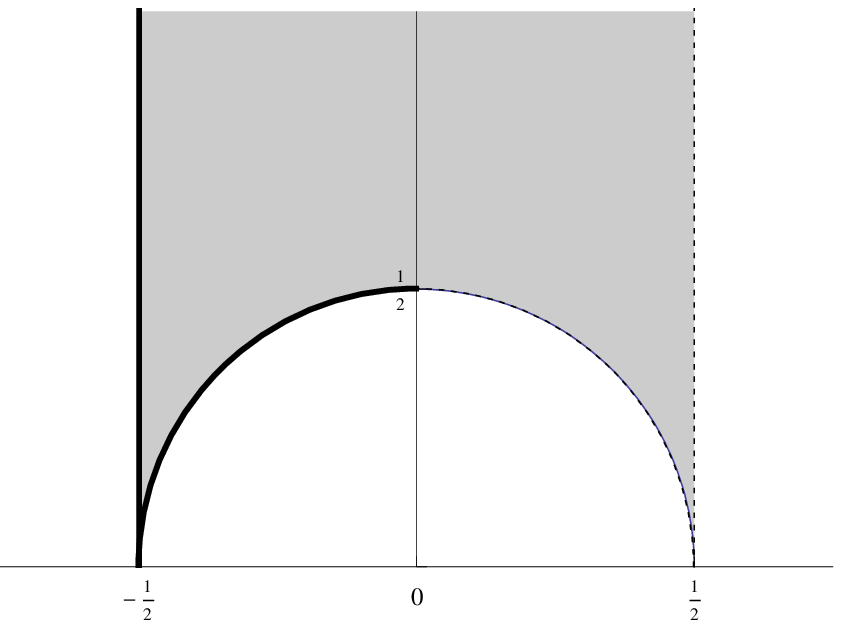}
\end{center}
\caption{$\Gamma_0^{*}(4)$}
\end{figure}

\paragraph{\bf Valence formula}
The cusps of $\Gamma_0^{*}(4)$ are $\infty$, $0$, and $-1/2$, and it has no elliptic point. Let $f$ be a modular function of weight $k$ for $\Gamma_0^{*}(4)$, which is not identically zero. We have
\begin{equation}
v_{\infty}(f) + v_{-1/2} (f) + \frac{1}{2} v_{i/2} + \sum_{\begin{subarray}{c} p \in \Gamma_0^{*}(4) \setminus \mathbb{H} \\ p \ne i/2\end{subarray}} v_p(f) = \frac{k}{4}.
\end{equation}\quad

\paragraph{\bf For the cusp $\infty$}
We have $\Gamma_{\infty} = \left\{ \pm \left( \begin{smallmatrix}1 & n \\ 0 & 1\end{smallmatrix} \right) \: ; \: n \in \mathbb{Z} \right\}$, and we have the Eisenstein series for the cusp $\infty$ associated with $\Gamma_0^{*}(4)$:
\begin{equation}
E_{k, 4+4}^{\infty}(z) := \frac{2^k E_k(4 z) - 2 E_k(2 z) + E_k(z)}{2^k - 1} \quad \text{for} \; k \geqslant 4.
\end{equation}\quad

\paragraph{\bf For the cusp $-1/2$}
We have $\Gamma_{-1/2} = \left\{ \pm \left( \begin{smallmatrix}n+1 & n/2 \\ -2n & - n+1\end{smallmatrix} \right) \: ; \: n \in \mathbb{Z} \right\}$ and $\gamma_{-1/2} = W_{4+, 2}$ and we have the Eisenstein series for the cusp $-1/2$ associated with $\Gamma_0^{*}(4)$:
\begin{equation}
E_{k, 4+4}^{-1/2}(z) := \frac{- 2^{k/2} (2^k E_k(4 z) - (2^k + 1) E_k(2 z) + E_k(z))}{2^k - 1} \quad \text{for} \; k \geqslant 4.
\end{equation}
We also have $\gamma_{-1/2}^{-1} \ \Gamma_0^{*}(4) \ \gamma_{-1/2} = \Gamma_0^{*}(4)$.\\

\paragraph{\bf The space of modular forms}
We have
\begin{equation*}
\Delta_{4+4}^{\infty} = \Delta_4^{\infty} \Delta_4^0,
 \quad \Delta_{4+4}^{-1/2} = (\Delta_4^{-1/2})^2.
\end{equation*}

Now, we have $M_k(\Gamma_0^{*}(4)) = \mathbb{C} E_{k, 4+4}^{\infty} \oplus \mathbb{C} E_{k, 4}^0 \oplus S_k(\Gamma_0^{*}(4))$ and $S_k(\Gamma_0^{*}(4)) = \Delta_{4+4} M_{k - 6}(\Gamma_0^{*}(4))$ for every even integer $k \geqslant 4$. Then, we have $M_{4 n + 2}(\Gamma_0^{*}(4)) = {E_{2, 4+4}}' M_{4 n}(\Gamma_0^{*}(4))$ and
\begin{allowdisplaybreaks}
\begin{align*}
M_{8 n}(\Gamma_0^{*}(4)) = &\mathbb{C} ((E_{4, 4+4}^{\infty})^2)^n \oplus \mathbb{C} ((E_{4, 4+4}^{\infty})^2)^{n-1} \Delta_{4+4} \oplus \cdots \oplus \mathbb{C} (E_{4, 4+4}^{\infty})^2 (\Delta_{4+4})^{n-1}\\
 &\oplus \mathbb{C} ((E_{4, 4+4}^{-1/2})^2)^n \oplus \mathbb{C} ((E_{4, 4+4}^{-1/2})^2)^{n-1} \Delta_{4+4} \oplus \cdots \oplus \mathbb{C} (E_{4, 4+4}^{-1/2})^2 (\Delta_{4+4})^{n-1} \oplus \mathbb{C} (\Delta_{4+4})^n,\\
M_{8 n + 4}(\Gamma_0^{*}(4)) = &E_{4, 4+4}^{\infty} (\mathbb{C} ((E_{4, 4+4}^{\infty})^2)^n \oplus \mathbb{C} ((E_{4, 4+4}^{\infty})^2)^{n-1} \Delta_{4+4} \oplus \cdots \oplus \mathbb{C} (\Delta_{4+4})^n)\\
 &\oplus E_{4, 4+4}^{-1/2} (\mathbb{C} ((E_{4, 4+4}^{-1/2})^2)^n \oplus \mathbb{C} ((E_{4, 4+4}^{-1/2})^2)^{n-1} \Delta_{4+4} \oplus \cdots \oplus \mathbb{C} (\Delta_{4+4})^n).
\end{align*}
\end{allowdisplaybreaks}

Here, since we have $E_{4, 4+4}^{\infty} = \Delta_{4+4}^{-1/2}$ and $E_{4, 4+4}^{-1/2} = \Delta_{4+4}^{\infty}$, we can write
\begin{equation*}
M_{4 n}(\Gamma_0^{*}(4)) = \mathbb{C} (\Delta_{4+4}^{\infty})^n \oplus \mathbb{C} (\Delta_{4+4}^{\infty})^{n-1} \Delta_{4+4}^{-1/2} \oplus \cdots \oplus \mathbb{C} (\Delta_{4+4}^{-1/2})^n.
\end{equation*}\quad

\paragraph{\bf Hauptmodul}
We define the {\it hauptmodul} of $\Gamma_0^{*}(4)$:
\begin{equation}
J_{4+4} := \Delta_{4+4}^{-1/2} / \Delta_{4+4}^{\infty} \: (= \eta^{48}(2 z) / (\eta^{24}(z) \eta^{24}(4 z))) = \frac{1}{q} + 24 + 276 q + 2048 q^2 + 11202 q^3 + \cdots,
\end{equation}
where $v_{\infty}(J_{4+4}) = -1$ and $v_{-1/2}(J_{4+4}) = 1$. Then, we have
\begin{equation}
J_{4+4} : \partial \mathbb{F}_{4+4} \setminus \{z \in \mathbb{H} \: ; \: Re(z) = \pm 1/2\} \to [0, 64] \subset \mathbb{R}.
\end{equation}\quad

\subsection{$\Gamma_0(4)$} (see \cite{SJ1})

We have $\Gamma_0(4) = V_2^{-1} \Gamma(2) V_2$.\\

\paragraph{\bf Fundamental domain}
We have a fundamental domain for $\Gamma_0(4)$ as follows:
\begin{equation}
\mathbb{F}_4 = \left\{|z + 1/4| \geqslant 1/4, \: - 1/2 \leqslant Re(z) \leqslant 0\right\}
 \bigcup \left\{|z - 1/4| > 1/4, \: 0 < Re(z) < 1/2 \right\},
\end{equation}
where $\left( \begin{smallmatrix} -1 & 0 \\ 4 & -1 \end{smallmatrix} \right) : (e^{i \theta} + 1) / 4 \rightarrow (e^{i (\pi - \theta)} - 1) / 4$. Then, we have
\begin{equation}
\Gamma_0(4) = \langle -I, \: \left( \begin{smallmatrix} 1 & 1 \\ 0 & 1 \end{smallmatrix} \right), \: \left( \begin{smallmatrix} 1 & 0 \\ 4 & 1 \end{smallmatrix} \right) \rangle.
\end{equation}
\begin{figure}[hbtp]
\begin{center}
\includegraphics[width=1.5in]{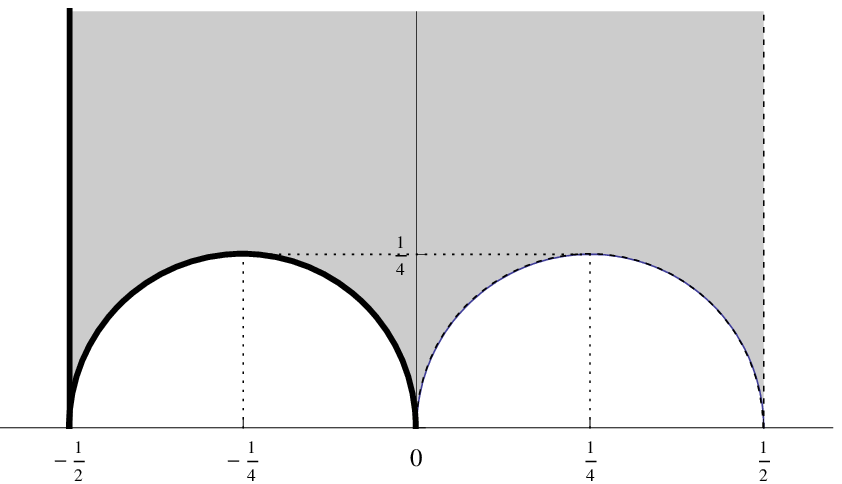}
\end{center}
\caption{$\Gamma_0(4)$}
\end{figure}

\paragraph{\bf Valence formula}
The cusps of $\Gamma_0(4)$ are $\infty$, $0$, and $-1/2$, and it has no elliptic point. Let $f$ be a modular function of weight $k$ for $\Gamma_0(4)$, which is not identically zero. We have
\begin{equation}
v_{\infty}(f) + v_0(f) + v_{-1/2} (f) + \sum_{p \in \Gamma_0(4) \setminus \mathbb{H}} v_p(f) = \frac{k}{2}.
\end{equation}\quad

\paragraph{\bf For the cusp $\infty$}
We have $\Gamma_{\infty} = \left\{ \pm \left( \begin{smallmatrix}1 & n \\ 0 & 1\end{smallmatrix} \right) \: ; \: n \in \mathbb{Z} \right\}$, and we have the Eisenstein series for the cusp $\infty$ associated with $\Gamma_0(4)$:
\begin{equation}
E_{k, 4}^{\infty}(z) := \frac{2^k E_k(4 z) - E_k(2 z)}{2^k - 1} \quad \text{for} \; k \geqslant 4.
\end{equation}
Note that we have $E_{k, 4}^{\infty}(z) = E_{k, 2}^{\infty}(2 z)$.\\

\paragraph{\bf For the cusp $0$}
We have $\Gamma_0 = \left\{ \pm \left( \begin{smallmatrix}1 & 0 \\ 4 n & 1\end{smallmatrix} \right) \: ; \: n \in \mathbb{Z} \right\}$ and $\gamma_0 = W_4$, and we have the Eisenstein series for the cusp $0$ associated with $\Gamma_0(4)$:
\begin{equation}
E_{k, 4}^0(z) := \frac{- (E_k(2 z) - E_k(z))}{2^k - 1} \quad \text{for} \; k \geqslant 4.
\end{equation}
Note that we have $E_{k, 4}^0 = 2^{-k/2} E_{k, 2}^0$. We also have $\gamma_0^{-1} \ \Gamma_0(4) \ \gamma_0 = \Gamma_0(4)$.\\

\paragraph{\bf For the cusp $-1/2$}
We have $\Gamma_{-1/2} = \left\{ \pm \left( \begin{smallmatrix}2 n + 1 & n \\ - 4 n & - 2 n + 1\end{smallmatrix} \right) \: ; \: n \in \mathbb{Z} \right\}$ and $\gamma_{-1/2} = W_{4-, 2}$ and we have the Eisenstein series for the cusp $-1/2$ associated with $\Gamma_0(4)$:
\begin{equation}
E_{k, 4}^{-1/2}(z) := \frac{- (2^k E_k(4 z) - (2^k + 1) E_k(2 z) + E_k(z))}{2^k - 1} \quad \text{for} \; k \geqslant 4.
\end{equation}
Note that we have $E_{k, 4}^{-1/2} = 2^{-k/2} E_{k, 2+}^{-1/2}$. We also have $\gamma_{-1/2}^{-1} \ \Gamma_0(4) \ \gamma_{-1/2} = \Gamma_0(4)$.\\

\paragraph{\bf The space of modular forms}
We have ${E_{2, 2}}'$ as a modular form for $\Gamma_0(4) \subset \Gamma_0(2)$ of weight $2$ such that $v_{-1/2 + i/2}({E_{2, 2}}') =1$. Moreover, we have $M_2(\Gamma_0(4)) = \mathbb{C} {E_{2, 2}}' \oplus \mathbb{C} {E_{2, 4}}'$.

Now, we have $M_k(\Gamma_0(4)) = \mathbb{C} E_{k, 4}^{\infty} \oplus \mathbb{C} E_{k, 4}^0 \oplus \mathbb{C} E_{k, 4}^{-1/2} \oplus S_k(\Gamma_0(4))$ and $S_k(\Gamma_0(4)) = \Delta_4 M_{k - 6}(\Gamma_0(4))$ for every even integer $k \geqslant 4$. Then, we have $M_{6 n + 2}(\Gamma_0(4)) = {E_{2, 4}}' M_{6 n}(\Gamma_0(4)) \oplus \mathbb{C} {E_{2, 2}}' (\Delta_4)^n$ and
\begin{allowdisplaybreaks}
\begin{align*}
M_{6 n}(\Gamma_0(4)) = &\mathbb{C} (E_{6, 4}^{\infty})^n \oplus \mathbb{C} (E_{6, 4}^{\infty})^{n-1} \Delta_4 \oplus \cdots \oplus \mathbb{C} E_{6, 4}^{\infty} (\Delta_4)^{n-1}\\
 &\oplus \mathbb{C} (E_{6, 4}^0)^n \oplus \mathbb{C} (E_{6, 4}^0)^{n-1} \Delta_4 \oplus \cdots \oplus \mathbb{C} E_{6, 4}^0 (\Delta_4)^{n-1}\\
 &\oplus \mathbb{C} (E_{6, 4}^{-1/2})^n \oplus \mathbb{C} (E_{6, 4}^{-1/2})^{n-1} \Delta_4 \oplus \cdots \oplus \mathbb{C} E_{6, 4}^{-1/2} (\Delta_4)^{n-1}\oplus \mathbb{C} (\Delta_4)^n,\\
M_{6 n + 4}(\Gamma_0(4)) = &E_{4, 4}^{\infty} (\mathbb{C} (E_{6, 4}^{\infty})^n \oplus \mathbb{C} (E_{6, 4}^{\infty})^{n-1} \Delta_4 \oplus \cdots \oplus \mathbb{C} (\Delta_4)^n)\\
 &\oplus E_{4, 4}^0 (\mathbb{C} (E_{6, 4}^0)^n \oplus \mathbb{C} (E_{6, 4}^0)^{n-1} \Delta_4 \oplus \cdots \oplus \mathbb{C} (\Delta_4)^n)\\
 &\oplus E_{4, 4}^{-1/2} (\mathbb{C} (E_{6, 4}^{-1/2})^n \oplus \mathbb{C} (E_{6, 4}^{-1/2})^{n-1} \Delta_4 \oplus \cdots \oplus \mathbb{C} (\Delta_4)^n).
\end{align*}
\end{allowdisplaybreaks}

Here, we have ${E_{2, 2}}' = 32 \Delta_4^{\infty} + \Delta_4^0$, ${E_{2, 4}}' = 16 \Delta_4^{\infty} + \Delta_4^0$, $E_{4, 4}^{\infty} = 16 \Delta_4^{\infty} \Delta_4^0 + (\Delta_4^0)^2$, $(1/16) E_{4, 4}^0 = 16 (\Delta_4^{\infty})^2 + \Delta_4^{\infty} \Delta_4^0$, $(-1/16) E_{4, 4}^{-1/2} = \Delta_4^{\infty} \Delta_4^0$, $E_{6, 4}^{\infty} = 128 (\Delta_4^{\infty})^2 \Delta_4^0 + 24 \Delta_4^{\infty} (\Delta_4^0)^2 + (\Delta_4^0)^3$, $(-1/8) E_{6, 4}^0 = 512 (\Delta_4^{\infty})^3 + 48 (\Delta_4^{\infty})^2 \Delta_4^0 + \Delta_4^{\infty} (\Delta_4^0)^2$, and $(1/8) E_{6, 4}^{-1/2} = -16 (\Delta_4^{\infty})^2 \Delta_4^0 + \Delta_4^{\infty} (\Delta_4^0)^2$. Now, we can write
\begin{equation*}
M_{2 n}(\Gamma_0(4)) = \mathbb{C} (\Delta_4^{\infty})^n \oplus \mathbb{C} (\Delta_4^{\infty})^{n-1} \Delta_4^0 \oplus \cdots \oplus \mathbb{C} (\Delta_4^0)^n.
\end{equation*}\quad

\paragraph{\bf Hauptmodul}
We define the {\it hauptmodul} of $\Gamma_0(4)$:
\begin{equation}
J_4 := \Delta_4^0 / \Delta_4^{\infty} \: (= \eta^8(z) / \eta^8(4 z)) = \frac{1}{q} - 8 + 20 q - 62 q^3 + 216 q^5 + \cdots,
\end{equation}
where $v_{\infty}(J_4) = -1$ and $v_0(J_4) = 1$. Then, we have
\begin{equation}
J_4 : \partial \mathbb{F}_4 \setminus \{z \in \mathbb{H} \: ; \: Re(z) = \pm 1/2\} \to [-16, 0] \subset \mathbb{R}.
\end{equation}

\clearpage

\section{Level $5$}

We have $\Gamma_0(5)+=\Gamma_0^{*}(5)$ and $\Gamma_0(5)-=\Gamma_0(5)$.

We have $W_5 = \left(\begin{smallmatrix}0&-1 / \sqrt{5}\\ \sqrt{5}&0\end{smallmatrix}\right)$, and denote $\rho_{5, 1} := - 1/2 + i / (2 \sqrt{5})$, $\rho_{5, 2} := - 2/5 + i / 5$, and $\rho_{5, 3} := 2/5 + i / 5$. We define
\begin{equation}
\begin{split}
&\Delta_5^{\infty}(z) := \eta^5(5 z) / \eta(z), \quad \Delta_5^0(z) := \eta^5(z) / \eta(5 z),\\
&\Delta_5(z) := \Delta_5^{\infty}(z) \Delta_5^0(z) = \eta^4(z) \eta^4(5 z),
\end{split}
\end{equation}
where $\Delta_5^{\infty}$ and $\Delta_5^0$ are $2$nd semimodular forms for $\Gamma_0(5)$ of weight $2$ such that $v_{\infty}(\Delta_5^{\infty}) = v_0(\Delta_5^0) = 1$, and $\Delta_5$ is a cusp form for $\Gamma_0(5)$ and $\Gamma_0^{*}(5)$ of weight $4$. Furthermore, we define
\begin{equation}
{E_{2, 5}}'(z) := (5 E_2(5 z) - E_2(z)) / 4,
\end{equation}
which is a modular form for $\Gamma_0(5)$ and $2$nd semimodular form for $\Gamma_0^{*}(5)$ of weight $2$ such that $v_{\rho_{5, 2}}({E_{2, 5}}') = 1$.\\

\subsection{$\Gamma_0^{*}(5)$} (see \cite{SJ2})

\paragraph{\bf Fundamental domain}
We have a fundamental domain for $\Gamma_0^{*}(5)$ as follows:
{\small \begin{equation}
\begin{split}
\mathbb{F}_{5+} = &\left\{|z + 1/2| \geqslant 1 / (2 \sqrt{5}), \: - 1/2 \leqslant Re(z) < - 2/5 \right\}
 \bigcup \left\{|z| \geqslant 1 / \sqrt{5}, \: - 2/5 \leqslant Re(z) \leqslant 0 \right\}\\
 &\bigcup \left\{|z| > 1 / \sqrt{5}, \: 0 < Re(z) \leqslant 2/5 \right\}
 \bigcup \left\{|z - 1/2| > 1 / (2 \sqrt{5}), \: 2/5 < Re(z) < 1/2 \right\},
\end{split}
\end{equation}
}where $W_5 :  e^{i \theta} / \sqrt{5} \rightarrow e^{i  (\pi - \theta)} / \sqrt{5}$ and $\left(\begin{smallmatrix} -2 & -1 \\ 5 & 2 \end{smallmatrix}\right) W_5 : e^{i \theta} / (2 \sqrt{5}) + 1/2 \rightarrow e^{i  (\pi - \theta)} / (2 \sqrt{5}) - 1/2$. Then, we have
\begin{equation}
\Gamma_0^{*}(5) = \langle \left( \begin{smallmatrix} 1 & 1 \\ 0 & 1 \end{smallmatrix} \right), \: W_5, \: \left( \begin{smallmatrix} 3 & 1 \\ 5 & 2 \end{smallmatrix} \right) \rangle.
\end{equation}
\begin{figure}[hbtp]
\begin{center}
{{$\mathbb{F}_{5+}$}\includegraphics[width=1.5in]{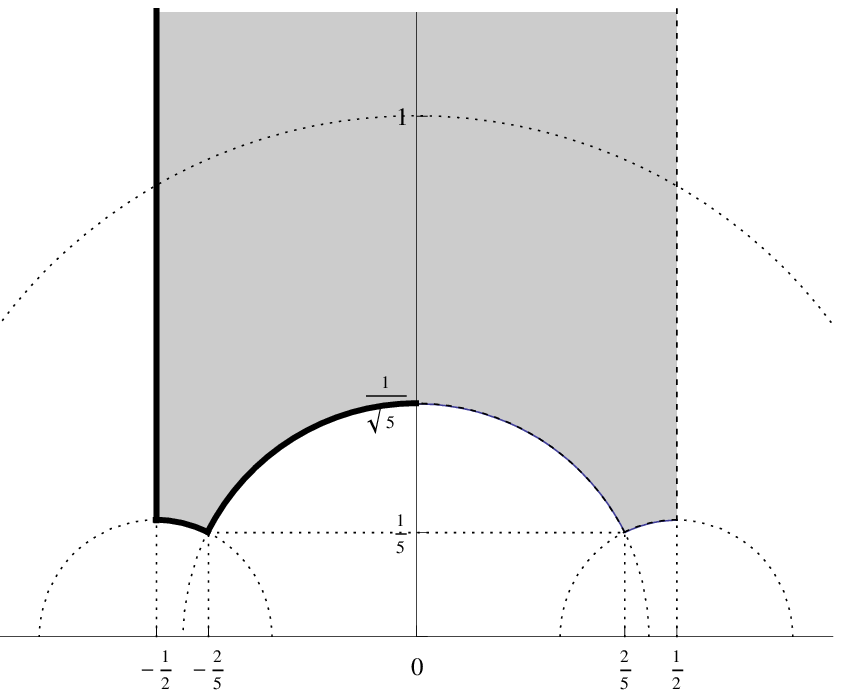}}
\end{center}
\caption{$\Gamma_0^{*}(5)$}
\end{figure}

\paragraph{\bf Valence formula}
The cusp of $\Gamma_0^{*}(5)$ is $\infty$, and the elliptic points are $i / \sqrt{5}$, $\rho_{5, 1} = - 1/2 + i \sqrt{5} / 10$ and $\rho_{5, 2} = - 2/5 + i / 5$. Let $f$ be a modular function of weight $k$ for $\Gamma_0^{*}(5)$, which is not identically zero. We have
\begin{equation}
v_{\infty}(f) + \frac{1}{2} v_{i / \sqrt{5}}(f) + \frac{1}{2} v_{\rho_{5, 1}} (f) + \frac{1}{2} v_{\rho_{5, 2}} (f) + \sum_{\begin{subarray}{c} p \in \Gamma_0^{*}(5) \setminus \mathbb{H} \\ p \ne i / \sqrt{5}, \; \rho_{5, 1}, \; \rho_{5, 2}\end{subarray}} v_p(f) = \frac{k}{4}.
\end{equation}

Furthermore, the stabilizer of the elliptic point $i / \sqrt{5}$ (resp. $\rho_{5, 1}$, $\rho_{5, 2}$) is $\left\{ \pm I, \pm W_5 \right\}$\\ (resp. $\left\{ \pm I, \pm \left( \begin{smallmatrix}3 & -1 \\ -5 & 2\end{smallmatrix} \right) W_5 \right\}$, $\left\{ \pm I, \pm \left( \begin{smallmatrix}-2 & -1 \\ 5 & 2\end{smallmatrix} \right) \right\}$)\\

\paragraph{\bf For the cusp $\infty$}
We have $\Gamma_{\infty} = \left\{ \pm \left( \begin{smallmatrix}1 & n \\ 0 & 1\end{smallmatrix} \right) \: ; \: n \in \mathbb{Z} \right\}$, and we have the Eisenstein series associated with $\Gamma_0^{*}(5)$:
\begin{equation}
E_{k, 5+}(z) := \frac{5^{k/2} E_k(5 z) + E_k(z)}{5^{k/2} + 1} \quad \text{for} \; k \geqslant 4.
\end{equation}\quad

\paragraph{\bf The space of modular forms}
We have $M_k(\Gamma_0^{*}(5)) = \mathbb{C} E_{k, 5+} \oplus S_k(\Gamma_0^{*}(5))$ and $S_k(\Gamma_0^{*}(5)) = \Delta_5 M_{k - 4}(\Gamma_0^{*}(5))$ for every even integer $k \geqslant 4$. Then, we have $M_{4 n + 6}(\Gamma_0^{*}(5)) = E_{6, 5+} M_{4 n}(\Gamma_0^{*}(5))$ and
\begin{equation*}
M_{4 n}(\Gamma_0^{*}(5)) = \mathbb{C} (({E_{2, 5}}')^2)^n \oplus \mathbb{C} (({E_{2, 5}}')^2)^{n-1} \Delta_5 \oplus \cdots \oplus \mathbb{C} (\Delta_5)^n.
\end{equation*}\quad

\paragraph{\bf Hauptmodul}
We define the {\it hauptmodul} of $\Gamma_0^{*}(5)$:
\begin{equation}
J_{5+} := ({E_{2, 5}}')^2 / \Delta_5 = \frac{1}{q} + 16 + 134 q + 760 q^2 + 3345 q^3 + \cdots,
\end{equation}
where $v_{\infty}(J_{5+}) = -1$ and $v_{\rho_{5, 2}}(J_{5+}) = 2$. Then, we have
\begin{equation}
J_{5+} : \partial \mathbb{F}_{5+} \setminus \{z \in \mathbb{H} \: ; \: Re(z) = \pm 1/2\} \to [22 - 10 \sqrt{5}, 22 + 10 \sqrt{5}] \subset \mathbb{R}.
\end{equation}\quad

\subsection{$\Gamma_0(5)$} (see \cite{SJ1})

\paragraph{\bf Fundamental domain}
We have a fundamental domain for $\Gamma_0(5)$ as follows:
{\small \begin{equation}
\begin{split}
\mathbb{F}_5 = &\left\{|z + 2/5| \geqslant 1/5, \: - 1/2 \leqslant Re(z) \leqslant - 2/5\right\}
 \bigcup \left\{|z + 2/5| > 1/5, \: - 2/5 < Re(z) < - 3 / 10 \right\}\\
&\bigcup \left\{|z + 1/5| \geqslant 1/5, \: - 3 / 10 \leqslant Re(z) \leqslant 0\right\}
 \bigcup \left\{|z - 1/5| > 1/5, \: 0 < Re(z) < 3 / 10 \right\}\\
&\bigcup \left\{|z - 2/5| \geqslant 1/5, \: 3 / 10 \leqslant Re(z) \leqslant 2/5 \right\}
 \bigcup \left\{|z - 2/5| > 1/5, \: 2/5 < Re(z) < 1/2 \right\},
\end{split}
\end{equation}
}where $\left( \begin{smallmatrix} -1 & 0 \\ 5 & -1 \end{smallmatrix} \right) : (e^{i \theta} + 1) / 5 \rightarrow (e^{i (\pi - \theta)} - 1) / 5$, $\left( \begin{smallmatrix} 2 & -1 \\ 5 & -2 \end{smallmatrix} \right) : (e^{i \theta} + 2) / 5 \rightarrow (e^{i (\pi - \theta)} + 2) / 5$, and  $\left( \begin{smallmatrix} -2 & -1 \\ 5 & 2 \end{smallmatrix} \right) : (e^{i \theta} - 2) / 5 \rightarrow (e^{i (\pi - \theta)} - 2) / 5$. Then, we have
\begin{equation}
\Gamma_0(5) = \langle \left( \begin{smallmatrix} 1 & 1 \\ 0 & 1 \end{smallmatrix} \right), \: \left( \begin{smallmatrix} 1 & 0 \\ 5 & 1 \end{smallmatrix} \right), \: \left( \begin{smallmatrix} 3 & 1 \\ 5 & 2 \end{smallmatrix} \right) \rangle.
\end{equation}
\begin{figure}[hbtp]
\begin{center}
\includegraphics[width=1.5in]{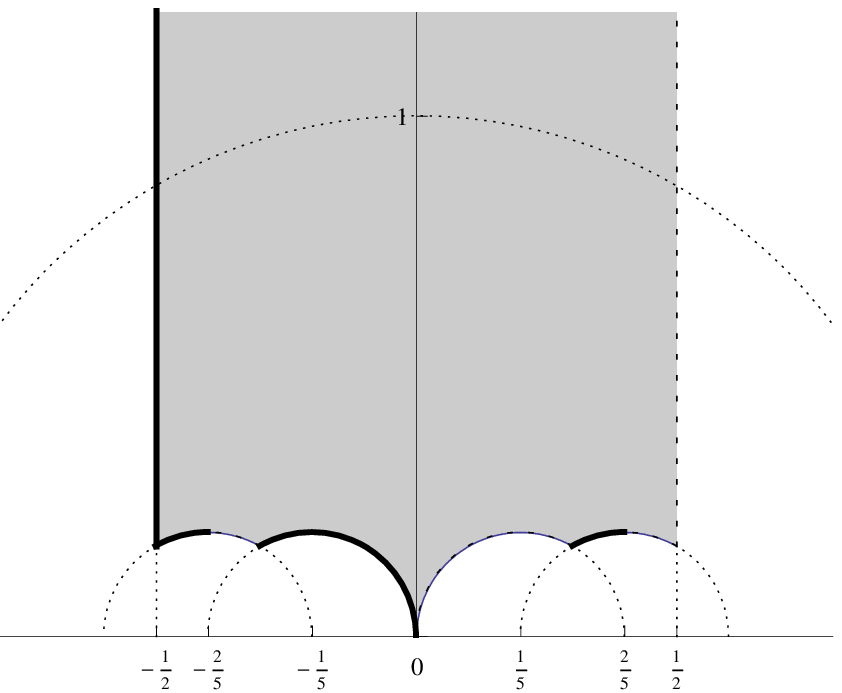}
\end{center}
\caption{$\Gamma_0(5)$}
\end{figure}

\paragraph{\bf Valence formula}
The cusps of $\Gamma_0(5)$ are $\infty$ and $0$, and the elliptic points are $\rho_{5, 2}$ and $\rho_{5, 3} = 2/5 + i / 5$. Let $f$ be a modular function of weight $k$ for $\Gamma_0(5)$, which is not identically zero. We have
\begin{equation}
v_{\infty}(f) + v_0(f) + \frac{1}{2} v_{\rho_{5, 2}} (f) + \frac{1}{2} v_{\rho_{5, 3}} (f) + \sum_{\begin{subarray}{c} p \in \Gamma_0(5) \setminus \mathbb{H} \\ p \ne \rho_{5, 2}, \rho_{5, 3}\end{subarray}} v_p(f) = \frac{k}{2}.
\end{equation}

Furthermore, the stabilizer of the elliptic point $\rho_{5, 2}$ (resp. $\rho_{5, 3}$) is $\left\{ \pm I, \pm \left( \begin{smallmatrix}-2 & -1 \\ 5 & 2\end{smallmatrix} \right) \right\}$ (resp. $\left\{ \pm I, \pm \left( \begin{smallmatrix}2 & -1 \\ 5 & -2\end{smallmatrix} \right) \right\}$).\\

\paragraph{\bf For the cusp $\infty$}
We have $\Gamma_{\infty} = \left\{ \pm \left( \begin{smallmatrix}1 & n \\ 0 & 1\end{smallmatrix} \right) \: ; \: n \in \mathbb{Z} \right\}$, and we have the Eisenstein series for the cusp $\infty$ associated with $\Gamma_0(5)$:
\begin{equation}
E_{k, 5}^{\infty}(z) := \frac{5^k E_k(5 z) - E_k(z)}{5^k - 1} \quad \text{for} \; k \geqslant 4.
\end{equation}\quad

\paragraph{\bf For the cusp $0$}
We have $\Gamma_0 = \left\{ \pm \left( \begin{smallmatrix}1 & 0 \\ 5 n & 1\end{smallmatrix} \right) \: ; \: n \in \mathbb{Z} \right\}$ and $\gamma_0 = W_5$, and we have the Eisenstein series for the cusp $0$ associated with $\Gamma_0(5)$:
\begin{equation}
E_{k, 5}^0(z) := \frac{- 5^{k/2} (E_k(5 z) - E_k(z))}{5^k - 1} \quad \text{for} \; k \geqslant 4.
\end{equation}
We also have $\gamma_0^{-1} \ \Gamma_0(5) \ \gamma_0 = \Gamma_0(5)$.\\

\paragraph{\bf The space of modular forms}
We have $M_k(\Gamma_0(5)) = \mathbb{C} E_{k, 5}^{\infty} \oplus \mathbb{C} E_{k, 5}^0 \oplus S_k(\Gamma_0(5))$ and $S_k(\Gamma_0(5)) = \Delta_5 M_{k - 4}(\Gamma_0(5))$ for every even integer $k \geqslant 4$. Then, we have $M_{4 n + 2}(\Gamma_0(5)) = {E_{2, 5}}' M_{4 n}(\Gamma_0(5))$ and
\begin{align*}
M_{4 n}(\Gamma_0(5)) &= \mathbb{C} (E_{4, 5}^{\infty})^n \oplus \mathbb{C} (E_{4, 5}^{\infty})^{n-1} \Delta_5 \oplus \cdots \oplus \mathbb{C} E_{4, 5}^{\infty} (\Delta_5)^{n-1}\\
 &\oplus \mathbb{C} (E_{4, 5}^0)^n \oplus \mathbb{C} (E_{4, 5}^0)^{n-1} \Delta_5 \oplus \cdots \oplus \mathbb{C} E_{4, 5}^0 (\Delta_5)^{n-1} \oplus \mathbb{C} (\Delta_5)^n.
\end{align*}

Here, we have $E_{4, 5}^{\infty} = (125/13) \Delta_5^{\infty} \Delta_5^0 + (\Delta_5^0)^2$ and $(13/125) E_{4, 5}^0 = 13 (\Delta_5^{\infty})^2 + \Delta_5^{\infty} \Delta_5^0$, then we can write
\begin{equation*}
M_{4 n}(\Gamma_0(5)) = \mathbb{C} (\Delta_5^{\infty})^{2n} \oplus \mathbb{C} (\Delta_5^{\infty})^{2n-1} \Delta_5^0 \oplus \cdots \oplus \mathbb{C} (\Delta_5^0)^{2n}.
\end{equation*}\quad

\paragraph{\bf Hauptmodul}
We define the {\it hauptmodul} of $\Gamma_0(5)$:
\begin{equation}
J_5 := \Delta_5^0 / \Delta_5^{\infty} \: (= \eta^6(z) / \eta^6(5 z)) = \frac{1}{q} - 6 + 9 q + 10 q^2 - 30 q^3 + \cdots,
\end{equation}
where $v_{\infty}(J_5) = -1$ and $v_0(J_5) = 1$. Then, we have
{\small \begin{equation}
\begin{split}
J_5 : \qquad \left\{|z + 1/5| = 1/5, \: - 3 / 10 \leqslant Re(z) \leqslant 0\right\} &\to [-5 - 2 \sqrt{5}, 0] \subset \mathbb{R},\\
\left\{|z + 2/5| = 1/5, \: - 1/2 \leqslant Re(z) \leqslant - 2/5\right\} &\to \{- 11 \leqslant Re(z) \leqslant -5 - 2 \sqrt{5}, \: 0 \leqslant Im(z) \leqslant 2\},\\
\left\{|z - 2/5| = 1/5, \: 3 / 10 \leqslant Re(z) \leqslant 2/5 \right\} &\to \{- 11 \leqslant Re(z) \leqslant -5 - 2 \sqrt{5}, \: - 2 \leqslant Im(z) \leqslant 0\}.
\end{split}
\end{equation}}

\begin{figure}[hbtp]
\begin{center}
{{\small Lower arcs of $\partial \mathbb{F}_5$}\includegraphics[width=2.5in]{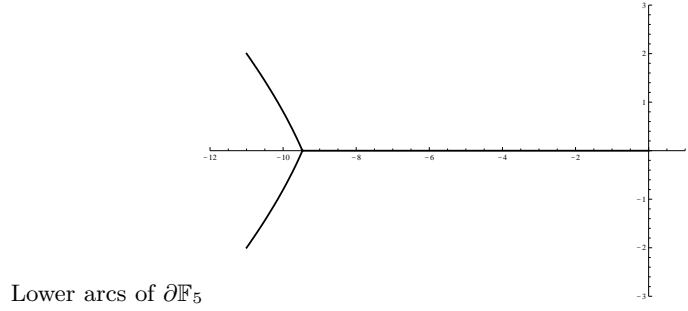}}
\end{center}
\caption{Image by $J_5$}\label{Im-J5B}
\end{figure}

%% file: report06-08.tex
\section{Level $6$}

We have $\Gamma_0(6)+$, $\Gamma_0(6)+6=\Gamma_0^{*}(6)$, $\Gamma_0(6)+3$, $\Gamma_0(6)+2$, and $\Gamma_0(6)-=\Gamma_0(6)$.

We have $W_6 = \left(\begin{smallmatrix}0 & - 1 / \sqrt{6}\\ \sqrt{6} & 0\end{smallmatrix}\right)$, $W_{6, 2} := \left(\begin{smallmatrix} -\sqrt{2} & -1/\sqrt{2}\\ 3\sqrt{2} & \sqrt{2}\end{smallmatrix}\right)$, and $W_{6, 3} := \left(\begin{smallmatrix} -\sqrt{3} & -2/\sqrt{3}\\ 2\sqrt{3} & \sqrt{3}\end{smallmatrix}\right)$, and we denote $\rho_{6, 1} := \rho_3 = - 1/2 + i / (2 \sqrt{3})$, $\rho_{6, 2} := - 1/3 + i / (3 \sqrt{2})$, $\rho_{6, 3} := - 2/5 + i / (5 \sqrt{6})$, $\rho_{6, 4} := - 1/4 + i / (4 \sqrt{3})$, and $\rho_{6, 5} := 1/3 + i / (3 \sqrt{2})$. We define
\begin{equation}
\begin{aligned}
\Delta_6^{\infty}(z) &:= \eta(z) \eta^{-2}(2 z) \eta^{-3}(3 z) \eta^6(6 z),\\
\Delta_6^{-1/2}(z) &:= \eta^{-3}(z) \eta^6(2 z) \eta(3 z) \eta^{-2}(6 z),\\
\Delta_6(z) &:= \Delta_6^{\infty}(z) \Delta_6^0(z) \Delta_6^{-1/3}(z) \Delta_6^{-1/2}(z)
\end{aligned}
\begin{aligned}
\Delta_6^0(z) &:= \eta^6(z) \eta^{-3}(2 z) \eta^{-2}(3 z) \eta(6 z), \\
\Delta_6^{-1/3}(z) &:= \eta^{-2}(z) \eta(2 z) \eta^6(3 z) \eta^{-3}(6 z),\\
= \eta^2(z) \eta^2(2 z) &\eta^2(3 z) \eta^2(6 z),
\end{aligned}
\end{equation}
where $\Delta_6^{\infty}$, $\Delta_6^0$, $\Delta_6^{-1/2}$, and $\Delta_6^{-1/3}$ are $2$nd semimodular forms for $\Gamma_0(6)$ of weight $1$ such that $v_{\infty}(\Delta_6^{\infty}) = v_0(\Delta_6^0) = v_{-1/2}(\Delta_6^{-1/2}) = v_{-1/3}(\Delta_6^{-1/3}) = 1$, and $\Delta_6$ is a cusp form for $\Gamma_0(6)$ of weight $4$. Furthermore, we define
\begin{equation}
\begin{split}
{E_{2, 6+6}}'(z) &:= (6 E_2(6 z) - 3 E_2(3 z) - 2 E_2(2 z) + E_2(z)) / 2,\\
{E_{2, 6+3}}'(z) &:= (6 E_2(6 z) - 3 E_2(3 z) + 2 E_2(2 z) - E_2(z)) / 4,\\
{E_{2, 6+2}}'(z) &:= (6 E_2(6 z) + 3 E_2(3 z) - 2 E_2(2 z) - E_2(z)) / 6,
\end{split}
\end{equation}
which are modular forms for $\Gamma_0(6)$ of weight $2$, and we have $v_{i / \sqrt{6}}({E_{2, 6+6}}') = v_{\rho_{6, 3}}({E_{2, 6+6}}') = 1$, $v_{\rho_{6, 1}}({E_{2, 6+3}}') = v_{\rho_{6, 4}}({E_{2, 6+3}}') = 1$, and $v_{\rho_{6, 2}}({E_{2, 6+2}}') = v_{\rho_{6, 5}}({E_{2, 6+2}}') = 1$.\\

\subsection{$\Gamma_0(6)+$}\quad

We have
\begin{equation*}
\Gamma_0(6)+ = \Gamma_0(6)+2,3,6 = \Gamma_0(6) \cup \Gamma_0(6) W_{6, 2} \cup \Gamma_0(6) W_{6, 3} \cup \Gamma_0(6) W_6.
\end{equation*}

\paragraph{\bf Fundamental domain}
We have a fundamental domain for $\Gamma_0(6)+$ as follows:
{\small \begin{equation}
\begin{split}
\mathbb{F}_{6+} = &\left\{|z + 1/2| \geqslant 1 / (2 \sqrt{3}), \: - 1/2 \leqslant Re(z) < - 1/3 \right\}
 \bigcup \left\{|z| \geqslant 1 / \sqrt{6}, \: - 1/3 \leqslant Re(z) \leqslant 0 \right\}\\
 &\bigcup \left\{|z| > 1 / \sqrt{6}, \: 0 < Re(z) \leqslant 1/3 \right\}
 \bigcup \left\{|z - 1/2| > 1 / (2 \sqrt{3}), \: 1/3 < Re(z) < 1/2 \right\},
\end{split}
\end{equation}
}where, $W_6 :  e^{i \theta} / \sqrt{6} \rightarrow e^{i  (\pi - \theta)} / \sqrt{6}$ and $\left(\begin{smallmatrix} -5 & -3 \\ 12 & 7 \end{smallmatrix}\right) W_{6, 3} : e^{i \theta} / (2 \sqrt{3}) + 1/2 \rightarrow e^{i  (\pi - \theta)} / (2 \sqrt{3}) - 1/2$. Then, we have
\begin{equation}
\Gamma_0(6)+ = \langle \left( \begin{smallmatrix} 1 & 1 \\ 0 & 1 \end{smallmatrix} \right), \: W_6, \: W_{6, 3} \rangle.
\end{equation}
\begin{figure}[hbtp]
\begin{center}
{{$\mathbb{F}_{6+}$}\includegraphics[width=1.5in]{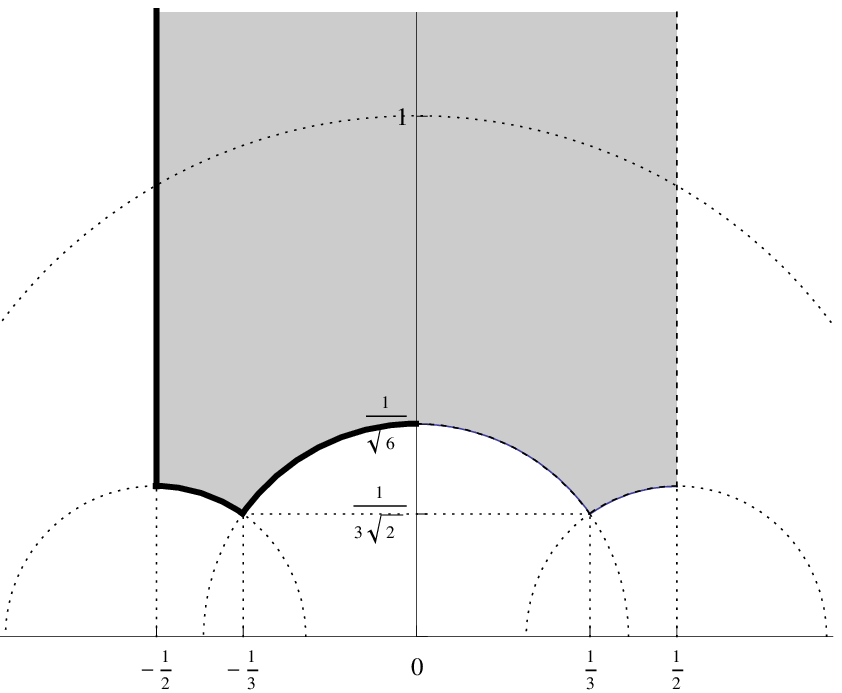}}
\end{center}
\caption{$\Gamma_0(6)+$}
\end{figure}

\paragraph{\bf Valence formula}
The cusp of $\Gamma_0(6)+$ is $\infty$, and the elliptic points are $i / \sqrt{6}$, $\rho_{6, 1} = - 1/2 + i / (2 \sqrt{3})$, and $\rho_{6, 2} = - 1/3 + i / (3 \sqrt{2})$. Let $f$ be a modular function of weight $k$ for $\Gamma_0(6)+$, which is not identically zero. We have
\begin{equation}
v_{\infty}(f) + \frac{1}{2} v_{i / \sqrt{6}}(f) + \frac{1}{2} v_{\rho_{6, 1}} (f) + \frac{1}{2} v_{\rho_{6, 2}} (f) + \sum_{\begin{subarray}{c} p \in \Gamma_0(6)+ \setminus \mathbb{H} \\ p \ne i / \sqrt{6}, \; \rho_{6, 1}, \; \rho_{6, 2}\end{subarray}} v_p(f) = \frac{k}{4}.
\end{equation}

Furthermore, the stabilizer of the elliptic point $i / \sqrt{6}$ (resp. $\rho_{6, 1}$, $\rho_{6, 2}$) is $\left\{ \pm I, \pm W_6 \right\}$\\ (resp. $\left\{ \pm I, \pm W_{6, 3} \right\}$, $\left\{ \pm I, \pm W_{6, 2} \right\}$)\\

\paragraph{\bf For the cusp $\infty$}
We have $\Gamma_{\infty} = \left\{ \pm \left( \begin{smallmatrix}1 & n \\ 0 & 1\end{smallmatrix} \right) \: ; \: n \in \mathbb{Z} \right\}$, and we have the Eisenstein series associated with $\Gamma_0(6)+$:
\begin{equation}
E_{k, 6+}(z) := \frac{6^{k/2} E_k(6 z) + 3^{k/2} E_k(3 z) + 2^{k/2} E_k(2 z) + E_k(z)}{(3^{k/2} + 1) (2^{k/2} + 1)} \quad \text{for} \; k \geqslant 4.
\end{equation}\quad

\paragraph{\bf The space of modular forms}
We have $M_k(\Gamma_0(6)+) = \mathbb{C} E_{k, 6+} \oplus S_k(\Gamma_0(6)+)$ and $S_k(\Gamma_0(6)+) = \Delta_6 M_{k - 4}(\Gamma_0(6)+)$ for every even integer $k \geqslant 4$. Then, we have $M_{4 n + 6}(\Gamma_0(6)+) = E_{6, 6+} M_{4 n}(\Gamma_0(6)+)$ and
\begin{equation*}
M_{4 n}(\Gamma_0(6)+) = \mathbb{C} (({E_{2, 6+2}}')^2)^n \oplus \mathbb{C} (({E_{2, 6+2}}')^2)^{n-1} \Delta_6 \oplus \cdots \oplus \mathbb{C} (\Delta_6)^n.
\end{equation*}\quad

\paragraph{\bf Hauptmodul}
We define the {\it hauptmodul} of $\Gamma_0(6)+$:
\begin{equation}
J_{6+} := ({E_{2, 6+2}}')^2 / \Delta_6 = \frac{1}{q} + 10 + 79 q + 352 q^2 + 1431 q^3 + \cdots,
\end{equation}
where $v_{\infty}(J_{6+}) = -1$ and $v_{\rho_{6, 2}}(J_{6+}) = 2$. Then, we have
\begin{equation}
J_{6+} : \partial \mathbb{F}_{6+} \setminus \{z \in \mathbb{H} \: ; \: Re(z) = \pm 1/2\} \to [- 4, 32] \subset \mathbb{R}.
\end{equation}\quad

\subsection{$\Gamma_0(6)+6 = \Gamma_0^{*}(6)$}

\paragraph{\bf Fundamental domain}
We have a fundamental domain for $\Gamma_0^{*}(6)$ as follows:
{\small \begin{equation}
\begin{split}
\mathbb{F}_{6+6} = &\left\{|z + 5 / 12| \geqslant 1 / 12, \: - 1/2 \leqslant Re(z) < - 2/5 \right\}
 \bigcup \left\{|z| \geqslant 1 / \sqrt{6}, \: - 2/5 \leqslant Re(z) \leqslant 0 \right\}\\
 &\bigcup \left\{|z| > 1 / \sqrt{6}, \: 0 < Re(z) \leqslant 2/5 \right\}
 \bigcup \left\{|z - 5 / 12| > 1 / 12, \: 2/5 < Re(z) < 1/2 \right\},
\end{split}
\end{equation}
}where $W_6 :  e^{i \theta} / \sqrt{6} \rightarrow e^{i  (\pi - \theta)} / \sqrt{6}$ and $\left(\begin{smallmatrix} -5 & 2 \\ 12 & -5 \end{smallmatrix} \right) : (e^{i \theta} + 5) / 12 \rightarrow (e^{i (\pi - \theta)} - 5) / 12$. Then, we have
\begin{equation}
\Gamma_0^{*}(6) = \langle \left( \begin{smallmatrix} 1 & 1 \\ 0 & 1 \end{smallmatrix} \right), \: W_6, \: \left( \begin{smallmatrix} 5 & 2 \\ 12 & 5 \end{smallmatrix} \right) \rangle.
\end{equation}
\begin{figure}[hbtp]
\begin{center}
\includegraphics[width=1.5in]{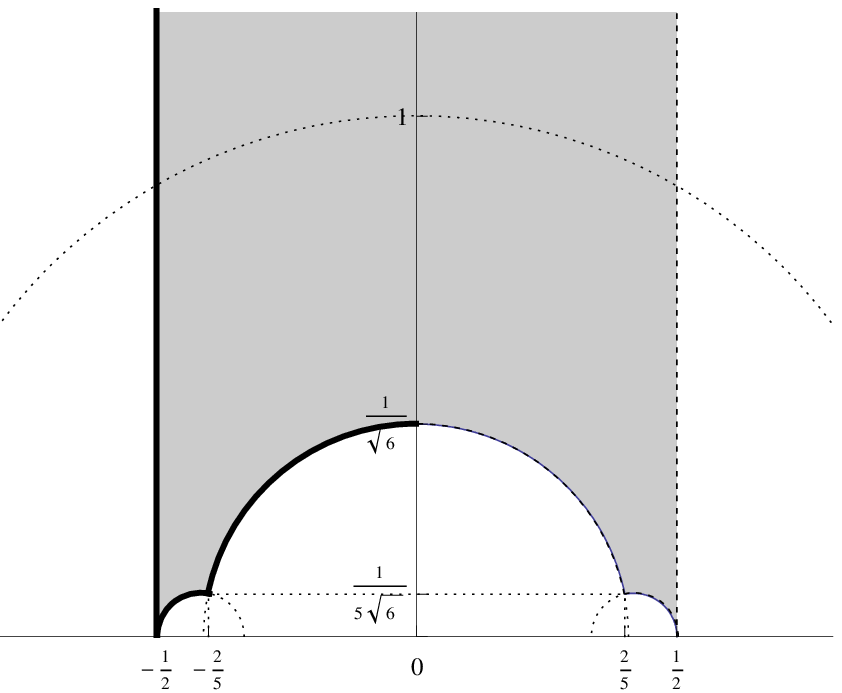}
\end{center}
\caption{$\Gamma_0^{*}(6)$}
\end{figure}

\paragraph{\bf Valence formula}
The cusps of $\Gamma_0^{*}(6)$ are $\infty$ and $-1/2$, and the elliptic points are $i / \sqrt{6}$ and $\rho_{6, 3} = - 2/5 + i / (5 \sqrt{6})$. Let $f$ be a modular function of weight $k$ for $\Gamma_0^{*}(6)$, which is not identically zero. We have
\begin{equation}
v_{\infty}(f) + v_{-1/2}(f) + \frac{1}{2} v_{i / \sqrt{6}} (f) + \frac{1}{2} v_{\rho_{6, 3}} (f) + \sum_{\begin{subarray}{c} p \in \Gamma_0^{*}(6) \setminus \mathbb{H} \\ p \ne i / \sqrt{6}, \rho_{6, 3}\end{subarray}} v_p(f) = \frac{k}{2}.
\end{equation}

Furthermore, the stabilizer of the elliptic point $i / \sqrt{6}$ (resp. $\rho_{6, 3}$) is $\left\{ \pm I, \pm W_6 \right\}$ (resp. $\left\{ \pm I, \pm \left( \begin{smallmatrix}5 & - 2 \\ - 12 & 5\end{smallmatrix} \right) W_6 \right\}$).\\

\paragraph{\bf For the cusp $\infty$}
We have $\Gamma_{\infty} = \left\{ \pm \left( \begin{smallmatrix}1 & n \\ 0 & 1\end{smallmatrix} \right) \: ; \: n \in \mathbb{Z} \right\}$, and we have the Eisenstein series for the cusp $\infty$ associated with $\Gamma_0^{*}(6)$:
\begin{equation}
E_{k, 6+6}^{\infty}(z) := \frac{(6^{k/2} + 1) (6^{k/2} E_k(6 z) + E_k(z)) - (3^{k/2} + 2^{k/2}) (3^{k/2} E_k(3 z) + 2^{k/2} E_k(2 z))}{(3^k - 1)(2^k - 1)} \quad \text{for} \; k \geqslant 4.
\end{equation}\quad

\paragraph{\bf For the cusp $-1/2$}
We have $\Gamma_{-1/2} = \left\{ \pm \left( \begin{smallmatrix}6 n + 1 & 3 n \\ - 12 n & - 6 n + 1\end{smallmatrix} \right) \: ; \: n \in \mathbb{Z} \right\}$ and $\gamma_{-1/2} = W_{6, 3}$, and we have the Eisenstein series for the cusp $-1/2$ associated with $\Gamma_0^{*}(6)$:
\begin{equation}
E_{k, 6+6}^{-1/2}(z) := \frac{- (3^{k/2} + 2^{k/2}) (6^{k/2} E_k(6 z) + E_k(z)) + (6^{k/2} + 1) (3^{k/2} E_k(3 z) + 2^{k/2} E_k(2 z))}{(3^k - 1)(2^k - 1)} \quad \text{for} \; k \geqslant 4.
\end{equation}
We also have $\gamma_{-1/2}^{-1} \ \Gamma_0^{*}(6) \ \gamma_{-1/2} = \Gamma_0^{*}(6)$.\\

\paragraph{\bf The space of modular forms}
We define
\begin{equation*}
\Delta_{6+6}^{\infty} := \Delta_6^{\infty} \Delta_6^0, \quad \Delta_{6+6}^{-1/2} := \Delta_6^{-1/2} \Delta_6^{-1/3},
\end{equation*}
which are $2$nd semimodular forms for $\Gamma_0^{*}(6)$ of weight $2$.

Now, we have $M_k(\Gamma_0^{*}(6)) = \mathbb{C} E_{k, 6+6}^{\infty} \oplus \mathbb{C} E_{k, 6+6}^{-1/2} \oplus S_k(\Gamma_0^{*}(6))$ and $S_k(\Gamma_0^{*}(6)) = \Delta_6 M_{k - 4}(\Gamma_0^{*}(6))$ for every even integer $k \geqslant 4$. Then, we have $M_{4 n + 2}(\Gamma_0^{*}(6)) = {E_{2, 6+6}}' M_{4 n}(\Gamma_0^{*}(6))$ and
\begin{align*}
M_{4 n}(\Gamma_0^{*}(6)) &= \mathbb{C} (E_{4, 6+6}^{\infty})^n \oplus \mathbb{C} (E_{4, 6+6}^{\infty})^{n-1} \Delta_6 \oplus \cdots \oplus \mathbb{C} E_{4, 6+6}^{\infty} (\Delta_6)^{n-1}\\
 &\oplus \mathbb{C} (E_{4, 6+6}^{-1/2})^n \oplus \mathbb{C} (E_{4, 6+6}^{-1/2})^{n-1} \Delta_6 \oplus \cdots \oplus \mathbb{C} E_{4, 6+6}^{-1/2} (\Delta_6)^{n-1} \oplus \mathbb{C} (\Delta_6)^n.
\end{align*}

Here, we have $E_{4, 6+6}^{\infty} = (-13/5) \Delta_{6+6}^{\infty} \Delta_{6+6}^{-1/2} + (\Delta_{6+6}^{-1/2})^2$ and $(-5/13) E_{4, 6+6}^{-1/2} = (-5/13) (\Delta_{6+6}^{\infty})^2 + \Delta_{6+6}^{\infty} \Delta_{6+6}^{-1/2}$, then we can write
\begin{equation*}
M_{4 n}(\Gamma_0^{*}(6)) = \mathbb{C} (\Delta_{6+6}^{\infty})^{2n} \oplus \mathbb{C} (\Delta_{6+6}^{\infty})^{2n-1} \Delta_{6+6}^{-1/2} \oplus \cdots \oplus \mathbb{C} (\Delta_{6+6}^{-1/2})^{2n}.
\end{equation*}\quad

\paragraph{\bf Hauptmodul}
We define the {\it hauptmodul} of $\Gamma_0^{*}(6)$:
\begin{equation}
J_{6+6} := \Delta_{6+6}^{-1/2} / \Delta_{6+6}^{\infty} \: (= \eta^5(z) \eta^{-1}(2 z) \eta(3 z) \eta^{-5}(6 z)) = \frac{1}{q} + 12 + 78 q + 364 q^2 + 1365 q^3 + \cdots,
\end{equation}
where $v_{\infty}(J_{6+6}) = -1$ and $v_{-1/2}(J_{6+6}) = 1$. Then, we have
\begin{equation}
J_{6+6} : \partial \mathbb{F}_{6+6} \setminus \{z \in \mathbb{H} \: ; \: Re(z) = \pm 1/2\} \to [0, 17 + 12 \sqrt{2}] \subset \mathbb{R}.
\end{equation}\quad

\subsection{$\Gamma_0(6)+3$}

\paragraph{\bf Fundamental domain}
We have a fundamental domain for $\Gamma_0(6)+3$ as follows:
{\small \begin{equation}
\begin{split}
\mathbb{F}_{6+3} = &\left\{|z + 1/2| \geqslant 1 / (2 \sqrt{3}), \: - 1/2 \leqslant Re(z) < -1/4 \right\}
 \bigcup \left\{|z + 1/6| \geqslant 1/6, \: - 1/4 \leqslant Re(z) \leqslant 0 \right\}\\
 &\bigcup \left\{|z - 1/6| > 1/6, \: 0 < Re(z) \leqslant 1/4 \right\}
 \bigcup \left\{|z - 1/2| > 1 / (2 \sqrt{3}), \: 1/4 < Re(z) < 1/2 \right\},
\end{split}
\end{equation}
}where $\left(\begin{smallmatrix} -1 & 0 \\ 6 & -1 \end{smallmatrix}\right) : (e^{i \theta} + 1) / 6 \rightarrow (e^{i (\pi - \theta)} - 1) / 6$ and $\left(\begin{smallmatrix} -5 & -3 \\ 12 & 7 \end{smallmatrix}\right) W_{6, 3} : e^{i \theta} / (2 \sqrt{3}) + 1/2 \rightarrow e^{i  (\pi - \theta)} / (2 \sqrt{3}) - 1/2$. Then, we have
\begin{equation}
\Gamma_0(6)+3 = \langle \left( \begin{smallmatrix} 1 & 1 \\ 0 & 1 \end{smallmatrix} \right), \: \left( \begin{smallmatrix} 1 & 0 \\ 6 & 1 \end{smallmatrix} \right), \: W_{6, 3} \rangle.
\end{equation}
\begin{figure}[hbtp]
\begin{center}
\includegraphics[width=1.5in]{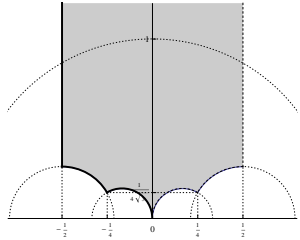}
\end{center}
\caption{$\Gamma_0(6)+3$}
\end{figure}

\paragraph{\bf Valence formula}
The cusps of $\Gamma_0(6)+3$ are $\infty$ and $0$, and the elliptic points are $\rho_{6, 1} = - 1/2 + i / (2 \sqrt{3})$ and $\rho_{6, 4} = - 1/4 + i / (4 \sqrt{3})$. Let $f$ be a modular function of weight $k$ for $\Gamma_0(6)+3$, which is not identically zero. We have
\begin{equation}
v_{\infty}(f) + v_0(f) + \frac{1}{2} v_{\rho_{6, 1}} (f) + \frac{1}{2} v_{\rho_{6, 4}} (f) + \sum_{\begin{subarray}{c} p \in \Gamma_0(6)+3 \setminus \mathbb{H} \\ p \ne \rho_{6, 1}, \rho_{6, 4}\end{subarray}} v_p(f) = \frac{k}{2}.
\end{equation}

Furthermore, the stabilizer of the elliptic point $\rho_{6, 1}$ (resp. $\rho_{6, 4}$) is $\left\{ \pm I, \pm W_{6, 3} \right\}$ (resp. $\left\{ \pm I, \pm \left( \begin{smallmatrix}-1 & -1 \\ 6 & 5\end{smallmatrix} \right) W_{6, 3} \right\}$).\\

\paragraph{\bf For the cusp $\infty$}
We have $\Gamma_{\infty} = \left\{ \pm \left( \begin{smallmatrix}1 & n \\ 0 & 1\end{smallmatrix} \right) \: ; \: n \in \mathbb{Z} \right\}$, and we have the Eisenstein series for the cusp $\infty$ associated with $\Gamma_0(6)+3$:
\begin{equation}
E_{k, 6+3}^{\infty}(z) := \frac{2^k 3^{k/2} E_k(6 z) - 3^{k/2} E_k(3 z) + 2^k E_k(2 z) - E_k(z)}{(3^{k/2} + 1)(2^k - 1)} \quad \text{for} \; k \geqslant 4.
\end{equation}\quad

\paragraph{\bf For the cusp $0$}
We have $\Gamma_0 = \left\{ \pm \left( \begin{smallmatrix}1 & 0 \\ 6 n & 1\end{smallmatrix} \right) \: ; \: n \in \mathbb{Z} \right\}$ and $\gamma_0 = W_6$, and we have the Eisenstein series for the cusp $0$ associated with $\Gamma_0(6)+3$:
\begin{equation}
E_{k, 6+3}^0(z) := \frac{- 2^{k/2} (3^{k/2} E_k(6 z) - 3^{k/2} E_k(3 z) + E_k(2 z) - E_k(z))}{(3^{k/2} + 1)(2^k - 1)} \quad \text{for} \; k \geqslant 4.
\end{equation}
We also have $\gamma_0^{-1} (\Gamma_0(6)+3) \gamma_0 = \Gamma_0(6)+3$.\\

\paragraph{\bf The space of modular forms}
We define
\begin{equation*}
\Delta_{6+3}^{\infty} := \Delta_6^{\infty} \Delta_6^{-1/2}, \quad \Delta_{6+3}^0 := \Delta_6^0 \Delta_6^{-1/3},
\end{equation*}
which are $2$nd semimodular forms for $\Gamma_0(6)+3$ of weight $2$.

Now, we have $M_k(\Gamma_0(6)+3) = \mathbb{C} E_{k, 6+3}^{\infty} \oplus \mathbb{C} E_{k, 6+3}^0 \oplus S_k(\Gamma_0(6)+3)$ and $S_k(\Gamma_0(6)+3) = \Delta_6 M_{k - 4}(\Gamma_0(6)+3)$ for every even integer $k \geqslant 4$. Then, we have $M_{4 n + 2}(\Gamma_0(6)+3) = {E_{2, 6+3}}' M_{4 n}(\Gamma_0(6)+3)$ and
\begin{align*}
M_{4 n}(\Gamma_0(6)+3) &= \mathbb{C} (E_{4, 6+3}^{\infty})^n \oplus \mathbb{C} (E_{4, 6+3}^{\infty})^{n-1} \Delta_6 \oplus \cdots \oplus \mathbb{C} E_{4, 6+3}^{\infty} (\Delta_6)^{n-1}\\
 &\oplus \mathbb{C} (E_{4, 6+3}^0)^n \oplus \mathbb{C} (E_{4, 6+3}^0)^{n-1} \Delta_6 \oplus \cdots \oplus \mathbb{C} E_{4, 6+3}^0 (\Delta_6)^{n-1} \oplus \mathbb{C} (\Delta_6)^n.
\end{align*}

Here, we have $E_{4, 6+3}^{\infty} = (32/5) \Delta_{6+3}^{\infty} \Delta_{6+3}^0 + (\Delta_{6+3}^0)^2$ and $(5/32) E_{4, 6+3}^0 = 10 (\Delta_{6+3}^{\infty})^2 + \Delta_{6+3}^{\infty} \Delta_{6+3}^0$, then we can write
\begin{equation*}
M_{4 n}(\Gamma_0(6)+3) = \mathbb{C} (\Delta_{6+3}^{\infty})^{2n} \oplus \mathbb{C} (\Delta_{6+3}^{\infty})^{2n-1} \Delta_{6+3}^0 \oplus \cdots \oplus \mathbb{C} (\Delta_{6+3}^0)^{2n}.
\end{equation*}\quad

\paragraph{\bf Hauptmodul}
We define the {\it hauptmodul} of $\Gamma_0(6)+3$:
\begin{equation}
J_{6+3} := \Delta_{6+3}^0 / \Delta_{6+3}^{\infty} \: (= \eta^6(z) \eta^{-6}(2 z) \eta^6(3 z) \eta^{-6}(6 z)) = \frac{1}{q} - 6 + 15 q - 32 q^2 + 87 q^3 - \cdots,
\end{equation}
where $v_{\infty}(J_{6+3}) = -1$ and $v_0(J_{6+3}) = 1$. Then, we have
\begin{equation}
J_{6+3} : \partial \mathbb{F}_{6+3} \setminus \{z \in \mathbb{H} \: ; \: Re(z) = \pm 1/2\} \to [-16, 0] \subset \mathbb{R}.
\end{equation}\quad

\subsection{$\Gamma_0(6)+2$}

\paragraph{\bf Fundamental domain}
We have a fundamental domain for $\Gamma_0(6)+2$ as follows:
{\small \begin{equation}
\begin{split}
\mathbb{F}_{6+2} = &\left\{|z + 1/3| \geqslant 1 / (3 \sqrt{2}), \: - 1/2 \leqslant Re(z) \leqslant -1/3 \right\}
 \bigcup \left\{|z + 1/3| > 1 / (3 \sqrt{2}), \: - 1/3 < Re(z) < -1/6 \right\}\\
 &\bigcup \left\{|z + 1/6| \geqslant 1/6, \: 0 \leqslant Re(z) \leqslant 0 \right\}
 \bigcup \left\{|z - 1/6| > 1/6, \: 0 < Re(z) < 1/6 \right\}\\
 &\bigcup \left\{|z - 1/3| \geqslant 1 / (3 \sqrt{2}), \: 1/6 \leqslant Re(z) \leqslant 1/3 \right\}
 \bigcup \left\{|z - 1/3| > 1 / (3 \sqrt{2}), \: 1/3 < Re(z) < 1/2 \right\},
\end{split}
\end{equation}
}where $\left(\begin{smallmatrix} -1 & 0 \\ 6 & -1 \end{smallmatrix}\right) : (e^{i \theta} + 1) / 6 \rightarrow (e^{i (\pi - \theta)} - 1) / 6$, $\left(\begin{smallmatrix} 5 & 2 \\ 12 & 5 \end{smallmatrix}\right) W_{6, 3} : e^{i \theta} / (3 \sqrt{2}) + 1/3 \rightarrow e^{i  (\pi - \theta)} / (3 \sqrt{2}) + 1/3$, and $W_{6, 3} : e^{i \theta} / (3 \sqrt{2}) - 1/3 \rightarrow e^{i  (\pi - \theta)} / (3 \sqrt{2}) - 1/3$. Then, we have
\begin{equation}
\Gamma_0(6)+2 = \langle \left( \begin{smallmatrix} 1 & 1 \\ 0 & 1 \end{smallmatrix} \right), \: \left( \begin{smallmatrix} 1 & 0 \\ 6 & 1 \end{smallmatrix} \right), \: W_{6, 2} \rangle.
\end{equation}
\begin{figure}[hbtp]
\begin{center}
\includegraphics[width=1.5in]{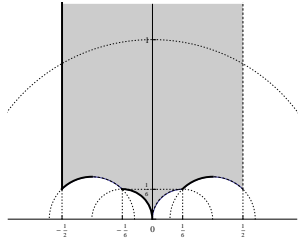}
\end{center}
\caption{$\Gamma_0(6)+2$}
\end{figure}

\paragraph{\bf Valence formula}
The cusps of $\Gamma_0(6)+2$ are $\infty$ and $0$, and the elliptic points are $\rho_{6, 2} = - 1/3 + i / (3 \sqrt{2})$ and $\rho_{6, 5} = 1/3 + i / (3 \sqrt{2})$. Let $f$ be a modular function of weight $k$ for $\Gamma_0(6)+2$, which is not identically zero. We have
\begin{equation}
v_{\infty}(f) + v_0(f) + \frac{1}{2} v_{\rho_{6, 2}} (f) + \frac{1}{2} v_{\rho_{6, 5}} (f) + \sum_{\begin{subarray}{c} p \in \Gamma_0(6)+2 \setminus \mathbb{H} \\ p \ne \rho_{6, 2}, \rho_{6, 5}\end{subarray}} v_p(f) = \frac{k}{2}.
\end{equation}

Furthermore, the stabilizer of the elliptic point $\rho_{6, 2}$ (resp. $\rho_{6, 5}$) is $\left\{ \pm I, \pm W_{6, 2} \right\}$ (resp. $\left\{ \pm I, \pm \left( \begin{smallmatrix}5 & 2 \\ 12 & 5\end{smallmatrix} \right) W_{6, 2} \right\}$).\\

\paragraph{\bf For the cusp $\infty$}
We have $\Gamma_{\infty} = \left\{ \pm \left( \begin{smallmatrix}1 & n \\ 0 & 1\end{smallmatrix} \right) \: ; \: n \in \mathbb{Z} \right\}$, and we have the Eisenstein series for the cusp $\infty$ associated with $\Gamma_0(6)+2$:
\begin{equation}
E_{k, 6+2}^{\infty}(z) := \frac{2^{k/2} 3^k E_k(6 z) + 3^k E_k(3 z) - 2^{k/2} E_k(2 z) - E_k(z)}{(3^k - 1)(2^{k/2} + 1)} \quad \text{for} \; k \geqslant 4.
\end{equation}\quad

\paragraph{\bf For the cusp $0$}
We have $\Gamma_0 = \left\{ \pm \left( \begin{smallmatrix}1 & 0 \\ 6 n & 1\end{smallmatrix} \right) \: ; \: n \in \mathbb{Z} \right\}$ and $\gamma_0 = W_6$, and we have the Eisenstein series for the cusp $0$ associated with $\Gamma_0(6)+2$:
\begin{equation}
E_{k, 6+2}^0(z) := \frac{- 3^{k/2} (2^{k/2} E_k(6 z) + E_k(3 z) - 2^{k/2} E_k(2 z) - E_k(z))}{(3^k - 1)(2^{k/2} + 1)} \quad \text{for} \; k \geqslant 4.
\end{equation}
We also have $\gamma_0^{-1} (\Gamma_0(6)+2) \gamma_0 = \Gamma_0(6)+2$.\\

\paragraph{\bf The space of modular forms}
We define
\begin{equation*}
\Delta_{6+2}^{\infty} := \Delta_6^{\infty} \Delta_6^{-1/3}, \quad \Delta_{6+2}^0 := \Delta_6^0 \Delta_6^{-1/2},
\end{equation*}
which are $2$nd semimodular forms for $\Gamma_0(6)+2$ of weight $2$.

Now, we have $M_k(\Gamma_0(6)+2) = \mathbb{C} E_{k, 6+2}^{\infty} \oplus \mathbb{C} E_{k, 6+2}^0 \oplus S_k(\Gamma_0(6)+2)$ and $S_k(\Gamma_0(6)+2) = \Delta_6 M_{k - 4}(\Gamma_0(6)+2)$ for every even integer $k \geqslant 4$. Then, we have $M_{4 n + 2}(\Gamma_0(6)+2) = {E_{2, 6+2}}' M_{4 n}(\Gamma_0(6)+2)$ and
\begin{align*}
M_{4 n}(\Gamma_0(6)+2) &= \mathbb{C} (E_{4, 6+2}^{\infty})^n \oplus \mathbb{C} (E_{4, 6+2}^{\infty})^{n-1} \Delta_6 \oplus \cdots \oplus \mathbb{C} E_{4, 6+2}^{\infty} (\Delta_6)^{n-1}\\
 &\oplus \mathbb{C} (E_{4, 6+2}^0)^n \oplus \mathbb{C} (E_{4, 6+2}^0)^{n-1} \Delta_6 \oplus \cdots \oplus \mathbb{C} E_{4, 6+2}^0 (\Delta_6)^{n-1} \oplus \mathbb{C} (\Delta_6)^n.
\end{align*}

Here, we have $E_{4, 6+2}^{\infty} = (27/5) \Delta_{6+2}^{\infty} \Delta_{6+2}^0 + (\Delta_{6+2}^0)^2$ and $(5/27) E_{4, 6+2}^0 = 7 (\Delta_{6+2}^{\infty})^2 + \Delta_{6+2}^{\infty} \Delta_{6+2}^0$, then we can write
\begin{equation*}
M_{4 n}(\Gamma_0(6)+2) = \mathbb{C} (\Delta_{6+2}^{\infty})^{2n} \oplus \mathbb{C} (\Delta_{6+2}^{\infty})^{2n-1} \Delta_{6+2}^0 \oplus \cdots \oplus \mathbb{C} (\Delta_{6+2}^0)^{2n}.
\end{equation*}\quad

\paragraph{\bf Hauptmodul}
We define the {\it hauptmodul} of $\Gamma_0(6)+2$:
\begin{equation}
J_{6+2} := \Delta_{6+2}^0 / \Delta_{6+2}^{\infty} \: (= \eta^4(z) \eta^4(2 z) \eta^{-4}(3 z) \eta^{-4}(6 z)) = \frac{1}{q} - 4 - 2 q + 28 q^2 - 27 q^3 - \cdots,
\end{equation}
where $v_{\infty}(J_{6+2}) = -1$ and $v_0(J_{6+2}) = 1$. Then, we have
{\small \begin{equation}
\begin{split}
J_{6+2} : \qquad \qquad \left\{|z + 1/6| = 1/6, \: - 1/6 \leqslant Re(z) \leqslant 0\right\} &\to [-3, 0] \subset \mathbb{R},\\
\left\{|z + 1/3| = 1 / (3 \sqrt{2}), \: -1/2 \leqslant Re(z) \leqslant -1/3\right\} &\to \{-7 \leqslant Re(z) \leqslant -3, \: 0 \leqslant Im(z) \leqslant 4 \sqrt{2}\},\\
\left\{|z - 1/3| = 1 / (3 \sqrt{2}), \: 1/6 \leqslant Re(z) \leqslant 1/3 \right\} &\to \{-7 \leqslant Re(z) \leqslant -3, \: - 4 \sqrt{2} \leqslant Im(z) \leqslant 0\}.
\end{split}
\end{equation}
}Thus, $J_{6+2}$ does not take real value on some arcs of $\partial \mathbb{F}_{6+2}$.

\begin{figure}[hbtp]
\begin{center}
{{\small Lower arcs of $\partial \mathbb{F}_{6+2}$}\includegraphics[width=2.5in]{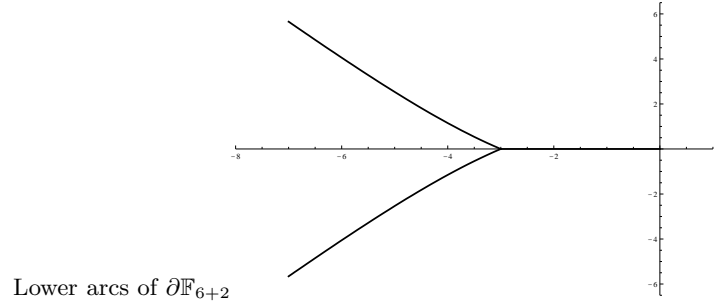}}
\end{center}
\caption{Image by $J_{6+2}$}\label{Im-J6D}
\end{figure}

\subsection{$\Gamma_0(6)$}

\paragraph{\bf Fundamental domain}
We have a fundamental domain for $\Gamma_0(6)$ as follows:
{\small \begin{equation}
\begin{split}
\mathbb{F}_6 = &\left\{|z + 5/12| \geqslant 1/12, \: - 1/2 \leqslant Re(z) < -1/3 \right\}
 \bigcup \left\{|z + 1/6| \geqslant 1/6, \: - 1/3 \leqslant Re(z) \leqslant 0 \right\}\\
 &\bigcup \left\{|z - 1/6| > 1/6, \: 0 < Re(z) \leqslant 1/3 \right\}
 \bigcup \left\{|z - 5/12| > 1/12, \: 1/3 < Re(z) < 1/2 \right\},
\end{split}
\end{equation}
}where $\left(\begin{smallmatrix} -1 & 0 \\ 6 & -1 \end{smallmatrix}\right) : (e^{i \theta} + 1) / 6 \rightarrow (e^{i (\pi - \theta)} - 1) / 6$ and $\left(\begin{smallmatrix} -5 & 2 \\ 12 & -5 \end{smallmatrix} \right) : (e^{i \theta} + 5) / 12 \rightarrow (e^{i (\pi - \theta)} - 5) / 12$. Then, we have
\begin{equation}
\Gamma_0(6) = \langle - I, \: \left( \begin{smallmatrix} 1 & 1 \\ 0 & 1 \end{smallmatrix} \right), \: \left( \begin{smallmatrix} 1 & 0 \\ 6 & 1 \end{smallmatrix} \right), \: \left( \begin{smallmatrix} 5 & 2 \\ 12 & 5 \end{smallmatrix} \right) \rangle.
\end{equation}
\begin{figure}[hbtp]
\begin{center}
\includegraphics[width=1.5in]{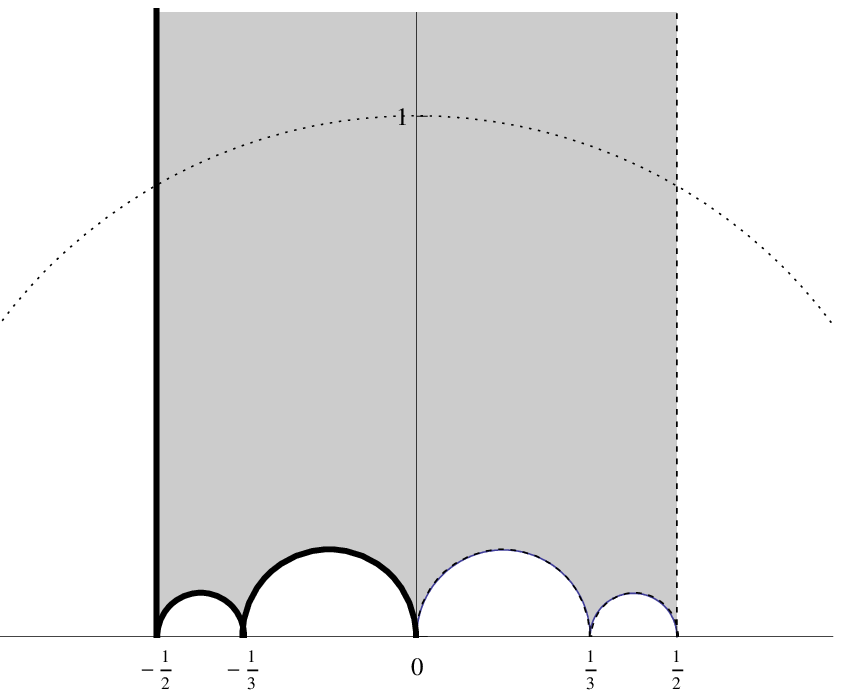}
\end{center}
\caption{$\Gamma_0(6)$}
\end{figure}

\paragraph{\bf Valence formula}
The cusps of $\Gamma_0(6)$ are $\infty$, $0$, $-1/2$, and $-1/3$. Let $f$ be a modular function of weight $k$ for $\Gamma_0(6)$, which is not identically zero. We have
\begin{equation}
v_{\infty}(f) + v_0(f) + v_{-1/2} (f) + v_{-1/3} (f) + \sum_{p \in \Gamma_0(6) \setminus \mathbb{H}} v_p(f) = k.
\end{equation}\quad

\paragraph{\bf For the cusp $\infty$}
We have $\Gamma_{\infty} = \left\{ \pm \left( \begin{smallmatrix}1 & n \\ 0 & 1\end{smallmatrix} \right) \: ; \: n \in \mathbb{Z} \right\}$, and we have the Eisenstein series for the cusp $\infty$ associated with $\Gamma_0(6)$:
\begin{equation}
E_{k, 6}^{\infty}(z) := \frac{6^k E_k(6 z) - 3^k E_k(3 z) - 2^k E_k(2 z) + E_k(z)}{(3^k - 1)(2^k - 1)} \quad \text{for} \; k \geqslant 4.
\end{equation}\quad

\paragraph{\bf For the cusp $0$}
We have $\Gamma_0 = \left\{ \pm \left( \begin{smallmatrix}1 & 0 \\ 6 n & 1\end{smallmatrix} \right) \: ; \: n \in \mathbb{Z} \right\}$ and $\gamma_0 = W_6$, and we have the Eisenstein series for the cusp $0$ associated with $\Gamma_0(6)$:
\begin{equation}
E_{k, 6}^0(z) := \frac{6^{k/2} (E_k(6 z) - E_k(3 z) - E_k(2 z) + E_k(z))}{(3^k - 1)(2^k - 1)} \quad \text{for} \; k \geqslant 4.
\end{equation}
We also have $\gamma_0^{-1} (\Gamma_0(6)) \gamma_0 = \Gamma_0(6)$.\\

\paragraph{\bf For the cusp $-1/2$}
We have $\Gamma_{-1/2} = \left\{ \pm \left( \begin{smallmatrix}6 n + 1 & 3 n \\ - 12 n & - 6 n + 1\end{smallmatrix} \right) \: ; \: n \in \mathbb{Z} \right\}$ and $\gamma_{-1/2} = W_{6, 3}$, and we have the Eisenstein series for the cusp $-1/2$ associated with $\Gamma_0(6)$:
\begin{equation}
E_{k, 6}^{-1/2}(z) := \frac{- 3^{k/2} (2^k E_k(6 z) - E_k(3 z) - 2^k E_k(2 z) + E_k(z))}{(3^k - 1)(2^k - 1)} \quad \text{for} \; k \geqslant 4.
\end{equation}
We also have $\gamma_{-1/2}^{-1} (\Gamma_0(6)) \gamma_{-1/2} = \Gamma_0(6)$.\\

\paragraph{\bf For the cusp $-1/3$}
We have $\Gamma_{-1/3} = \left\{ \pm \left( \begin{smallmatrix}6 n + 1 & 2 n \\ - 18 n & - 6 n + 1\end{smallmatrix} \right) \: ; \: n \in \mathbb{Z} \right\}$ and $\gamma_{-1/3} = W_{6, 2}$, and we have the Eisenstein series for the cusp $-1/3$ associated with $\Gamma_0(6)$:
\begin{equation}
E_{k, 6}^{-1/3}(z) := \frac{- 2^{k/2} (3^k E_k(6 z) - 3^k E_k(3 z) - E_k(2 z) + E_k(z))}{(3^k - 1)(2^k - 1)} \quad \text{for} \; k \geqslant 4.
\end{equation}
We also have $\gamma_{-1/3}^{-1} (\Gamma_0(6)) \gamma_{-1/3} = \Gamma_0(6)$.\\

\paragraph{\bf The space of modular forms}
We have $M_k(\Gamma_0(6)) = \mathbb{C} E_{k, 6}^{\infty} \oplus \mathbb{C} E_{k, 6}^0 \oplus \mathbb{C} E_{k, 6}^{-1/2} \oplus \mathbb{C} E_{k, 6}^{-1/3} \oplus S_k(\Gamma_0(6))$ and $S_k(\Gamma_0(6)) = \Delta_6 M_{k - 4}(\Gamma_0(6))$ for every even integer $k \geqslant 4$. Then, we have $M_{4 n + 2}(\Gamma_0(6)) = {E_{2, 6+6}}' M_{4 n}(\Gamma_0(6)) \oplus \mathbb{C} {E_{2, 6+3}}' (\Delta_6)^n \oplus \mathbb{C} {E_{2, 6+2}}' (\Delta_6)^n$ and
\begin{align*}
M_{4 n}(\Gamma_0(6)) = &\mathbb{C} (E_{4, 6}^{\infty})^n \oplus \mathbb{C} (E_{4, 6}^{\infty})^{n-1} \Delta_6 \oplus \cdots \oplus \mathbb{C} E_{4, 6}^{\infty} (\Delta_6)^{n-1}\\
 &\oplus \mathbb{C} (E_{4, 6}^0)^n \oplus \mathbb{C} (E_{4, 6}^0)^{n-1} \Delta_6 \oplus \cdots \oplus \mathbb{C} E_{4, 6}^0 (\Delta_6)^{n-1}\\
 &\oplus \mathbb{C} (E_{4, 6}^{-1/2})^n \oplus \mathbb{C} (E_{4, 6}^{-1/2})^{n-1} \Delta_6 \oplus \cdots \oplus \mathbb{C} E_{4, 6}^{-1/2} (\Delta_6)^{n-1}\\
 &\oplus \mathbb{C} (E_{4, 6}^{-1/3})^n \oplus \mathbb{C} (E_{4, 6}^{-1/3})^{n-1} \Delta_6 \oplus \cdots \oplus \mathbb{C} E_{4, 6}^{-1/3} (\Delta_6)^{n-1} \oplus \mathbb{C} (\Delta_6)^n.
\end{align*}

Here, we have ${E_{2, 6+6}}' = -72 (\Delta_6^{\infty})^2 + (\Delta_6^0)^2$, ${E_{2, 6+3}}' = 72 (\Delta_6^{\infty})^2 + 18 \Delta_6^{\infty} \Delta_6^0 + (\Delta_6^0)^2$, ${E_{2, 6+2}}' = 72 (\Delta_6^{\infty})^2 + 16 \Delta_6^{\infty} \Delta_6^0 + (\Delta_6^0)^2$, $E_{4, 6}^{\infty} = (2592/5) (\Delta_6^{\infty})^3 \Delta_6^0 + (972/5) (\Delta_6^{\infty})^2 (\Delta_6^0)^2 + (121/5) \Delta_6^{\infty} (\Delta_6^0)^3 + (\Delta_6^0)^4$, $(5/36)E_{4, 6}^0 = 720 (\Delta_6^{\infty})^4 + 242 (\Delta_6^{\infty})^3 \Delta_6^0 + 27 (\Delta_6^{\infty})^2 (\Delta_6^0)^2 + \Delta_6^{\infty} (\Delta_6^0)^3$, $(-5/9) E_{4, 6}^{-1/2} = 32 (\Delta_6^{\infty})^3 \Delta_6^0 + 12 (\Delta_6^{\infty})^2 (\Delta_6^0)^2 + \Delta_6^{\infty} (\Delta_6^0)^3$, and $(-5/4) E_{4, 6}^{-1/3} = 162 (\Delta_6^{\infty})^3 \Delta_6^0 + 27 (\Delta_6^{\infty})^2 (\Delta_6^0)^2 + \Delta_6^{\infty} (\Delta_6^0)^3$. Now, we can write
\begin{equation*}
M_{2 n}(\Gamma_0(6)) = \mathbb{C} (\Delta_6^{\infty})^{2n} \oplus \mathbb{C} (\Delta_6^{\infty})^{2n-1} \Delta_6^0 \oplus \cdots \oplus \mathbb{C} (\Delta_6^0)^{2n}.
\end{equation*}\quad

\paragraph{\bf Hauptmodul}
We define the {\it hauptmodul} of $\Gamma_0(6)$:
\begin{equation}
J_6 := \Delta_{6}^0 / \Delta_{6}^{\infty} \: (= \eta^5(z) \eta^{-1}(2 z) \eta(3 z) \eta^{-5}(6 z)) = \frac{1}{q} - 5 + 6 q + 4 q^2 - 3 q^3 - \cdots,
\end{equation}
where $v_{\infty}(J_6) = -1$ and $v_0(J_6) = 1$. Then, we have
\begin{equation}
J_6 : \partial \mathbb{F}_6 \setminus \{z \in \mathbb{H} \: ; \: Re(z) = \pm 1/2\} \to [-9, 0] \subset \mathbb{R}.
\end{equation}

\clearpage

\section{Level $7$}

We have $\Gamma_0(7)+=\Gamma_0^{*}(7)$ and $\Gamma_0(7)-=\Gamma_0(7)$.

We have $W_7 = \left(\begin{smallmatrix}0&-1 / \sqrt{7}\\ \sqrt{7}&0\end{smallmatrix}\right)$, and denote $\rho_{7, 1} := - 1/2 + i / (2 \sqrt{7})$, $\rho_{7, 2} := - 5/14 + \sqrt{3} i / 14$, and $\rho_{7, 3} := 5/14 + \sqrt{3} i / 14$. We define
\begin{equation}
\begin{split}
&\Delta_7^{\infty}(z) := \sqrt{\eta^7(7 z) / \eta(z)}, \quad \Delta_7^0(z) := \sqrt{\eta^7(z) / \eta(7 z)},\\
&\Delta_7(z) := \Delta_7^{\infty}(z) \Delta_7^0(z) = \eta^3(z) \eta^3(7 z),
\end{split}
\end{equation}
where $\Delta_7^{\infty}$ and $\Delta_7^0$ are $4$th semimodular forms for $\Gamma_0(7)$ of weight $3/2$ such that $v_{\infty}(\Delta_7^{\infty}) = v_0(\Delta_7^0) = 1$ and $\Delta_7^{\infty} (\Delta_7^0)^3$ is a modular form for $\Gamma_0(7)$, and $\Delta_7$ is a cusp form for $\Gamma_0(7)$ and $2$nd semimodular form for $\Gamma_0^{*}(7)$ of weight $6$. Furthermore, we define
\begin{equation}
{E_{2, 7}}'(z) := (7 E_2(7 z) - E_2(z)) / 6, \quad {E_{1, 7}}' := \sqrt{{E_{2, 7}}'},
\end{equation}
where $\sqrt{ \ \cdot \ }$ is selected so that constant term of its Fourier expansion is positive. Here, ${E_{1, 7}}'$ is a $2$nd semimodular form for $\Gamma_0(7)$ and $4$th semimodular form for $\Gamma_0^{*}(7)$ such that $v_{\rho_{7, 2}}({E_{1, 7}}') = 1$, and ${E_{1, 7}}' \Delta_7$ is a modular form. Then, ${E_{2, 7}}'$ is a square of it.\\

\subsection{$\Gamma_0^{*}(7)$} (see \cite{SJ2})

\paragraph{\bf Fundamental domain}
We have a fundamental domain for $\Gamma_0^{*}(7)$ as follows:
{\small \begin{equation}
\begin{split}
\mathbb{F}_{7+} = &\left\{|z + 1/2| \geqslant 1 / (2 \sqrt{7}), \: - 1/2 \leqslant Re(z) < - 5/14 \right\}
 \bigcup \left\{|z| \geqslant 1 / \sqrt{7}, \: - 5/14 \leqslant Re(z) \leqslant 0 \right\}\\
 &\bigcup \left\{|z| > 1 / \sqrt{7}, \: 0 < Re(z) \leqslant 5/14 \right\}
 \bigcup \left\{|z - 1/2| > 1 / (2 \sqrt{7}), \: 5/14 < Re(z) < 1/2 \right\}.
\end{split}
\end{equation}
}where $W_7 :  e^{i \theta} / \sqrt{7} \rightarrow e^{i  (\pi - \theta)} / \sqrt{7}$ and $\left(\begin{smallmatrix} -3 & -1 \\ 7 & 2 \end{smallmatrix}\right) W_7 : e^{i \theta} / (2 \sqrt{7}) + 1/2 \rightarrow e^{i  (\pi - \theta)} / (2 \sqrt{7}) - 1/2$. Then, we have
\begin{equation}
\Gamma_0^{*}(7) = \langle \left( \begin{smallmatrix} 1 & 1 \\ 0 & 1 \end{smallmatrix} \right), \: W_7, \: \left( \begin{smallmatrix} 3 & 1 \\ 7 & 2 \end{smallmatrix} \right) \rangle.
\end{equation}
\begin{figure}[hbtp]
\begin{center}
{{$\mathbb{F}_{7+}$}\includegraphics[width=1.5in]{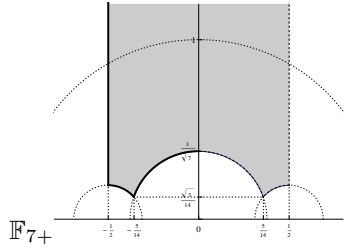}}
\end{center}
\caption{$\Gamma_0^{*}(7)$}
\end{figure}

\paragraph{\bf Valence formula}
The cusp of $\Gamma_0^{*}(7)$ is $\infty$, and the elliptic points are $i / \sqrt{7}$, $\rho_{7, 1} = - 1/2 + i \sqrt{7} / 10$ and $\rho_{7, 2} = - 5/14 + \sqrt{3} i / 14$. Let $f$ be a modular function of weight $k$ for $\Gamma_0^{*}(7)$, which is not identically zero. We have
\begin{equation}
v_{\infty}(f) + \frac{1}{2} v_{i / \sqrt{7}}(f) + \frac{1}{2} v_{\rho_{7, 1}} (f) + \frac{1}{3} v_{\rho_{7, 2}} (f) + \sum_{\begin{subarray}{c} p \in \Gamma_0^{*}(7) \setminus \mathbb{H} \\ p \ne i / \sqrt{7}, \; \rho_{7, 1}, \; \rho_{7, 2}\end{subarray}} v_p(f) = \frac{k}{3}.
\end{equation}

Furthermore, the stabilizer of the elliptic point $i / \sqrt{7}$ (resp. $\rho_{7, 1}$, $\rho_{7, 2}$) is $\left\{ \pm I, \pm W_7 \right\}$\\ (resp. $\left\{ \pm I, \pm \left( \begin{smallmatrix} 4 & -1 \\ -7 & 2 \end{smallmatrix} \right) W_7 \right\}$, $\left\{ \pm I, \pm \left( \begin{smallmatrix} -3 & -1 \\ 7 & 2 \end{smallmatrix} \right), \pm \left( \begin{smallmatrix} -2 & -1 \\ 7 & 3 \end{smallmatrix} \right) \right\}$)\\

\paragraph{\bf For the cusp $\infty$}
We have $\Gamma_{\infty} = \left\{ \pm \left( \begin{smallmatrix}1 & n \\ 0 & 1\end{smallmatrix} \right) \: ; \: n \in \mathbb{Z} \right\}$, and we have the Eisenstein series associated with $\Gamma_0^{*}(7)$:
\begin{equation}
E_{k, 7+}(z) := \frac{7^{k/2} E_k(7 z) + E_k(z)}{7^{k/2} + 1} \quad \text{for} \; k \geqslant 4.
\end{equation}\quad

\paragraph{\bf The space of modular forms}
We define the following functions:
\begin{gather*}
\Delta_{4, 7} := (5 / 16) (({E_{2, 7}}')^2 - E_{4, 7+}),\\
\Delta_{10, 0, 7+} := (559 / 690) ((41065 / 137592) (E_{4, 7+} E_{6, 7+} - E_{10, 7+}) - E_{6, 7+} \Delta_{4, 7}),
\end{gather*}
Here, $\Delta_{4, 7}$ and $\Delta_{10, 0, 7+}$ are cusp forms for $\Gamma_0^{*}(7)$ and $\Gamma_0(7)$ of weight $4$ and $10$ such that $v_{\infty} (\Delta_{4, 7}) = v_{\rho_{7, 2}}(\Delta_{4, 7}) = 1$ and $v_{i / \sqrt{7}}(\Delta_{10, 0, 7+}) = v_{\rho_{7, 1}}(\Delta_{10, 0, 7+}) = 1$, $v_{\infty}(\Delta_{10, 0, 7+}) = 2$, respectively.

Furthermore, we have $M_k(\Gamma_0^{*}(7)) = \mathbb{C} E_{k, 7+} \oplus S_k(\Gamma_0^{*}(7))$ and $S_k(\Gamma_0^{*}(7)) = S_{12}(\Gamma_0^{*}(7)) M_{k - 12}(\Gamma_0^{*}(7))$, where
\begin{equation*}
S_{12}(\Gamma_0^{*}(7)) = \mathbb{C} ({E_{2, 7}}')^4 \Delta_{7, 4} \oplus \mathbb{C} ({E_{2, 7}}')^2 (\Delta_{7, 4})^2 \oplus \mathbb{C} (\Delta_{7, 4})^3 \oplus \mathbb{C} (\Delta_7)^4
\end{equation*}
Then, when we write $k = 12 + l$ for $l = 0, 4, 6, 8, 10, 14$, we have
\begin{align*}
M_{12+l}(\Gamma_0^{*}(7))
 &= E_{l, 7+} (\mathbb{C} ({E_{2, 7}}')^{6 n} \oplus ({E_{2, 7}}')^{6 (n-1)} S_{12}(\Gamma_0^{*}(7)) \oplus ({E_{2, 7}}')^{6 (n-2)} (\Delta_7)^4 S_{12}(\Gamma_0^{*}(7))\\
 &\hspace{0.7in} \oplus \cdots \oplus (\Delta_7)^{4(n-1)} S_{12}(\Gamma_0^{*}(7)))\\
 &\quad \oplus (\Delta_7)^{4n} S_l(\Gamma_0^{*}(7)).
\end{align*}
where we denote $\Delta_{6, 7+} := \Delta_{10, 0, 7+} / \Delta_{4, 7}$, and we have $S_4(\Gamma_0^{*}(7)) = \mathbb{C} \Delta_{7, 4}$, $S_6(\Gamma_0^{*}(7)) = \mathbb{C} \Delta_{6, 7+}$,
\begin{gather*}
S_8(\Gamma_0^{*}(7)) = \mathbb{C} ({E_{2, 7}}')^2 \Delta_{7, 4} \oplus \mathbb{C} (\Delta_{7, 4})^2,
 \quad S_{10}(\Gamma_0^{*}(7)) = \mathbb{C} ({E_{2, 7}}')^2 \Delta_{6, 7+} \oplus \mathbb{C} \Delta_{10, 0, 7+},\\
S_{14}(\Gamma_0^{*}(7)) = \mathbb{C} ({E_{2, 7}}')^4 \Delta_{6, 7+} \oplus \mathbb{C} ({E_{2, 7}}')^2 \Delta_{10, 0, 7+} \oplus \mathbb{C} \Delta_{7, 4} \Delta_{10, 0, 7+}.
\end{gather*}\quad

Here, we define
\begin{equation*}
{E_{3, 7}}' := \Delta_{10, 0, 7+} / (\Delta_{4, 7} \Delta_7),
\end{equation*}
which is $4$th semimodular form such that $v_{i / \sqrt{7}}({E_{3, 7}}') = v_{\rho_{7, 1}}({E_{3, 7}}') = 1$, and $({E_{1, 7}}')^3 {E_{3, 7}}'$ is a modular form.

Then, we have $\Delta_{4, 7} = {E_{1, 7}}' \Delta_7$, $\Delta_{6, 7+} = {E_{3, 7}}' \Delta_7$, and $\Delta_{10, 0, 7+} = {E_{1, 7}}' {E_{3, 7}}' (\Delta_7)^2$, and
\begin{equation*}
S_k(\Gamma_0^{*}(7)) = (\mathbb{C} ({E_{1, 7}}')^9 \Delta_7 \oplus \mathbb{C} ({E_{1, 7}}')^6 (\Delta_7)^2 \oplus \mathbb{C} ({E_{1, 7}}')^3 (\Delta_7)^3 \oplus \mathbb{C} (\Delta_7)^4) M_{k - 12}(\Gamma_0^{*}(7)).
\end{equation*}
Furthermore, we can write
\begin{equation*}
M_k(\Gamma_0^{*}(7)) = {E_{\overline{k}, 7+}}' (\mathbb{C} (({E_{1, 7}}')^3)^n \oplus \mathbb{C} (({E_{1, 7}}')^3)^{n-1} \Delta_7 \oplus \cdots \oplus \mathbb{C} (\Delta_7)^n),
\end{equation*}
where $n = \dim(M_k(\Gamma_0^{*}(7))) - 1 = \lfloor k/3 - 2 (k/4 - \lfloor k/4 \rfloor)\rfloor$, and where ${E_{\overline{k}, 7+}}' := 1$, $({E_{1, 7}}')^2 {E_{3, 7}}'$, ${E_{1, 7}}'$, ${E_{3, 7}}'$, $({E_{1, 7}}')^2$, and ${E_{1, 7}}' {E_{3, 7}}'$, when $k \equiv 0$, $2$, $4$, $6$, $8$, and $10 \pmod{12}$, respectively.\\

\paragraph{\bf Hauptmodul}
We define the {\it hauptmodul} of $\Gamma_0^{*}(7)$:
\begin{equation}
J_{7+} := ({E_{1, 7}}')^3 / \Delta_7 = \frac{1}{q} + 9 + 51 q + 204 q^2 + 681 q^3 + \cdots,
\end{equation}
where $v_{\infty}(J_{7+}) = -1$ and $v_{\rho_{7, 2}}(J_{7+}) = 3$. Then, we have
\begin{equation}
J_{7+} : \partial \mathbb{F}_{7+} \setminus \{z \in \mathbb{H} \: ; \: Re(z) = \pm 1/2\} \to [-1, 27] \subset \mathbb{R}.
\end{equation}\quad

\subsection{$\Gamma_0(7)$} (see \cite{SJ1})

\paragraph{\bf Fundamental domain}
We have a fundamental domain for $\Gamma_0(7)$ as follows:
{\small \begin{equation}
\begin{split}
\mathbb{F}_7 = &\left\{|z + 3/7| \geqslant 1/7, \: - 1/2 \leqslant Re(z) \leqslant - 5/14\right\}
 \bigcup \left\{|z + 2/7| > 1/7, \: - 5/14 < Re(z) < - 3/14 \right\}\\
&\bigcup \left\{|z + 1/7| \geqslant 1/7, \: - 3/14 \leqslant Re(z) \leqslant 0\right\}
 \bigcup \left\{|z - 1/7| > 1/7, \: 0 < Re(z) < 3/14 \right\}\\
&\bigcup \left\{|z - 2/7| \geqslant 1/7, \: 3/14 \leqslant Re(z) \leqslant 5/14 \right\}
 \bigcup \left\{|z - 3/7| > 1/7, \: 5/14 < Re(z) < 1/2 \right\}.
\end{split}
\end{equation}
}where $\left( \begin{smallmatrix} -1 & 0 \\ 7 & -1 \end{smallmatrix} \right) : (e^{i \theta} + 1) / 7 \rightarrow (e^{i (\pi - \theta)} - 1) / 7$, $\left( \begin{smallmatrix} -3 & -1 \\ 7 & 2 \end{smallmatrix} \right) : (e^{i \theta} - 2) / 7 \rightarrow (e^{i (\pi - \theta)} - 3) / 7$, and $\left( \begin{smallmatrix} 2 & -1 \\ 7 & -3 \end{smallmatrix} \right) : (e^{i \theta} + 3) / 7 \rightarrow (e^{i (\pi - \theta)} + 2) / 7$. Then, we have
\begin{equation}
\Gamma_0(7) = \langle \left( \begin{smallmatrix} 1 & 1 \\ 0 & 1 \end{smallmatrix} \right), \: - \left( \begin{smallmatrix} 1 & 0 \\ 7 & 1 \end{smallmatrix} \right), \: \left( \begin{smallmatrix} 4 & 1 \\ 7 & 2 \end{smallmatrix} \right) \rangle.
\end{equation}
\begin{figure}[hbtp]
\begin{center}
\includegraphics[width=1.5in]{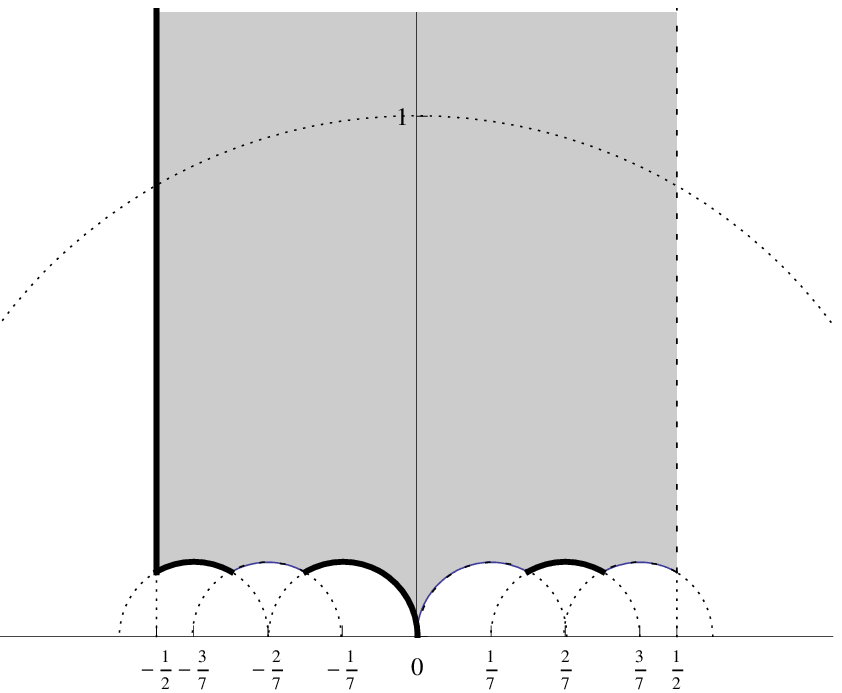}
\end{center}
\caption{$\Gamma_0(7)$}
\end{figure}

\paragraph{\bf Valence formula}
The cusps of $\Gamma_0(7)$ are $\infty$ and $0$, and the elliptic points are $\rho_{7, 2}$ and $\rho_{7, 3} = 5/14 + \sqrt{3} i / 14$. Let $f$ be a modular function of weight $k$ for $\Gamma_0(7)$, which is not identically zero. We have
\begin{equation}
v_{\infty}(f) + v_0(f) + \frac{1}{3} v_{\rho_{7, 2}} (f) + \frac{1}{3} v_{\rho_{7, 3}} (f) + \sum_{\begin{subarray}{c} p \in \Gamma_0(7) \setminus \mathbb{H} \\ p \ne \rho_{7, 2}, \rho_{7, 3}\end{subarray}} v_p(f) = \frac{2 k}{3}.
\end{equation}

Furthermore, the stabilizer of the elliptic point $\rho_{7, 2}$ (resp. $\rho_{7, 3}$) is $\left\{ \pm I, \pm \left( \begin{smallmatrix} -3 & -1 \\ 7 & 2 \end{smallmatrix} \right), \pm \left( \begin{smallmatrix} -2 & -1 \\ 7 & 3 \end{smallmatrix} \right) \right\}$ (resp. $\left\{ \pm I, \pm \left( \begin{smallmatrix} 3 & -1 \\ 7 & -2 \end{smallmatrix} \right), \pm \left( \begin{smallmatrix} 2 & -1 \\ 7 & -3 \end{smallmatrix} \right) \right\}$).\\

\paragraph{\bf For the cusp $\infty$}
We have $\Gamma_{\infty} = \left\{ \pm \left( \begin{smallmatrix}1 & n \\ 0 & 1\end{smallmatrix} \right) \: ; \: n \in \mathbb{Z} \right\}$, and we have the Eisenstein series for the cusp $\infty$ associated with $\Gamma_0(7)$:
\begin{equation}
E_{k, 7}^{\infty}(z) := \frac{7^k E_k(7 z) - E_k(z)}{7^k - 1} \quad \text{for} \; k \geqslant 4.
\end{equation}\quad

\paragraph{\bf For the cusp $0$}
We have $\Gamma_0 = \left\{ \pm \left( \begin{smallmatrix}1 & 0 \\ 7 n & 1\end{smallmatrix} \right) \: ; \: n \in \mathbb{Z} \right\}$ and $\gamma_0 = W_7$, and we have the Eisenstein series for the cusp $0$ associated with $\Gamma_0(7)$:
\begin{equation}
E_{k, 7}^0(z) := \frac{- 7^{k/2} (E_k(7 z) - E_k(z))}{7^k - 1} \quad \text{for} \; k \geqslant 4.
\end{equation}
We also have $\gamma_0^{-1} \ \Gamma_0(7) \ \gamma_0 = \Gamma_0(7)$.\\

\paragraph{\bf The space of modular forms}
We define the following functions:
\begin{equation*}
\Delta_{6, 0, 7} := \Delta_7^{\infty} (\Delta_7^0)^3, \quad \Delta_{6, 1, 7} := (\Delta_7^{\infty})^3 \Delta_7^0,
\end{equation*}
where they are cusp forms of weight $6$.

Let $k$ be an even integer $k \geqslant 4$. We have $M_k(\Gamma_0(7)) = \mathbb{C} E_{k, 7}^{\infty} \oplus \mathbb{C} E_{k, 7}^0 \oplus S_k(\Gamma_0(7))$ and $S_k(\Gamma_0(7)) = (\mathbb{C} \Delta_{6, 0, 7} \oplus \mathbb{C} \Delta_{6, 1, 7} \oplus \mathbb{C} (\Delta_7)^2) M_{k - 6}(\Gamma_0(7))$. Then, we have $M_{6 n + 2}(\Gamma_0(7)) = {E_{2, 7}}' M_{6 n}(\Gamma_0(7))$ and
\begin{align*}
M_{6 n}(\Gamma_0(7)) &= \mathbb{C} (E_{6, 7}^{\infty})^n \oplus \mathbb{C} (E_{6, 7}^{\infty})^{n-1} \Delta_{6, 0, 7} \oplus \mathbb{C} (E_{6, 7}^{\infty})^{n-1} (\Delta_7)^2 \oplus \cdots \oplus \mathbb{C} \Delta_{6, 0, 7} (\Delta_7)^{2(n-1)}\\
 &\oplus \mathbb{C} (E_{6, 7}^0)^n \oplus \mathbb{C} (E_{6, 7}^0)^{n-1} \Delta_{6, 1, 7} \oplus \mathbb{C} (E_{6, 7}^0)^{n-1} (\Delta_7)^2 \oplus \cdots \oplus \mathbb{C} \Delta_{6, 1, 7} (\Delta_7)^{2(n-1)} \oplus \mathbb{C} (\Delta_7)^{2n},\\
M_{6 n + 4}(\Gamma_0(7)) &= E_{4, 7}^{\infty} (\mathbb{C} (E_{6, 7}^{\infty})^n \oplus \mathbb{C} (E_{6, 7}^{\infty})^{n-1} \Delta_{6, 0, 7} \oplus \mathbb{C} (E_{6, 7}^{\infty})^{n-1} (\Delta_7)^2 \oplus \cdots \oplus \mathbb{C} (\Delta_7)^{2n})\\
 &\oplus E_{4, 7}^0 (\oplus \mathbb{C} (E_{6, 7}^0)^n \oplus \mathbb{C} (E_{6, 7}^0)^{n-1} \Delta_{6, 1, 7} \oplus \mathbb{C} (E_{6, 7}^0)^{n-1} (\Delta_7)^2 \oplus \cdots \oplus \mathbb{C} (\Delta_7)^{2n}) \oplus \mathbb{C} \Delta_{4, 7} (\Delta_7)^{2n}.
\end{align*}

Furthermore, we can write
\begin{equation*}
M_k(\Gamma_0(7)) = {E_{k - 3n/2, 7}}' (\mathbb{C} (\Delta_7^{\infty})^n \oplus \mathbb{C} (\Delta_7^{\infty})^{n-1} \Delta_7^0 \oplus \cdots \oplus \mathbb{C} (\Delta_7^0)^n),
\end{equation*}
where $n = \dim(M_k(\Gamma_0(7))) - 1 = \lfloor 2k/3 -2 (k/3 - \lfloor k/3 \rfloor) \rfloor$ and where ${E_{0, 7}}' := 1$.\\

\paragraph{\bf Hauptmodul}
We define the {\it hauptmodul} of $\Gamma_0(7)$:
\begin{equation}
J_7 := \Delta_7^0 / \Delta_7^{\infty} \: (= \eta^4(z) / \eta^4(7 z)) = \frac{1}{q} - 4 + 2 q + 8 q^2 - 5 q^3 - \cdots,
\end{equation}
where $v_{\infty}(J_7) = -1$ and $v_0(J_7) = 1$. Then, we have
{\small \begin{equation}
\begin{split}
J_7 : \qquad \left\{|z + 1/7| = 1/7, \: - 3/14 \leqslant Re(z) \leqslant 0\right\} &\to [- 2 \sqrt{7} + 1/2, 0] \subset \mathbb{R},\\
\left\{|z + 3/7| = 1/7, \: - 1/2 \leqslant Re(z) \leqslant - 5/14\right\} &\to \{- 13/2 \leqslant Re(z) \leqslant - 2 \sqrt{7} + 1/2, \: 0 \leqslant Im(z) \leqslant 3 \sqrt{3} / 2\},\\
\left\{|z - 2/7| = 1/7, \: 3/14 \leqslant Re(z) \leqslant 5/14 \right\} &\to \{- 13/2 \leqslant Re(z) \leqslant - 2 \sqrt{7} + 1/2, \: - 3 \sqrt{3} / 2 \leqslant Im(z) \leqslant 0\}.
\end{split}
\end{equation}
}Thus, $J_7$ does not take real value on some arcs of $\partial \mathbb{F}_7$.

\begin{figure}[hbtp]
\begin{center}
{{\small Lower arcs of $\partial \mathbb{F}_7$}\includegraphics[width=2.5in]{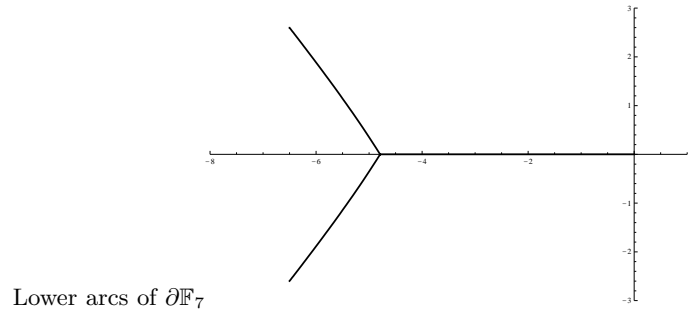}}
\end{center}
\caption{Image by $J_7$}\label{Im-J7B}
\end{figure}

\clearpage

\section{Level $8$}

We have $\Gamma_0(8)+ = \Gamma_0(8)+8 = \Gamma_0^{*}(8)$ and $\Gamma_0(8)-=\Gamma_0(8)$.

We have $W_8 = \left(\begin{smallmatrix} 0 & -1 / (2 \sqrt{2}) \\ 2 \sqrt{2} & 0 \end{smallmatrix}\right)$, $W_{8-, 2} := \left(\begin{smallmatrix} -1 & -1/2 \\ 4 & 1 \end{smallmatrix}\right)$, and $W_{8-, 4} := \left(\begin{smallmatrix} -\sqrt{2} & -3/(2 \sqrt{2}) \\ 2 \sqrt{2} & \sqrt{2} \end{smallmatrix}\right)$, and we denote $\rho_8 :=  -1/3 + i / (6 \sqrt{2})$. We define
\begin{equation}
\begin{aligned}
\Delta_8^{\infty}(z) &:= \eta^4(8 z) / \eta^2(4 z),\\
\Delta_8^{-1/2}(z) &:= \eta^{10}(2 z) / (\eta^4(z) \eta^4(4 z)),\\
\Delta_8(z) &:= \Delta_8^{\infty}(z) \Delta_8^0(z) \Delta_8^{-1/2}(z) \Delta_8^{-1/4}(z)
\end{aligned}
\begin{aligned}
\Delta_8^0(z) &:= \eta^4(z) / \eta^2(2 z), \\
\Delta_8^{-1/4}(z) &:= \eta^{10}(4 z) / (\eta^4(2 z) \eta^4(8 z)),\\
= \eta^4(2 z) &\eta^4(4 z),
\end{aligned}
\end{equation}
where $\Delta_8^{\infty}$, $\Delta_8^0$, $\Delta_8^{-1/2}$, and $\Delta_8^{-1/4}$ are $2$nd semimodular forms for $\Gamma_0(8)$ of weight $1$ such that $v_{\infty}(\Delta_8^{\infty}) = v_0(\Delta_8^0) = v_{-1/2}(\Delta_8^{-1/2}) = v_{-1/4}(\Delta_8^{-1/4}) = 1$, and $\Delta_8$ is a cusp form for $\Gamma_0(8)$ and $\Gamma_0^{*}(8)$ of weight $4$. Furthermore, we define
\begin{equation}
\begin{split}
{E_{2, 8}}'(z) &:= 2 E_2(8 z) - E_2(4 z),\\
{E_{2, 8+8}}'(z) &:= (8 E_2(8 z) - 4 E_2(4 z) - 2 E_2(2 z) + E_2(z)) / 3,
\end{split}
\end{equation}
which are modular forms for $\Gamma_0(8)$ of weight $2$, and we have $v_{-1/8 + i/8}({E_{2, 8}}') = v_{-3/8 + i/8}({E_{2, 8}}') = 1$ and $v_{i / (2 \sqrt{2})}({E_{2, 8+8}}') = v_{\rho_8}({E_{2, 8+8}}') = 1$.\\

\subsection{$\Gamma_0(8)+8 = \Gamma_0^{*}(8)$}

\paragraph{\bf Fundamental domain}
We have a fundamental domain for $\Gamma_0^{*}(8)$ as follows:
{\small \begin{equation}
\begin{split}
\mathbb{F}_{8+8} = &\left\{|z + 3/8| \geqslant 1/8, \: - 1/2 \leqslant Re(z) < - 1/4 \right\}
 \bigcup \left\{|z| \geqslant 1 / (2 \sqrt{2}), \: - 1/4 \leqslant Re(z) \leqslant 0 \right\}\\
 &\bigcup \left\{|z| > 1 / (2 \sqrt{2}), \: 0 < Re(z) \leqslant 1/4 \right\}
 \bigcup \left\{|z - 3/8| > 1/8, \: 1/4 < Re(z) < 1/2 \right\}.
\end{split}
\end{equation}
}where $W_8 :  e^{i \theta} / (2 \sqrt{2}) \rightarrow e^{i  (\pi - \theta)} / (2 \sqrt{2})$ and $\left(\begin{smallmatrix} -3 & 1 \\ 8 & -3 \end{smallmatrix} \right) : (e^{i \theta} + 3) / 8 \rightarrow (e^{i (\pi - \theta)} - 3) / 8$. Then, we have
\begin{equation}
\Gamma_0^{*}(8) = \langle \left( \begin{smallmatrix} 1 & 1 \\ 0 & 1 \end{smallmatrix} \right), \: W_8, \: \left( \begin{smallmatrix} 3 & 1 \\ 8 & 3 \end{smallmatrix} \right) \rangle.
\end{equation}
\begin{figure}[hbtp]
\begin{center}
\includegraphics[width=1.5in]{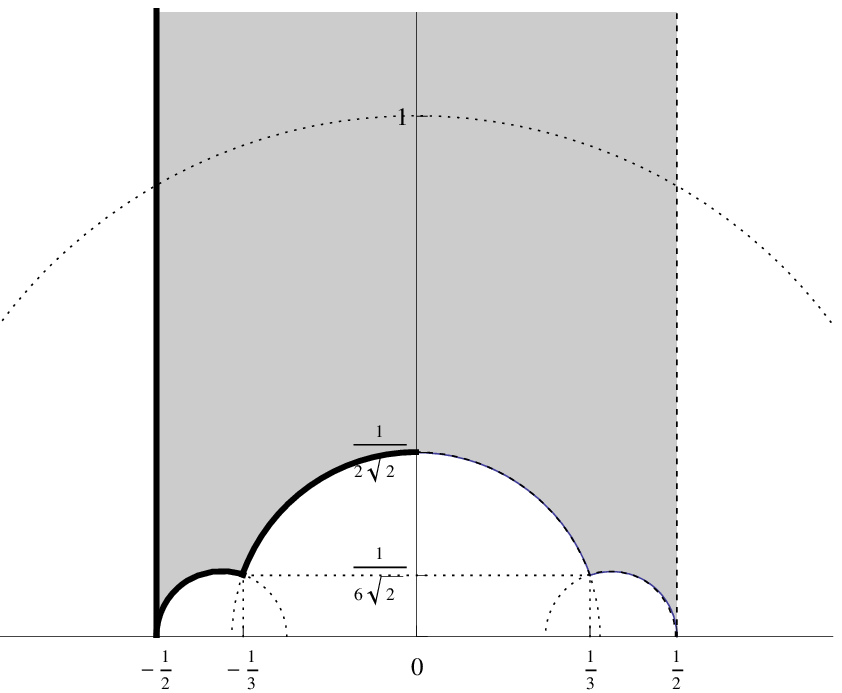}
\end{center}
\caption{$\Gamma_0^{*}(8)$}
\end{figure}

\paragraph{\bf Valence formula}
The cusps of $\Gamma_0^{*}(8)$ are $\infty$ and $-1/2$, and the elliptic points are $i / (2 \sqrt{2})$ and $\rho_8 = - 1/4 + i / (5 \sqrt{6})$. Let $f$ be a modular function of weight $k$ for $\Gamma_0^{*}(8)$, which is not identically zero. We have
\begin{equation}
v_{\infty}(f) + v_{-1/2}(f) + \frac{1}{2} v_{i / (2 \sqrt{2})} (f) + \frac{1}{2} v_{\rho_8} (f) + \sum_{\begin{subarray}{c} p \in \Gamma_0^{*}(8) \setminus \mathbb{H} \\ p \ne i / (2 \sqrt{2}), \rho_8\end{subarray}} v_p(f) = \frac{k}{2}.
\end{equation}

Furthermore, the stabilizer of the elliptic point $i / (2 \sqrt{2})$ (resp. $\rho_8$) is $\left\{ \pm I, \pm W_8 \right\}$ (resp. $\left\{ \pm I, \pm \left( \begin{smallmatrix} 3 & -1 \\ -8 & 3 \end{smallmatrix} \right) W_8 \right\}$).\\

\paragraph{\bf For the cusp $\infty$}
We have $\Gamma_{\infty} = \left\{ \pm \left( \begin{smallmatrix}1 & n \\ 0 & 1\end{smallmatrix} \right) \: ; \: n \in \mathbb{Z} \right\}$, and we have the Eisenstein series for the cusp $\infty$ associated with $\Gamma_0^{*}(8)$:
\begin{equation}
E_{k, 8+8}^{\infty}(z) := \frac{2^{3k/2} E_k(8 z) - 2^{k/2} E_k(4 z) - E_k(2 z) + E_k(z)}{2^{k/2}(2^k - 1)} \quad \text{for} \; k \geqslant 4.
\end{equation}\quad

\paragraph{\bf For the cusp $-1/2$}
We have $\Gamma_{-1/2} = \left\{ \pm \left( \begin{smallmatrix} 4 n + 1 & 2 n \\ - 8 n & - 4 n + 1\end{smallmatrix} \right) \: ; \: n \in \mathbb{Z} \right\}$ and $\gamma_{-1/2} = W_{8-, 4}$, and we have the Eisenstein series for the cusp $-1/2$ associated with $\Gamma_0^{*}(8)$:
\begin{equation}
E_{k, 8+8}^{-1/2}(z) := \frac{- (2^{3k/2} E_k(8 z) - 2^{k/2} (2^k - 2^{k/2} + 1) E_k(4 z) - (2^k - 2^{k/2} + 1) E_k(2 z) + E_k(z))}{2^{k/2}(2^k - 1)} \quad \text{for} \; k \geqslant 4.
\end{equation}
We also have $\gamma_{-1/2}^{-1} \ \Gamma_0^{*}(8) \ \gamma_{-1/2} = \Gamma_0^{*}(8)$.\\

\paragraph{\bf The space of modular forms}
We define
\begin{equation*}
\Delta_{8+8}^{\infty} := \Delta_8^{\infty} \Delta_8^0, \quad \Delta_{8+8}^{-1/2} := \Delta_8^{-1/2} \Delta_8^{-1/4},
\end{equation*}
which are $2$nd semimodular forms for $\Gamma_0^{*}(8)$ of weight $2$.

Let $k$ be an even integer $k \geqslant 4$. We have $M_k(\Gamma_0^{*}(8)) = \mathbb{C} E_{k, 8+8}^{\infty} \oplus \mathbb{C} E_{k, 8+8}^{-1/2} \oplus S_k(\Gamma_0^{*}(8))$ and $S_k(\Gamma_0^{*}(8)) = \Delta_8 M_{k - 4}(\Gamma_0^{*}(8))$. Then, we have $M_{4 n + 2}(\Gamma_0^{*}(8)) = {E_{2, 8+8}}' M_{4 n}(\Gamma_0^{*}(8))$ and
\begin{align*}
M_{4 n}(\Gamma_0^{*}(8)) &= \mathbb{C} (E_{4, 8+8}^{\infty})^n \oplus \mathbb{C} (E_{4, 8+8}^{\infty})^{n-1} \Delta_8 \oplus \cdots \oplus \mathbb{C} E_{4, 8+8}^{\infty} (\Delta_8)^{n-1}\\
 &\oplus \mathbb{C} (E_{4, 8+8}^{-1/2})^n \oplus \mathbb{C} (E_{4, 8+8}^{-1/2})^{n-1} \Delta_8 \oplus \cdots \oplus \mathbb{C} E_{4, 8+8}^{-1/2} (\Delta_8)^{n-1} \oplus \mathbb{C} (\Delta_8)^n.
\end{align*}

Furthremore, we can write
\begin{equation*}
M_{4 n}(\Gamma_0^{*}(8)) = \mathbb{C} (\Delta_{8+8}^{\infty})^{2n} \oplus \mathbb{C} (\Delta_{8+8}^{\infty})^{2n-1} \Delta_{8+8}^{-1/2} \oplus \cdots \oplus \mathbb{C} (\Delta_{8+8}^{-1/2})^{2n}.
\end{equation*}\quad

\paragraph{\bf Hauptmodul}
We define the {\it hauptmodul} of $\Gamma_0^{*}(8)$:
\begin{equation}
J_{8+8} := \Delta_{8+8}^{-1/2} / \Delta_{8+8}^{\infty} \: (= \eta^{-8}(z) \eta^8(2 z) \eta^8(4 z) \eta^{-8}(8 z)) = \frac{1}{q} + 8 + 36 q + 128 q^2 + 386 q^3 + \cdots,
\end{equation}
where $v_{\infty}(J_{8+8}) = -1$ and $v_{-1/2}(J_{8+8}) = 1$. Then, we have
\begin{equation}
J_{8+8} : \partial \mathbb{F}_{8+8} \setminus \{z \in \mathbb{H} \: ; \: Re(z) = \pm 1/2\} \to [0, 12 + 8 \sqrt{2}] \subset \mathbb{R}.
\end{equation}\quad

\subsection{$\Gamma_0(8)$}

\paragraph{\bf Fundamental domain}
We have a fundamental domain for $\Gamma_0(8)$ as follows:
{\small \begin{equation}
\begin{split}
\mathbb{F}_8 = &\left\{|z + 3/8| \geqslant 1/8, \: - 1/2 \leqslant Re(z) < -1/4 \right\}
 \bigcup \left\{|z + 1/8| \geqslant 1/8, \: - 1/4 \leqslant Re(z) \leqslant 0 \right\}\\
 &\bigcup \left\{|z - 1/8| > 1/8, \: 0 < Re(z) \leqslant 1/4 \right\}
 \bigcup \left\{|z - 3/8| > 1/8, \: 1/4 < Re(z) < 1/2 \right\}.
\end{split}
\end{equation}
}where $\left(\begin{smallmatrix} -1 & 0 \\ 8 & -1 \end{smallmatrix}\right) : (e^{i \theta} + 1) / 8 \rightarrow (e^{i (\pi - \theta)} - 1) / 8$ and $\left(\begin{smallmatrix} -3 & 1 \\ 8 & -3 \end{smallmatrix} \right) : (e^{i \theta} + 3) / 8 \rightarrow (e^{i (\pi - \theta)} - 3) / 8$.
\begin{equation}
\Gamma_0(8) = \langle - I, \: \left( \begin{smallmatrix} 1 & 1 \\ 0 & 1 \end{smallmatrix} \right), \: \left( \begin{smallmatrix} 1 & 0 \\ 8 & 1 \end{smallmatrix} \right), \: \left( \begin{smallmatrix} 3 & 1 \\ 8 & 3 \end{smallmatrix} \right) \rangle.
\end{equation}
\begin{figure}[hbtp]
\begin{center}
\includegraphics[width=1.5in]{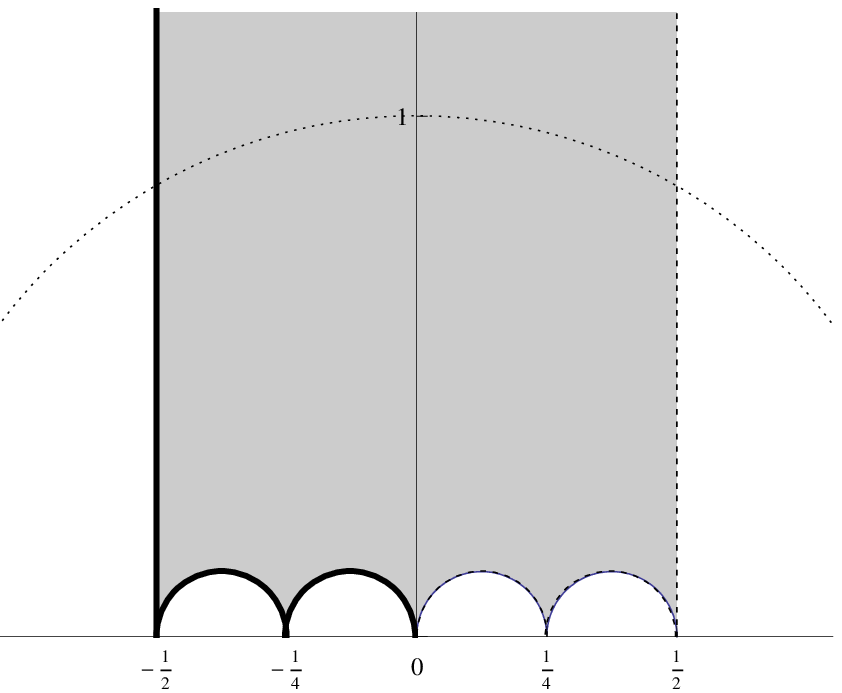}
\end{center}
\caption{$\Gamma_0(8)$}
\end{figure}

\paragraph{\bf Valence formula}
The cusps of $\Gamma_0(8)$ are $\infty$, $0$, $-1/2$, and $-1/4$. Let $f$ be a modular function of weight $k$ for $\Gamma_0(8)$, which is not identically zero. We have
\begin{equation}
v_{\infty}(f) + v_0(f) + v_{-1/2} (f) + v_{-1/4} (f) + \sum_{p \in \Gamma_0(8) \setminus \mathbb{H}} v_p(f) = k.
\end{equation}\quad

\paragraph{\bf For the cusp $\infty$}
We have $\Gamma_{\infty} = \left\{ \pm \left( \begin{smallmatrix}1 & n \\ 0 & 1\end{smallmatrix} \right) \: ; \: n \in \mathbb{Z} \right\}$, and we have the Eisenstein series for the cusp $\infty$ associated with $\Gamma_0(8)$:
\begin{equation}
E_{k, 8}^{\infty}(z) := \frac{2^k E_k(8 z) - E_k(4 z)}{2^k - 1} \quad \text{for} \; k \geqslant 4.
\end{equation}
Note that we have $E_{k, 8}^{\infty}(z) = E_{k, 2}^{\infty}(4 z)$.\\

\paragraph{\bf For the cusp $0$}
We have $\Gamma_0 = \left\{ \pm \left( \begin{smallmatrix}1 & 0 \\ 8 n & 1\end{smallmatrix} \right) \: ; \: n \in \mathbb{Z} \right\}$ and $\gamma_0 = W_8$, and we have the Eisenstein series for the cusp $0$ associated with $\Gamma_0(8)$:
\begin{equation}
E_{k, 8}^0(z) := \frac{- (E_k(2 z) - E_k(z))}{2^{k/2} (2^k - 1)} \quad \text{for} \; k \geqslant 4.
\end{equation}
Note that we have $E_{k, 4}^0(z) = 2^{-k} E_{k, 2}^0(z)$. We also have $\gamma_0^{-1} \ \Gamma_0(8) \ \gamma_0 = \Gamma_0(8)$.\\

\paragraph{\bf For the cusp $-1/2$}
We have $\Gamma_{-1/2} = \left\{ \pm \left( \begin{smallmatrix} 4 n + 1 & 2 n \\ - 8 n & - 4 n + 1 \end{smallmatrix} \right) \: ; \: n \in \mathbb{Z} \right\}$ and $\gamma_{-1/2} = W_{8-, 4}$, and we have the Eisenstein series for the cusp $-1/2$ associated with $\Gamma_0(8)$:
\begin{equation}
E_{k, 8}^{-1/2}(z) := \frac{- (2^k E_k(4 z) - (2^k + 1) E_k(2 z) + E_k(z))}{2^{k/2} (2^k - 1)} \quad \text{for} \; k \geqslant 4.
\end{equation}
Note that we have $E_{k, 4}^{-1/2}(z) = 2^{-k/2} E_{k, 4+4}^{-1/2}(z)$. We also have $\gamma_{-1/2}^{-1} \ \Gamma_0(8) \ \gamma_{-1/2} = \Gamma_0(8)$.\\

\paragraph{\bf For the cusp $-1/4$}
We have $\Gamma_{-1/4} = \left\{ \pm \left( \begin{smallmatrix} 4 n + 1 & n \\ - 16 n & - 4 n + 1 \end{smallmatrix} \right) \: ; \: n \in \mathbb{Z} \right\}$ and $\gamma_{-1/4} = W_{8-, 2}$, and we have the Eisenstein series for the cusp $-1/4$ associated with $\Gamma_0(8)$:
\begin{equation}
E_{k, 8}^{-1/4}(z) := \frac{- (2^k E_k(8 z) - (2^k + 1) E_k(4 z) + E_k(2 z))}{2^k - 1} \quad \text{for} \; k \geqslant 4.
\end{equation}
Note that we have $E_{k, 4}^{-1/4}(z) = E_{k, 4+4}^{-1/2}(2 z)$. We also have $\gamma_{-1/4}^{-1} \ \Gamma_0(8) \ \gamma_{-1/4} = \Gamma_0(8)$.\\

\paragraph{\bf The space of modular forms}
Let $k$ be an even integer $k \geqslant 4$. We have $M_k(\Gamma_0(8)) = \mathbb{C} E_{k, 8}^{\infty} \oplus \mathbb{C} E_{k, 8}^0 \oplus \mathbb{C} E_{k, 8}^{-1/2} \oplus \mathbb{C} E_{k, 8}^{-1/4} \oplus S_k(\Gamma_0(8))$ and $S_k(\Gamma_0(8)) = \Delta_8 M_{k - 4}(\Gamma_0(8))$. Then, we have $M_{4 n + 2}(\Gamma_0(8)) = {E_{2, 8}}' M_{4 n}(\Gamma_0(8)) \oplus \mathbb{C} {E_{2, 4}}' (\Delta_8)^n \oplus \mathbb{C} {E_{2, 2}}' (\Delta_8)^n$ and
\begin{align*}
M_{4 n}(\Gamma_0(8)) = &\mathbb{C} (E_{4, 8}^{\infty})^n \oplus \mathbb{C} (E_{4, 8}^{\infty})^{n-1} \Delta_8 \oplus \cdots \oplus \mathbb{C} E_{4, 8}^{\infty} (\Delta_8)^{n-1}\\
 &\oplus \mathbb{C} (E_{4, 8}^0)^n \oplus \mathbb{C} (E_{4, 8}^0)^{n-1} \Delta_8 \oplus \cdots \oplus \mathbb{C} E_{4, 8}^0 (\Delta_8)^{n-1}\\
 &\oplus \mathbb{C} (E_{4, 8}^{-1/2})^n \oplus \mathbb{C} (E_{4, 8}^{-1/2})^{n-1} \Delta_8 \oplus \cdots \oplus \mathbb{C} E_{4, 8}^{-1/2} (\Delta_8)^{n-1}\\
 &\oplus \mathbb{C} (E_{4, 8}^{-1/4})^n \oplus \mathbb{C} (E_{4, 8}^{-1/4})^{n-1} \Delta_8 \oplus \cdots \oplus \mathbb{C} E_{4, 8}^{-1/4} (\Delta_8)^{n-1} \oplus \mathbb{C} (\Delta_8)^n.
\end{align*}

Furthermore, we can write
\begin{equation*}
M_{2 n}(\Gamma_0(8)) = \mathbb{C} (\Delta_8^{\infty})^{2n} \oplus \mathbb{C} (\Delta_8^{\infty})^{2n-1} \Delta_8^0 \oplus \cdots \oplus \mathbb{C} (\Delta_8^0)^{2n}.
\end{equation*}\quad

\paragraph{\bf Hauptmodul}
We define the {\it hauptmodul} of $\Gamma_0(8)$:
\begin{equation}
J_8 := \Delta_8^0 / \Delta_8^{\infty} \: (= \eta^4(z) \eta^{-2}(2 z) \eta^2(4 z) \eta^{-4}(8 z)) = \frac{1}{q} -4 + 4 q + 2 q^3 - 8 q^5 - \cdots,
\end{equation}
where $v_{\infty}(J_8) = -1$ and $v_0(J_8) = 1$. Then, we have
\begin{equation}
J_8 : \partial \mathbb{F}_8 \setminus \{z \in \mathbb{H} \: ; \: Re(z) = \pm 1/2\} \to [-8, 0] \subset \mathbb{R}.
\end{equation}

%% file: report09-11.tex
\section{Level $9$}

We have $\Gamma_0(9)+ = \Gamma_0(9)+9 = \Gamma_0^{*}(9)$ and $\Gamma_0(9)-=\Gamma_0(9)$.

We have $W_9 = \left(\begin{smallmatrix} 0 & -1/3 \\ 3 & 0 \end{smallmatrix}\right)$, $W_{9-, 3} := \left(\begin{smallmatrix} -1 & -2/3 \\ 3 & 1 \end{smallmatrix}\right)$, and $W_{9-, -3} := \left(\begin{smallmatrix} 1 & -2/3 \\ 3 & -1 \end{smallmatrix}\right)$, and we define $\rho_9 := - 1/2 + i/6$. We define
\begin{equation}
\begin{aligned}
\Delta_9^{\infty}(z) &:= \eta^3(9 z) / \eta(3 z),\\
\Delta_9^{-1/3}(z) &:= \eta^3(z + 1/3) / \eta(3 z),\\
\Delta_9(z) &:= \Delta_9^{\infty}(z) \Delta_9^0(z) \Delta_9^{-1/3}(z) \Delta_9^{1/3}(z)
\end{aligned}
\begin{aligned}
\Delta_9^0(z) &:= \eta^3(z) / \eta(3 z), \\
\Delta_9^{1/3}(z) &:= \eta^3(z - 1/3) / \eta(3 z),\\
= \eta^3(8 z)&,
\end{aligned}
\end{equation}
where $\Delta_9^{\infty}$, $\Delta_9^0$, $\Delta_9^{-1/2}$, and $\Delta_9^{-1/3}$ are $2$nd semimodular forms for $\Gamma_0(9)$ of weight $1$ such that $v_{\infty}(\Delta_9^{\infty}) = v_0(\Delta_9^0) = v_{-1/3}(\Delta_9^{-1/3}) = v_{1/3}(\Delta_9^{1/3}) = 1$, and $\Delta_9$ is a cusp form for $\Gamma_0(9)$ and $\Gamma_0^{*}(9)$ of weight $4$. Furthermore, we define
\begin{equation}
\begin{split}
{E_{2, 9}}'(z) &:= (3 E_2(9 z) - E_2(3 z)) / 2,\\
{E_{2, 9+3}}'(z) &:= (3 E_2(3 z) - E_2(z + 1/3)) / 2,\\
{E_{2, 9+9}}'(z) &:= (9 E_2(9 z) - 6 E_2(3 z) + E_2(z)) / 3,
\end{split}
\end{equation}
which are modular forms for $\Gamma_0(9)$ of weight $2$, and we have $v_{-1/6 + \sqrt{3} i/18}({E_{2, 9}}') = 2$, $v_{1/6 + \sqrt{3} i/6}({E_{2, 9+3}}') = 2$ and $v_{i/3}({E_{2, 9+9}}') = v_{\rho_9}({E_{2, 9+9}}') = 1$.\\

\subsection{$\Gamma_0(9)+9 = \Gamma_0^{*}(9)$}

\paragraph{\bf Fundamental domain}
We have a fundamental domain for $\Gamma_0^{*}(9)$ as follows:
{\small \begin{equation}
\begin{split}
\mathbb{F}_{9+9} = &\left\{|z + 5/9| \geqslant 1/9, \: - 1/2 \leqslant Re(z) < - 1/3 \right\}
 \bigcup \left\{|z| \geqslant 1/3, \: - 1/3 \leqslant Re(z) \leqslant 0 \right\}\\
 &\bigcup \left\{|z| > 1/3, \: 0 < Re(z) \leqslant 1/3 \right\}
 \bigcup \left\{|z - 5/9| > 1/9, \: 1/3 < Re(z) < 1/2 \right\}.
\end{split}
\end{equation}
}where $W_9 :  e^{i \theta}/3 \rightarrow e^{i  (\pi - \theta)}/3$ and $\left(\begin{smallmatrix} -4 & -1 \\ 9 & 2 \end{smallmatrix} \right) W_9 : (e^{i \theta} + 3) / 6 \rightarrow (e^{i (\pi - \theta)} - 3) / 6$. Then, we have
\begin{equation}
\Gamma_0^{*}(9) = \langle \left( \begin{smallmatrix} 1 & 1 \\ 0 & 1 \end{smallmatrix} \right), \: W_9, \: \left( \begin{smallmatrix} 5 & 1 \\ 9 & 2 \end{smallmatrix} \right) \rangle.
\end{equation}
\begin{figure}[hbtp]
\begin{center}
\includegraphics[width=1.5in]{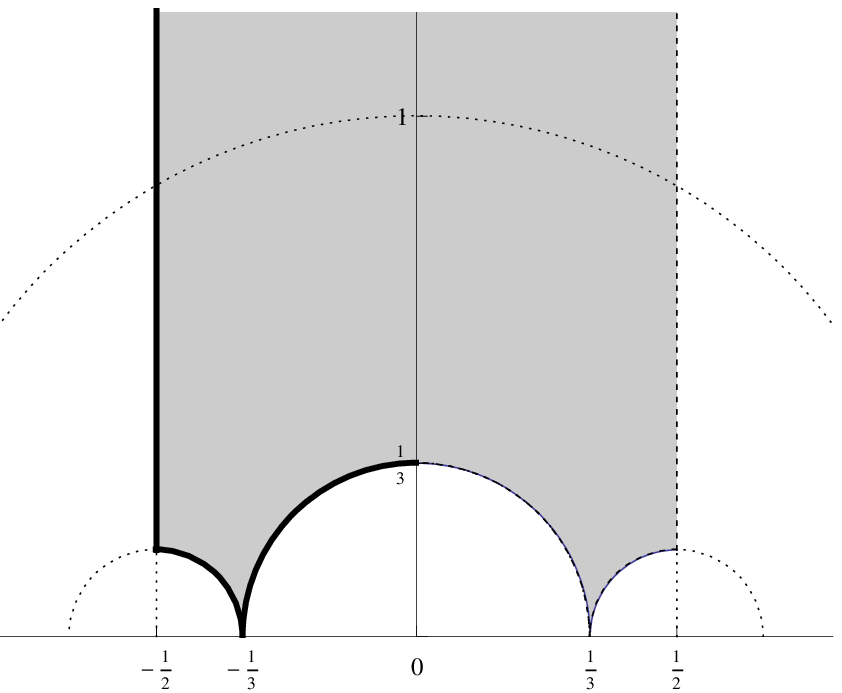}
\end{center}
\caption{$\Gamma_0^{*}(9)$}
\end{figure}

\paragraph{\bf Valence formula}
The cusps of $\Gamma_0^{*}(9)$ are $\infty$ and $-1/3$, and the elliptic points are $i/3$ and $\rho_9 = - 1/2 + i/6$. Let $f$ be a modular function of weight $k$ for $\Gamma_0^{*}(9)$, which is not identically zero. We have
\begin{equation}
v_{\infty}(f) + v_{-1/3}(f) + \frac{1}{2} v_{i/3} (f) + \frac{1}{2} v_{\rho_9} (f) + \sum_{\begin{subarray}{c} p \in \Gamma_0^{*}(9) \setminus \mathbb{H} \\ p \ne i/3, \rho_9\end{subarray}} v_p(f) = \frac{k}{2}.
\end{equation}

Furthermore, the stabilizer of the elliptic point $i/3$ (resp. $\rho_9$) is $\left\{ \pm I, \pm W_9 \right\}$ (resp. $\left\{ \pm I, \pm \left(\begin{smallmatrix} 5 & -1 \\ -9 & 2 \end{smallmatrix} \right) W_9 \right\}$).\\

\paragraph{\bf For the cusp $\infty$}
We have $\Gamma_{\infty} = \left\{ \pm \left( \begin{smallmatrix}1 & n \\ 0 & 1\end{smallmatrix} \right) \: ; \: n \in \mathbb{Z} \right\}$, and we have the Eisenstein series for the cusp $\infty$ associated with $\Gamma_0^{*}(9)$:
\begin{equation}
E_{k, 9+9}^{\infty}(z) := \frac{3^k E_k(9 z) - 2 E_k(3 z) + E_k(z)}{3^k - 1} \quad \text{for} \; k \geqslant 4.
\end{equation}\quad

\paragraph{\bf For the cusp $-1/3$}
We have $\Gamma_{-1/3} = \left\{ \pm \left( \begin{smallmatrix} 3 n + 1 & n \\ - 9 n & - 3 n + 1\end{smallmatrix} \right) \: ; \: n \in \mathbb{Z} \right\}$ and $\gamma_{-1/3} = W_{9-, 3}$, and we have the Eisenstein series for the cusp $-1/3$ associated with $\Gamma_0^{*}(9)$:
\begin{equation}
E_{k, 9+9}^{-1/3}(z) := \frac{- (3^k E_k(9 z) - (3^k + 1) E_k(3 z) + E_k(z))}{(3^k - 1)} \quad \text{for} \; k \geqslant 4.
\end{equation}
We also have $\gamma_{-1/3}^{-1} \ \Gamma_0^{*}(9) \ \gamma_{-1/3} = \Gamma_0^{*}(9)$.\\

\paragraph{\bf The space of modular forms}
We define
\begin{equation*}
\Delta_{9+9}^{\infty} := \Delta_9^{\infty} \Delta_9^0, \quad \Delta_{9+9}^{-1/3} := \Delta_9^{-1/3} \Delta_9^{1/3},
\end{equation*}
which are $2$nd semimodular forms for $\Gamma_0^{*}(9)$ of weight $2$.

Let $k$ be an even integer $k \geqslant 4$. We have $M_k(\Gamma_0^{*}(9)) = \mathbb{C} E_{k, 9+9}^{\infty} \oplus \mathbb{C} E_{k, 9+9}^{-1/3} \oplus S_k(\Gamma_0^{*}(9))$ and $S_k(\Gamma_0^{*}(9)) = \Delta_9 M_{k - 4}(\Gamma_0^{*}(9))$. Then, we have $M_{4 n + 2}(\Gamma_0^{*}(9)) = {E_{2, 9+9}}' M_{4 n}(\Gamma_0^{*}(9))$ and
\begin{align*}
M_{4 n}(\Gamma_0^{*}(9)) &= \mathbb{C} (E_{4, 9+9}^{\infty})^n \oplus \mathbb{C} (E_{4, 9+9}^{\infty})^{n - 1} \Delta_9 \oplus \cdots \oplus \mathbb{C} E_{4, 9+9}^{\infty} \Delta_9^{n - 1}\\
 &\oplus \mathbb{C} (E_{4, 9+9}^{-1/3})^n \oplus \mathbb{C} (E_{4, 9+9}^{-1/3})^{n - 1} \Delta_9 \oplus \cdots \oplus \mathbb{C} E_{4, 9+9}^{-1/3} \Delta_9^{n - 1} \oplus \mathbb{C} (\Delta_9)^n.
\end{align*}

Furthermore, we can write
\begin{equation*}
M_{4 n}(\Gamma_0^{*}(9)) = \mathbb{C} (\Delta_{9+9}^{\infty})^{2n} \oplus \mathbb{C} (\Delta_{9+9}^{\infty})^{2n-1} \Delta_{9+9}^{-1/3} \oplus \cdots \oplus \mathbb{C} (\Delta_{9+9}^{-1/3})^{2n}.
\end{equation*}\quad

\paragraph{\bf Hauptmodul}
We define the {\it hauptmodul} of $\Gamma_0^{*}(9)$:
\begin{equation}
J_{9+9} := \Delta_{9+9}^{-1/3} / \Delta_{9+9}^{\infty} \: (= \eta^{12}(3 z) / (\eta^6(z) \eta^6(9 z))) = \frac{1}{q} + 6 + 27 q + 86 q^2 + 243 q^3 + \cdots,
\end{equation}
where $v_{\infty}(J_{9+9}) = -1$ and $v_{-1/3}(J_{9+9}) = 1$. Then, we have
\begin{equation}
J_{9+9} : \partial \mathbb{F}_{9+9} \setminus \{z \in \mathbb{H} \: ; \: Re(z) = \pm 1/2\} \to [9 - 6 \sqrt{3}, 9 + 6 \sqrt{3}] \subset \mathbb{R}.
\end{equation}\quad

\subsection{$\Gamma_0(9)$}

\paragraph{\bf Fundamental domain}
We have a fundamental domain for $\Gamma_0(9)$ as follows:
{\small \begin{equation}
\begin{split}
\mathbb{F}_9 = &\left\{|z + 4/9| \geqslant 1/9, \: - 1/2 \leqslant Re(z) \leqslant -1/3 \right\}
 \bigcup \left\{|z + 2/9| >1/9, \: - 1/3 < Re(z) < -1/6 \right\}\\
 &\bigcup \left\{|z + 1/9| \geqslant 1/9, \: -1/6 \leqslant Re(z) \leqslant 0 \right\}
 \bigcup \left\{|z - 1/9| > 1/9, \: 0 < Re(z) < 1/6 \right\}\\
 &\bigcup \left\{|z - 2/9| \geqslant 1/9, \: 1/6 \leqslant Re(z) \leqslant 1/3 \right\}
 \bigcup \left\{|z - 4/9| > 1/9, \: 1/3 < Re(z) < 1/2 \right\}.
\end{split}
\end{equation}
}where $\left(\begin{smallmatrix} -1 & 0 \\ 9 & -1 \end{smallmatrix}\right) : (e^{i \theta} + 1) / 9 \rightarrow (e^{i (\pi - \theta)} - 1) / 9$, $\left(\begin{smallmatrix} -4 & -1 \\ 9 & 2 \end{smallmatrix} \right) : (e^{i \theta} - 2) / 9 \rightarrow (e^{i (\pi - \theta)} - 4) / 9$, and $\left(\begin{smallmatrix} 2 & -1 \\ 9 & -4 \end{smallmatrix} \right) : (e^{i \theta} + 4) / 9 \rightarrow (e^{i (\pi - \theta)} + 2) / 9$.
\begin{equation}
\Gamma_0(9) = \langle - I, \: \left( \begin{smallmatrix} 1 & 1 \\ 0 & 1 \end{smallmatrix} \right), \: \left( \begin{smallmatrix} 1 & 0 \\ 9 & 1 \end{smallmatrix} \right), \: \left( \begin{smallmatrix} 5 & 1 \\ 9 & 2 \end{smallmatrix} \right) \rangle.
\end{equation}
\begin{figure}[hbtp]
\begin{center}
\includegraphics[width=1.5in]{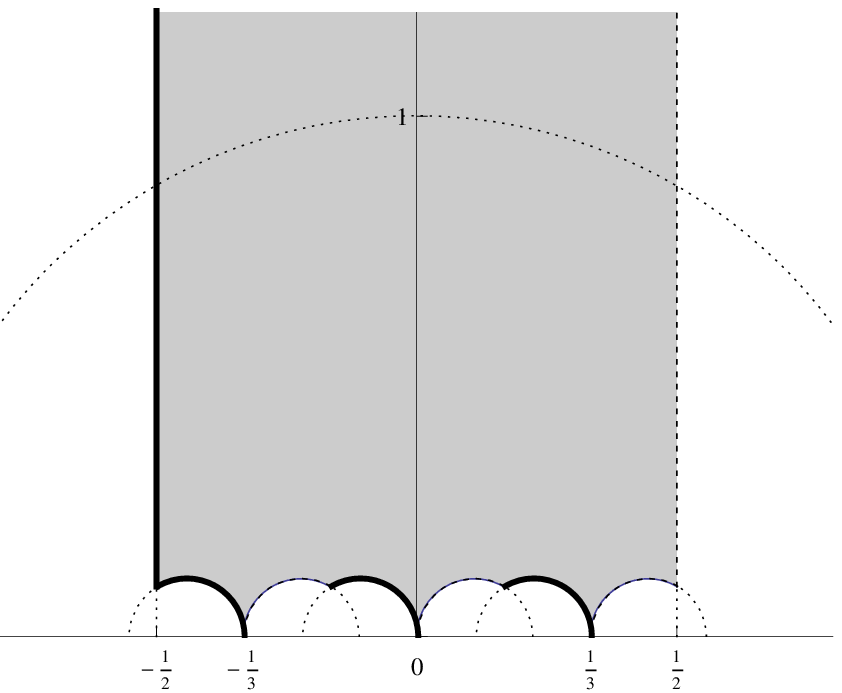}
\end{center}
\caption{$\Gamma_0(9)$}
\end{figure}

\paragraph{\bf Valence formula}
The cusps of $\Gamma_0(9)$ are $\infty$, $0$, $-1/3$, and $1/3$. Let $f$ be a modular function of weight $k$ for $\Gamma_0(9)$, which is not identically zero. We have
\begin{equation}
v_{\infty}(f) + v_0(f) + v_{-1/3} (f) + v_{1/3} (f) + \sum_{p \in \Gamma_0(9) \setminus \mathbb{H}} v_p(f) = k.
\end{equation}\quad

\paragraph{\bf For the cusp $\infty$}
We have $\Gamma_{\infty} = \left\{ \pm \left( \begin{smallmatrix}1 & n \\ 0 & 1\end{smallmatrix} \right) \: ; \: n \in \mathbb{Z} \right\}$, and we have the Eisenstein series for the cusp $\infty$ associated with $\Gamma_0(9)$:
\begin{equation}
E_{k, 9}^{\infty}(z) := \frac{3^k E_k(9 z) - E_k(3 z)}{3^k - 1} \quad \text{for} \; k \geqslant 4.
\end{equation}
Note that we have $E_{k, 9}^{\infty}(z) = E_{k, 3}^{\infty}(3 z)$.\\

\paragraph{\bf For the cusp $0$}
We have $\Gamma_0 = \left\{ \pm \left( \begin{smallmatrix}1 & 0 \\ 9 n & 1\end{smallmatrix} \right) \: ; \: n \in \mathbb{Z} \right\}$ and $\gamma_0 = W_9$, and we have the Eisenstein series for the cusp $0$ associated with $\Gamma_0(9)$:
\begin{equation}
E_{k, 9}^0(z) := \frac{- (E_k(3 z) - E_k(z))}{3^k - 1} \quad \text{for} \; k \geqslant 4.
\end{equation}
Note that we have $E_{k, 9}^0(z) = 3^{-k/2} E_{k, 3}^0(z)$. We also have $\gamma_0^{-1} \ \Gamma_0(9) \ \gamma_0 = \Gamma_0(9)$.\\

\paragraph{\bf For the cusp $-1/3$}
We have $\Gamma_{-1/3} = \left\{ \pm \left( \begin{smallmatrix} 3 n + 1 & 2 n \\ - 9 n & -3 n + 1 \end{smallmatrix} \right) \: ; \: n \in \mathbb{Z} \right\}$ and $\gamma_{-1/3} = W_{9-, 3}$, and we have the Eisenstein series for the cusp $-1/3$ associated with $\Gamma_0(9)$:
\begin{equation}
E_{k, 9}^{-1/3}(z) := \frac{- (E_k(3 z) - E_k(z + 1/3))}{3^k - 1} \quad \text{for} \; k \geqslant 4.
\end{equation}
Note that we have $E_{k, 9}^{-1/3}(z) = 3^{-k/2} E_{k, 3}^0(z + 1/3)$. We also have $\gamma_{-1/3}^{-1} \ \Gamma_0(9) \ \gamma_{-1/3} = \Gamma_0(9)$.\\

\paragraph{\bf For the cusp $1/3$}
We have $\Gamma_{1/3} = \left\{ \pm \left( \begin{smallmatrix} - 3 n + 1 & 2 n \\ - 9 n & 3 n + 1 \end{smallmatrix} \right) \: ; \: n \in \mathbb{Z} \right\}$ and $\gamma_{1/3} = W_{9-, -3}$, and we have the Eisenstein series for the cusp $1/3$ associated with $\Gamma_0(9)$:
\begin{equation}
E_{k, 9}^{1/3}(z) := \frac{- (E_k(3 z) - E_k(z - 1/3))}{3^k - 1} \quad \text{for} \; k \geqslant 4.
\end{equation}
Note that we have $E_{k, 9}^{1/3}(z) = 3^{-k/2} E_{k, 3}^0(z - 1/3)$. We also have $\gamma_{1/3}^{-1} \ \Gamma_0(9) \ \gamma_{1/3} = \Gamma_0(9)$.\\

\paragraph{\bf The space of modular forms}
Let $k$ be an even integer $k \geqslant 4$. We have $M_k(\Gamma_0(9)) = \mathbb{C} E_{k, 9}^{\infty} \oplus \mathbb{C} E_{k, 9}^0 \oplus \mathbb{C} E_{k, 9}^{-1/3} \oplus \mathbb{C} E_{k, 9}^{1/3} \oplus S_k(\Gamma_0(9))$ and $S_k(\Gamma_0(9)) = \Delta_9 M_{k - 4}(\Gamma_0(9))$. Then, we have $M_{4 n + 2}(\Gamma_0(9)) = {E_{2, 9}}' M_{4 n}(\Gamma_0(9)) \oplus \mathbb{C} {E_{2, 3}}' (\Delta_9)^n \oplus \mathbb{C} {E_{2, 9+3}}' (\Delta_9)^n$ and
\begin{align*}
M_{4 n}(\Gamma_0(9)) = &\mathbb{C} (E_{4, 9}^{\infty})^n \oplus \mathbb{C} (E_{4, 9}^{\infty})^{n-1} \Delta_9 \oplus \cdots \oplus \mathbb{C} E_{4, 9}^{\infty} (\Delta_9)^{n-1}\\
 &\oplus \mathbb{C} (E_{4, 9}^0)^n \oplus \mathbb{C} (E_{4, 9}^0)^{n-1} \Delta_9 \oplus \cdots \oplus \mathbb{C} E_{4, 9}^0 (\Delta_9)^{n-1}\\
 &\oplus \mathbb{C} (E_{4, 9}^{-1/3})^n \oplus \mathbb{C} (E_{4, 9}^{-1/3})^{n-1} \Delta_9 \oplus \cdots \oplus \mathbb{C} E_{4, 9}^{-1/3} (\Delta_9)^{n-1}\\
 &\oplus \mathbb{C} (E_{4, 9}^{1/3})^n \oplus \mathbb{C} (E_{4, 9}^{1/3})^{n-1} \Delta_9 \oplus \cdots \oplus \mathbb{C} E_{4, 9}^{1/3} (\Delta_9)^{n-1} \oplus \mathbb{C} (\Delta_9)^n,
\end{align*}

Furthermore, we can write
\begin{equation*}
M_{2 n}(\Gamma_0(9)) = \mathbb{C} (\Delta_9^{\infty})^{2n} \oplus \mathbb{C} (\Delta_9^{\infty})^{2n-1} \Delta_9^0 \oplus \cdots \oplus \mathbb{C} (\Delta_9^0)^{2n}.
\end{equation*}\quad

\paragraph{\bf Hauptmodul}
We define the {\it hauptmodul} of $\Gamma_0(9)$:
\begin{equation}
J_9 := \Delta_9^0 / \Delta_9^{\infty} \: (= \eta^3(z) / \eta^3(9 z)) = \frac{1}{q} - 3 + 5 q^2 - 7 q^5 + 3 q^8 + \cdots,
\end{equation}
where $v_{\infty}(J_9) = -1$ and $v_0(J_9) = 1$. Then, we have
{\small \begin{equation}
\begin{split}
J_9 : \qquad \left\{|z + 1/9| = 1/9, \: -1/6 \leqslant Re(z) \leqslant 0\right\} &\to -3 + [0, 3] \subset \mathbb{R},\\
\left\{|z + 4/9| = 1 / 9, \: -1/2 \leqslant Re(z) \leqslant -1/3\right\} &\to -3 + e^{2 \pi i / 3} [0, 3],\\
\left\{|z - 2/9| = 1 / 9, \: 1/6 \leqslant Re(z) \leqslant 1/3 \right\} &\to -3 + e^{- 2 \pi i / 3} [0, 3].
\end{split}
\end{equation}}

\begin{figure}[hbtp]
\begin{center}
{{\small Lower arcs of $\partial \mathbb{F}_9$}\includegraphics[width=2.5in]{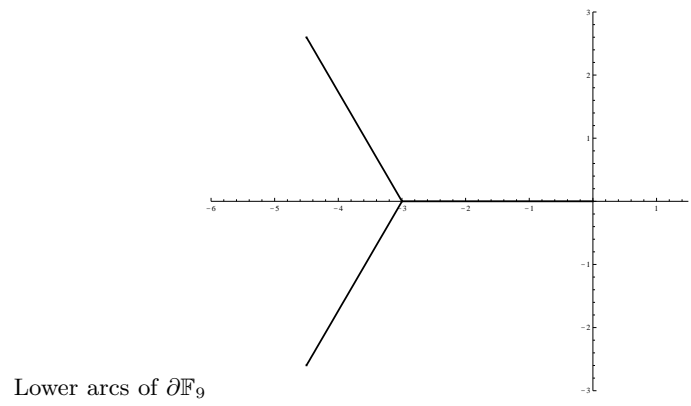}}
\end{center}
\caption{Image by $J_9$}\label{Im-J9B}
\end{figure}

\clearpage

\section{Level $10$}

We have $\Gamma_0(10)+$, $\Gamma_0(10)+10=\Gamma_0^{*}(10)$, $\Gamma_0(10)+5$, $\Gamma_0(10)+2$, and $\Gamma_0(10)-=\Gamma_0(10)$.

We have $W_{10} = \left(\begin{smallmatrix} 0 & - 1 / \sqrt{10} \\ \sqrt{10} & 0 \end{smallmatrix}\right)$, $W_{10, 2} := \left(\begin{smallmatrix} -\sqrt{2} & -1/\sqrt{2} \\ 5\sqrt{2} & 2 \sqrt{2} \end{smallmatrix}\right)$, and $W_{10, 5} := \left(\begin{smallmatrix} -\sqrt{5} & -3/\sqrt{5}\\ 2\sqrt{5} & \sqrt{5}\end{smallmatrix}\right)$, and we denote $\rho_{10, 1} := -3/10 + i/10$, $\rho_{10, 2} := -1/2 + i/(2 \sqrt{5})$, $\rho_{10, 3} := -1/7 + i/(7 \sqrt{10})$, $\rho_{10, 4} := -1/6 + i / (6 \sqrt{5})$, and $\rho_{10, 5} := 3/10 + i /10$. We define
\begin{equation}
\begin{aligned}
\Delta_{10}^{\infty}(z) &:= \sqrt[3]{\eta(z) \eta^{-2}(2 z) \eta^{-5}(5 z) \eta^{10}(10 z)},\\
\Delta_{10}^{-1/2}(z) &:= \sqrt[3]{\eta^{-5}(z) \eta^{10}(2 z) \eta(5 z) \eta^{-2}(10 z)},\\
\Delta_{10}(z) \quad &:= \quad \Delta_{10}^{\infty}(z) \Delta_{10}^0(z) \Delta_{10}^{-1/5}(z) \Delta_{10}^{-1/2}(z)
\end{aligned}
\begin{aligned}
\Delta_{10}^0&(z) := \sqrt[3]{\eta^{10}(z) \eta^{-5}(2 z) \eta^{-2}(5 z) \eta(10 z)}, \\
\Delta_{10}^{-1/5}&(z) := \sqrt[3]{\eta^{-2}(z) \eta(2 z) \eta^{10}(5 z) \eta^{-5}(10 z)},\\
= \quad &\sqrt[3]{\eta^4(z) \eta^4(2 z) \eta^4(5 z) \eta^4(10 z)},
\end{aligned}
\end{equation}
where $\Delta_{10}^{\infty}$, $\Delta_{10}^0$, $\Delta_{10}^{-1/2}$, and $\Delta_{10}^{-1/5}$ are $6$th semimodular forms for $\Gamma_0(10)$ of weight $2/3$ such that $v_{\infty}(\Delta_{10}^{\infty}) = v_0(\Delta_{10}^0) = v_{-1/2}(\Delta_{10}^{-1/2}) = v_{-1/5}(\Delta_{10}^{-1/5}) = 1$, and $\Delta_{10}$ is a $3$rd semimodular form for $\Gamma_0(10)$ of weight $8/3$. Furthermore, we define
\begin{equation}
\begin{split}
{E_{2, 10+10}}'(z) &:= (10 E_2(10 z) - 5 E_2(5 z) - 2 E_2(2 z) + E_2(z)) / 4,\\
{E_{2, 10+5}}'(z) &:= (10 E_2(10 z) - 5 E_2(5 z) + 2 E_2(2 z) - E_2(z)) / 6,\\
{E_{2, 10+2}}'(z) &:= (10 E_2(10 z) + 5 E_2(5 z) - 2 E_2(2 z) - E_2(z)) / 12,
\end{split}
\end{equation}
which are modular forms for $\Gamma_0(10)$ of weight $2$, and we have $v_{i / \sqrt{10}}({E_{2, 10+10}}') = v_{\rho_{10, 1}}({E_{2, 10+10}}') = v_{\rho_{10, 3}}({E_{2, 10+10}}') = 1$, $v_{\rho_{10, 1}}({E_{2, 10+5}}') = v_{\rho_{10, 2}}({E_{2, 10+5}}') = v_{\rho_{10, 4}}({E_{2, 10+5}}') = 1$, and $v_{\rho_{10, 1}}({E_{2, 10+2}}') = v_{\rho_{10, 5}}({E_{2, 10+2}}') = 3$. In addition, we define
\begin{equation*}
{E_{2/3, 10}}' := \sqrt[3]{{E_{2, 10+2}}'}, \quad {E_{8/3, 10}}' := {E_{2, 10+10}}' {E_{2, 10+5}}' / ({E_{2/3, 10}}')^2,
\end{equation*}
which are $3$rd semimodular forms for $\Gamma_0(10)$ of weight $2/3$ and $8/3$, respectively. We also have $v_{\rho_{10, 1}}({E_{2/3, 10}}') = v_{\rho_{10, 5}}({E_{2/3, 10}}') = 3$ and $v_{i / \sqrt{10}}({E_{8/3, 10}}') = v_{\rho_{10, 2}}({E_{8/3, 10}}') = v_{\rho_{10, 3}}({E_{8/3, 10}}') = v_{\rho_{10, 4}}({E_{8/3, 10}}') = 1$.

Furthermore, $\prod_{i = 1}^6 \Delta_{10}^{\kappa_i}$, ${E_{2/3, 10}}' \Delta_{10}^{\kappa_1} \Delta_{10}^{\kappa_2}$, and ${E_{2/3, 10}}' ({E_{8/3, 10}}')^2$ are modular form for $\Gamma_0(10)$ for $\kappa_i \in \{ \infty, 0, -1/2, -1/5 \}$.

\subsection{$\Gamma_0(10)+$}\quad

We have
\begin{equation*}
\Gamma_0(10)+ = \Gamma_0(10)+2,5,10 = \Gamma_0(10) \cup \Gamma_0(10) W_{10, 2} \cup \Gamma_0(10) W_{10, 5} \cup \Gamma_0(10) W_{10}.
\end{equation*}

\paragraph{\bf Fundamental domain}
We have a fundamental domain for $\Gamma_0(10)+$ as follows:
{\small \begin{equation}
\begin{split}
\mathbb{F}_{10+} = &\left\{|z + 1/2| \geqslant 1 / (2 \sqrt{5}), \: - 1/2 \leqslant Re(z) < - 3/10 \right\}
 \bigcup \left\{|z| \geqslant 1 / \sqrt{10}, \: - 3/10 \leqslant Re(z) \leqslant 0 \right\}\\
 &\bigcup \left\{|z| > 1 / \sqrt{10}, \: 0 < Re(z) \leqslant 3/10 \right\}
 \bigcup \left\{|z - 1/2| > 1 / (2 \sqrt{5}), \: 3/10 < Re(z) < 1/2 \right\},
\end{split}
\end{equation}
}where, $W_{10} :  e^{i \theta} / \sqrt{10} \rightarrow e^{i  (\pi - \theta)} / \sqrt{10}$ and $\left(\begin{smallmatrix} -9 & -5 \\ 20 & 11 \end{smallmatrix}\right) W_{10, 5} : e^{i \theta} / (2 \sqrt{5}) + 1/2 \rightarrow e^{i  (\pi - \theta)} / (2 \sqrt{5}) - 1/2$. Then, we have
\begin{equation}
\Gamma_0(10)+ = \langle \left( \begin{smallmatrix} 1 & 1 \\ 0 & 1 \end{smallmatrix} \right), \: W_{10}, \: W_{10, 5} \rangle.
\end{equation}
\begin{figure}[hbtp]
\begin{center}
{{$\mathbb{F}_{10+}$}\includegraphics[width=1.5in]{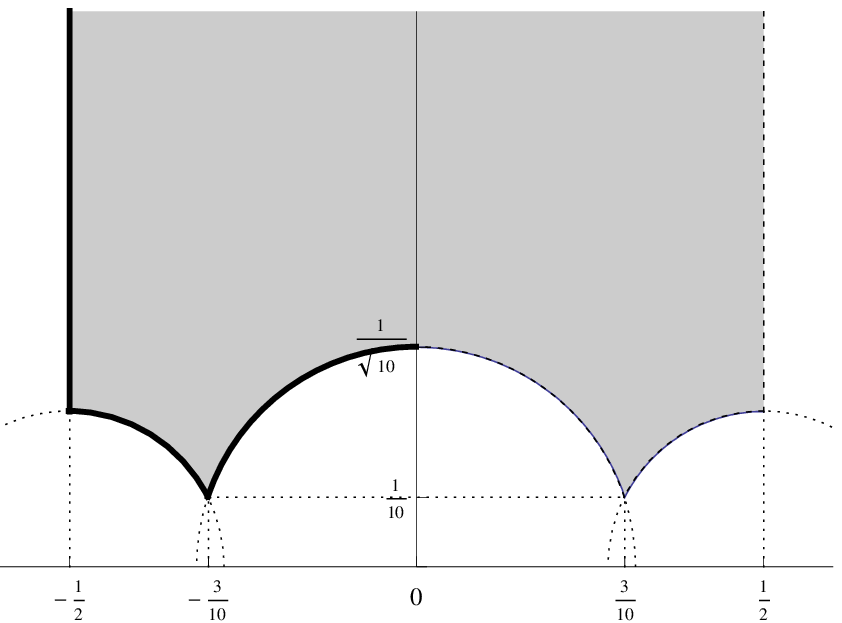}}
\end{center}
\caption{$\Gamma_0(10)+$}
\end{figure}

\paragraph{\bf Valence formula}
The cusp of $\Gamma_0(10)+$ is $\infty$, and the elliptic points are $i / \sqrt{10}$, $\rho_{10, 1} = -3/10 + i/10$, and $\rho_{10, 2} = -1/2 + i/(2 \sqrt{5})$. Let $f$ be a modular function of weight $k$ for $\Gamma_0(10)+$, which is not identically zero. We have
\begin{equation}
v_{\infty}(f) + \frac{1}{2} v_{i / \sqrt{10}}(f) + \frac{1}{4} v_{\rho_{10, 1}} (f) + \frac{1}{2} v_{\rho_{10, 2}} (f) + \sum_{\begin{subarray}{c} p \in \Gamma_0(10)+ \setminus \mathbb{H} \\ p \ne i / \sqrt{10}, \; \rho_{10, 1}, \; \rho_{10, 2}\end{subarray}} v_p(f) = \frac{3 k}{8}.
\end{equation}

Furthermore, the stabilizer of the elliptic point $i / \sqrt{10}$ (resp. $\rho_{10, 1}$, $\rho_{10, 2}$) is $\left\{ \pm I, \pm W_{10} \right\}$\\ (resp. $\left\{ \pm I, \pm \left(\begin{smallmatrix} -3 & -1 \\ 10 & 3 \end{smallmatrix}\right), \pm W_{10, 2}, \pm \left(\begin{smallmatrix} -3 & -1 \\ 10 & 3 \end{smallmatrix}\right) W_{10, 2} \right\}$, $\left\{ \pm I, \pm W_{10, 5} \right\}$)\\

\paragraph{\bf For the cusp $\infty$}
We have $\Gamma_{\infty} = \left\{ \pm \left( \begin{smallmatrix}1 & n \\ 0 & 1\end{smallmatrix} \right) \: ; \: n \in \mathbb{Z} \right\}$, and we have the Eisenstein series associated with $\Gamma_0(10)+$:
\begin{equation}
E_{k, 10+}(z) := \frac{10^{k/2} E_k(10 z) + 5^{k/2} E_k(5 z) + 2^{k/2} E_k(2 z) + E_k(z)}{(5^{k/2} + 1) (2^{k/2} + 1)} \quad \text{for} \; k \geqslant 4.
\end{equation}\quad

\paragraph{\bf The space of modular forms}
Let $k$ be an even integer $k \geqslant 4$. We have $M_k(\Gamma_0(10)+) = \mathbb{C} E_{k, 10+} \oplus S_k(\Gamma_0(10)+)$ and $S_k(\Gamma_0(10)+) = (\mathbb{C} ({E_{2/3, 10}}')^8 \Delta_{10} \oplus \mathbb{C} ({E_{2/3, 10}}')^4 (\Delta_{10})^2 \oplus \mathbb{C} (\Delta_{10})^3) M_{k - 8}(\Gamma_0(10)+)$. Then, we have
\begin{equation*}
M_k(\Gamma_0(10)+) = {E_{\overline{k}, 10+}}' (\mathbb{C} (({E_{2/3, 10}}')^4)^n \oplus \mathbb{C} (({E_{2/3, 10}}')^4)^{n - 1} \Delta_{10} \oplus \cdots \oplus \mathbb{C} (\Delta_{10})^n),
\end{equation*}
where $n = \dim(M_k(\Gamma_0(10)+)) - 1 = \lfloor 3k/8 - 2 (k/4 - \lfloor k/4 \rfloor)\rfloor$, and where ${E_{\overline{k}, 10+}}' := 1$, $({E_{2/3, 10}}')^3 {E_{8/3, 10}}'$, $({E_{2/3, 10}}')^2$, and ${E_{2/3, 10}}' {E_{8/3, 10}}'$, when $k \equiv 0$, $2$, $4$, and $6 \pmod{8}$, respectively.\\

\paragraph{\bf Hauptmodul}
We define the {\it hauptmodul} of $\Gamma_0(10)+$:
\begin{equation}
J_{10+} := ({E_{2/3, 10}}')^4 / \Delta_{10} = \frac{1}{q} + 4 + 22 q + 56 q^2 + 177 q^3 + \cdots,
\end{equation}
where $v_{\infty}(J_{10+}) = -1$ and $v_{\rho_{10, 1}}(J_{10+}) = 4$. Then, we have
\begin{equation}
J_{10+} : \partial \mathbb{F}_{10+} \setminus \{z \in \mathbb{H} \: ; \: Re(z) = \pm 1/2\} \to [- 4, 16] \subset \mathbb{R}.
\end{equation}\quad

\subsection{$\Gamma_0(10)+10 = \Gamma_0^{*}(10)$}

\paragraph{\bf Fundamental domain}
We have a fundamental domain for $\Gamma_0^{*}(10)$ as follows:
{\small \begin{equation}
\begin{split}
\mathbb{F}_{10+10} = &\left\{|z + 9/20| \geqslant 1/20, \: - 1/2 \leqslant Re(z) < - 3/7 \right\}
 \left\{|z + 1/3| \geqslant 1/(3 \sqrt{10}), \: - 3/7 \leqslant Re(z) < - 3/10 \right\}\\
&\bigcup \left\{|z| \geqslant 1/\sqrt{10}, \: - 3/10 \leqslant Re(z) \leqslant 0 \right\}
 \bigcup \left\{|z| > 1/\sqrt{10}, \: 0 < Re(z) \leqslant 3/10 \right\}\\
&\bigcup \left\{|z - 1/3| > 1/(3 \sqrt{10}), \: 3/10 < Re(z) \leqslant 3/7 \right\}
 \bigcup \left\{|z - 9/20| > 1/20, \: 3/7 < Re(z) < 1/2 \right\},
\end{split}
\end{equation}
}where $W_{10} :  e^{i \theta} / \sqrt{10} \rightarrow e^{i  (\pi - \theta)} / \sqrt{10}$, $\left(\begin{smallmatrix} -3 & -1 \\ 10 & 3 \end{smallmatrix} \right) W_{10} :  e^{i \theta} / (3 \sqrt{10}) + 1/3 \rightarrow e^{i  (\pi - \theta)} / (3 \sqrt{10}) - 1/3$, and $\left(\begin{smallmatrix} -9 & 4 \\ 20 & -9 \end{smallmatrix} \right) : (e^{i \theta} + 5) / 12 \rightarrow (e^{i (\pi - \theta)} - 5) / 12$. Then, we have
\begin{equation}
\Gamma_0^{*}(10) = \langle \left( \begin{smallmatrix} 1 & 1 \\ 0 & 1 \end{smallmatrix} \right), \: W_{10}, \: \left(\begin{smallmatrix} -3 & -1 \\ 10 & 3 \end{smallmatrix} \right), \: \left( \begin{smallmatrix} 9 & 4 \\ 20 & 9 \end{smallmatrix} \right) \rangle.
\end{equation}
\begin{figure}[hbtp]
\begin{center}
\includegraphics[width=1.5in]{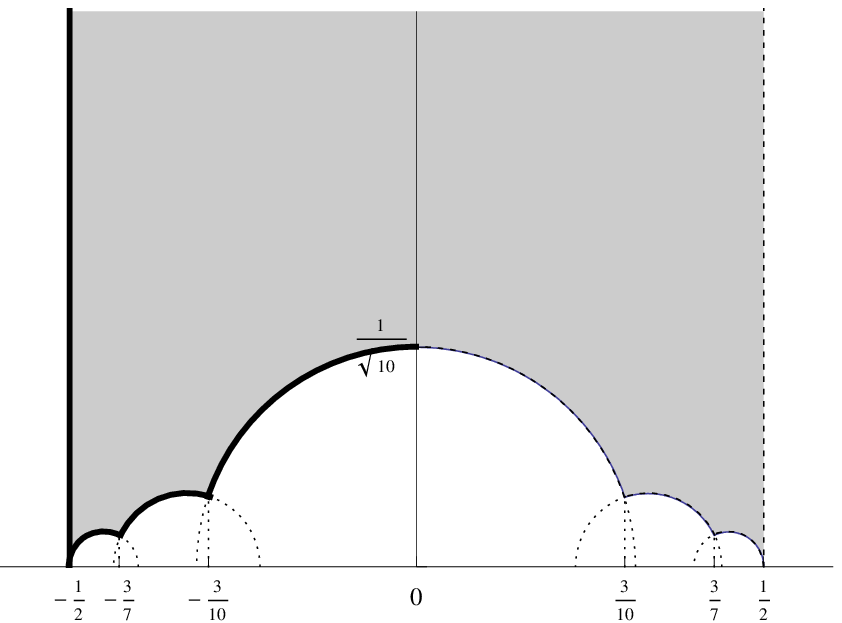}
\end{center}
\caption{$\Gamma_0^{*}(10)$}
\end{figure}

\paragraph{\bf Valence formula}
The cusps of $\Gamma_0^{*}(10)$ are $\infty$ and $-1/2$, and the elliptic points are $i / \sqrt{10}$ and $\rho_{10, 5} = - 3/7 + i / (5 \sqrt{10})$. Let $f$ be a modular function of weight $k$ for $\Gamma_0^{*}(10)$, which is not identically zero. We have
\begin{equation}
v_{\infty}(f) + v_{-1/2}(f) + \frac{1}{2} v_{i / \sqrt{10}} (f) + \frac{1}{2} v_{\rho_{10, 1}} (f) + \frac{1}{2} v_{\rho_{10, 3}} (f) + \sum_{\begin{subarray}{c} p \in \Gamma_0^{*}(10) \setminus \mathbb{H} \\ p \ne i / \sqrt{10}, \rho_{10, 1}, \rho_{10, 3}\end{subarray}} v_p(f) = \frac{3 k}{4}.
\end{equation}

Furthermore, the stabilizer of the elliptic point $i / \sqrt{10}$ (resp. $\rho_{10, 1}$, $\rho_{10, 3}$) is $\left\{ \pm I, \pm W_{10} \right\}$ (resp. $\left\{ \pm I, \pm \left(\begin{smallmatrix} -3 & -1 \\ 10 & 3 \end{smallmatrix} \right) \right\}$, $\left\{ \pm I, \pm \left( \begin{smallmatrix} -13 & 3 \\ 30 & -7\end{smallmatrix} \right) W_{10} \right\}$).\\

\paragraph{\bf For the cusp $\infty$}
We have $\Gamma_{\infty} = \left\{ \pm \left( \begin{smallmatrix}1 & n \\ 0 & 1\end{smallmatrix} \right) \: ; \: n \in \mathbb{Z} \right\}$, and we have the Eisenstein series for the cusp $\infty$ associated with $\Gamma_0^{*}(10)$:
\begin{equation}
E_{k, 10+10}^{\infty}(z) := \frac{(10^{k/2} + 1) (10^{k/2} E_k(10 z) + E_k(z)) - (5^{k/2} + 2^{k/2}) (5^{k/2} E_k(5 z) + 2^{k/2} E_k(2 z))}{(5^k - 1)(2^k - 1)} \quad \text{for} \; k \geqslant 4.
\end{equation}\quad

\paragraph{\bf For the cusp $-1/2$}
We have $\Gamma_{-1/2} = \left\{ \pm \left( \begin{smallmatrix}10 n + 1 & 5 n \\ - 20 n & - 10 n + 1\end{smallmatrix} \right) \: ; \: n \in \mathbb{Z} \right\}$ and $\gamma_{-1/2} = W_{10, 5}$, and we have the Eisenstein series for the cusp $-1/2$ associated with $\Gamma_0^{*}(10)$:
\begin{equation}
E_{k, 10+10}^{-1/2}(z) := \frac{- (5^{k/2} + 2^{k/2}) (10^{k/2} E_k(10 z) + E_k(z)) + (10^{k/2} + 1) (5^{k/2} E_k(5 z) + 2^{k/2} E_k(2 z))}{(5^k - 1)(2^k - 1)} \quad \text{for} \; k \geqslant 4.
\end{equation}
We also have $\gamma_{-1/2}^{-1} \ \Gamma_0^{*}(10) \ \gamma_{-1/2} = \Gamma_0^{*}(10)$.\\

\paragraph{\bf The space of modular forms}
We define
\begin{equation*}
\Delta_{10+10}^{\infty} := \Delta_{10}^{\infty} \Delta_{10}^0, \quad \Delta_{10+10}^{-1/2} := \Delta_{10}^{-1/2} \Delta_{10}^{-1/5},
\end{equation*}
which are $2$nd semimodular forms for $\Gamma_0^{*}(10)$ of weight $2$. Furthermore, we define
\begin{align*}
\Delta_{8, 1, 10+10} &:= \Delta_{10+10}^{\infty} (\Delta_{10+10}^{-1/2})^5,
&\Delta_{8, 2, 10+10} &:= (\Delta_{10+10}^{\infty})^5 \Delta_{10+10}^{-1/2},\\
\Delta_{8, 3, 10+10} &:= (\Delta_{10+10}^{\infty})^2 (\Delta_{10+10}^{-1/2})^4,
&\Delta_{8, 4, 10+10} &:= (\Delta_{10+10}^{\infty})^4 (\Delta_{10+10}^{-1/2})^2.
\end{align*}

Now, we have
\begin{gather*}
M_k(\Gamma_0^{*}(10)) = \mathbb{C} E_{k, 10+10}^{\infty} \oplus \mathbb{C} E_{k, 10+10}^{-1/2} \oplus S_k(\Gamma_0^{*}(10)),\\
S_k(\Gamma_0^{*}(10)) = (\mathbb{C} \Delta_{8, 1, 10+10} \oplus \mathbb{C} \Delta_{8, 2, 10+10} \oplus \mathbb{C} \Delta_{8, 3, 10+10} \oplus \mathbb{C} \Delta_{8, 4, 10+10} \oplus \mathbb{C} (\Delta_{10})^3) M_{k - 8}(\Gamma_0^{*}(10))
\end{gather*}
for every even integer $k \geqslant 4$. Then, we have $M_{4 n + 2}(\Gamma_0^{*}(10)) = {E_{2, 10+10}}' M_{4 n}(\Gamma_0^{*}(10))$ and
\begin{align*}
M_{8 n}(\Gamma_0^{*}(10)) &= \mathbb{C} (E_{8, 10+10}^{\infty})^n \oplus \mathbb{C} (E_{8, 10+10}^{\infty})^{n-1} \Delta_{8, 1, 10+10} \oplus \mathbb{C} (E_{8, 10+10}^{\infty})^{n-1} \Delta_{8, 3, 10+10}\\
&\qquad \qquad \qquad \qquad \oplus \mathbb{C} (E_{8, 10+10}^{\infty})^{n-1} (\Delta_{10})^3 \oplus \cdots \oplus \mathbb{C} \Delta_{8, 3, 10+10} (\Delta_{10})^{3(n-1)}\\
&\oplus \mathbb{C} (E_{8, 10+10}^{-1/2})^n \oplus \mathbb{C} (E_{8, 10+10}^{-1/2})^{n-1} \Delta_{8, 2, 10+10} \oplus \mathbb{C} (E_{8, 10+10}^{-1/2})^{n-1} \Delta_{8, 4, 10+10}\\
&\qquad \qquad \qquad \qquad \oplus \mathbb{C} (E_{8, 10+10}^{-1/2})^{n-1} (\Delta_{10})^3 \oplus \cdots \oplus \mathbb{C} \Delta_{8, 4, 10+10} (\Delta_{10})^{3(n-1)} \oplus \mathbb{C} (\Delta_{10})^{3n},\\
M_{8 n + 4}(\Gamma_0^{*}(10)) &= E_{4, 10+10}^{\infty} (\mathbb{C} (E_{8, 10+10}^{\infty})^n \oplus \mathbb{C} (E_{8, 10+10}^{\infty})^{n-1} \Delta_{8, 1, 10+10} \oplus \mathbb{C} (E_{8, 10+10}^{\infty})^{n-1} \Delta_{8, 3, 10+10}\\
&\qquad \qquad \qquad \qquad \oplus \mathbb{C} (E_{8, 10+10}^{\infty})^{n-1} (\Delta_{10})^3 \oplus \cdots \oplus \mathbb{C} \Delta_{8, 3, 10+10} (\Delta_{10})^{3(n-1)} \oplus \mathbb{C} (\Delta_{10})^{3n})\\
&\oplus E_{4, 10+10}^{-1/2} (\mathbb{C} (E_{8, 10+10}^{-1/2})^n \oplus \mathbb{C} (E_{8, 10+10}^{-1/2})^{n-1} \Delta_{8, 2, 10+10} \oplus \mathbb{C} (E_{8, 10+10}^{-1/2})^{n-1} \Delta_{8, 4, 10+10}\\
&\qquad \qquad \qquad \qquad \oplus \mathbb{C} (E_{8, 10+10}^{-1/2})^{n-1} (\Delta_{10})^3 \oplus \cdots \oplus \mathbb{C} \Delta_{8, 4, 10+10} (\Delta_{10})^{3(n-1)} \oplus \mathbb{C} (\Delta_{10})^{3n})\\
&\qquad \qquad \oplus (\mathbb{C} \Delta_{10+10}^{\infty} (\Delta_{10+10}^{-1/2})^2 \oplus \mathbb{C} (\Delta_{10+10}^{\infty})^2 \Delta_{10+10}^{-1/2}) (\Delta_{10})^{3n}.
\end{align*}

Furthremore, we can write
\begin{equation*}
M_{4 n}(\Gamma_0^{*}(10)) = \mathbb{C} (\Delta_{10+10}^{\infty})^{3n} \oplus \mathbb{C} (\Delta_{10+10}^{\infty})^{3n-1} \Delta_{10+10}^{-1/2} \oplus \cdots \oplus \mathbb{C} (\Delta_{10+10}^{-1/2})^{3n}.
\end{equation*}\quad

\paragraph{\bf Hauptmodul}
We define the {\it hauptmodul} of $\Gamma_0^{*}(10)$:
\begin{equation}
J_{10+10} := \Delta_{10+10}^{-1/2} / \Delta_{10+10}^{\infty} \: (= \eta^{-6}(z) \eta^6(2 z) \eta^{-6}(5 z) \eta^6(10 z)) = \frac{1}{q} + 6 + 21 q + 62 q^2 + 162 q^3 + \cdots,
\end{equation}
where $v_{\infty}(J_{10+10}) = -1$ and $v_{-1/2}(J_{10+10}) = 1$. Then, we have
\begin{equation}
J_{10+10} : \partial \mathbb{F}_{10+10} \setminus \{z \in \mathbb{H} \: ; \: Re(z) = \pm 1/2\} \to [0, 9 + 4 \sqrt{5}] \subset \mathbb{R}.
\end{equation}\quad

\subsection{$\Gamma_0(10)+5$}

\paragraph{\bf Fundamental domain}
We have a fundamental domain for $\Gamma_0(10)+5$ as follows:
{\small \begin{equation}
\begin{split}
\mathbb{F}_{10+5} = &\left\{|z + 1/2| \geqslant 1 / (2 \sqrt{5}), \: - 1/2 \leqslant Re(z) < -3/10 \right\}
 \bigcup \left\{|z + 1/4| \geqslant 1 / (4 \sqrt{5}), \: - 3/10 \leqslant Re(z) < -1/6 \right\}\\
&\bigcup \left\{|z + 1/10| \geqslant 1 / 10, \: - 1/6 \leqslant Re(z) \leqslant 0 \right\}
 \bigcup \left\{|z - 1/10| > 1 / 10, \: 0 < Re(z) \leqslant 1/6 \right\}\\
&\bigcup \left\{|z - 1/4| > 1 / (4 \sqrt{5}), \: 1/6 < Re(z) \leqslant 3/10 \right\}
 \bigcup \left\{|z - 1/2| > 1 / (2 \sqrt{5}), \: 3/10 < Re(z) < 1/2 \right\},
\end{split}
\end{equation}
}where $\left(\begin{smallmatrix} -1 & 0 \\ 10 & -1 \end{smallmatrix}\right) : (e^{i \theta} + 1) / 10 \rightarrow (e^{i (\pi - \theta)} - 1) / 10$, $\left(\begin{smallmatrix} -7 & -4 \\ 30 & 17 \end{smallmatrix}\right) W_{10, 5} : e^{i \theta} / (4 \sqrt{5}) + 1/4 \rightarrow e^{i  (\pi - \theta)} / (4 \sqrt{5}) - 1/4$, and $\left(\begin{smallmatrix} -9 & -5 \\ 20 & 11 \end{smallmatrix}\right) W_{10, 5} : e^{i \theta} / (2 \sqrt{5}) + 1/2 \rightarrow e^{i  (\pi - \theta)} / (2 \sqrt{5}) - 1/2$. Then, we have
\begin{equation}
\Gamma_0(10)+5 = \langle \left( \begin{smallmatrix} 1 & 1 \\ 0 & 1 \end{smallmatrix} \right), \: W_{10, 5}, \: \left( \begin{smallmatrix} 1 & 0 \\ 10 & 1 \end{smallmatrix} \right), \: \left(\begin{smallmatrix} -3 & -1 \\ 10 & 3 \end{smallmatrix} \right) \rangle.
\end{equation}
\begin{figure}[hbtp]
\begin{center}
\includegraphics[width=1.5in]{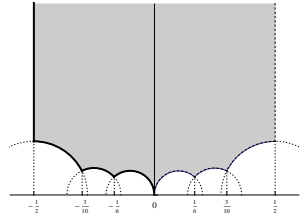}
\end{center}
\caption{$\Gamma_0(10)+5$}
\end{figure}

\paragraph{\bf Valence formula}
The cusps of $\Gamma_0(10)+5$ are $\infty$ and $0$, and the elliptic points are $\rho_{10, 1} = -3/10 + i/10$, $\rho_{10, 2} = -1/2 + i/(2 \sqrt{5})$,  and $\rho_{10, 4} = -1/6 + i / (6 \sqrt{5})$. Let $f$ be a modular function of weight $k$ for $\Gamma_0(10)+5$, which is not identically zero. We have
\begin{equation}
v_{\infty}(f) + v_0(f) + \frac{1}{2} v_{\rho_{10, 1}} (f) + \frac{1}{2} v_{\rho_{10, 2}} (f) + \frac{1}{2} v_{\rho_{10, 4}} (f) + \sum_{\begin{subarray}{c} p \in \Gamma_0(10)+5 \setminus \mathbb{H} \\ p \ne \rho_{10, 1}, \rho_{10, 2}, \rho_{10, 4}\end{subarray}} v_p(f) = \frac{3k}{4}.
\end{equation}

Furthermore, the stabilizer of the elliptic point $\rho_{10, 1}$ (resp. $\rho_{10, 2}$, $\rho_{10, 4}$) is $\left\{ \pm I, \pm \left( \begin{smallmatrix} -3 & -1 \\ 10 & 3 \end{smallmatrix} \right) \right\}$ (resp. $\left\{ \pm I, \pm W_{10, 5} \right\}$, $\left\{ \pm I, \pm \left( \begin{smallmatrix} -3 & -2 \\ 20 & 13 \end{smallmatrix} \right) W_{10, 5} \right\}$).\\

\paragraph{\bf For the cusp $\infty$}
We have $\Gamma_{\infty} = \left\{ \pm \left( \begin{smallmatrix}1 & n \\ 0 & 1\end{smallmatrix} \right) \: ; \: n \in \mathbb{Z} \right\}$, and we have the Eisenstein series for the cusp $\infty$ associated with $\Gamma_0(10)+5$:
\begin{equation}
E_{k, 10+5}^{\infty}(z) := \frac{2^k 5^{k/2} E_k(10 z) - 5^{k/2} E_k(5 z) + 2^k E_k(2 z) - E_k(z)}{(5^{k/2} + 1)(2^k - 1)} \quad \text{for} \; k \geqslant 4.
\end{equation}\quad

\paragraph{\bf For the cusp $0$}
We have $\Gamma_0 = \left\{ \pm \left( \begin{smallmatrix}1 & 0 \\ 10 n & 1\end{smallmatrix} \right) \: ; \: n \in \mathbb{Z} \right\}$ and $\gamma_0 = W_{10}$, and we have the Eisenstein series for the cusp $0$ associated with $\Gamma_0(10)+5$:
\begin{equation}
E_{k, 10+5}^0(z) := \frac{- 2^{k/2} (5^{k/2} E_k(10 z) - 5^{k/2} E_k(5 z) + E_k(2 z) - E_k(z))}{(5^{k/2} + 1)(2^k - 1)} \quad \text{for} \; k \geqslant 4.
\end{equation}
We also have $\gamma_0^{-1} \ (\Gamma_0(10)+5) \ \gamma_0 = \Gamma_0(10)+5$.\\

\paragraph{\bf The space of modular forms}
We define
\begin{equation*}
\Delta_{10+5}^{\infty} := \Delta_{10}^{\infty} \Delta_{10}^{-1/2}, \quad \Delta_{10+5}^0 := \Delta_{10}^0 \Delta_{10}^{-1/5},
\end{equation*}
which are $2$nd semimodular forms for $\Gamma_0(10)+5$ of weight $2$. Furthermore, we define
\begin{align*}
\Delta_{8, 1, 10+5} &:= \Delta_{10+5}^{\infty} (\Delta_{10+5}^0)^5,
&\Delta_{8, 2, 10+5} &:= (\Delta_{10+5}^{\infty})^5 \Delta_{10+5}^0,\\
\Delta_{8, 3, 10+5} &:= (\Delta_{10+5}^{\infty})^2 (\Delta_{10+5}^0)^4,
&\Delta_{8, 4, 10+5} &:= (\Delta_{10+5}^{\infty})^4 (\Delta_{10+5}^0)^2.
\end{align*}

Now, we have
\begin{gather*}
M_k(\Gamma_0(10)+5) = \mathbb{C} E_{k, 10+5}^{\infty} \oplus \mathbb{C} E_{k, 10+5}^0 \oplus S_k(\Gamma_0(10)+5),\\
S_k(\Gamma_0(10)+5) = (\mathbb{C} \Delta_{8, 1, 10+5} \oplus \mathbb{C} \Delta_{8, 2, 10+5} \oplus \mathbb{C} \Delta_{8, 3, 10+5} \oplus \mathbb{C} \Delta_{8, 4, 10+5} \oplus \mathbb{C} (\Delta_{10})^3) M_{k - 8}(\Gamma_0(10)+5)
\end{gather*}
for every even integer $k \geqslant 4$. Then, we have $M_{4 n + 2}(\Gamma_0(10)+5) = {E_{2, 10+5}}' M_{4 n}(\Gamma_0(10)+5)$ and
\begin{align*}
M_{8 n}(\Gamma_0(10)+5) &= \mathbb{C} (E_{8, 10+5}^{\infty})^n \oplus \mathbb{C} (E_{8, 10+5}^{\infty})^{n-1} \Delta_{8, 1, 10+5} \oplus \mathbb{C} (E_{8, 10+5}^{\infty})^{n-1} \Delta_{8, 3, 10+5}\\
&\qquad \qquad \qquad \qquad \oplus \mathbb{C} (E_{8, 10+5}^{\infty})^{n-1} (\Delta_{10})^3 \oplus \cdots \oplus \mathbb{C} \Delta_{8, 3, 10+5} (\Delta_{10})^{3(n-1)}\\
&\oplus \mathbb{C} (E_{8, 10+5}^0)^n \oplus \mathbb{C} (E_{8, 10+5}^0)^{n-1} \Delta_{8, 2, 10+5} \oplus \mathbb{C} (E_{8, 10+5}^0)^{n-1} \Delta_{8, 4, 10+5}\\
&\qquad \qquad \qquad \qquad \oplus \mathbb{C} (E_{8, 10+5}^0)^{n-1} (\Delta_{10})^3 \oplus \cdots \oplus \mathbb{C} \Delta_{8, 4, 10+5} (\Delta_{10})^{3(n-1)} \oplus \mathbb{C} (\Delta_{10})^{3n},\\
M_{8 n + 4}(\Gamma_0(10)+5) &= E_{4, 10+5}^{\infty} (\mathbb{C} (E_{8, 10+5}^{\infty})^n \oplus \mathbb{C} (E_{8, 10+5}^{\infty})^{n-1} \Delta_{8, 1, 10+5} \oplus \mathbb{C} (E_{8, 10+5}^{\infty})^{n-1} \Delta_{8, 3, 10+5}\\
&\qquad \qquad \qquad \qquad \oplus \mathbb{C} (E_{8, 10+5}^{\infty})^{n-1} (\Delta_{10})^3 \oplus \cdots \oplus \mathbb{C} \Delta_{8, 3, 10+5} (\Delta_{10})^{3(n-1)} \oplus \mathbb{C} (\Delta_{10})^{3n})\\
&\oplus E_{4, 10+5}^0 (\mathbb{C} (E_{8, 10+5}^0)^n \oplus \mathbb{C} (E_{8, 10+5}^0)^{n-1} \Delta_{8, 2, 10+5} \oplus \mathbb{C} (E_{8, 10+5}^0)^{n-1} \Delta_{8, 4, 10+5}\\
&\qquad \qquad \qquad \qquad \oplus \mathbb{C} (E_{8, 10+5}^0)^{n-1} (\Delta_{10})^3 \oplus \cdots \oplus \mathbb{C} \Delta_{8, 4, 10+5} (\Delta_{10})^{3(n-1)} \oplus \mathbb{C} (\Delta_{10})^{3n})\\
&\qquad \qquad \oplus (\mathbb{C} \Delta_{10+5}^{\infty} (\Delta_{10+5}^0)^2 \oplus \mathbb{C} (\Delta_{10+5}^{\infty})^2 \Delta_{10+5}^0) (\Delta_{10})^{3n}.
\end{align*}

Furthremore, we can write
\begin{equation*}
M_{4 n}(\Gamma_0(10)+5) = \mathbb{C} (\Delta_{10+5}^{\infty})^{3n} \oplus \mathbb{C} (\Delta_{10+5}^{\infty})^{3n-1} \Delta_{10+5}^0 \oplus \cdots \oplus \mathbb{C} (\Delta_{10+5}^0)^{3n}.
\end{equation*}\quad

\paragraph{\bf Hauptmodul}
We define the {\it hauptmodul} of $\Gamma_0(10)+5$:
\begin{equation}
J_{10+5} := \Delta_{10+5}^0 / \Delta_{10+5}^{\infty} \: (= \eta^4(z) \eta^{-4}(2 z) \eta^4(5 z) \eta^{-4}(10 z)) = \frac{1}{q} - 4 + 6 q - 8 q^2 + 17 q^3 - \cdots,
\end{equation}
where $v_{\infty}(J_{10+5}) = -1$ and $v_0(J_{10+5}) = 1$. Then, we have
\begin{equation}
J_{10+5} : \partial \mathbb{F}_{10+5} \setminus \{z \in \mathbb{H} \: ; \: Re(z) = \pm 1 / 2\} \to [-6 - 2 \sqrt{5}, 0] \subset \mathbb{R}.
\end{equation}\quad

\subsection{$\Gamma_0(10)+2$}

\paragraph{\bf Fundamental domain}
We have a fundamental domain for $\Gamma_0(10)+2$ as follows:
{\small \begin{equation}
\begin{split}
\mathbb{F}_{10+2} = &\left\{|z + 2/5| \geqslant 1 / (5 \sqrt{2}), \: - 1/2 \leqslant Re(z) \leqslant -3/10 \right\}
 \bigcup \left\{|z + 1/5| > 1 / (5 \sqrt{2}), \: -3/10 < Re(z) < -1/10 \right\}\\
 &\bigcup \left\{|z + 1/10| \geqslant 1/10, \: 0 \leqslant Re(z) \leqslant 0 \right\}
 \bigcup \left\{|z - 1/10| > 1/10, \: 0 < Re(z) < 1/10 \right\}\\
 &\bigcup \left\{|z - 1/5| \geqslant 1 / (5 \sqrt{2}), \: 1/10 \leqslant Re(z) \leqslant 3/10 \right\}
 \bigcup \left\{|z - 2/5| > 1 / (5 \sqrt{2}), \: 3/10 < Re(z) < 1/2 \right\},
\end{split}
\end{equation}
}where $\left(\begin{smallmatrix} -1 & 0 \\ 10 & -1 \end{smallmatrix}\right) : (e^{i \theta} + 1) / 10 \rightarrow (e^{i (\pi - \theta)} - 1) / 10$, $\left(\begin{smallmatrix} -3 & -1 \\ 10 & 3 \end{smallmatrix}\right) W_{10, 5} : e^{i \theta} / (5 \sqrt{2}) - 1/5 \rightarrow e^{i  (\pi - \theta)} / (5 \sqrt{2}) - 2/5$, and $\left(\begin{smallmatrix} 9 & 2 \\ 40 & 9 \end{smallmatrix}\right) W_{10, 5} : e^{i \theta} / (5 \sqrt{2}) + 2/5 \rightarrow e^{i  (\pi - \theta)} / (5 \sqrt{2}) + 1/5$. Then, we have
\begin{equation}
\Gamma_0(10)+2 = \langle \left( \begin{smallmatrix} 1 & 1 \\ 0 & 1 \end{smallmatrix} \right), \: W_{10, 2}, \: \left( \begin{smallmatrix} 1 & 0 \\ 10 & 1 \end{smallmatrix} \right), \: \left(\begin{smallmatrix} 9 & 2 \\ 40 & 9 \end{smallmatrix}\right) \rangle.
\end{equation}
\begin{figure}[hbtp]
\begin{center}
\includegraphics[width=1.5in]{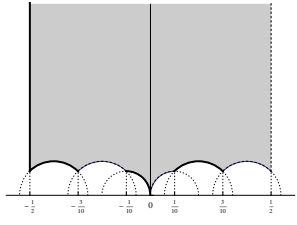}
\end{center}
\caption{$\Gamma_0(10)+2$}
\end{figure}

\paragraph{\bf Valence formula}
The cusps of $\Gamma_0(10)+2$ are $\infty$ and $0$, and the elliptic points are $\rho_{10, 1} = -3/10 + i/10$ and $\rho_{10, 5} = 3/10 + i/10$. Let $f$ be a modular function of weight $k$ for $\Gamma_0(10)+2$, which is not identically zero. We have
\begin{equation}
v_{\infty}(f) + v_0(f) + \frac{1}{4} v_{\rho_{10, 1}} (f) + \frac{1}{4} v_{\rho_{10, 5}} (f) + \sum_{\begin{subarray}{c} p \in \Gamma_0(10)+2 \setminus \mathbb{H} \\ p \ne \rho_{10, 1}, \rho_{10, 5}\end{subarray}} v_p(f) = \frac{3k}{4}.
\end{equation}

Furthermore, the stabilizer of the elliptic point $\rho_{10, 1}$ (resp. $\rho_{10, 5}$) is $\left\{ \pm I, \pm \left(\begin{smallmatrix} -3 & -1 \\ 10 & 3 \end{smallmatrix}\right), \pm W_{10, 2}, \pm \left(\begin{smallmatrix} -3 & -1 \\ 10 & 3 \end{smallmatrix}\right) W_{10, 2} \right\}$ (resp. $\left\{ \pm I, \pm \left(\begin{smallmatrix} 3 & -1 \\ 10 & -3 \end{smallmatrix}\right), \pm \left(\begin{smallmatrix} 13 & 3 \\ 30 & 7 \end{smallmatrix}\right) W_{10, 2}, \pm \left(\begin{smallmatrix} 9 & 2 \\ 40 & 9 \end{smallmatrix}\right) W_{10, 2} \right\}$).\\

\paragraph{\bf For the cusp $\infty$}
We have $\Gamma_{\infty} = \left\{ \pm \left( \begin{smallmatrix}1 & n \\ 0 & 1\end{smallmatrix} \right) \: ; \: n \in \mathbb{Z} \right\}$, and we have the Eisenstein series for the cusp $\infty$ associated with $\Gamma_0(10)+2$:
\begin{equation}
E_{k, 10+2}^{\infty}(z) := \frac{2^{k/2} 5^k E_k(10 z) + 5^k E_k(5 z) - 2^{k/2} E_k(2 z) - E_k(z)}{(5^k - 1)(2^{k/2} + 1)} \quad \text{for} \; k \geqslant 4.
\end{equation}\quad

\paragraph{\bf For the cusp $0$}
We have $\Gamma_0 = \left\{ \pm \left( \begin{smallmatrix}1 & 0 \\ 10 n & 1\end{smallmatrix} \right) \: ; \: n \in \mathbb{Z} \right\}$ and $\gamma_0 = W_{10}$, and we have the Eisenstein series for the cusp $0$ associated with $\Gamma_0(10)+2$:
\begin{equation}
E_{k, 10+2}^0(z) := \frac{- 5^{k/2} (2^{k/2} E_k(10 z) + E_k(5 z) - 2^{k/2} E_k(2 z) - E_k(z))}{(5^k - 1)(2^{k/2} + 1)} \quad \text{for} \; k \geqslant 4.
\end{equation}
We also have $\gamma_0^{-1} \ (\Gamma_0(10)+2) \ \gamma_0 = \Gamma_0(10)+2$.\\

\paragraph{\bf The space of modular forms}
We define
\begin{equation*}
\Delta_{10+2}^{\infty} := \Delta_{10}^{\infty} \Delta_{10}^{-1/5}, \quad \Delta_{10+2}^0 := \Delta_{10}^0 \Delta_{10}^{-1/2},
\end{equation*}
which are $2$nd semimodular forms for $\Gamma_0(10)+2$ of weight $2$. Furthermore, we define
\begin{align*}
\Delta_{8, 1, 10+2} &:= \Delta_{10+2}^{\infty} (\Delta_{10+2}^0)^5,
&\Delta_{8, 2, 10+2} &:= (\Delta_{10+2}^{\infty})^5 \Delta_{10+2}^0,\\
\Delta_{8, 3, 10+2} &:= (\Delta_{10+2}^{\infty})^2 (\Delta_{10+2}^0)^4,
&\Delta_{8, 4, 10+2} &:= (\Delta_{10+2}^{\infty})^4 (\Delta_{10+2}^0)^2.
\end{align*}

Now, we have
\begin{gather*}
M_k(\Gamma_0(10)+2) = \mathbb{C} E_{k, 10+2}^{\infty} \oplus \mathbb{C} E_{k, 10+2}^0 \oplus S_k(\Gamma_0(10)+2),\\
S_k(\Gamma_0(10)+2) = (\mathbb{C} \Delta_{8, 1, 10+2} \oplus \mathbb{C} \Delta_{8, 2, 10+2} \oplus \mathbb{C} \Delta_{8, 3, 10+2} \oplus \mathbb{C} \Delta_{8, 4, 10+2} \oplus \mathbb{C} (\Delta_{10})^3) M_{k - 8}(\Gamma_0(10)+2)
\end{gather*}
for every even integer $k \geqslant 4$. Then, we have $M_{8 n + 2}(\Gamma_0(10)+2) = {E_{2, 10+2}}' M_{8 n}(\Gamma_0(10)+2)$ and
\begin{align*}
M_{8 n}(\Gamma_0(10)+2) &= \mathbb{C} (E_{8, 10+2}^{\infty})^n \oplus \mathbb{C} (E_{8, 10+2}^{\infty})^{n-1} \Delta_{8, 1, 10+2} \oplus \mathbb{C} (E_{8, 10+2}^{\infty})^{n-1} \Delta_{8, 3, 10+2}\\
&\qquad \qquad \qquad \qquad \oplus \mathbb{C} (E_{8, 10+2}^{\infty})^{n-1} (\Delta_{10})^3 \oplus \cdots \oplus \mathbb{C} \Delta_{8, 3, 10+2} (\Delta_{10})^{3(n-1)}\\
&\oplus \mathbb{C} (E_{8, 10+2}^0)^n \oplus \mathbb{C} (E_{8, 10+2}^0)^{n-1} \Delta_{8, 2, 10+2} \oplus \mathbb{C} (E_{8, 10+2}^0)^{n-1} \Delta_{8, 4, 10+2}\\
&\qquad \qquad \qquad \qquad \oplus \mathbb{C} (E_{8, 10+2}^0)^{n-1} (\Delta_{10})^3 \oplus \cdots \oplus \mathbb{C} \Delta_{8, 4, 10+2} (\Delta_{10})^{3(n-1)} \oplus \mathbb{C} (\Delta_{10})^{3n},
\end{align*}
\begin{align*}
M_{8 n + 4}(\Gamma_0(10)+2) &= E_{4, 10+2}^{\infty} (\mathbb{C} (E_{8, 10+2}^{\infty})^n \oplus \mathbb{C} (E_{8, 10+2}^{\infty})^{n-1} \Delta_{8, 1, 10+2} \oplus \mathbb{C} (E_{8, 10+2}^{\infty})^{n-1} \Delta_{8, 3, 10+2}\\
&\qquad \qquad \qquad \qquad \oplus \mathbb{C} (E_{8, 10+2}^{\infty})^{n-1} (\Delta_{10})^3 \oplus \cdots \oplus \mathbb{C} \Delta_{8, 3, 10+2} (\Delta_{10})^{3(n-1)} \oplus \mathbb{C} (\Delta_{10})^{3n})\\
&\oplus E_{4, 10+2}^0 (\mathbb{C} (E_{8, 10+2}^0)^n \oplus \mathbb{C} (E_{8, 10+2}^0)^{n-1} \Delta_{8, 2, 10+2} \oplus \mathbb{C} (E_{8, 10+2}^0)^{n-1} \Delta_{8, 4, 10+2}\\
&\qquad \qquad \qquad \qquad \oplus \mathbb{C} (E_{8, 10+2}^0)^{n-1} (\Delta_{10})^3 \oplus \cdots \oplus \mathbb{C} \Delta_{8, 4, 10+2} (\Delta_{10})^{3(n-1)} \oplus \mathbb{C} (\Delta_{10})^{3n})\\
&\qquad \oplus ({E_{2/3, 10}}')^2 (\mathbb{C} (\Delta_{10+2}^{\infty})^2 \oplus \mathbb{C} (\Delta_{10+2}^0)^2 \oplus \mathbb{C} \Delta_{10}) (\Delta_{10})^{3n},\\
M_{8 n + 6}(\Gamma_0(10)+2) &= E_{6, 10+2}^{\infty} (\mathbb{C} (E_{8, 10+2}^{\infty})^n \oplus \mathbb{C} (E_{8, 10+2}^{\infty})^{n-1} \Delta_{8, 1, 10+2} \oplus \mathbb{C} (E_{8, 10+2}^{\infty})^{n-1} \Delta_{8, 3, 10+2}\\
&\qquad \qquad \qquad \qquad \oplus \mathbb{C} (E_{8, 10+2}^{\infty})^{n-1} (\Delta_{10})^3 \oplus \cdots \oplus \mathbb{C} \Delta_{8, 3, 10+2} (\Delta_{10})^{3(n-1)} \oplus \mathbb{C} (\Delta_{10})^{3n})\\
&\oplus E_{6, 10+2}^0 (\mathbb{C} (E_{8, 10+2}^0)^n \oplus \mathbb{C} (E_{8, 10+2}^0)^{n-1} \Delta_{8, 2, 10+2} \oplus \mathbb{C} (E_{8, 10+2}^0)^{n-1} \Delta_{8, 4, 10+2}\\
&\qquad \qquad \qquad \qquad \oplus \mathbb{C} (E_{8, 10+2}^0)^{n-1} (\Delta_{10})^3 \oplus \cdots \oplus \mathbb{C} \Delta_{8, 4, 10+2} (\Delta_{10})^{3(n-1)} \oplus \mathbb{C} (\Delta_{10})^{3n})\\
&\qquad \oplus {E_{2/3, 10}}' (\mathbb{C} (\Delta_{10+2}^{\infty})^4 \oplus \mathbb{C} (\Delta_{10+2}^0)^4 \oplus \mathbb{C} (\Delta_{10+2}^{\infty})^3 \Delta_{10+2}^0\\
&\qquad \qquad \qquad \qquad \qquad \oplus \mathbb{C} \Delta_{10+2}^{\infty} (\Delta_{10+2}^0)^3 \oplus \mathbb{C} (\Delta_{10})^2) (\Delta_{10})^{3n}.
\end{align*}

Here, we define ${E_{2/3, 10+2}}' := \sqrt[3]{{E_{2, 10+2}}'}$ where $v_{\rho_{10, 1}}({E_{2/3, 10+2}}') = v_{\rho_{10, 5}}({E_{2/3, 10+2}}') = 1$, and we can write
\begin{equation*}
M_k(\Gamma_0(10)+2) = {E_{\overline{k}, 10+2}}' (\mathbb{C} (\Delta_{10+2}^{\infty})^{n} \oplus \mathbb{C} (\Delta_{10+2}^{\infty})^{n-1} \Delta_{10+2}^0 \oplus \cdots \oplus \mathbb{C} (\Delta_{10+2}^0)^{n}),
\end{equation*}
where $n = \dim(M_k(\Gamma_0(10)+2)) - 1 = \lfloor 3k/4 - 2 (3k/8 -\lfloor 3k/8 \rfloor) \rfloor$, and where ${E_{\overline{k}, 10+2}}' := 1$, $({E_{2/3, 10+2}}')^3$, $({E_{2/3, 10+2}}')^2$, and ${E_{2/3, 10+2}}'$, when $k \equiv 0$, $2$, $4$, and $6 \pmod{8}$, respectively.\\

\paragraph{\bf Hauptmodul}
We define the {\it hauptmodul} of $\Gamma_0(10)+2$:
\begin{equation}
J_{10+2} := \Delta_{10+2}^0 / \Delta_{10+2}^{\infty} \: (= \eta^2(z) \eta^2(2 z) \eta^{-2}(5 z) \eta^{-2}(10 z)) = \frac{1}{q} - 2 - 3 q + 6 q^2 + 2 q^3 - \cdots,
\end{equation}
where $v_{\infty}(J_{10+2}) = -1$ and $v_0(J_{10+2}) = 1$. Then, we have
{\small \begin{equation}
\begin{split}
J_{10+2} : \qquad \qquad \left\{|z + 1/10| = 1/10, \: - 1/10 \leqslant Re(z) \leqslant 0\right\} &\to [-1, 0] \subset \mathbb{R},\\
\left\{|z + 2/5| = 1 / (5 \sqrt{2}), \: -1/2 \leqslant Re(z) \leqslant -3/10\right\} &\to \{-3 \leqslant Re(z) \leqslant -1, \: 0 \leqslant Im(z) \leqslant 4\},\\
\left\{|z - 1/5| = 1 / (5 \sqrt{2}), \: 1/10 \leqslant Re(z) \leqslant 3/10 \right\} &\to \{-3 \leqslant Re(z) \leqslant -1, \: - 4 \leqslant Im(z) \leqslant 0\}.
\end{split}
\end{equation}
}Thus, $J_{10+2}$ does not take real value on some arcs of $\partial \mathbb{F}_{10+2}$.

\begin{figure}[hbtp]
\begin{center}
{{\small Lower arcs of $\partial \mathbb{F}_{10+2}$}\includegraphics[width=2.5in]{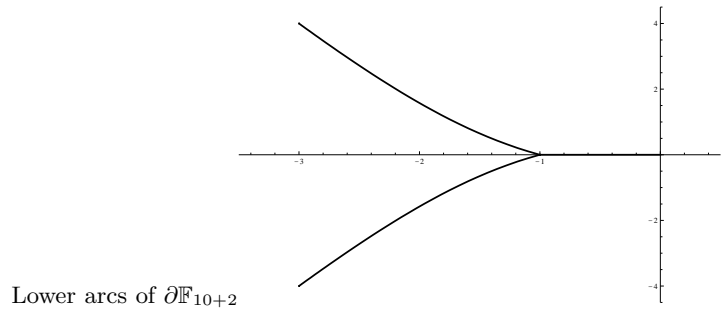}}
\end{center}
\caption{Image by $J_{10+2}$}\label{Im-J10C}
\end{figure}

\subsection{$\Gamma_0(10)$}

\paragraph{\bf Fundamental domain}
We have a fundamental domain for $\Gamma_0(10)$ as follows:
{\small \begin{equation}
\begin{split}
\mathbb{F}_{10} = &\left\{|z + 9/20| \geqslant 1/20, \: - 1/2 \leqslant Re(z) < -2/5 \right\}
 \bigcup \left\{|z + 3/10| \geqslant 1/10, \: -2/5 \leqslant Re(z) \leqslant -3/10 \right\}\\
&\bigcup \left\{|z + 3/10| > 1/10, \: -3/10 < Re(z) < -1/5 \right\}
 \bigcup \left\{|z + 1/10| \geqslant 1/10, \: - 1/5 \leqslant Re(z) \leqslant 0 \right\}\\
&\bigcup \left\{|z - 1/10| > 1/10, \: 0 < Re(z) < 1/5 \right\}
 \bigcup \left\{|z - 3/10| \leqslant 1/10, \: 1/5 \leqslant Re(z) \leqslant 3/10 \right\}\\
&\bigcup \left\{|z - 3/10| > 1/10, \: 3/10 < Re(z) \leqslant 2/5 \right\}
 \bigcup \left\{|z - 9/20| > 1/20, \: 2/5 < Re(z) < 1/2 \right\},
\end{split}
\end{equation}
}where $\left(\begin{smallmatrix} -1 & 0 \\ 10 & -1 \end{smallmatrix}\right) : (e^{i \theta} + 1) / 10 \rightarrow (e^{i (\pi - \theta)} - 1) / 10$, $\left(\begin{smallmatrix} -3 & -1 \\ 10 & 3 \end{smallmatrix}\right) : (e^{i \theta} - 3) / 10 \rightarrow (e^{i (\pi - \theta)} - 3) / 10$, $\left(\begin{smallmatrix} 3 & -1 \\ 10 & -3 \end{smallmatrix}\right) : (e^{i \theta} + 3) / 10 \rightarrow (e^{i (\pi - \theta)} + 3) / 10$, and $\left(\begin{smallmatrix} -9 & 4 \\ 20 & -9 \end{smallmatrix} \right) : (e^{i \theta} + 9) / 20 \rightarrow (e^{i (\pi - \theta)} - 9) / 20$. Then, we have
\begin{equation}
\Gamma_0(10) = \langle \left( \begin{smallmatrix} 1 & 1 \\ 0 & 1 \end{smallmatrix} \right), \: \left( \begin{smallmatrix} 1 & 0 \\ 10 & 1 \end{smallmatrix} \right), \: \left( \begin{smallmatrix} -3 & -1 \\ 10 & 3 \end{smallmatrix} \right), \: \left( \begin{smallmatrix} 3 & -1 \\ 10 & -3 \end{smallmatrix} \right), \: \left( \begin{smallmatrix} 9 & 4 \\ 20 & 9 \end{smallmatrix} \right) \rangle.
\end{equation}
\begin{figure}[hbtp]
\begin{center}
\includegraphics[width=1.5in]{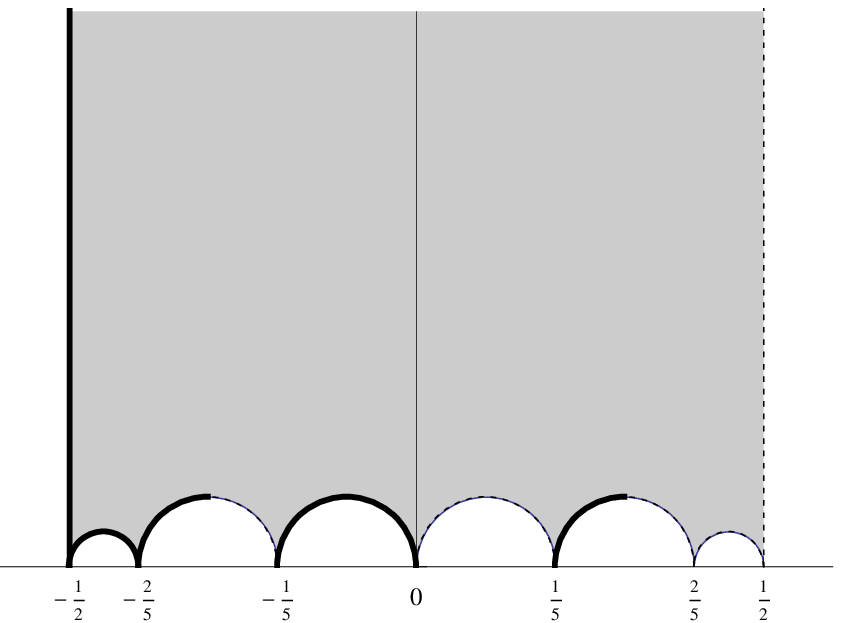}
\end{center}
\caption{$\Gamma_0(10)$}
\end{figure}

\paragraph{\bf Valence formula}
The cusps of $\Gamma_0(10)$ are $\infty$, $0$, $-1/2$, and $-1/5$. Let $f$ be a modular function of weight $k$ for $\Gamma_0(10)$, which is not identically zero. We have
\begin{equation}
v_{\infty}(f) + v_0(f) + v_{-1/2} (f) + v_{-1/5} (f)  + \frac{1}{2} v_{\rho_{10, 1}} (f) + \frac{1}{2} v_{\rho_{10, 5}} (f) + \sum_{\begin{subarray}{c} p \in \Gamma_0(10) \setminus \mathbb{H} \\ p \ne \rho_{10, 1}, \rho_{10, 5}\end{subarray}} v_p(f) = \frac{3k}{2}.
\end{equation}

Furthermore, the stabilizer of the elliptic point $\rho_{10, 1}$ (resp. $\rho_{10, 5}$) is $\left\{ \pm I, \pm \left(\begin{smallmatrix} -3 & -1 \\ 10 & 3 \end{smallmatrix}\right) \right\}$ (resp. $\left\{ \pm I, \pm \left(\begin{smallmatrix} 3 & -1 \\ 10 & -3 \end{smallmatrix}\right) \right\}$).\\

\paragraph{\bf For the cusp $\infty$}
We have $\Gamma_{\infty} = \left\{ \pm \left( \begin{smallmatrix}1 & n \\ 0 & 1\end{smallmatrix} \right) \: ; \: n \in \mathbb{Z} \right\}$, and we have the Eisenstein series for the cusp $\infty$ associated with $\Gamma_0(10)$:
\begin{equation}
E_{k, 10}^{\infty}(z) := \frac{10^k E_k(10 z) - 5^k E_k(5 z) - 2^k E_k(2 z) + E_k(z)}{(5^k - 1)(2^k - 1)} \quad \text{for} \; k \geqslant 4.
\end{equation}\quad

\paragraph{\bf For the cusp $0$}
We have $\Gamma_0 = \left\{ \pm \left( \begin{smallmatrix}1 & 0 \\ 10 n & 1\end{smallmatrix} \right) \: ; \: n \in \mathbb{Z} \right\}$ and $\gamma_0 = W_{10}$, and we have the Eisenstein series for the cusp $0$ associated with $\Gamma_0(10)$:
\begin{equation}
E_{k, 10}^0(z) := \frac{10^{k/2} (E_k(10 z) - E_k(5 z) - E_k(2 z) + E_k(z))}{(5^k - 1)(2^k - 1)} \quad \text{for} \; k \geqslant 4.
\end{equation}
We also have $\gamma_0^{-1} \ \Gamma_0(10) \ \gamma_0 = \Gamma_0(10)$.\\

\paragraph{\bf For the cusp $-1/2$}
We have $\Gamma_{-1/2} = \left\{ \pm \left( \begin{smallmatrix}10 n + 1 & 5 n \\ - 20 n & - 10 n + 1\end{smallmatrix} \right) \: ; \: n \in \mathbb{Z} \right\}$ and $\gamma_{-1/2} = W_{10, 5}$, and we have the Eisenstein series for the cusp $-1/2$ associated with $\Gamma_0(10)$:
\begin{equation}
E_{k, 10}^{-1/2}(z) := \frac{- 5^{k/2} (2^k E_k(10 z) - E_k(5 z) - 2^k E_k(2 z) + E_k(z))}{(5^k - 1)(2^k - 1)} \quad \text{for} \; k \geqslant 4.
\end{equation}
We also have $\gamma_{-1/2}^{-1} \ \Gamma_0(10) \ \gamma_{-1/2} = \Gamma_0(10)$.\\

\paragraph{\bf For the cusp $-1/5$}
We have $\Gamma_{-1/5} = \left\{ \pm \left( \begin{smallmatrix}10 n + 1 & 2 n \\ - 50 n & - 10 n + 1\end{smallmatrix} \right) \: ; \: n \in \mathbb{Z} \right\}$ and $\gamma_{-1/5} = W_{10, 2}$, and we have the Eisenstein series for the cusp $-1/5$ associated with $\Gamma_0(10)$:
\begin{equation}
E_{k, 10}^{-1/5}(z) := \frac{- 2^{k/2} (5^k E_k(10 z) - 5^k E_k(5 z) - E_k(2 z) + E_k(z))}{(5^k - 1)(2^k - 1)} \quad \text{for} \; k \geqslant 4.
\end{equation}
We also have $\gamma_{-1/5}^{-1} \ \Gamma_0(10) \ \gamma_{-1/5} = \Gamma_0(10)$.\\

\paragraph{\bf The space of modular forms}
We define
\begin{align*}
\Delta_{8, 1, 10} &:= ({E_{2/3, 10}}')^2 (\Delta_{10}^0)^2 (\Delta_{10}^{-1/2})^2 (\Delta_{10}^{-1/5})^2 \Delta_{10},
&\Delta_{8, 2, 10} &:= ({E_{2/3, 10}}')^2 (\Delta_{10}^{\infty})^2 (\Delta_{10}^{-1/2})^2 (\Delta_{10}^{-1/5})^2 \Delta_{10},\\
\Delta_{8, 3, 10} &:= ({E_{2/3, 10}}')^2 (\Delta_{10}^{\infty})^2 (\Delta_{10}^0)^2 (\Delta_{10}^{-1/5})^2 \Delta_{10},
&\Delta_{8, 4, 10} &:= ({E_{2/3, 10}}')^2 (\Delta_{10}^{\infty})^2 (\Delta_{10}^0)^2 (\Delta_{10}^{-1/2})^2 \Delta_{10},\\
\Delta_{8, 5, 10} &:= (\Delta_{10}^0)^2 \Delta_{10}^{-1/2} \Delta_{10}^{-1/5} (\Delta_{10})^2,
&\Delta_{8, 6, 10} &:= (\Delta_{10}^{\infty})^2 \Delta_{10}^{-1/2} \Delta_{10}^{-1/5} (\Delta_{10})^2,\\
\Delta_{8, 7, 10} &:= \Delta_{10}^{\infty} \Delta_{10}^0 (\Delta_{10}^{-1/5})^2 (\Delta_{10})^2,
&\Delta_{8, 8, 10} &:= \Delta_{10}^{\infty} \Delta_{10}^0 (\Delta_{10}^{-1/2})^2 (\Delta_{10})^2.
\end{align*}

Now, we have
\begin{align*}
M_k(\Gamma_0(10)) &= \mathbb{C} E_{k, 10}^{\infty} \oplus \mathbb{C} E_{k, 10}^0 \oplus \mathbb{C} E_{k, 10}^{-1/2} \oplus \mathbb{C} E_{k, 10}^{-1/5} \oplus S_k(\Gamma_0(10)),\\
S_k(\Gamma_0(10)) &= (\mathbb{C} \Delta_{8, 1, 10} \oplus \mathbb{C} \Delta_{8, 2, 10} \oplus \mathbb{C} \Delta_{8, 3, 10} \oplus \mathbb{C} \Delta_{8, 4, 10} \oplus \mathbb{C} \Delta_{8, 5, 10}\\
&\qquad \oplus \mathbb{C} \Delta_{8, 6, 10} \oplus \mathbb{C} \Delta_{8, 7, 10} \oplus \mathbb{C} \Delta_{8, 8, 10} \oplus \mathbb{C} (\Delta_{10})^3) M_{k - 8}(\Gamma_0(10))
\end{align*}
for every even integer $k \geqslant 4$. Then, we have $M_{8 n + 2}(\Gamma_0(10)) = {E_{2, 10+2}}' M_{8 n}(\Gamma_0(10)) \oplus \mathbb{C} {E_{2, 10+10}}' (\Delta_{10})^{3n} \oplus \mathbb{C} {E_{2, 10+5}}' (\Delta_{10})^{3n}$, $M_{8 n + 6}(\Gamma_0(10)) = {E_{2, 10+2}}' M_{8 n + 4}(\Gamma_0(10)) \oplus \mathbb{C} ({E_{2, 10+10}}')^3 (\Delta_{10})^{3n} \oplus \mathbb{C} ({E_{2, 10+5}}')^3 (\Delta_{10})^{3n}$, and
\begin{allowdisplaybreaks}
\begin{align*}
M_{8 n}(\Gamma_0(10)) &= \mathbb{C} (E_{8, 10}^{\infty})^n \oplus \mathbb{C} (E_{8, 10}^{\infty})^{n-1} \Delta_{8, 1, 10} \oplus \mathbb{C} (E_{8, 10}^{\infty})^{n-1} \Delta_{8, 5, 10}\\
&\qquad \qquad \qquad \qquad \oplus \mathbb{C} (E_{8, 10}^{\infty})^{n-1} (\Delta_{10})^3 \oplus \cdots \oplus \mathbb{C} \Delta_{8, 5, 10} (\Delta_{10})^{3(n-1)}\\
&\oplus \mathbb{C} (E_{8, 10}^0)^n \oplus \mathbb{C} (E_{8, 10}^0)^{n-1} \Delta_{8, 2, 10} \oplus \mathbb{C} (E_{8, 10}^0)^{n-1} \Delta_{8, 6, 10}\\
&\qquad \qquad \qquad \qquad \oplus \mathbb{C} (E_{8, 10}^0)^{n-1} (\Delta_{10})^3 \oplus \cdots \oplus \mathbb{C} \Delta_{8, 6, 10} (\Delta_{10})^{3(n-1)}\\
&\oplus \mathbb{C} (E_{8, 10}^{-1/2})^n \oplus \mathbb{C} (E_{8, 10}^{-1/2})^{n-1} \Delta_{8, 3, 10} \oplus \mathbb{C} (E_{8, 10}^{-1/2})^{n-1} \Delta_{8, 7, 10}\\
&\qquad \qquad \qquad \qquad \oplus \mathbb{C} (E_{8, 10}^{-1/2})^{n-1} (\Delta_{10})^3 \oplus \cdots \oplus \mathbb{C} \Delta_{8, 7, 10} (\Delta_{10})^{3(n-1)}\\
&\oplus \mathbb{C} (E_{8, 10}^{-1/5})^n \oplus \mathbb{C} (E_{8, 10}^{-1/5})^{n-1} \Delta_{8, 4, 10} \oplus \mathbb{C} (E_{8, 10}^{-1/5})^{n-1} \Delta_{8, 8, 10}\\
&\qquad \qquad \qquad \qquad \oplus \mathbb{C} (E_{8, 10}^{-1/5})^{n-1} (\Delta_{10})^3 \oplus \cdots \oplus \mathbb{C} \Delta_{8, 8, 10} (\Delta_{10})^{3(n-1)} \oplus \mathbb{C} (\Delta_{10})^{3n},\\
M_{8 n + 4}(\Gamma_0(10)) &= E_{4, 10}^{\infty} (\mathbb{C} (E_{8, 10}^{\infty})^n \oplus \mathbb{C} (E_{8, 10}^{\infty})^{n-1} \Delta_{8, 1, 10} \oplus \mathbb{C} (E_{8, 10}^{\infty})^{n-1} \Delta_{8, 5, 10}\\
&\qquad \qquad \qquad \qquad \oplus \mathbb{C} (E_{8, 10}^{\infty})^{n-1} (\Delta_{10})^3 \oplus \cdots \oplus \mathbb{C} \Delta_{8, 5, 10} (\Delta_{10})^{3(n-1)} \oplus \mathbb{C} (\Delta_{10})^{3n})\\
&\oplus E_{4, 10}^0 (\mathbb{C} (E_{8, 10}^0)^n \oplus \mathbb{C} (E_{8, 10}^0)^{n-1} \Delta_{8, 2, 10} \oplus \mathbb{C} (E_{8, 10}^0)^{n-1} \Delta_{8, 6, 10}\\
&\qquad \qquad \qquad \qquad \oplus \mathbb{C} (E_{8, 10}^0)^{n-1} (\Delta_{10})^3 \oplus \cdots \oplus \mathbb{C} \Delta_{8, 6, 10} (\Delta_{10})^{3(n-1)} \oplus \mathbb{C} (\Delta_{10})^{3n})\\
&\oplus E_{4, 10}^{-1/2} (\mathbb{C} (E_{8, 10}^{-1/2})^n \oplus \mathbb{C} (E_{8, 10}^{-1/2})^{n-1} \Delta_{8, 3, 10} \oplus \mathbb{C} (E_{8, 10}^{-1/2})^{n-1} \Delta_{8, 7, 10}\\
&\qquad \qquad \qquad \qquad \oplus \mathbb{C} (E_{8, 10}^{-1/2})^{n-1} (\Delta_{10})^3 \oplus \cdots \oplus \mathbb{C} \Delta_{8, 7, 10} (\Delta_{10})^{3(n-1)} \oplus \mathbb{C} (\Delta_{10})^{3n})\\
&\oplus E_{4, 10}^{-1/5} (\mathbb{C} (E_{8, 10}^{-1/5})^n \oplus \mathbb{C} (E_{8, 10}^{-1/5})^{n-1} \Delta_{8, 4, 10} \oplus \mathbb{C} (E_{8, 10}^{-1/5})^{n-1} \Delta_{8, 8, 10}\\
&\qquad \qquad \qquad \qquad \oplus \mathbb{C} (E_{8, 10}^{-1/5})^{n-1} (\Delta_{10})^3 \oplus \cdots \oplus \mathbb{C} \Delta_{8, 8, 10} (\Delta_{10})^{3(n-1)} \oplus \mathbb{C} (\Delta_{10})^{3n})\\
&\qquad \qquad \oplus (\mathbb{C} \Delta_{10}^{\infty} \Delta_{10}^0 \oplus \mathbb{C} \Delta_{10}^0 \Delta_{10}^{-1/2} \oplus \mathbb{C} \Delta_{10}^{-1/2} \Delta_{10}^{-1/5}) (\Delta_{10})^{3n+1}.
\end{align*}
\end{allowdisplaybreaks}

Furthermore, we can write
\begin{equation*}
M_k(\Gamma_0(10)) = {E_{\overline{k}, 10}}' (\mathbb{C} (\Delta_{10}^{\infty})^{n} \oplus \mathbb{C} (\Delta_{10}^{\infty})^{n-1} \Delta_{10}^0 \oplus \cdots \oplus \mathbb{C} (\Delta_{10}^0)^{n}),
\end{equation*}
where $n = \dim(M_k(\Gamma_0(10))) - 1 = \lfloor 3k/2 - 2 (k/4 -\lfloor k/4 \rfloor) \rfloor$, and where ${E_{\overline{k}, 10}}' := 1$ and ${E_{2/3, 10+2}}'$, when $k \equiv 0$ and $2 \pmod{4}$, respectively.\\

\paragraph{\bf Hauptmodul}
We define the {\it hauptmodul} of $\Gamma_0(10)$:
\begin{equation}
J_{10} := \Delta_{10}^0 / \Delta_{10}^{\infty} \: (= \eta^3(z) \eta^{-1}(2 z) \eta(5 z) \eta^{-3}(10 z)) = \frac{1}{q} - 3 + q + 2 q^2 + 2 q^3 - \cdots,
\end{equation}
where $v_{\infty}(J_{10}) = -1$ and $v_0(J_{10}) = 1$. Then, we have
{\small \begin{equation}
\begin{split}
J_{10} : \qquad  \left\{|z + 1/10| = 1/10, \: -1/10 \leqslant Re(z) \leqslant 0 \right\} &\to [-4, 0] \subset \mathbb{R},\\
\left\{|z + 3/10| = 1/10, \: -2/5 \leqslant Re(z) \leqslant -3/10 \right\} &\to \{-4 \leqslant Re(z) < -3.8 , \: 0 \leqslant Im(z) \leqslant 2\},\\
\left\{|z - 3/10| = 1/10, \: 1/5 \leqslant Re(z) \leqslant 3/10 \right\} &\to \{-4 \leqslant Re(z) < -3.8 , \: -2  \leqslant Im(z) \leqslant 0\},\\
\left\{|z + 9/20| = 1/20, \: -1/2 \leqslant Re(z) \leqslant -2/5 \right\} &\to [-5, -4] \subset \mathbb{R}.
\end{split}
\end{equation}
}Thus, $J_{10}$ does not take real value on some arcs of $\partial \mathbb{F}_{10}$.

\begin{figure}[hbtp]
\begin{center}
{{\small Lower arcs of $\partial \mathbb{F}_{10}$}\includegraphics[width=2.5in]{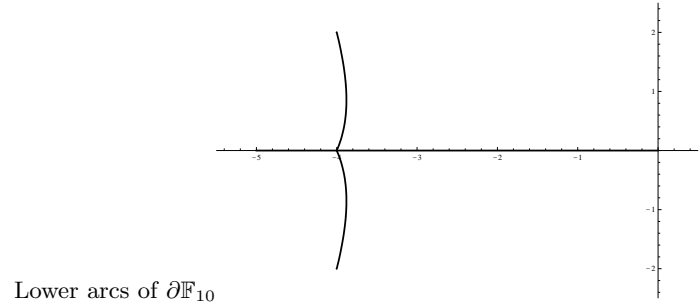}}
\end{center}
\caption{Image by $J_{10}$}\label{Im-J10E}
\end{figure}

\clearpage

\section{Level $11$}

We have $\Gamma_0(11)+=\Gamma_0^{*}(11)$ and $\Gamma_0(11)-=\Gamma_0(11)$, but $\Gamma_0(11)$ is of genus $1$.

We have $W_{11} = \left(\begin{smallmatrix}0&-1 / \sqrt{11}\\ \sqrt{11}&0\end{smallmatrix}\right)$, and denote $\rho_{11, 1} := - 1/2 + i / (2 \sqrt{11})$, $\rho_{11, 2} := - 1/3 + i / (3 \sqrt{11})$, and $\rho_{11, 3} := 1/3 + i / (3 \sqrt{11})$. We define
\begin{equation}
\begin{split}
&\Delta_{11}^{\infty}(z) := \sqrt[5]{\eta^{11}(11 z) / \eta(z)}, \quad \Delta_{11}^0(z) := \sqrt[5]{\eta^{11}(z) / \eta(11 z)},\\
&\Delta_{11}(z) := \Delta_{11}^{\infty}(z) \Delta_{11}^0(z) = \eta^2(z) \eta^2(11 z),
\end{split}
\end{equation}
where $\Delta_{11}^{\infty}$ and $\Delta_{11}^0$ are $2$nd semimodular forms for $\Gamma_0(11)$ of weight $1$ such that $v_{\infty}(\Delta_{11}^{\infty}) = v_0(\Delta_{11}^0) = 1$, and $\Delta_{11}$ is a cusp form for $\Gamma_0(11)$ and $2$nd semimodular form for $\Gamma_0^{*}(11)$ of weight $2$. Furthermore, we define
\begin{equation}
{E_{2, 11}}'(z) := (11 E_2(11 z) - E_2(z)) / 10,
\end{equation}
which is a modular form for $\Gamma_0(11)$ and $2$nd semimodular form for $\Gamma_0^{*}(11)$ of weight $2$.\\

\subsection{$\Gamma_0^{*}(11)$}

\paragraph{\bf Fundamental domain}
We have a fundamental domain for $\Gamma_0^{*}(11)$ as follows:
{\small \begin{equation}
\begin{split}
\mathbb{F}_{11+} = &\left\{|z + 1/2| \geqslant 1 / (2 \sqrt{11}), \: - 1/2 \leqslant Re(z) < - 25/66 \right\}
 \bigcup \left\{|z + 1/3| \geqslant 1 / (3 \sqrt{11}), \: - 25/66 \leqslant Re(z) \leqslant - 1/3 \right\}\\
 &\bigcup \left\{|z + 1/3| > 1 / (3 \sqrt{11}), \: - 1/3 < Re(z) < - 19/66 \right\}
 \bigcup \left\{|z| \geqslant 1 / \sqrt{11}, \: - 19/66 \leqslant Re(z) \leqslant 0 \right\}\\
 &\bigcup \left\{|z| > 1 / \sqrt{11}, \: 0 < Re(z) <19/66 \right\}
 \bigcup \left\{|z - 1/3| \geqslant 1 / (3 \sqrt{11}), \: 19/66 \leqslant Re(z) \leqslant 1/3 \right\}\\
 &\bigcup \left\{|z - 1/3| > 1 / (3 \sqrt{11}), \: 1/3 < Re(z) \leqslant 25/66 \right\}
 \bigcup \left\{|z - 1/2| > 1 / (2 \sqrt{11}), \: 25/66 < Re(z) < 1/2 \right\}.
\end{split}
\end{equation}
}where $W_{11} :  e^{i \theta} / \sqrt{11} \rightarrow e^{i  (\pi - \theta)} / \sqrt{11}$, $\left(\begin{smallmatrix} -5 & -1 \\ 11 & 2 \end{smallmatrix}\right) W_{11} : e^{i \theta} / (2 \sqrt{11}) + 1/2 \rightarrow e^{i  (\pi - \theta)} / (2 \sqrt{11}) - 1/2$, $\left(\begin{smallmatrix} 4 & -1 \\ -11 & 3 \end{smallmatrix}\right) W_{11} : e^{i \theta} / (3 \sqrt{11}) - 1/3 \rightarrow e^{i  (\pi - \theta)} / (3 \sqrt{11}) - 1/3$, and $\left(\begin{smallmatrix} 4 & 1 \\ 11 & 3 \end{smallmatrix}\right) W_{11} : e^{i \theta} / (3 \sqrt{11}) + 1/3 \rightarrow e^{i  (\pi - \theta)} / (3 \sqrt{11}) + 1/3$. Then, we have
\begin{equation}
\Gamma_0^{*}(11) = \langle \left( \begin{smallmatrix} 1 & 1 \\ 0 & 1 \end{smallmatrix} \right), \: W_{11}, \: \left( \begin{smallmatrix} 4 & 1 \\ 11 & 3 \end{smallmatrix} \right), \: \left( \begin{smallmatrix} 3 & 1 \\ 11 & 4 \end{smallmatrix} \right) \rangle.
\end{equation}
\begin{figure}[hbtp]
\begin{center}
\includegraphics[width=1.5in]{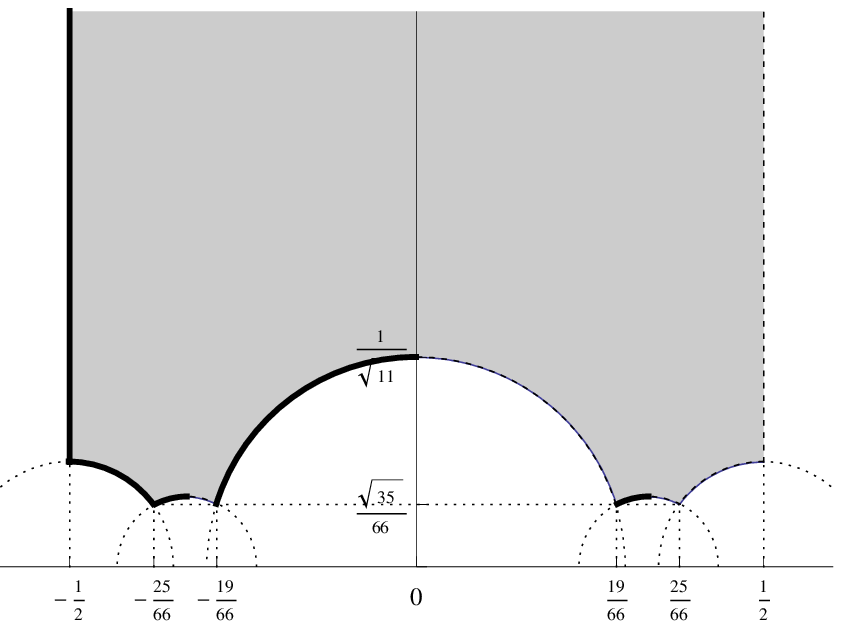}
\end{center}
\caption{$\Gamma_0^{*}(11)$}
\end{figure}

\paragraph{\bf Valence formula}
The cusp of $\Gamma_0^{*}(11)$ is $\infty$, and the elliptic points are $i / \sqrt{11}$, $\rho_{11, 1} = - 1/2 + i \sqrt{11} / 10$, $\rho_{11, 2} = - 1/3 + i / (3 \sqrt{11})$, and $\rho_{11, 3} = 1/3 + i / (3 \sqrt{11})$. Let $f$ be a modular function of weight $k$ for $\Gamma_0^{*}(11)$, which is not identically zero. We have
\begin{equation}
v_{\infty}(f) + \frac{1}{2} v_{i / \sqrt{11}}(f) + \frac{1}{2} v_{\rho_{11, 1}} (f) + \frac{1}{2} v_{\rho_{11, 2}} (f) + \frac{1}{2} v_{\rho_{11, 3}} (f) + \sum_{\begin{subarray}{c} p \in \Gamma_0^{*}(11) \setminus \mathbb{H} \\ p \ne i / \sqrt{11}, \; \rho_{11, 1}, \; \rho_{11, 2}, \; \rho_{11, 3} \end{subarray}} v_p(f) = \frac{k}{2}.
\end{equation}

Furthermore, the stabilizer of the elliptic point $i / \sqrt{11}$ (resp. $\rho_{11, 1}$, $\rho_{11, 2}$, $\rho_{11, 3}$) is $\left\{ \pm I, \pm W_{11} \right\}$\\ (resp. $\left\{ \pm I, \pm \left( \begin{smallmatrix} 6 & -1 \\ -11 & 2 \end{smallmatrix} \right) W_{11} \right\}$, $\left\{ \pm I, \pm \left( \begin{smallmatrix} 4 & -1 \\ -11 & 3 \end{smallmatrix} \right) W_{11} \right\}$, $\left\{ \pm I, \pm \left( \begin{smallmatrix} 4 & 1 \\ 11 & 3 \end{smallmatrix} \right) W_{11} \right\}$)\\

\paragraph{\bf For the cusp $\infty$}
We have $\Gamma_{\infty} = \left\{ \pm \left( \begin{smallmatrix}1 & n \\ 0 & 1\end{smallmatrix} \right) \: ; \: n \in \mathbb{Z} \right\}$, and we have the Eisenstein series associated with $\Gamma_0^{*}(11)$:
\begin{equation}
E_{k, 11+}(z) := \frac{11^{k/2} E_k(11 z) + E_k(z)}{11^{k/2} + 1} \quad \text{for} \; k \geqslant 4.
\end{equation}\quad

\paragraph{\bf The space of modular forms}
We define the following functions:
\begin{equation*}
{E_{4, 11}}' := (- 1525 ({E_{2, 11}}')^2 + 4320 {E_{2, 11}}' \Delta_{11} + 2016 (\Delta_{11})^2 + 3050 E_{4, 11}^{\infty}) / 1525,
\end{equation*}
where $E_{4, 11}^{\infty}$ is the Eisenstein series of weight $4$ for the cusp $\infty$ for $\Gamma_0(11)$. Then,  ${E_{4, 11}}'$ is a $2$nd semimodular form for $\Gamma_0^{*}(11)$ of weight $4$ such that $v_{i / \sqrt{11}}({E_{4, 11}}') = v_{\rho_{11, 1}}({E_{4, 11}}') = v_{\rho_{11, 2}}({E_{4, 11}}') = v_{\rho_{11, 3}}({E_{4, 11}}') = 1$.

Let $k$ be an even integer $k \geqslant 4$. We have $M_k(\Gamma_0^{*}(11)) = \mathbb{C} E_{k, 11+} \oplus S_k(\Gamma_0^{*}(11))$ and $S_k(\Gamma_0^{*}(11)) = (\mathbb{C} {E_{2, 11}}' \Delta_{11} \oplus \mathbb{C} (\Delta_{11})^2) M_{k - 4}(\Gamma_0^{*}(11))$.
Then, we have
\begin{align*}
M_{4 n}(\Gamma_0^{*}(11)) &= \mathbb{C} ({E_{2, 11}}')^{2 n} \oplus \mathbb{C} ({E_{2, 11}}')^{2 n - 1} \Delta_{11} \oplus \cdots \oplus \mathbb{C} (\Delta_{11})^{2 n},\\
M_{4 n + 6}(\Gamma_0^{*}(11)) &= {E_{4, 11}}' (\mathbb{C} ({E_{2, 11}}')^{2 n + 1} \oplus \mathbb{C} ({E_{2, 11}}')^{2 n} \Delta_{11} \oplus \cdots \oplus \mathbb{C} (\Delta_{11})^{2 n + 1}).
\end{align*}\quad

\paragraph{\bf Hauptmodul}
We define the {\it hauptmodul} of $\Gamma_0^{*}(11)$:
\begin{equation}
J_{11+} := {E_{2, 11}}' / \Delta_{11} = \frac{1}{q} + \frac{22}{5} + 17 q + 46 q^2 + 116 q^3 + \cdots,
\end{equation}
where $v_{\infty}(J_{11+}) = -1$. Then, we have
{\small \begin{equation}
\begin{split}
J_{11+} : \qquad \qquad \left\{|z| = 1 / \sqrt{11}, \: -19/66 \leqslant Re(z) \leqslant 0\right\} &\to [22/5 - 2 \sqrt{5}, 15.22750...] \subset \mathbb{R},\\
\left\{|z + 1/3| = 1 / (3 \sqrt{11}), \: -25/66 \leqslant Re(z) \leqslant -1/3\right\} &\to \{ 22/5 - 2 \sqrt{5} \leqslant Re(z) \leqslant -0.013750..., \: 0 \leqslant Im(z) \leqslant 0.31397... \},\\
\left\{|z - 1/3| = 1 / (3 \sqrt{11}), \: 19/66 \leqslant Re(z) \leqslant 1/3 \right\} &\to \{ 22/5 - 2 \sqrt{5} \leqslant Re(z) \leqslant -0.013750..., \: - 0.31397... \leqslant Im(z) \leqslant 0\}\\
\left\{|z + 1/2| = 1 / (2 \sqrt{11}), \: -1/2 \leqslant Re(z) \leqslant -25/66 \right\} &\to [-8/5, 22/5 - 2 \sqrt{5}] \subset \mathbb{R}.
\end{split}
\end{equation}
}Thus, $J_{11+}$ does not take real value on some arcs of $\partial \mathbb{F}_{11+}$.

\begin{figure}[hbtp]
\begin{center}
{{\small Lower arcs of $\partial \mathbb{F}_{11+}$}\includegraphics[width=2.5in]{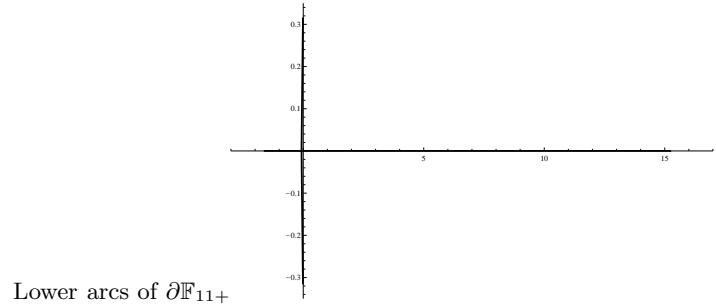}}
\end{center}
\caption{Image by $J_{11+}$}\label{Im-J11A}
\end{figure}

\subsection{$\Gamma_0(11)$}

\paragraph{\bf Fundamental domain}
We have a fundamental domain for $\Gamma_0(11)$ as follows:
{\small \begin{equation}
\begin{split}
\mathbb{F}_{11} = &\left\{|z + 5/11| \geqslant 1/11, \: - 1/2 \leqslant Re(z) < - 9/22\right\}
 \bigcup \left\{|z + 4/11| \geqslant 1/11, \: - 9/22 \leqslant Re(z) < - 7/22 \right\}\\
&\left\{|z + 3/11| \geqslant 1/11, \: - 7/22 \leqslant Re(z) \leqslant - 5/22\right\}
 \bigcup \left\{|z + 2/11| > 1/11, \: - 5/22 < Re(z) < - 3/22 \right\}\\
&\bigcup \left\{|z + 1/11| \geqslant 1/11, \: - 3/22 \leqslant Re(z) \leqslant 0\right\}
 \bigcup \left\{|z - 1/11| > 1/11, \: 0 < Re(z) < 3/22 \right\}\\
&\bigcup \left\{|z - 2/11| \geqslant 1/11, \: 3/22 \leqslant Re(z) \leqslant 5/22 \right\}
 \bigcup \left\{|z - 3/11| > 1/11, \: 5/22 < Re(z) < 7/22 \right\}\\
&\bigcup \left\{|z - 4/11| > 1/11, \: 7/22 \leqslant Re(z) < 9/22 \right\}
 \bigcup \left\{|z - 5/11| > 1/11, \: 9/22 \leqslant Re(z) < 1/2 \right\}.
\end{split}
\end{equation}
}where $\left( \begin{smallmatrix} -1 & 0 \\ 11 & -1 \end{smallmatrix} \right) : (e^{i \theta} + 1) / 11 \rightarrow (e^{i (\pi - \theta)} - 1) / 11$, $\left( \begin{smallmatrix} -5 & -1 \\ 11 & 2 \end{smallmatrix} \right) : (e^{i \theta} - 2) / 11 \rightarrow (e^{i (\pi - \theta)} - 5) / 11$, $\left( \begin{smallmatrix} 2 & -1 \\ 11 & -5 \end{smallmatrix} \right) : (e^{i \theta} + 5) / 11 \rightarrow (e^{i (\pi - \theta)} + 2) / 11$, $\left( \begin{smallmatrix} 3 & -1 \\ -11 & 4 \end{smallmatrix} \right) : (e^{i \theta} + 4) / 11 \rightarrow (e^{i (\pi - \theta)} - 3) / 11$, and $\left( \begin{smallmatrix} 4 & -1 \\ -11 & 3 \end{smallmatrix} \right) : (e^{i \theta} + 3) / 11 \rightarrow (e^{i (\pi - \theta)} - 4) / 11$. Then, we have
\begin{equation}
\Gamma_0(11) = \langle - I, \: \left( \begin{smallmatrix} 1 & 1 \\ 0 & 1 \end{smallmatrix} \right), \: \left( \begin{smallmatrix} 4 & 1 \\ 11 & 3 \end{smallmatrix} \right), \: \left( \begin{smallmatrix} 3 & 1 \\ 11 & 4 \end{smallmatrix} \right) \rangle.
\end{equation}
\begin{figure}[hbtp]
\begin{center}
\includegraphics[width=1.5in]{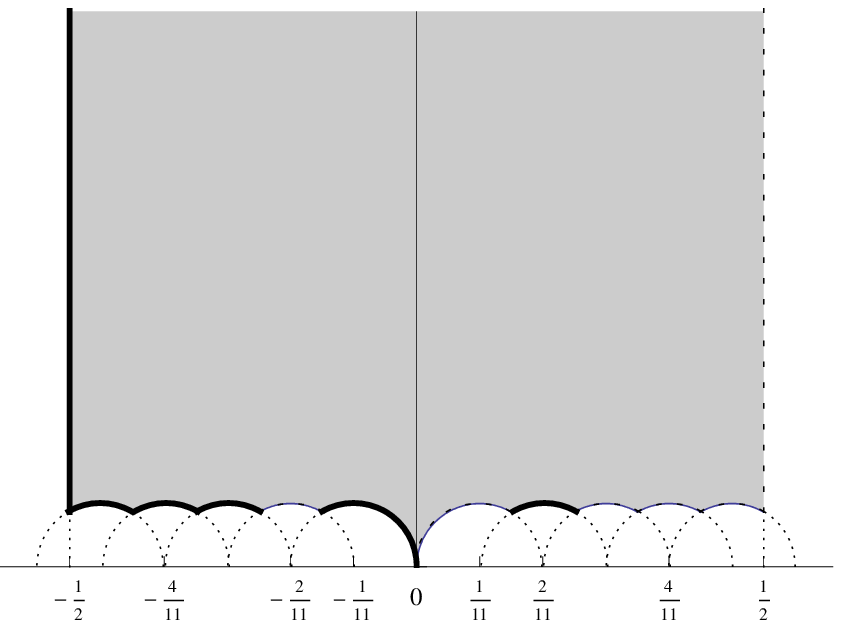}
\end{center}
\caption{$\Gamma_0(11)$}
\end{figure}

\paragraph{\bf Valence formula}
The cusps of $\Gamma_0(11)$ are $\infty$ and $0$. Let $f$ be a modular function of weight $k$ for $\Gamma_0(11)$, which is not identically zero. We have
\begin{equation}
v_{\infty}(f) + v_0(f) + \sum_{p \in \Gamma_0(11) \setminus \mathbb{H}} v_p(f) = k.
\end{equation}\quad

\paragraph{\bf For the cusp $\infty$}
We have $\Gamma_{\infty} = \left\{ \pm \left( \begin{smallmatrix}1 & n \\ 0 & 1\end{smallmatrix} \right) \: ; \: n \in \mathbb{Z} \right\}$, and we have the Eisenstein series for the cusp $\infty$ associated with $\Gamma_0(11)$:
\begin{equation}
E_{k, 11}^{\infty}(z) := \frac{11^k E_k(11 z) - E_k(z)}{11^k - 1} \quad \text{for} \; k \geqslant 4.
\end{equation}\quad

\paragraph{\bf For the cusp $0$}
We have $\Gamma_0 = \left\{ \pm \left( \begin{smallmatrix}1 & 0 \\ 11 n & 1\end{smallmatrix} \right) \: ; \: n \in \mathbb{Z} \right\}$ and $\gamma_0 = W_{11}$, and we have the Eisenstein series for the cusp $0$ associated with $\Gamma_0(11)$:
\begin{equation}
E_{k, 11}^0(z) := \frac{- 11^{k/2} (E_k(11 z) - E_k(z))}{11^k - 1} \quad \text{for} \; k \geqslant 4.
\end{equation}
We also have $\gamma_0^{-1} \ \Gamma_0(11) \ \gamma_0 = \Gamma_0(11)$.\\

%% file: report12-13.tex
\section{Level $12$}

We have $\Gamma_0(12)+$, $\Gamma_0(12)+12=\Gamma_0^{*}(12)$, $\Gamma_0(12)+4$, $\Gamma_0(12)+3$, and $\Gamma_0(12)-=\Gamma_0(12)$.

We have $W_{12} = \left(\begin{smallmatrix}0 & - 1 / (2 \sqrt{3})\\ 2 \sqrt{3} & 0\end{smallmatrix}\right)$, $W_{12, 3} := \left(\begin{smallmatrix} - \sqrt{3} & -1 / \sqrt{3} \\ 4 \sqrt{3} & \sqrt{3} \end{smallmatrix}\right)$, $W_{12, 4} := \left(\begin{smallmatrix} -2 & 1/2 \\ 6 & -2 \end{smallmatrix}\right)$, $W_{12-, 2} := \left(\begin{smallmatrix} -1 & 0 \\ 6 & -1 \end{smallmatrix}\right)$, $W_{12+, 2} := \left(\begin{smallmatrix} -1 / \sqrt{2} & -1 / (2 \sqrt{2}) \\ 3 \sqrt{2} & 1 / \sqrt{2} \end{smallmatrix}\right)$, $W_{12-, 6} := \left(\begin{smallmatrix} - \sqrt{3} & -2 / \sqrt{3} \\ 2 \sqrt{3} & \sqrt{3} \end{smallmatrix}\right)$, and $W_{12+, 6} := \left(\begin{smallmatrix} - \sqrt{6} / 2 & -5 / (2 \sqrt{6}) \\ \sqrt{6} & \sqrt{6} / 2 \end{smallmatrix}\right)$, and we denote $\rho_{12, 1} := -1/4 + i / (4 \sqrt{3})$, $\rho_{12, 2} := -2/7 + i / (14 \sqrt{3})$, and $\rho_{12, 3} := 1/4 + i / (4 \sqrt{3})$. We define
\begin{equation}
\begin{aligned}
\Delta_{12}^{\infty}(z) &:= \sqrt{\eta(2 z) \eta^{-2}(4 z) \eta^{-3}(6 z) \eta^6(12 z)},\\
\Delta_{12}^{-1/3}(z) &:= \sqrt{\eta^{-2}(z) \eta(2 z) \eta^6(3 z) \eta^{-3}(6 z)},\\
\Delta_{12}^{-1/2}(z) &:= \sqrt{\eta^{-6}(z) \eta^{15}(2 z) \eta^2(3 z) \eta^{-6} (4 z)\eta^{-5}(6 z) \eta^2(12 z)},\\
\Delta_{12}^{-1/6}(z) &:= \sqrt{\eta^2(z) \eta^{-5}(2 z) \eta^{-6}(3 z) \eta^2 (4 z)\eta^{15}(6 z) \eta^{-2}(12 z)},\\
\Delta_{12} &:= \Delta_{12}^{\infty} \Delta_{12}^0 \Delta_{12}^{-1/3} \Delta_{12}^{-1/4} \Delta_{12}^{-1/2} \Delta_{12}^{-1/6},
\end{aligned}
\hspace{-0.5in}
\begin{aligned}
\Delta_{12}^0(z) &:= \sqrt{\eta^6(z) \eta^{-3}(2 z) \eta^{-2}(3 z) \eta(6 z)}, \\
\Delta_{12}^{-1/4}(z) &:= \sqrt{\eta^{-3}(2 z) \eta^6(4 z) \eta(6 z) \eta^{-2}(12 z)},\\
& \quad\\ & \quad\\
\Delta_{12+} := &\Delta_{12}^{\infty} \Delta_{12}^0 \Delta_{12}^{-1/3} \Delta_{12}^{-1/4} (\Delta_{12}^{-1/2})^2 (\Delta_{12}^{-1/6})^2,
\end{aligned}
\end{equation}
where $\Delta_{12}^{\infty}$, $\Delta_{12}^0$, $\Delta_{12}^{-1/3}$, $\Delta_{12}^{-1/4}$, $\Delta_{12}^{-1/2}$, and $\Delta_{12}^{-1/6}$ are $4$th semimodular forms for $\Gamma_0(12)$ of weight $1$ such that $v_{\infty}(\Delta_{12}^{\infty}) = v_0(\Delta_{12}^0) = v_{-1/3}(\Delta_{12}^{-1/3}) = v_{-1/4}(\Delta_{12}^{-1/4}) = v_{-1/2}(\Delta_{12}^{-1/2}) = v_{-1/6}(\Delta_{12}^{-1/6}) = 1$. Furthermore, we define
\begin{equation}
\begin{split}
{E_{2, 12+}}'(z) &:= (12 E_2(12 z) - 12 E_2(6 z) + 4 E_2(4 z) + 3 E_2(3 z) - 4 E_2(2 z) + E_2(z)) / 4,\\
{E_{1, 12+12}}'(z) &:= \sqrt{- (12 E_2(12 z) - 18 E_2(6 z) - 4 E_2(4 z) + 3 E_2(3 z) + 6 E_2(2 z) - E_2(z)) / 2},\\
{E_{1, 12+3}}'(z) &:= \sqrt{(3 E_2(6 z) - E_2(z)) / 2},
\end{split}
\end{equation}
where ${E_{2, 12+}}'$ is a modular form of weight $2$ and ${E_{1, 12+12}}'$, ${E_{1, 12+3}}'$ are $2$nd semimodular forms of weight $1$ for $\Gamma_0(12)$, and we have $v_{i / (2 \sqrt{3})}({E_{2, 12+}}') = v_{\rho_{12, 1}}({E_{2, 12+}}') = v_{\rho_{12, 2}}({E_{2, 12+}}') = v_{\rho_{12, 3}}({E_{2, 12+}}') = 1$, $v_{i / (2 \sqrt{3})}({E_{1, 12+12}}') = v_{\rho_{12, 2}}({E_{1, 12+12}}') = 1$, and $v_{\rho_{12, 1}}({E_{1, 12+3}}') = v_{\rho_{12, 3}}({E_{1, 12+3}}') = 1$.\\

\subsection{$\Gamma_0(12)+$}\quad

We have
\begin{equation*}
\Gamma_0(12)+ = \Gamma_0(12)+3,4,12 = \Gamma_0(12) \cup \Gamma_0(12) W_{12, 3} \cup \Gamma_0(12) W_{12, 4} \cup \Gamma_0(12) W_{12},
\end{equation*}
and $\Gamma_0(12)+ = T_{1/2}^{-1} (\Gamma_0(6)+3) T_{1/2}$.\\

\paragraph{\bf Fundamental domain}
We have a fundamental domain for $\Gamma_0(12)+$ as follows:
{\small \begin{equation}
\begin{split}
\mathbb{F}_{12+} = &\left\{|z + 1/3| \geqslant 1/6, \: - 1/2 \leqslant Re(z) \leqslant -1/4 \right\}
 \bigcup \left\{|z| \geqslant 1 / (2 \sqrt{3}), \: - 1/4 \leqslant Re(z) \leqslant 0 \right\}\\
 &\bigcup \left\{|z| > 1 / (2 \sqrt{3}), \: 0 < Re(z) \leqslant 1/4 \right\}
 \bigcup \left\{|z - 1/3| > 1/6, \: 1/4 < Re(z) \leqslant 1/2 \right\},
\end{split}
\end{equation}
}where $W_{12} : e^{i \theta} / (2 \sqrt{3}) \rightarrow e^{i  (\pi - \theta)} / (2 \sqrt{3})$ and $W_{12, 4} : e^{i \theta} / 6 + 1/3\rightarrow e^{i (\pi - \theta)} / 6 - 1/3$. Then, we have
\begin{equation}
\Gamma_0(12)+ = \langle \left( \begin{smallmatrix} 1 & 1 \\ 0 & 1 \end{smallmatrix} \right), \: W_{12}, \: W_{12, 4} \rangle.
\end{equation}
\begin{figure}[hbtp]
\begin{center}
\includegraphics[width=1.5in]{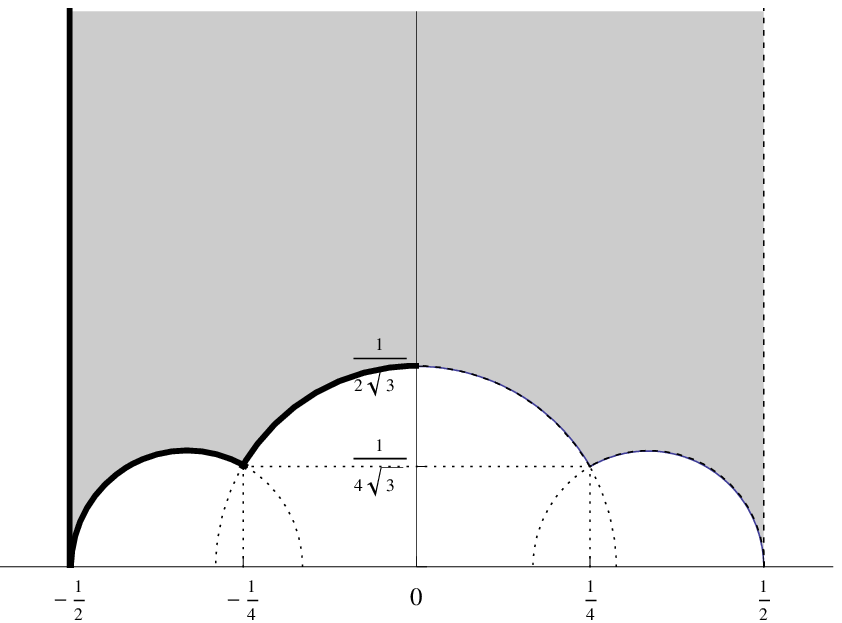}
\end{center}
\caption{$\Gamma_0(12)+$}
\end{figure}

\paragraph{\bf Valence formula}
The cusps of $\Gamma_0(12)+$ are $\infty$ and $0$, and the elliptic points are $i / (2 \sqrt{3})$ and $\rho_{12, 1} = - 1/4 + i / (4 \sqrt{3})$. Let $f$ be a modular function of weight $k$ for $\Gamma_0(12)+$, which is not identically zero. We have
\begin{equation}
v_{\infty}(f) + v_{-1/2}(f) + \frac{1}{2} v_{i / (2 \sqrt{3})} (f) + \frac{1}{2} v_{\rho_{12, 1}} (f) + \sum_{\begin{subarray}{c} p \in \Gamma_0(12)+ \setminus \mathbb{H} \\ p \ne i / (2 \sqrt{3}), \rho_{12, 1}\end{subarray}} v_p(f) = \frac{k}{2}.
\end{equation}

Furthermore, the stabilizer of the elliptic point $i / (2 \sqrt{3})$ (resp. $\rho_{12, 1}$) is $\left\{ \pm I, \pm W_{12} \right\}$ (resp. $\left\{ \pm I, \pm W_{12, 3} \right\}$).\\

\paragraph{\bf For the cusp $\infty$}
We have $\Gamma_{\infty} = \left\{ \pm \left( \begin{smallmatrix}1 & n \\ 0 & 1\end{smallmatrix} \right) \: ; \: n \in \mathbb{Z} \right\}$, and we have the Eisenstein series for the cusp $\infty$ associated with $\Gamma_0(12)+$ of weight $k \geqslant 4$:
\begin{equation}
E_{k, 12+}^{\infty}(z) := \frac{2^k 3^{k/2} E_k(12 z) - 2 \ 3^{k/2} E_k(6 z) + 2^k E_k(4 z) + 3^{k/2} E_k(3 z) - 2 E_k(2 z) + E_k(z)}{(3^{k/2} + 1)(2^k - 1)}.
\end{equation}\quad

\paragraph{\bf For the cusp $-1/2$}
We have $\Gamma_{-1/2} = \left\{ \pm \left( \begin{smallmatrix} 3n+1 & 3n/2 \\ -6n & -3n+1 \end{smallmatrix} \right) \: ; \: n \in \mathbb{Z} \right\}$ and $\gamma_{-1/2} = W_{12+, 6}$, and we have the Eisenstein series for the cusp $-1/2$ associated with $\Gamma_0(12)+$ of weight $k \geqslant 4$:
\begin{multline}
E_{k, 12+}^{-1/2}(z) :=\\ \frac{- 2 \ 2^{k/2} (2^k 3^{k/2} E_k(12 z) - 3^{k/2} (2^k + 1) E_k(6 z) + 2^k E_k(4 z) + 3^{k/2} E_k(3 z) - (2^k + 1) E_k(2 z) + E_k(z))}{(3^{k/2} + 1)(2^k - 1)}.
\end{multline}
We also have $\gamma_{-1/2}^{-1} \ \Gamma_0(12)+ \ \gamma_{-1/2} = \Gamma_0(12)+$.\\

\paragraph{\bf The space of modular forms}
We define
\begin{equation*}
\Delta_{12+}^{\infty} := \Delta_{12}^{\infty} \Delta_{12}^0 \Delta_{12}^{-1/3} \Delta_{12}^{-1/4}, \quad \Delta_{12+}^{-1/2} := (\Delta_{12}^{-1/2})^2 (\Delta_{12}^{-1/6})^2,
\end{equation*}
which are $2$nd semimodular forms for $\Gamma_0(12)+$ of weight $2$.

Now, we have $M_k(\Gamma_0(12)+) = \mathbb{C} E_{k, 12+}^{\infty} \oplus \mathbb{C} E_{k, 12+}^{-1/2} \oplus S_k(\Gamma_0(12)+)$ and $S_k(\Gamma_0(12)+) = \Delta_{12+} M_{k - 4}(\Gamma_0(12)+)$ for every even integer $k \geqslant 4$. Then, we have $M_{4 n + 2}(\Gamma_0(12)+) = {E_{2, 12+}}' M_{4 n}(\Gamma_0(12)+)$ and
\begin{align*}
M_{4 n}(\Gamma_0(12)+) &= \mathbb{C} (E_{4, 12+}^{\infty})^n \oplus \mathbb{C} (E_{4, 12+}^{\infty})^{n-1} \Delta_{12+} \oplus \cdots \oplus \mathbb{C} E_{4, 12+}^{\infty} (\Delta_{12+})^{n-1}\\
 &\oplus \mathbb{C} (E_{4, 12+}^{-1/2})^n \oplus \mathbb{C} (E_{4, 12+}^{-1/2})^{n-1} \Delta_{12+} \oplus \cdots \oplus \mathbb{C} E_{4, 12+}^{-1/2} (\Delta_{12+})^{n-1} \oplus \mathbb{C} (\Delta_{12+})^n.
\end{align*}

Furthremore, we can write
\begin{equation*}
M_{4 n}(\Gamma_0(12)+) = \mathbb{C} (\Delta_{12+}^{\infty})^{2n} \oplus \mathbb{C} (\Delta_{12+}^{\infty})^{2n-1} \Delta_{12+}^{-1/2} \oplus \cdots \oplus \mathbb{C} (\Delta_{12+}^{-1/2})^{2n}.
\end{equation*}\quad

\paragraph{\bf Hauptmodul}
We define the {\it hauptmodul} of $\Gamma_0(12)+$:
\begin{equation}
\begin{split}
J_{12+} &:= \Delta_{12+}^{-1/2} / \Delta_{12+}^{\infty} \: (= \eta^{-6}(z) \eta^{12}(2 z) \eta^{-6}(3 z) \eta^{-6}(4 z) \eta^{12}(6 z) \eta^{-6}(12 z))\\
 &= \frac{1}{q} + 6 + 15 q + 32 q^2 + 87 q^3 + \cdots,
\end{split}
\end{equation}
where $v_{\infty}(J_{12+}) = -1$ and $v_{-1/2}(J_{12+}) = 1$. Then, we have
\begin{equation}
J_{12+} : \partial \mathbb{F}_{12+} \setminus \{z \in \mathbb{H} \: ; \: Re(z) = \pm 1/2\} \to [0, 16] \subset \mathbb{R}.
\end{equation}\quad

\subsection{$\Gamma_0(12)+12 = \Gamma_0^{*}(12)$}

\paragraph{\bf Fundamental domain}
We have a fundamental domain for $\Gamma_0^{*}(12)$ as follows:
{\small \begin{equation}
\begin{split}
\mathbb{F}_{12+12} = &\left\{|z + 5/12| \geqslant 1/12, \: - 1/2 \leqslant Re(z) < - 1/3 \right\}
 \bigcup \left\{|z + 7/24| \geqslant 1/24, \: - 1/3 \leqslant Re(z) < - 2/7 \right\}\\
 &\bigcup \left\{|z| \geqslant 1 / (2 \sqrt{3}), \: - 2/7 \leqslant Re(z) \leqslant 0 \right\}
 \bigcup \left\{|z| > 1 / (2 \sqrt{3}), \: 0 < Re(z) \leqslant 2/7 \right\}\\
 &\bigcup \left\{|z - 7/24| > 1/24, \: 2/7 < Re(z) \leqslant 1/3 \right\}
 \bigcup \left\{|z - 5/12| > 1/12, \: 1/3 < Re(z) < 1/2 \right\},
\end{split}
\end{equation}
}where $W_{12} :  e^{i \theta} / (2 \sqrt{3}) \rightarrow e^{i  (\pi - \theta)} / (2 \sqrt{3})$, $\left(\begin{smallmatrix} -7 & 2 \\ 24 & -7 \end{smallmatrix} \right) : (e^{i \theta} + 7)/24 \rightarrow (e^{i (\pi - \theta)} - 7)/24$, and $\left(\begin{smallmatrix} -5 & 2 \\ 12 & -5 \end{smallmatrix} \right) : (e^{i \theta} + 5)/12 \rightarrow (e^{i (\pi - \theta)} - 5)/12$. Then, we have
\begin{equation}
\Gamma_0^{*}(12) = \langle \left( \begin{smallmatrix} 1 & 1 \\ 0 & 1 \end{smallmatrix} \right), \: W_{12}, \: \left( \begin{smallmatrix} 5 & 2 \\ 12 & 5 \end{smallmatrix} \right), \: \left( \begin{smallmatrix} 7 & 2 \\ 24 & 7 \end{smallmatrix} \right) \rangle.
\end{equation}
\begin{figure}[hbtp]
\begin{center}
\includegraphics[width=1.5in]{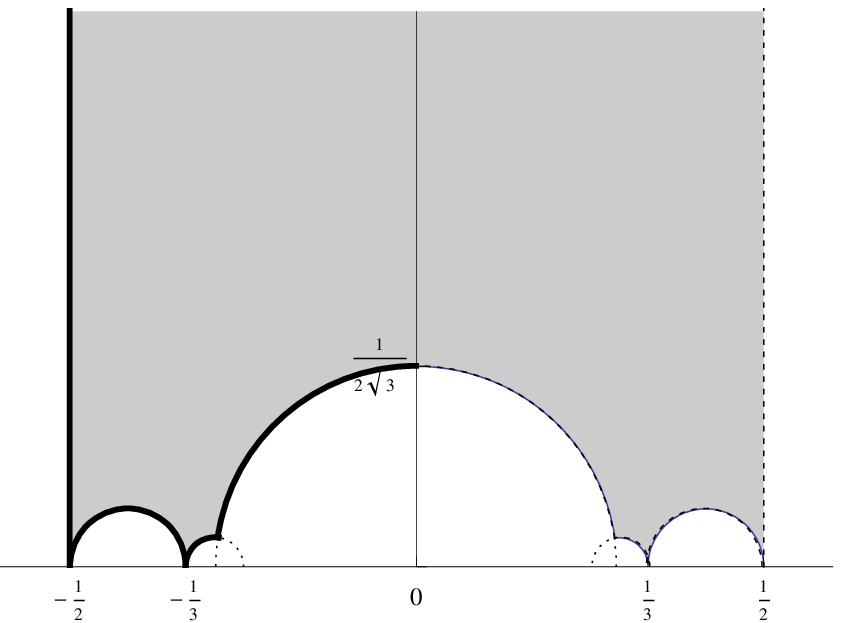}
\end{center}
\caption{$\Gamma_0^{*}(12)$}
\end{figure}

\paragraph{\bf Valence formula}
The cusps of $\Gamma_0^{*}(12)$ are $\infty$, $-1/3$, and $-1/2$, and the elliptic points are $i / (2 \sqrt{3})$ and $\rho_{12, 2} = - 2/7 + i / (14 \sqrt{3})$. Let $f$ be a modular function of weight $k$ for $\Gamma_0^{*}(12)$, which is not identically zero. We have
\begin{equation}
v_{\infty}(f) + v_{-1/3}(f) + v_{-1/2}(f) + \frac{1}{2} v_{i / (2 \sqrt{3})} (f) + \frac{1}{2} v_{\rho_{12, 2}} (f) + \sum_{\begin{subarray}{c} p \in \Gamma_0^{*}(12) \setminus \mathbb{H} \\ p \ne i / (2 \sqrt{3}), \rho_{12, 2}\end{subarray}} v_p(f) = k.
\end{equation}

Furthermore, the stabilizer of the elliptic point $i / (2 \sqrt{3})$ (resp. $\rho_{12, 2}$) is $\left\{ \pm I, \pm W_{12} \right\}$ (resp. $\left\{ \pm I, \pm \left( \begin{smallmatrix} 7 & - 2 \\ - 24 & 7 \end{smallmatrix} \right) W_{12} \right\}$).\\

\paragraph{\bf For the cusp $\infty$}
We have $\Gamma_{\infty} = \left\{ \pm \left( \begin{smallmatrix} 1 & n \\ 0 & 1 \end{smallmatrix} \right) \: ; \: n \in \mathbb{Z} \right\}$, and we have the Eisenstein series for the cusp $\infty$ associated with $\Gamma_0^{*}(12)$ of weight $k \geqslant 4$:
\begin{multline}
E_{k, 12+12}^{\infty}(z) :=\\ \frac{2^k 3^k E_k(12 z) + 3^{k/2} (3^{k/2} - 1) E_k(6 z) - 2^k 3^{k/2} E_k(4 z) - 3^k E_k(3 z) - (3^{k/2} - 1) E_k(2 z) + E_k(z)}{(3^k - 1)(2^k - 1)}.
\end{multline}\quad

\paragraph{\bf For the cusp $-1/3$}
We have $\Gamma_{-1/3} = \left\{ \pm \left( \begin{smallmatrix} 12n+1 & 4n \\ -36n & -12n+1\end{smallmatrix} \right) \: ; \: n \in \mathbb{Z} \right\}$ and $\gamma_{-1/3} = W_{12, 4}$, and we have the Eisenstein series for the cusp $-1/3$ associated with $\Gamma_0^{*}(12)$ of weight $k \geqslant 4$:
\begin{multline}
E_{k, 12+12}^{-1/3}(z) :=\\ \frac{- (2^k 3^{k/2} E_k(12 z) + 3^{k/2} (3^{k/2} - 1) E_k(6 z) - 2^k 3^{k/2} E_k(4 z) - 3^k E_k(3 z) + (3^{k/2} - 1) E_k(2 z) + E_k(z))}{(3^k - 1)(2^k - 1)}.
\end{multline}
We also have $\gamma_{-1/3}^{-1} \ \Gamma_0^{*}(12) \ \gamma_{-1/3} = \Gamma_0^{*}(12)$.\\

\paragraph{\bf For the cusp $-1/2$}
We have $\Gamma_{-1/2} = \left\{ \pm \left( \begin{smallmatrix} 6n+1 & 3n \\ -12n & -6n+1 \end{smallmatrix} \right) \: ; \: n \in \mathbb{Z} \right\}$ and $\gamma_{-1/2} = W_{12-, 6}$, and we have the Eisenstein series for the cusp $-1/2$ associated with $\Gamma_0^{*}(12)$ of weight $k \geqslant 4$:
\begin{multline}
E_{k, 12+12}^{-1/2}(z) :=\\ \frac{- (2^k 3^{k/2} E_k(12 z) - 3^{k/2} (2^k + 1) E_k(6 z) + 2^k E_k(4 z) + 3^{k/2} E_k(3 z) - (2^k + 1) E_k(2 z) + E_k(z))}{(3^{k/2} + 1)(2^k - 1)}.
\end{multline}
Note that we have $E_{k, 12+12}^{-1/2} = 2^{-1} 2^{-k/2} E_{k, 12+}^{-1/2}$. Note that $\gamma_{-1/2}^{-1} \ \Gamma_0^{*}(12) \ \gamma_{-1/2} \ne \Gamma_0^{*}(12)$.\\

\paragraph{\bf The space of modular forms}
We define
\begin{equation*}
\Delta_{12+12}^{\infty} := \Delta_{12}^{\infty} \Delta_{12}^0, \quad
 \Delta_{12+12}^{-1/3} := \Delta_{12}^{-1/3} \Delta_{12}^{-1/4}, \quad
 \Delta_{12+12}^{-1/2} := \Delta_{12}^{-1/2} \Delta_{12}^{-1/6},
\end{equation*}
which are $4$th semimodular forms for $\Gamma_0^{*}(12)$ of weight $1$.  Furthermore, we define
\begin{align*}
\Delta_{12, 1, 12+12} &:= ({E_{1, 12+12}}')^5 {E_{2, 12+}}' \Delta_{12+12}^{-1/3} \Delta_{12+12}^{-1/2} \Delta_{12},
&\Delta_{12, 2, 12+12} &:= ({E_{1, 12+12}}')^5 {E_{2, 12+}}' \Delta_{12+12}^{\infty} \Delta_{12+12}^{-1/2} \Delta_{12},\\
\Delta_{12, 3, 12+12} &:= ({E_{1, 12+12}}')^5 {E_{2, 12+}}' \Delta_{12+12}^{\infty} \Delta_{12+12}^{-1/3} \Delta_{12},
&\Delta_{12, 4, 12+12} &:= ({E_{1, 12+12}}')^4 \Delta_{12+12}^{-1/3} \Delta_{12+12}^{-1/2} (\Delta_{12})^2,\\
\Delta_{12, 5, 12+12} &:= ({E_{1, 12+12}}')^4 \Delta_{12+12}^{\infty} \Delta_{12+12}^{-1/2} (\Delta_{12})^2,
&\Delta_{12, 6, 12+12} &:= ({E_{1, 12+12}}')^4 \Delta_{12+12}^{\infty} \Delta_{12+12}^{-1/3} (\Delta_{12})^2,\\
\Delta_{12, 7, 12+12} &:= (\Delta_{12+12}^{-1/3})^2 \Delta_{12+12}^{-1/2} (\Delta_{12})^3,
&\Delta_{12, 8, 12+12} &:= \Delta_{12+12}^{\infty} (\Delta_{12+12}^{-1/2})^2 (\Delta_{12})^3,\\
\Delta_{12, 9, 12+12} &:= (\Delta_{12+12}^{\infty})^2 \Delta_{12+12}^{-1/3} (\Delta_{12})^3.
\end{align*}
In addition, we denote
\begin{allowdisplaybreaks}
\begin{align*}
A_1 &= {E_{1, 12+12}}' \Delta_{12} (\mathbb{C} \Delta_{12+12}^{\infty} \Delta_{12+12}^{-1/2} \oplus \mathbb{C} \Delta_{12+12}^{-1/3} \Delta_{12+12}^{-1/2} \oplus \mathbb{C} \Delta_{12+12}^{\infty} \Delta_{12+12}^{-1/3}),\\
A_2 &= (\Delta_{12})^2 (\mathbb{C} \Delta_{12+12}^{\infty} \Delta_{12+12}^{-1/2} \oplus \mathbb{C} \Delta_{12+12}^{-1/3} \Delta_{12+12}^{-1/2} \oplus \mathbb{C} \Delta_{12+12}^{\infty} \Delta_{12+12}^{-1/3}),\\
B_1 &= \mathbb{C} (E_{12, 12+12}^{\infty})^n \oplus \mathbb{C} (E_{12, 12+12}^{\infty})^{n-1} \Delta_{12, 1, 12+12} \oplus \mathbb{C} (E_{12, 12+12}^{\infty})^{n-1} \Delta_{12, 4, 12+12}\\
 &\qquad \oplus \mathbb{C} (E_{12, 12+12}^{\infty})^{n-1} \Delta_{12, 7, 12+12} \oplus \mathbb{C} (E_{12, 12+12}^{\infty})^{n-1} (\Delta_{12})^4 \oplus \mathbb{C} (E_{12, 12+12}^{\infty})^{n-2} \Delta_{12, 1, 12+12} (\Delta_{12})^4\\
&\qquad \oplus \cdots \oplus \mathbb{C} \Delta_{12, 7, 12+12} (\Delta_{12})^{4(n-1)} \oplus \mathbb{C} (\Delta_{12})^{4n},\\
B_2 &= \mathbb{C} (E_{12, 12+12}^{-1/3})^n \oplus \mathbb{C} (E_{12, 12+12}^{-1/3})^{n-1} \Delta_{12, 2, 12+12} \oplus \mathbb{C} (E_{12, 12+12}^{-1/3})^{n-1} \Delta_{12, 5, 12+12}\\
 &\qquad \oplus \mathbb{C} (E_{12, 12+12}^{-1/3})^{n-1} \Delta_{12, 8, 12+12} \oplus \mathbb{C} (E_{12, 12+12}^{-1/3})^{n-1} (\Delta_{12})^4 \oplus \mathbb{C} (E_{12, 12+12}^{-1/3})^{n-2} \Delta_{12, 2, 12+12} (\Delta_{12})^4\\
&\qquad \oplus \cdots \oplus \mathbb{C} \Delta_{12, 8, 12+12} (\Delta_{12})^{4(n-1)} \oplus \mathbb{C} (\Delta_{12})^{4n},\\
B_3 &= \mathbb{C} (E_{12, 12+12}^{-1/2})^n \oplus \mathbb{C} (E_{12, 12+12}^{-1/2})^{n-1} \Delta_{12, 3, 12+12} \oplus \mathbb{C} (E_{12, 12+12}^{-1/2})^{n-1} \Delta_{12, 6, 12+12}\\
 &\qquad \oplus \mathbb{C} (E_{12, 12+12}^{-1/2})^{n-1} \Delta_{12, 9, 12+12} \oplus \mathbb{C} (E_{12, 12+12}^{-1/2})^{n-1} (\Delta_{12})^4 \oplus \mathbb{C} (E_{12, 12+12}^{-1/2})^{n-2} \Delta_{12, 3, 12+12} (\Delta_{12})^4\\
&\qquad \oplus \cdots \oplus \mathbb{C} \Delta_{12, 8, 12+12} (\Delta_{12})^{4(n-1)} \oplus \mathbb{C} (\Delta_{12})^{4n}.
\end{align*}
\end{allowdisplaybreaks}

Now, we have
\begin{gather*}
M_k(\Gamma_0^{*}(12)) = \mathbb{C} E_{k, 12+12}^{\infty} \oplus \mathbb{C} E_{k, 12+12}^{-1/3} \oplus \mathbb{C} E_{k, 12+12}^{-1/2} \oplus S_k(\Gamma_0^{*}(12)),\\
S_k(\Gamma_0^{*}(12)) = (\mathbb{C} \Delta_{12, 1, 12+12} \oplus \mathbb{C} \Delta_{12, 2, 12+12} \oplus \cdots \oplus \mathbb{C} \Delta_{12, 9, 12+12} \oplus \mathbb{C} (\Delta_{12})^4) M_{k - 12}(\Gamma_0^{*}(12))
\end{gather*}
for every even integer $k \geqslant 4$. Then, we have
\begin{allowdisplaybreaks}
\begin{align*}
M_{12 n}(\Gamma_0^{*}(12)) &= B_1 \oplus B_2 \oplus B_3 \oplus \mathbb{C} (\Delta_{12})^{4n},\\
M_{12 n + 2}(\Gamma_0^{*}(12)) &= {E_{2, 12+}}' M_{12 n}(\Gamma_0^{*}(12)) \oplus \mathbb{C} {E_{1, 12+12}}' \Delta_{12+12}^{\infty} (\Delta_{12})^{4n},\\
M_{12 n + 4}(\Gamma_0^{*}(12)) &= E_{4, 12+12}^{\infty} (B_1 \oplus \mathbb{C} (\Delta_{12})^{4n}) \oplus E_{4, 12+12}^{-1/3} (B_2 \oplus \mathbb{C} (\Delta_{12})^{4n}) \oplus E_{4, 12+12}^{-1/2} (B_3 \oplus \mathbb{C} (\Delta_{12})^{4n})\\
&\qquad \qquad \oplus (\mathbb{C} {E_{1, 12+12}}' \oplus \mathbb{C} \Delta_{12+12}^{\infty}) (\Delta_{12})^{4n+1},\\
M_{12 n + 6}(\Gamma_0^{*}(12)) &= E_{6, 12+12}^{\infty} (B_1 \oplus \mathbb{C} (\Delta_{12})^{4n}) \oplus E_{6, 12+12}^{-1/3} (B_2 \oplus \mathbb{C} (\Delta_{12})^{4n}) \oplus E_{6, 12+12}^{-1/2} (B_3 \oplus \mathbb{C} (\Delta_{12})^{4n})\\
&\qquad \qquad \oplus A_1 (\Delta_{12})^{4n},\\
M_{12 n + 8}(\Gamma_0^{*}(12)) &= E_{8, 12+12}^{\infty} (B_1 \oplus \mathbb{C} (\Delta_{12})^{4n}) \oplus E_{8, 12+12}^{-1/3} (B_2 \oplus \mathbb{C} (\Delta_{12})^{4n}) \oplus E_{8, 12+12}^{-1/2} (B_3 \oplus \mathbb{C} (\Delta_{12})^{4n})\\
&\qquad \qquad \oplus {E_{2, 12+}}' A_1 (\Delta_{12})^{4n} \oplus A_2 (\Delta_{12})^{4n},\\
M_{12 n + 10}(\Gamma_0^{*}(12)) &= E_{10, 12+12}^{\infty} (B_1 \oplus \mathbb{C} (\Delta_{12})^{4n}) \oplus E_{10, 12+12}^{-1/3} (B_2 \oplus \mathbb{C} (\Delta_{12})^{4n}) \oplus E_{10, 12+12}^{-1/2} (B_3 \oplus \mathbb{C} (\Delta_{12})^{4n})\\
&\qquad \qquad \oplus \mathbb{C} ({E_{1, 12+12}}')^4 A_1 (\Delta_{12})^{4n} \oplus {E_{2, 12+}}' A_2 (\Delta_{12})^{4n+1} \oplus \mathbb{C} {E_{1, 12+12}}' (\Delta_{12})^{4n+3}.
\end{align*}
\end{allowdisplaybreaks}

Furthermore, we can write
\begin{equation*}
M_k(\Gamma_0^{*}(12)) = {E_{\overline{k}, 12+12}}' (\mathbb{C} (\Delta_{12+12}^{\infty})^{n} \oplus \mathbb{C} (\Delta_{12+12}^{\infty})^{n-1} \Delta_{12+12}^{-1/3} \oplus \cdots \oplus \mathbb{C} (\Delta_{12+12}^{-1/3})^{n}),
\end{equation*}
where $n = \dim(M_k(\Gamma_0^{*}(12))) - 1 = \lfloor k - 2 (k/4 -\lfloor k/4 \rfloor) \rfloor$, and where ${E_{\overline{k}, 12+12}}' := 1$ and ${E_{1, 12+12}}'$, when $k \equiv 0$ and $2 \pmod{4}$, respectively.\\

\paragraph{\bf Hauptmodul}
We define the {\it hauptmodul} of $\Gamma_0^{*}(12)$:
\begin{equation}
J_{12+12} := \Delta_{12+12}^{-1/3} / \Delta_{12+12}^{\infty} \: (= \eta^{-4}(z) \eta^4(3 z) \eta^4(4 z) \eta^{-4}(12 z))
 = \frac{1}{q} + 4 + 14 q + 36 q^2 + 85 q^3 + \cdots,
\end{equation}
where $v_{\infty}(J_{12+12}) = -1$ and $v_{-1/3}(J_{12+12}) = 1$. Then, we have
\begin{equation}
J_{12+12} : \partial \mathbb{F}_{12+12} \setminus \{z \in \mathbb{H} \: ; \: Re(z) = \pm 1/2\} \to [-1, 7 + 4 \sqrt{3}] \subset \mathbb{R}.
\end{equation}
\newpage

\subsection{$\Gamma_0(12)+4$}\quad

We have $\Gamma_0(12)+4 = T_{1/2}^{-1} \Gamma_0(6) T_{1/2}$.\\

\paragraph{\bf Fundamental domain}
We have a fundamental domain for $\Gamma_0(12)+4$ as follows:
{\small \begin{equation}
\begin{split}
\mathbb{F}_{12+4} = &\left\{|z + 1/3| \geqslant 1/6, \: - 1/2 \leqslant Re(z) < -1/6 \right\}
 \bigcup \left\{|z + 1/12| \geqslant 1/12, \: - 1/6 \leqslant Re(z) \leqslant 0 \right\}\\
 &\bigcup \left\{|z - 1/12| > 1/12, \: 0 < Re(z) \leqslant 1/6 \right\}
 \bigcup \left\{|z - 1/3| > 1/6, \: 1/6 < Re(z) < 1/2 \right\},
\end{split}
\end{equation}
}where $\left(\begin{smallmatrix} -1 & 0 \\ 12 & -1 \end{smallmatrix} \right) : (e^{i \theta} + 1) / 12 \rightarrow (e^{i (\pi - \theta)} - 1) / 12$ and $W_{12,4} : 1/3 + e^{i \theta} / 6 \rightarrow -1/3 + e^{i (\pi - \theta)} / 6$. Then, we have
\begin{equation}
\Gamma_0(12)+4 = \langle - I, \: \left( \begin{smallmatrix} 1 & 1 \\ 0 & 1 \end{smallmatrix} \right), \: W_{12, 4}, \: \left( \begin{smallmatrix} 1 & 0 \\ 12 & 1 \end{smallmatrix} \right) \rangle.
\end{equation}
\begin{figure}[hbtp]
\begin{center}
\includegraphics[width=1.5in]{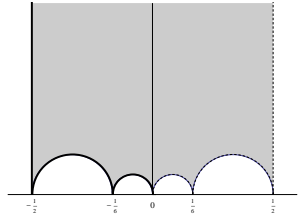}
\end{center}
\caption{$\Gamma_0(12)+4$}
\end{figure}

\paragraph{\bf Valence formula}
The cusps of $\Gamma_0(12)+4$ are $\infty$, $0$, $-1/2$, and $-1/6$. Let $f$ be a modular function of weight $k$ for $\Gamma_0(12)+4$, which is not identically zero. We have
\begin{equation}
v_{\infty}(f) + v_0(f) + v_{-1/2} (f) + v_{-1/6} (f) + \sum_{p \in \Gamma_0(12)+4 \setminus \mathbb{H}} v_p(f) = k.
\end{equation}\quad

\paragraph{\bf For the cusp $\infty$}
We have $\Gamma_{\infty} = \left\{ \pm \left( \begin{smallmatrix}1 & n \\ 0 & 1\end{smallmatrix} \right) \: ; \: n \in \mathbb{Z} \right\}$, and we have the Eisenstein series for the cusp $\infty$ associated with $\Gamma_0(12)+4$:
\begin{equation}
E_{k, 12+4}^{\infty}(z) := \frac{2^k 3^k E_k(12 z) - 2 \ 3^k E_k(6 z) - 2^k E_k(4 z) + 3^k E_k(3 z) + E_k(2 z) - E_k(z)}{(3^k - 1)(2^k - 1)} \quad \text{for} \; k \geqslant 4.
\end{equation}\quad

\paragraph{\bf For the cusp $0$}
We have $\Gamma_0 = \left\{ \pm \left( \begin{smallmatrix}1 & 0 \\ 12 n & 1\end{smallmatrix} \right) \: ; \: n \in \mathbb{Z} \right\}$ and $\gamma_0 = W_{12}$, and we have the Eisenstein series for the cusp $0$ associated with $\Gamma_0(12)+4$:
\begin{equation}
E_{k, 12+4}^0(z) := \frac{3^{k/2} (2^k  E_k(12 z) - 2 E_k(6 z) - 2^k E_k(4 z) + E_k(3 z) + 2 E_k(2 z) - E_k(z))}{(3^k - 1)(2^k - 1)} \quad \text{for} \; k \geqslant 4.
\end{equation}
We also have $\gamma_0^{-1} (\Gamma_0(12)+4) \gamma_0 = \Gamma_0(12)+4$.\\

\paragraph{\bf For the cusp $-1/2$}
We have $\Gamma_{-1/2} = \left\{ \pm \left( \begin{smallmatrix}3 n + 1 & 3 n / 2 \\ - 6 n & - 3 n + 1\end{smallmatrix} \right) \: ; \: n \in \mathbb{Z} \right\}$ and $\gamma_{-1/2} = W_{12+, 6}$, and we have the Eisenstein series for the cusp $-1/2$ associated with $\Gamma_0(12)+4$ of weight $k \geqslant 4$:
\begin{multline}
E_{k, 12+4}^{-1/2}(z) :=\\ \frac{2 \ 2^{k/2} 3^{k/2} (2^k 3^k E_k(12 z) - 3^k (2^k + 1) E_k(6 z) - 2^k E_k(4 z) + 3^k E_k(3 z) + (2^k + 1) E_k(2 z) - E_k(z))}{(3^k - 1)(2^k - 1)}.
\end{multline}
We also have $\gamma_{-1/2}^{-1} (\Gamma_0(12)+4) \ \gamma_{-1/2} = \Gamma_0(12)+4$.\\

\paragraph{\bf For the cusp $-1/6$}
We have $\Gamma_{-1/6} = \left\{ \pm \left( \begin{smallmatrix}3 n + 1 & n / 2 \\ - 18 n & - 3 n + 1\end{smallmatrix} \right) \: ; \: n \in \mathbb{Z} \right\}$ and $\gamma_{-1/6} = W_{12+, 2}$, and we have the Eisenstein series for the cusp $-1/6$ associated with $\Gamma_0(12)+4$ of weight $k \geqslant 4$:
\begin{multline}
E_{k, 12+4}^{-1/6}(z) :=\\ \frac{- 2 \ 2^{k/2} (2^k 3^{k/2} E_k(12 z) - 3^{k/2} (2^k + 1) E_k(6 z) + 2^k E_k(4 z) + 3^{k/2} E_k(3 z) - (2^k + 1) E_k(2 z) + E_k(z))}{(3^k - 1)(2^k - 1)}.
\end{multline}
We also have $\gamma_{-1/6}^{-1} (\Gamma_0(12)+4) \ \gamma_{-1/6} = \Gamma_0(12)+4$.\\

\paragraph{\bf The space of modular forms}
We define
\begin{align*}
\Delta_{12+4}^{\infty} &:= \Delta_{12}^{\infty} \Delta_{12}^{-1/3},
&\Delta_{12+4}^0 &:= \Delta_{12}^0 \Delta_{12}^{-1/4}, \\
\Delta_{12+4}^{-1/2} &:= (\Delta_{12}^{-1/2})^2,
&\Delta_{12+4}^{-1/6} &:= (\Delta_{12}^{-1/6})^2,
\end{align*}
which are $4$th semimodular forms for $\Gamma_0^{*}(12)$ of weight $1$.

We have $M_k(\Gamma_0(12)+4) = \mathbb{C} E_{k, 12+4}^{\infty} \oplus \mathbb{C} E_{k, 12+4}^0 \oplus \mathbb{C} E_{k, 12+4}^{-1/2} \oplus \mathbb{C} E_{k, 12+4}^{-1/6} \oplus S_k(\Gamma_0(12)+4)$ and $S_k(\Gamma_0(12)+4) = \Delta_{12+} M_{k - 4}(\Gamma_0(12)+4)$ for every even integer $k \geqslant 4$. Then, we have $M_{4 n + 2}(\Gamma_0(12)+4) = {E_{2, 12+}}' M_{4 n}(\Gamma_0(12)+4) \oplus \mathbb{C} ({E_{1, 12+12}}')^2 (\Delta_{12+})^n \oplus \mathbb{C} ({E_{1, 12+3}}')^2 (\Delta_{12+})^n$ and
\begin{align*}
M_{4 n}(\Gamma_0(12)+4) = &\mathbb{C} (E_{4, 12+4}^{\infty})^n \oplus \mathbb{C} (E_{4, 12+4}^{\infty})^{n-1} \Delta_{12+} \oplus \cdots \oplus \mathbb{C} E_{4, 12+4}^{\infty} (\Delta_{12+})^{n-1}\\
 &\oplus \mathbb{C} (E_{4, 12+4}^0)^n \oplus \mathbb{C} (E_{4, 12+4}^0)^{n-1} \Delta_{12+} \oplus \cdots \oplus \mathbb{C} E_{4, 12+4}^0 (\Delta_{12+})^{n-1}\\
 &\oplus \mathbb{C} (E_{4, 12+4}^{-1/2})^n \oplus \mathbb{C} (E_{4, 12+4}^{-1/2})^{n-1} \Delta_{12+} \oplus \cdots \oplus \mathbb{C} E_{4, 12+4}^{-1/2} (\Delta_{12+})^{n-1}\\
 &\oplus \mathbb{C} (E_{4, 12+4}^{-1/6})^n \oplus \mathbb{C} (E_{4, 12+4}^{-1/6})^{n-1} \Delta_{12+} \oplus \cdots \oplus \mathbb{C} E_{4, 12+4}^{-1/6} (\Delta_{12+})^{n-1} \oplus \mathbb{C} (\Delta_{12+})^n.
\end{align*}

Furthermore, we can write
\begin{equation*}
M_{2 n}(\Gamma_0(12)+4) = \mathbb{C} (\Delta_{12+4}^{\infty})^{2n} \oplus \mathbb{C} (\Delta_{12+4}^{\infty})^{2n-1} \Delta_{12+4}^0 \oplus \cdots \oplus \mathbb{C} (\Delta_{12+4}^0)^{2n}.
\end{equation*}\quad

\paragraph{\bf Hauptmodul}
We define the {\it hauptmodul} of $\Gamma_0(12)+4$:
\begin{equation}
\begin{split}
J_{12+4} &:= \Delta_{12+4}^0 / \Delta_{12+4}^{\infty} \: (= \eta^4(z) \eta^{-4}(2 z) \eta^{-4}(3 z) \eta^4(4 z) \eta^4(6 z) \eta^{-4}(12 z))\\
 &= \frac{1}{q} - 4 + 6 q - 4 q^2 - 3 q^3 + \cdots,
\end{split}
\end{equation}
where $v_{\infty}(J_{12+4}) = -1$ and $v_0(J_{12+4}) = 1$. Then, we have
\begin{equation}
J_{12+4} : \partial \mathbb{F}_{12+4} \setminus \{z \in \mathbb{H} \: ; \: Re(z) = \pm 1/2\} \to [-9, 0] \subset \mathbb{R}.
\end{equation}\quad

\subsection{$\Gamma_0(12)+3$}

\paragraph{\bf Fundamental domain}
We have a fundamental domain for $\Gamma_0(12)+3$ as follows:
{\small \begin{equation}
\begin{split}
\mathbb{F}_{12+3} = &\left\{|z + 5/12| \geqslant 1/12, \: - 1/2 \leqslant Re(z) \leqslant - 3/8 \right\}
 \bigcup \left\{|z + 1/4| \geqslant 1 / (4 \sqrt{3}), \: - 3/8 < Re(z) \leqslant - 1/4 \right\}\\
 &\bigcup \left\{|z + 1/4| > 1 / (4 \sqrt{3}), \: -  1/4 < Re(z) < - 1/8 \right\}
 \bigcup \left\{|z + 1/12| \geqslant 1/12, \: - 1/8 \leqslant Re(z) \leqslant 0 \right\}\\
 &\bigcup \left\{|z - 1/12| > 1/12, \: 0 < Re(z) < 1/8 \right\}
 \bigcup \left\{|z - 1/4| \geqslant 1 / (4 \sqrt{3}), \: 1/8 \leqslant Re(z) \leqslant 1/4 \right\}\\
 &\bigcup \left\{|z - 1/4| > 1 / (4 \sqrt{3}), \: 1/4 < Re(z) \leqslant 3/8 \right\}
 \bigcup \left\{|z - 5/12| > 1/12, \: 3/8 < Re(z) < 1/2 \right\},
\end{split}
\end{equation}
}where $W_{12, 3} :  -1/4 + e^{i \theta} / (4 \sqrt{3}) \rightarrow -1/4 + e^{i  (\pi - \theta)} / (4 \sqrt{3})$, $\left( \begin{smallmatrix} 7 & 2 \\ 24 & 7 \end{smallmatrix} \right) W_{12, 3} :  1/4 + e^{i \theta} / (4 \sqrt{3}) \rightarrow 1/4 + e^{i  (\pi - \theta)} / (4 \sqrt{3})$, $\left(\begin{smallmatrix} -1 & 0 \\ 12 & -1 \end{smallmatrix} \right) : (e^{i \theta} + 1)/12 \rightarrow (e^{i (\pi - \theta)} - 1)/12$, and $\left(\begin{smallmatrix} -5 & 2 \\ 12 & -5 \end{smallmatrix} \right) : (e^{i \theta} + 5)/12 \rightarrow (e^{i (\pi - \theta)} - 5)/12$. Then, we have
\begin{equation}
\Gamma_0(12)+3 = \langle \left( \begin{smallmatrix} 1 & 1 \\ 0 & 1 \end{smallmatrix} \right), \: W_{12, 3}, \: \left( \begin{smallmatrix} 1 & 0 \\ 12 & 1 \end{smallmatrix} \right), \: \left( \begin{smallmatrix} 5 & 2 \\ 12 & 5 \end{smallmatrix} \right) \rangle.
\end{equation}
\begin{figure}[hbtp]
\begin{center}
\includegraphics[width=1.5in]{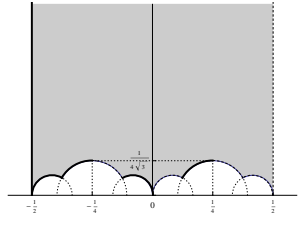}
\end{center}
\caption{$\Gamma_0(12)+3$}
\end{figure}

\paragraph{\bf Valence formula}
The cusps of $\Gamma_0(12)+3$ are $\infty$, $0$, and $-1/2$, and the elliptic points are $\rho_{12, 1}$ and $\rho_{12, 3} = 1/4 + i / (4 \sqrt{3})$. Let $f$ be a modular function of weight $k$ for $\Gamma_0(12)+3$, which is not identically zero. We have
\begin{equation}
v_{\infty}(f) + v_{-1/3}(f) + v_{-1/2}(f) + \frac{1}{2} v_{\rho_{12, 1}} (f) + \frac{1}{2} v_{\rho_{12, 3}} (f) + \sum_{\begin{subarray}{c} p \in \Gamma_0(12)+3 \setminus \mathbb{H} \\ p \ne \rho_{12, 1}, \rho_{12, 3}\end{subarray}} v_p(f) = k.
\end{equation}

Furthermore, the stabilizer of the elliptic point $\rho_{12, 1}$ (resp. $\rho_{12, 3}$) is $\left\{ \pm I, \pm W_{12, 3} \right\}$ (resp. $\left\{ \pm I, \pm \left( \begin{smallmatrix} 7 & 2 \\ 24 & 7 \end{smallmatrix} \right) W_{12, 3} \right\}$).\\

\paragraph{\bf For the cusp $\infty$}
We have $\Gamma_{\infty} = \left\{ \pm \left( \begin{smallmatrix} 1 & n \\ 0 & 1 \end{smallmatrix} \right) \: ; \: n \in \mathbb{Z} \right\}$, and we have the Eisenstein series for the cusp $\infty$ associated with $\Gamma_0(12)+3$:
\begin{equation}
E_{k, 12+3}^{\infty}(z) := \frac{2^k 3^{k/2} E_k(12 z) - 3^{k/2} E_k(6 z) + 2^k E_k(4 z) - E_k(2 z)}{(3^{k/2} + 1)(2^k - 1)} \quad \text{for} \; k \geqslant 4.
\end{equation}
Note that we have $E_{k, 12+3}^{\infty}(z) = E_{k, 6+3}^{\infty}(2 z)$.\\

\paragraph{\bf For the cusp $0$}
We have $\Gamma_0 = \left\{ \pm \left( \begin{smallmatrix} 1 & 0 \\ 12n & 1\end{smallmatrix} \right) \: ; \: n \in \mathbb{Z} \right\}$ and $\gamma_0 = W_{12}$, and we have the Eisenstein series for the cusp $0$ associated with $\Gamma_0(12)+3$:
\begin{equation}
E_{k, 12+3}^0(z) := \frac{- (3^{k/2} E_k(6 z) - 3^{k/2} E_k(3 z) + E_k(2 z) - E_k(z))}{(3^{k/2} + 1)(2^k - 1)} \quad \text{for} \; k \geqslant 4.
\end{equation}
Note that we have $E_{k, 12+3}^0(z) = 2^{-k/2} E_{k, 6+3}^0(z)$. We also have $\gamma_0^{-1} (\Gamma_0(12)+3) \gamma_0 = \Gamma_0(12)+3$.\\

\paragraph{\bf For the cusp $-1/2$}
We have $\Gamma_{-1/2} = \left\{ \pm \left( \begin{smallmatrix} 6n+1 & 3n \\ -12n & -6n+1 \end{smallmatrix} \right) \: ; \: n \in \mathbb{Z} \right\}$ and $\gamma_{-1/2} = W_{12-, 6}$, and we have the Eisenstein series for the cusp $-1/2$ associated with $\Gamma_0(12)+3$ of weight $k \geqslant 4$:
\begin{multline}
E_{k, 12+3}^{-1/2}(z) :=\\ \frac{- (2^k 3^{k/2} E_k(12 z) - 3^{k/2} (2^k + 1) E_k(6 z) + 2^k E_k(4 z) + 3^{k/2} E_k(3 z) - (2^k + 1) E_k(2 z) + E_k(z))}{(3^{k/2} + 1)(2^k - 1)}.
\end{multline}
Note that we have $E_{k, 12+3}^{-1/2}(z) = E_{k, 12+12}^{-1/2}(z) = 2^{-1} 2^{-k/2} E_{k, 12+}^{-1/2}(z)$. We also have $\gamma_{-1/2}^{-1} (\Gamma_0(12)+3) \ \gamma_{-1/2} = \Gamma_0(12)+3$.\\

\paragraph{\bf The space of modular forms}
We define
\begin{equation*}
\Delta_{12+3}^{\infty} := \Delta_{12}^{\infty} \Delta_{12}^{-1/4}, \quad
 \Delta_{12+3}^0 := \Delta_{12}^0 \Delta_{12}^{-1/3}, \quad
 \Delta_{12+3}^{-1/2} := \Delta_{12}^{-1/2} \Delta_{12}^{-1/6},
\end{equation*}
which are $4$th semimodular forms for $\Gamma_0(12)+3$ of weight $1$.  Furthermore, we define
\begin{align*}
\Delta_{12, 1, 12+3} &:= ({E_{1, 12+3}}')^5 {E_{2, 12+}}' \Delta_{12+3}^0 \Delta_{12+3}^{-1/2} \Delta_{12},
&\Delta_{12, 2, 12+3} &:= ({E_{1, 12+3}}')^5 {E_{2, 12+}}' \Delta_{12+3}^{\infty} \Delta_{12+3}^{-1/2} \Delta_{12},\\
\Delta_{12, 3, 12+3} &:= ({E_{1, 12+3}}')^5 {E_{2, 12+}}' \Delta_{12+3}^{\infty} \Delta_{12+3}^0 \Delta_{12},
&\Delta_{12, 4, 12+3} &:= ({E_{1, 12+3}}')^4 \Delta_{12+3}^0 \Delta_{12+3}^{-1/2} (\Delta_{12})^2,\\
\Delta_{12, 5, 12+3} &:= ({E_{1, 12+3}}')^4 \Delta_{12+3}^{\infty} \Delta_{12+3}^{-1/2} (\Delta_{12})^2,
&\Delta_{12, 6, 12+3} &:= ({E_{1, 12+3}}')^4 \Delta_{12+3}^{\infty} \Delta_{12+3}^0 (\Delta_{12})^2,\\
\Delta_{12, 7, 12+3} &:= (\Delta_{12+3}^0)^2 \Delta_{12+3}^{-1/2} (\Delta_{12})^3,
&\Delta_{12, 8, 12+3} &:= \Delta_{12+3}^{\infty} (\Delta_{12+3}^{-1/2})^2 (\Delta_{12})^3,\\
\Delta_{12, 9, 12+3} &:= (\Delta_{12+3}^{\infty})^2 \Delta_{12+3}^0 (\Delta_{12})^3.
\end{align*}
In addition, we denote
\begin{allowdisplaybreaks}
\begin{align*}
A_1 &= {E_{1, 12+3}}' \Delta_{12} (\mathbb{C} \Delta_{12+3}^{\infty} \Delta_{12+3}^{-1/2} \oplus \mathbb{C} \Delta_{12+3}^0 \Delta_{12+3}^{-1/2} \oplus \mathbb{C} \Delta_{12+3}^{\infty} \Delta_{12+3}^0),\\
A_2 &= (\Delta_{12})^2 (\mathbb{C} \Delta_{12+3}^{\infty} \Delta_{12+3}^{-1/2} \oplus \mathbb{C} \Delta_{12+3}^0 \Delta_{12+3}^{-1/2} \oplus \mathbb{C} \Delta_{12+3}^{\infty} \Delta_{12+3}^0),\\
B_1 &= \mathbb{C} (E_{12, 12+3}^{\infty})^n \oplus \mathbb{C} (E_{12, 12+3}^{\infty})^{n-1} \Delta_{12, 1, 12+3} \oplus \mathbb{C} (E_{12, 12+3}^{\infty})^{n-1} \Delta_{12, 4, 12+3}\\
 &\qquad \oplus \mathbb{C} (E_{12, 12+3}^{\infty})^{n-1} \Delta_{12, 7, 12+3} \oplus \mathbb{C} (E_{12, 12+3}^{\infty})^{n-1} (\Delta_{12})^4 \oplus \mathbb{C} (E_{12, 12+3}^{\infty})^{n-2} \Delta_{12, 1, 12+3} (\Delta_{12})^4\\
&\qquad \oplus \cdots \oplus \mathbb{C} \Delta_{12, 7, 12+3} (\Delta_{12})^{4(n-1)} \oplus \mathbb{C} (\Delta_{12})^{4n},\\
B_2 &= \mathbb{C} (E_{12, 12+3}^0)^n \oplus \mathbb{C} (E_{12, 12+3}^0)^{n-1} \Delta_{12, 2, 12+3} \oplus \mathbb{C} (E_{12, 12+3}^0)^{n-1} \Delta_{12, 5, 12+3}\\
 &\qquad \oplus \mathbb{C} (E_{12, 12+3}^0)^{n-1} \Delta_{12, 8, 12+3} \oplus \mathbb{C} (E_{12, 12+3}^0)^{n-1} (\Delta_{12})^4 \oplus \mathbb{C} (E_{12, 12+3}^0)^{n-2} \Delta_{12, 2, 12+3} (\Delta_{12})^4\\
&\qquad \oplus \cdots \oplus \mathbb{C} \Delta_{12, 8, 12+3} (\Delta_{12})^{4(n-1)} \oplus \mathbb{C} (\Delta_{12})^{4n},\\
B_3 &= \mathbb{C} (E_{12, 12+3}^{-1/2})^n \oplus \mathbb{C} (E_{12, 12+3}^{-1/2})^{n-1} \Delta_{12, 3, 12+3} \oplus \mathbb{C} (E_{12, 12+3}^{-1/2})^{n-1} \Delta_{12, 6, 12+3}\\
 &\qquad \oplus \mathbb{C} (E_{12, 12+3}^{-1/2})^{n-1} \Delta_{12, 9, 12+3} \oplus \mathbb{C} (E_{12, 12+3}^{-1/2})^{n-1} (\Delta_{12})^4 \oplus \mathbb{C} (E_{12, 12+3}^{-1/2})^{n-2} \Delta_{12, 3, 12+3} (\Delta_{12})^4\\
&\qquad \oplus \cdots \oplus \mathbb{C} \Delta_{12, 8, 12+3} (\Delta_{12})^{4(n-1)} \oplus \mathbb{C} (\Delta_{12})^{4n}.
\end{align*}
\end{allowdisplaybreaks}

Now, we have
\begin{gather*}
M_k(\Gamma_0(12)+3) = \mathbb{C} E_{k, 12+3}^{\infty} \oplus \mathbb{C} E_{k, 12+3}^0 \oplus \mathbb{C} E_{k, 12+3}^{-1/2} \oplus S_k(\Gamma_0(12)+3),\\
S_k(\Gamma_0(12)+3) = (\mathbb{C} \Delta_{12, 1, 12+3} \oplus \mathbb{C} \Delta_{12, 2, 12+3} \oplus \cdots \oplus \mathbb{C} \Delta_{12, 9, 12+3} \oplus \mathbb{C} (\Delta_{12})^4) M_{k - 12}(\Gamma_0(12)+3)
\end{gather*}
for every even integer $k \geqslant 4$. Then, we have
\begin{allowdisplaybreaks}
\begin{align*}
M_{12 n}(\Gamma_0(12)+3) &= B_1 \oplus B_2 \oplus B_3 \oplus \mathbb{C} (\Delta_{12})^{4n},\\
M_{12 n + 2}(\Gamma_0(12)+3) &= {E_{2, 12+}}' M_{12 n}(\Gamma_0(12)+3) \oplus \mathbb{C} {E_{1, 12+3}}' \Delta_{12+3}^{\infty} (\Delta_{12})^{4n},\\
M_{12 n + 4}(\Gamma_0(12)+3) &= E_{4, 12+3}^{\infty} (B_1 \oplus \mathbb{C} (\Delta_{12})^{4n}) \oplus E_{4, 12+3}^0 (B_2 \oplus \mathbb{C} (\Delta_{12})^{4n}) \oplus E_{4, 12+3}^{-1/2} (B_3 \oplus \mathbb{C} (\Delta_{12})^{4n})\\
&\qquad \qquad \oplus (\mathbb{C} {E_{1, 12+3}}' \oplus \mathbb{C} \Delta_{12+3}^{\infty}) (\Delta_{12})^{4n+1},\\
M_{12 n + 6}(\Gamma_0(12)+3) &= E_{6, 12+3}^{\infty} (B_1 \oplus \mathbb{C} (\Delta_{12})^{4n}) \oplus E_{6, 12+3}^0 (B_2 \oplus \mathbb{C} (\Delta_{12})^{4n}) \oplus E_{6, 12+3}^{-1/2} (B_3 \oplus \mathbb{C} (\Delta_{12})^{4n})\\
&\qquad \qquad \oplus A_1 (\Delta_{12})^{4n},\\
M_{12 n + 8}(\Gamma_0(12)+3) &= E_{8, 12+3}^{\infty} (B_1 \oplus \mathbb{C} (\Delta_{12})^{4n}) \oplus E_{8, 12+3}^0 (B_2 \oplus \mathbb{C} (\Delta_{12})^{4n}) \oplus E_{8, 12+3}^{-1/2} (B_3 \oplus \mathbb{C} (\Delta_{12})^{4n})\\
&\qquad \qquad \oplus {E_{2, 12+}}' A_1 (\Delta_{12})^{4n} \oplus A_2 (\Delta_{12})^{4n},\\
M_{12 n + 10}(\Gamma_0(12)+3) &= E_{10, 12+3}^{\infty} (B_1 \oplus \mathbb{C} (\Delta_{12})^{4n}) \oplus E_{10, 12+3}^0 (B_2 \oplus \mathbb{C} (\Delta_{12})^{4n}) \oplus E_{10, 12+3}^{-1/2} (B_3 \oplus \mathbb{C} (\Delta_{12})^{4n})\\
&\qquad \qquad \oplus \mathbb{C} ({E_{1, 12+3}}')^4 A_1 (\Delta_{12})^{4n} \oplus {E_{2, 12+}}' A_2 (\Delta_{12})^{4n+1} \oplus \mathbb{C} {E_{1, 12+3}}' (\Delta_{12})^{4n+3}.
\end{align*}
\end{allowdisplaybreaks}

Furthermore, we can write
\begin{equation*}
M_k(\Gamma_0(12)+3) = {E_{\overline{k}, 12+3}}' (\mathbb{C} (\Delta_{12+3}^{\infty})^{n} \oplus \mathbb{C} (\Delta_{12+3}^{\infty})^{n-1} \Delta_{12+3}^0 \oplus \cdots \oplus \mathbb{C} (\Delta_{12+3}^0)^{n}),
\end{equation*}
where $n = \dim(M_k(\Gamma_0(12)+3)) - 1 = \lfloor k - 2 (k/4 -\lfloor k/4 \rfloor) \rfloor$, and where ${E_{\overline{k}, 12+3}}' := 1$ and ${E_{1, 12+3}}'$, when $k \equiv 0$ and $2 \pmod{4}$, respectively.\\

\paragraph{\bf Hauptmodul}
We define the {\it hauptmodul} of $\Gamma_0(12)+3$:
\begin{equation}
J_{12+3} := \Delta_{12+3}^0 / \Delta_{12+3}^{\infty} \: (= \eta^2(z) \eta^2(3 z) \eta^{-2}(4 z) \eta^{-2}(12 z))
 = \frac{1}{q} - 2 - q + 7 q^3 - 9 q^5 + \cdots,
\end{equation}
where $v_{\infty}(J_{12+3}) = -1$ and $v_0(J_{12+3}) = 1$. Then, we have
{\small \begin{equation}
\begin{split}
J_{12+3} : \qquad \left\{|z + 1/12| = 1/12, \: -1/8 \leqslant Re(z) \leqslant 0\right\} &\to -2 + [0, 2] \subset \mathbb{R},\\
\left\{|z + 5/12| = 1/12, \: -1/2 \leqslant Re(z) \leqslant -3/8 \right\} &\to -2 - [0, 2] \subset \mathbb{R},\\
\left\{|z + 1/4| = 1 / (4 \sqrt{3}), \: -3/8 \leqslant Re(z) \leqslant -1/4\right\} &\to -2 + 2 \sqrt{3} [0, 1],\\
\left\{|z - 1/4| = 1 / (4 \sqrt{3}), \: 1/8 \leqslant Re(z) \leqslant 1/4\right\} &\to -2 - 2 \sqrt{3} [0, 1].
\end{split}
\end{equation}}

\begin{figure}[hbtp]
\begin{center}
{{\small Lower arcs of $\partial \mathbb{F}_{12+3}$}\includegraphics[width=2.5in]{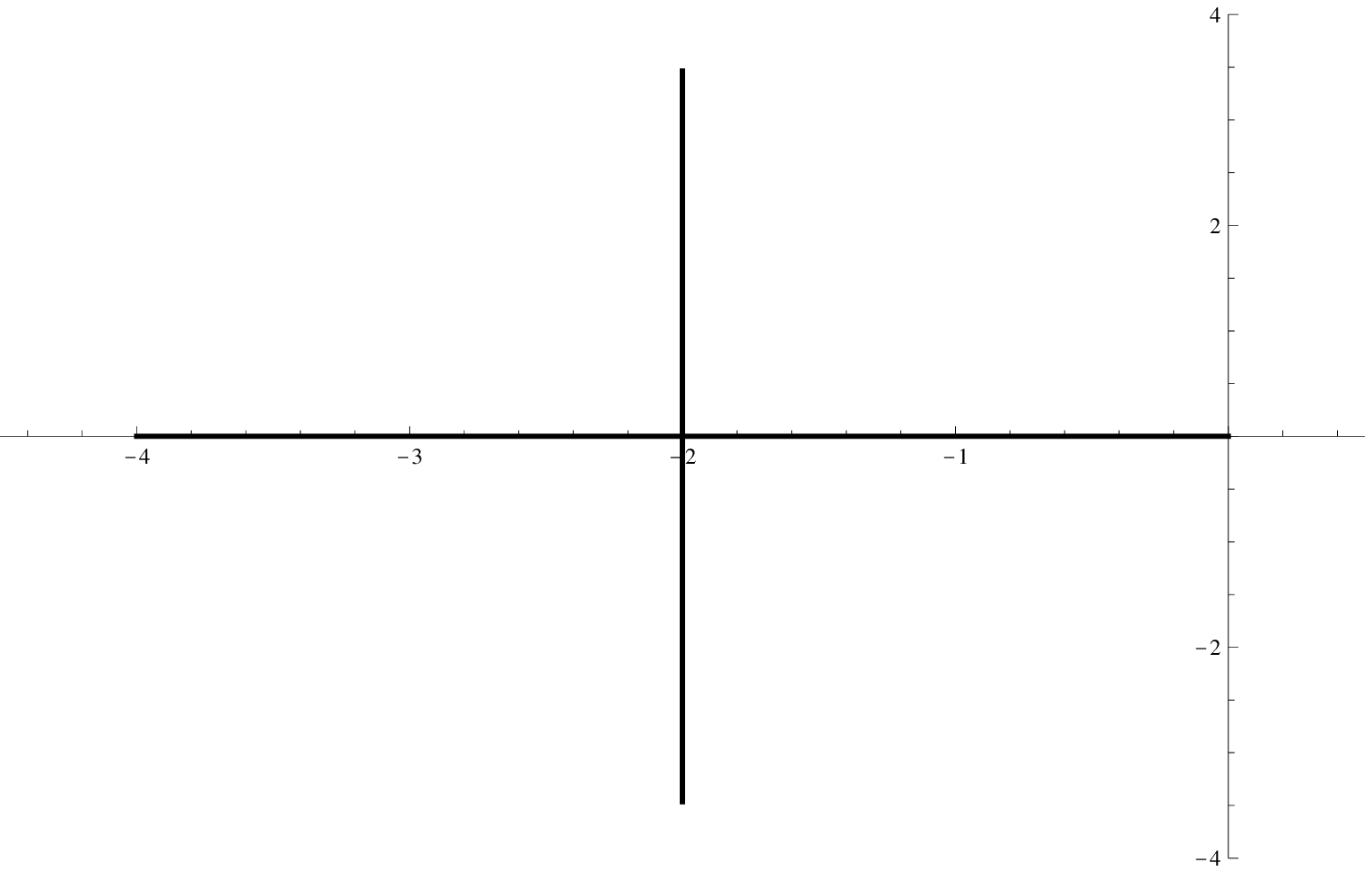}}
\end{center}
\caption{Image by $J_{12+3}$}\label{Im-J12E}
\end{figure}

\subsection{$\Gamma_0(12)$}

\paragraph{\bf Fundamental domain}
We have a fundamental domain for $\Gamma_0(12)$ as follows:
{\small \begin{equation}
\begin{split}
\mathbb{F}_{12} = &\left\{|z + 5/12| \geqslant 1/12, \: - 1/2 \leqslant Re(z) < -1/3 \right\}
 \bigcup \left\{|z + 7/24| \geqslant 1/24, \: - 1/3 \leqslant Re(z) < -1/4 \right\}\\
 &\bigcup \left\{|z + 5/24| \geqslant 1/24, \: - 1/4 \leqslant Re(z) < -1/6 \right\}
 \bigcup \left\{|z + 1/12| \geqslant 1/12, \: - 1/6 \leqslant Re(z) \leqslant 0 \right\}\\
 &\bigcup \left\{|z - 1/12| > 1/12, \: 0 < Re(z) \leqslant 1/6 \right\}
 \bigcup \left\{|z - 5/24| > 1/24, \: 1/6 < Re(z) \leqslant 1/4 \right\}\\
 &\bigcup \left\{|z - 7/24| > 1/24, \: 1/4 < Re(z) \leqslant 1/3 \right\}
 \bigcup \left\{|z - 5/12| > 1/12, \: 1/3 < Re(z) < 1/2 \right\},
\end{split}
\end{equation}
}where $\left(\begin{smallmatrix} -1 & 0 \\ 12 & -1 \end{smallmatrix} \right) : (e^{i \theta} + 1) / 12 \rightarrow (e^{i (\pi - \theta)} - 1) / 12$, $\left(\begin{smallmatrix} -5 & 1 \\ 24 & -5 \end{smallmatrix} \right) : (e^{i \theta} + 5) / 24 \rightarrow (e^{i (\pi - \theta)} - 5) / 24$, $\left(\begin{smallmatrix} -7 & 2 \\ 24 & -7 \end{smallmatrix} \right) : (e^{i \theta} + 7) / 24 \rightarrow (e^{i (\pi - \theta)} - 7) / 24$, and $\left(\begin{smallmatrix} -5 & 2 \\ 12 & -5 \end{smallmatrix} \right) : (e^{i \theta} + 5) / 12 \rightarrow (e^{i (\pi - \theta)} - 5) / 12$. Then, we have
\begin{equation}
\Gamma_0(12) = \langle - I, \: \left( \begin{smallmatrix} 1 & 1 \\ 0 & 1 \end{smallmatrix} \right), \: \left( \begin{smallmatrix} 1 & 0 \\ 12 & 1 \end{smallmatrix} \right), \: \left( \begin{smallmatrix} 5 & 2 \\ 12 & 5 \end{smallmatrix} \right), \: \left( \begin{smallmatrix} 5 & 1 \\ 24 & 5 \end{smallmatrix} \right), \: \left( \begin{smallmatrix} 7 & 2 \\ 24 & 7 \end{smallmatrix} \right) \rangle.
\end{equation}
\begin{figure}[hbtp]
\begin{center}
\includegraphics[width=1.5in]{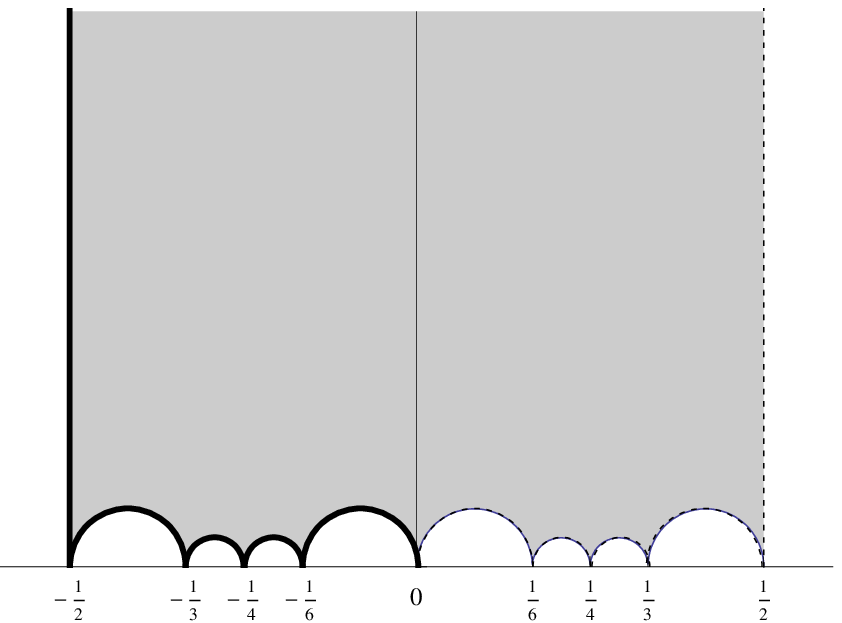}
\end{center}
\caption{$\Gamma_0(12)$}
\end{figure}

\paragraph{\bf Valence formula}
The cusps of $\Gamma_0(12)$ are $\infty$, $0$, $-1/3$, $-1/4$, $-1/2$, and $-1/6$. Let $f$ be a modular function of weight $k$ for $\Gamma_0(12)$, which is not identically zero. We have
\begin{equation}
v_{\infty}(f) + v_0(f) + v_{-1/3} (f) + v_{-1/4} (f) + v_{-1/2} (f) + v_{-1/6} (f) + \sum_{p \in \Gamma_0(12) \setminus \mathbb{H}} v_p(f) = 2 k.
\end{equation}\quad

\paragraph{\bf For the cusp $\infty$}
We have $\Gamma_{\infty} = \left\{ \pm \left( \begin{smallmatrix}1 & n \\ 0 & 1\end{smallmatrix} \right) \: ; \: n \in \mathbb{Z} \right\}$, and we have the Eisenstein series for the cusp $\infty$ associated with $\Gamma_0(12)$:
\begin{equation}
E_{k, 12}^{\infty}(z) := \frac{2^k 3^k E_k(12 z) - 3^k E_k(6 z) - 2^k E_k(4 z) + E_k(2 z)}{(3^k - 1)(2^k - 1)} \quad \text{for} \; k \geqslant 4.
\end{equation}
Note that we have $E_{k, 12}^{\infty}(z) = E_{k, 6}^{\infty}(2 z)$.\\

\paragraph{\bf For the cusp $0$}
We have $\Gamma_0 = \left\{ \pm \left( \begin{smallmatrix}1 & 0 \\ 12 n & 1\end{smallmatrix} \right) \: ; \: n \in \mathbb{Z} \right\}$ and $\gamma_0 = W_{12}$, and we have the Eisenstein series for the cusp $0$ associated with $\Gamma_0(12)$:
\begin{equation}
E_{k, 12}^0(z) := \frac{3^{k/2} (E_k(6 z) - E_k(3 z) - E_k(2 z) + E_k(z))}{(3^k - 1)(2^k - 1)} \quad \text{for} \; k \geqslant 4.
\end{equation}
Note that we have $E_{k, 12}^0(z) = 2^{-k/2} E_{k, 6}^0(z)$. We also have $\gamma_0^{-1} \ \Gamma_0(12) \ \gamma_0 = \Gamma_0(12)$.\\

\paragraph{\bf For the cusp $-1/3$}
We have $\Gamma_{-1/3} = \left\{ \pm \left( \begin{smallmatrix} 12n+1 & 4n \\ -36n & -12n+1 \end{smallmatrix} \right) \: ; \: n \in \mathbb{Z} \right\}$ and $\gamma_{-1/3} = W_{12, 4}$, and we have the Eisenstein series for the cusp $\infty$ associated with $\Gamma_0(12)$:
\begin{equation}
E_{k, 12}^{-1/3}(z) := \frac{- (3^k E_k(6 z) - 3^k E_k(3 z) - E_k(2 z) + E_k(z))}{(3^k - 1)(2^k - 1)} \quad \text{for} \; k \geqslant 4.
\end{equation}
Note that we have $E_{k, 12}^{-1/3}(z) = 2^{-k/2} E_{k, 6}^{-1/3}(z)$. We also have $\gamma_{-1/3}^{-1} \ \Gamma_0(12) \ \gamma_{-1/3} = \Gamma_0(12)$.\\

\paragraph{\bf For the cusp $-1/4$}
We have $\Gamma_{-1/4} = \left\{ \pm \left( \begin{smallmatrix} 12n+1 & 3n \\ -48n & -12n+1 \end{smallmatrix} \right) \: ; \: n \in \mathbb{Z} \right\}$ and $\gamma_{-1/4} = W_{12, 3}$, and we have the Eisenstein series for the cusp $-1/4$ associated with $\Gamma_0(12)$:
\begin{equation}
E_{k, 12}^{-1/4}(z) := \frac{- 3^{k/2} (2^k  E_k(12 z) - E_k(6 z) - 2^k E_k(4 z) + E_k(2 z))}{(3^k - 1)(2^k - 1)} \quad \text{for} \; k \geqslant 4.
\end{equation}
Note that we have $E_{k, 12}^{-1/4}(z) = E_{k, 6}^{-1/4}(2 z)$. We also have $\gamma_{-1/4}^{-1} \ \Gamma_0(12) \ \gamma_{-1/4} = \Gamma_0(12)$.\\

\paragraph{\bf For the cusp $-1/2$}
We have $\Gamma_{-1/2} = \left\{ \pm \left( \begin{smallmatrix} 6n+1 & 3n \\ -12n & -6n+1 \end{smallmatrix} \right) \: ; \: n \in \mathbb{Z} \right\}$ and $\gamma_{-1/2} = W_{12-, 6}$, and we have the Eisenstein series for the cusp $-1/2$ associated with $\Gamma_0(12)$ of weight $k \geqslant 4$:
\begin{multline}
E_{k, 12}^{-1/2}(z) :=\\ \frac{3^{k/2} (2^k 3^k E_k(12 z) - 3^k (2^k + 1) E_k(6 z) - 2^k E_k(4 z) + 3^k E_k(3 z) + (2^k + 1) E_k(2 z) - E_k(z))}{(3^k - 1)(2^k - 1)}.
\end{multline}
Note that we have $E_{k, 12}^{-1/2}(z) = 2^{-1} 2^{-k/2} E_{k, 12+4}^{-1/2}(z)$. We also have $\gamma_{-1/2}^{-1} \ \Gamma_0(12) \ \gamma_{-1/2} = \Gamma_0(12)$.\\

\paragraph{\bf For the cusp $-1/6$}
We have $\Gamma_{-1/6} = \left\{ \pm \left( \begin{smallmatrix} 6 n+1 & n \\ -36n & -6n + 1 \end{smallmatrix} \right) \: ; \: n \in \mathbb{Z} \right\}$ and $\gamma_{-1/6} = W_{12-, 2}$, and we have the Eisenstein series for the cusp $-1/6$ associated with $\Gamma_0(12)$ of weight $k \geqslant 4$:
\begin{multline}
E_{k, 12}^{-1/6}(z) :=\\ \frac{- (2^k 3^{k/2} E_k(12 z) - 3^{k/2} (2^k + 1) E_k(6 z) + 2^k E_k(4 z) + 3^{k/2} E_k(3 z) - (2^k + 1) E_k(2 z) + E_k(z))}{(3^k - 1)(2^k - 1)}.
\end{multline}
Note that we have $E_{k, 12}^{-1/6}(z) = 2^{-1} 2^{-k/2} E_{k, 12+4}^{-1/6}(z)$. We also have $\gamma_{-1/6}^{-1} \ \Gamma_0(12) \ \gamma_{-1/6} = \Gamma_0(12)$.\\

\paragraph{\bf The space of modular forms}
We define
\begin{align*}
\Delta_{6, 1, 12-} &:= (\Delta_{12}^0)^2 \Delta_{12}^{-1/3} \Delta_{12}^{-1/4} \Delta_{12}^{-1/2} \Delta_{12}^{-1/6},
&\Delta_{6, 2, 12-} &:= (\Delta_{12}^{\infty})^2 \Delta_{12}^{-1/3} \Delta_{12}^{-1/4} \Delta_{12}^{-1/2} \Delta_{12}^{-1/6},\\
\Delta_{6, 3, 12-} &:= \Delta_{12}^{\infty} \Delta_{12}^0 (\Delta_{12}^{-1/4})^2 \Delta_{12}^{-1/2} \Delta_{12}^{-1/6},
&\Delta_{6, 4, 12-} &:= \Delta_{12}^{\infty} \Delta_{12}^0 (\Delta_{12}^{-1/3})^2 \Delta_{12}^{-1/2} \Delta_{12}^{-1/6},\\
\Delta_{6, 5, 12-} &:= \Delta_{12}^{\infty} \Delta_{12}^0 \Delta_{12}^{-1/3} \Delta_{12}^{-1/4} (\Delta_{12}^{-1/6})^2,
&\Delta_{6, 6, 12-} &:= \Delta_{12}^{\infty} \Delta_{12}^0 \Delta_{12}^{-1/3} \Delta_{12}^{-1/4} (\Delta_{12}^{-1/2})^2.
\end{align*}

Now, we have
\begin{gather*}
M_k(\Gamma_0(12)) = \mathbb{C} E_{k, 12}^{\infty} \oplus \mathbb{C} E_{k, 12}^0 \oplus \mathbb{C} E_{k, 12}^{-1/3} \oplus \mathbb{C} E_{k, 12}^{-1/4} \oplus E_{k, 12}^{-1/2} \oplus \mathbb{C} E_{k, 12}^{-1/6} \oplus S_k(\Gamma_0(12)),\\
S_k(\Gamma_0(12)) = (\mathbb{C} \Delta_{6, 1, 12-} \oplus \mathbb{C} \Delta_{6, 2, 12-} \oplus \cdots \oplus \mathbb{C} \Delta_{6, 6, 12-} \oplus \mathbb{C} (\Delta_{12})^2) M_{k - 6}(\Gamma_0(12)),
\end{gather*}
for every even integer $k \geqslant 4$. Then, we have
\begin{allowdisplaybreaks}
\begin{align*}
M_{6 n}(\Gamma_0&(12)) = \mathbb{C} (E_{6, 12}^{\infty})^n \oplus \mathbb{C} (E_{6, 12}^{\infty})^{n-1} \Delta_{6, 1, 12-} \oplus \mathbb{C} (E_{6, 12}^{\infty})^{n-1} (\Delta_{12})^2 \oplus \cdots \oplus \mathbb{C} \Delta_{6, 1, 12-} (\Delta_{12})^{2(n-1)}\\
 &\oplus \mathbb{C} (E_{6, 12}^0)^n \oplus \mathbb{C} (E_{6, 12}^0)^{n-1} \Delta_{6, 2, 12-} \oplus \mathbb{C} (E_{6, 12}^0)^{n-1} (\Delta_{12})^2 \oplus \cdots \oplus \mathbb{C} \Delta_{6, 2, 12-} (\Delta_{12})^{2(n-1)}\\
 &\oplus \mathbb{C} (E_{6, 12}^{-1/3})^n \oplus \mathbb{C} (E_{6, 12}^{-1/3})^{n-1} \Delta_{6, 3, 12-} \oplus \mathbb{C} (E_{6, 12}^{-1/3})^{n-1} (\Delta_{12})^2 \oplus \cdots \oplus \mathbb{C} \Delta_{6, 3, 12-} (\Delta_{12})^{2(n-1)}\\
 &\oplus \mathbb{C} (E_{6, 12}^{-1/4})^n \oplus \mathbb{C} (E_{6, 12}^{-1/4})^{n-1} \Delta_{6, 4, 12-} \oplus \mathbb{C} (E_{6, 12}^{-1/4})^{n-1} (\Delta_{12})^2 \oplus \cdots \oplus \mathbb{C} \Delta_{6, 4, 12-} (\Delta_{12})^{2(n-1)}\\
 &\oplus \mathbb{C} (E_{6, 12}^{-1/2})^n \oplus \mathbb{C} (E_{6, 12}^{-1/2})^{n-1} \Delta_{6, 5, 12-} \oplus \mathbb{C} (E_{6, 12}^{-1/2})^{n-1} (\Delta_{12})^2 \oplus \cdots \oplus \mathbb{C} \Delta_{6, 5, 12-} (\Delta_{12})^{2(n-1)}\\
 &\oplus \mathbb{C} (E_{6, 12}^{-1/6})^n \oplus \mathbb{C} (E_{6, 12}^{-1/6})^{n-1} \Delta_{6, 6, 12-} \oplus \mathbb{C} (E_{6, 12}^{-1/6})^{n-1} (\Delta_{12})^2 \oplus \cdots \oplus \mathbb{C} \Delta_{6, 6, 12-} (\Delta_{12})^{2(n-1)}\\
 &\qquad \qquad \oplus \mathbb{C} (\Delta_{12})^{2n},\\
M_{6 n + 2}(\Gamma_0&(12)) = {E_{2, 12+}}' M_{6 n}(\Gamma_0(12))\\
 &\qquad \qquad \oplus (\mathbb{C} ({E_{1, 12+12}}')^2 \oplus \mathbb{C} ({E_{1, 12+3}}')^2 \oplus \mathbb{C} (\Delta_{12}^{\infty})^4 \oplus \mathbb{C} (\Delta_{12}^0)^4) (\Delta_{12})^{2n},\\
M_{6 n + 4}(\Gamma_0&(12)) = E_{4, 12}^{\infty}(\mathbb{C} (E_{6, 12}^{\infty})^n \oplus \mathbb{C} (E_{6, 12}^{\infty})^{n-1} \Delta_{6, 1, 12-} \oplus \mathbb{C} (E_{6, 12}^{\infty})^{n-1} (\Delta_{12})^2 \oplus \cdots \oplus \mathbb{C} (\Delta_{12})^{2n})\\
 &\oplus E_{4, 12}^0(\mathbb{C} (E_{6, 12}^0)^n \oplus \mathbb{C} (E_{6, 12}^0)^{n-1} \Delta_{6, 2, 12-} \oplus \mathbb{C} (E_{6, 12}^0)^{n-1} (\Delta_{12})^2 \oplus \cdots \oplus \mathbb{C} (\Delta_{12})^{2n})\\
 &\oplus E_{4, 12}^{-1/3}(\mathbb{C} (E_{6, 12}^{-1/3})^n \oplus \mathbb{C} (E_{6, 12}^{-1/3})^{n-1} \Delta_{6, 3, 12-} \oplus \mathbb{C} (E_{6, 12}^{-1/3})^{n-1} (\Delta_{12})^2 \oplus \cdots \oplus \mathbb{C} (\Delta_{12})^{2n})\\
 &\oplus E_{4, 12}^{-1/4}(\mathbb{C} (E_{6, 12}^{-1/4})^n \oplus \mathbb{C} (E_{6, 12}^{-1/4})^{n-1} \Delta_{6, 4, 12-} \oplus \mathbb{C} (E_{6, 12}^{-1/4})^{n-1} (\Delta_{12})^2 \oplus \cdots \oplus \mathbb{C} (\Delta_{12})^{2n})\\
 &\oplus E_{4, 12}^{-1/2}(\mathbb{C} (E_{6, 12}^{-1/2})^n \oplus \mathbb{C} (E_{6, 12}^{-1/2})^{n-1} \Delta_{6, 5, 12-} \oplus \mathbb{C} (E_{6, 12}^{-1/2})^{n-1} (\Delta_{12})^2 \oplus \cdots \oplus \mathbb{C} (\Delta_{12})^{2n})\\
 &\oplus E_{4, 12}^{-1/6}(\mathbb{C} (E_{6, 12}^{-1/6})^n \oplus \mathbb{C} (E_{6, 12}^{-1/6})^{n-1} \Delta_{6, 6, 12-} \oplus \mathbb{C} (E_{6, 12}^{-1/6})^{n-1} (\Delta_{12})^2 \oplus \cdots \oplus \mathbb{C} (\Delta_{12})^{2n})\\
 &\qquad \qquad \oplus (\mathbb{C} (\Delta_{12}^{\infty})^4 (\Delta_{12}^0)^4 \oplus \mathbb{C} (\Delta_{12}^{-1/3})^4 (\Delta_{12}^{-1/4})^4) (\Delta_{12})^{2n} \oplus \mathbb{C} {E_{1, 12+12}}' (\Delta_{12})^{2n+1}.
\end{align*}
\end{allowdisplaybreaks}

Furthermore, we can write
\begin{equation*}
M_{2 n}(\Gamma_0(12)) = \mathbb{C} (\Delta_{12}^{\infty})^{4n} \oplus \mathbb{C} (\Delta_{12}^{\infty})^{4n-1} \Delta_{12}^0 \oplus \cdots \oplus \mathbb{C} (\Delta_{12}^0)^{4n}.
\end{equation*}\quad

\paragraph{\bf Hauptmodul}
We define the {\it hauptmodul} of $\Gamma_0(12)$:
\begin{equation}
\begin{split}
J_{12} &:= \Delta_{12}^0 / \Delta_{12}^{\infty} \: (= \eta^3(z) \eta^{-2}(2 z) \eta^{-1}(3 z) \eta(4 z) \eta^2(6 z) \eta^{-3}(12 z))\\
 &= \frac{1}{q} - 3 + 2 q + q^3 - 2 q^7 - \cdots,
\end{split}
\end{equation}
where $v_{\infty}(J_{12}) = -1$ and $v_0(J_{12}) = 1$. Then, we have
\begin{equation}
J_{12} : \partial \mathbb{F}_{12} \setminus \{z \in \mathbb{H} \: ; \: Re(z) = \pm 1/2\} \to [-6, 0] \subset \mathbb{R}.
\end{equation}

%% file: report.bbl
\begin{thebibliography}{MNS}

\bibitem[ACMS]{ACMS}
D. Alexander, C. Cummins J. McKay, and C. Simons, {\it Completely replicable functions}. In: {\it Groups, combinatorics $\&$ geometry} (M. Liebeck and J. Saxl eds.), 87--98, London Math. Soc. Lecture Note Ser., No. 165, Cambridge Univ. Press, Cambridge, 1992. (Proceedings of the L.M.S. Symposium on Groups and Combinatorics, Durham, 1990.)

\bibitem[AKN]{AKN}
T. Asai, M. Kaneko, and H. Ninomiya, {\it Zeros of certain modular functions and an application}, Comment.
Math. Univ. St. Paul. {\bfseries 46} (1997), 93--101.

\bibitem[BKM]{BKM}
E. Bannai, K. Kojima, and T. Miezaki, {\it On the zeros of Hecke type Faber polynomial}, to appear in Kyushu J. Math.

\bibitem[CN]{CN}
J. H. Conway, S. P. Norton, {\it Monstrous moonshine}, Bull. London Math. Soc., {\bfseries 11}(1979), 308--339.

\bibitem[G]{G}
J. Getz, {\it A generalization of a theorem of Rankin and Swinnerton-Dyer on zeros of modular forms}, Proc. Amer. Math. Soc., {\bfseries 132}(2004), No. 8, 2221--2231.

\bibitem[H]{H}
H. Hahn, {\it On zeros of Eisenstein series for genus zero Fuchsian groups}, Proc. Amer. Math. Soc. {\bfseries 135}(2007), No. 8, 2391--2401.

\bibitem[Ko]{Ko}
N. Koblitz, {\it Introduction to Elliptic Curves and Modular Forms}, Graduate Texts in Mathematics, No. 97, Springer-Verlag, New York, 1984.

\bibitem[Kr]{Kr}
A. Krieg, {\it Modular Forms on the Fricke Group.}, Abh. Math. Sem. Univ. Hamburg, {\bfseries 65}(1995), 293--299.

\bibitem[MNS]{MNS}
T. Miezaki, H. Nozaki, and J. Shigezumi, {\it On the zeros of Eisenstein series for $\Gamma_0^* (2)$ and $\Gamma_0^* (3)$}, J. Math. Soc. Japan, {\bfseries 59}(2007), 693--706.

\bibitem[Q]{Q}
H. -G. Quebbemann, {\it Atkin-Lehner eigenforms and strongly modular lattices}, Enseign. Math. (2), {\bfseries 43}(1997), No. 1-2, 55--65.

\bibitem[RSD]{RSD}
F. K. C. Rankin, H. P. F. Swinnerton-Dyer, {\it On the zeros of Eisenstein Series}, Bull. London Math. Soc., {\bfseries 2}(1970), 169--170.

\bibitem[Se]{Se}
J. -P. Serre, {\it A Course in Arithmetic}, Graduate Texts in Mathematics, No. 7, Springer-Verlag, New York-Heidelberg, 1973. (Translation of {\it Cours d'arithm\'etique $($French$)$}, Presses Univ. France, Paris, 1970.)

\bibitem[SG]{SG}
G. Shimura, {\it On Eisenstein Series}, Duke Math. J., {\bfseries 50}(1983), No. 2, 417--476.

\bibitem[SH]{SH}
H. Shimizu, {\it Hokei kansu. I-III. $($Japanese$)$ $[$Automorphic functions. I-III$]$}, Iwanami Shoten Kiso Sugaku [Iwanami Lectures on Fundamental Mathematics] 8, Iwanami Shoten Publishers, Tokyo, 1977--1978.

\bibitem[SJ1]{SJ1}
J. Shigezumi, {\it On the zeros of Eisenstein series for $\Gamma_0^{*}(p)$ and $\Gamma_0(p)$ of low levels}, M.S. thesis, Kyushu University, 2006.

\bibitem[SJ2]{SJ2}
J. Shigezumi, {\it On the zeros of the Eisenstein series for $\Gamma_0^{*}(5)$ and $\Gamma_0^{*}(7)$}, Kyushu. J. Math. {\bfseries 61}(2007), 527--549.

\end{thebibliography}
